\documentclass[12pt]{article}
\usepackage{amsthm, amssymb}
\usepackage{amsfonts}
\usepackage{epsfig,multicol}
\title{Planar Algebras, I}
\author{V. F. R. Jones
\thanks{Research supported in part by NSF Grant DMS93--22675, the Marsden fund
UOA520, and the Swiss National Science Foundation.}
}
\date{Department of Mathematics\\University of
    California\\Berkeley, California 94720--3840}
\begin{document}
\maketitle

\newtheorem{prop}{Proposition}[section]
\newtheorem{lem}{Lemma}[section]
\newtheorem{thm}{Theorem}[section]
\newtheorem{cor}{Corollary}[section]
\newtheorem{rem}{Remark}[section]
\newtheorem{note}{Note}[section]
\newtheorem{exm}{Examples}
\newcommand{\dsize}{\displaystyle}


\def\mapleftright#1{\smash{\mathop{\longleftrightarrow}\limits^{#1}}}

\def\beye{\mathop{\setlength{\unitlength}{1pt}
\begin{picture}(12,15)(0,0)
\put(-4,-3.5){
\epsfig{file=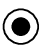}}
\end{picture}}}

\def\weye{\mathop{\setlength{\unitlength}{1pt}
\begin{picture}(12,15)(0,0)
\put(-4,-3.5){
\epsfig{file=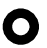}}
\end{picture}}}

\def\deye{\mathop{\setlength{\unitlength}{1pt}
\begin{picture}(12,15)(0,0)
\put(-4,-3.5){
\epsfig{file=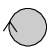}}
\end{picture}}}

\def\deyea{\mathop{\setlength{\unitlength}{1pt}
\begin{picture}(12,15)(0,0)
\put(-4,-3.5){
\epsfig{file=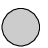}}
\end{picture}}}

\def\deyer{\mathop{\setlength{\unitlength}{1pt}
\begin{picture}(12,15)(0,0)
\put(-4,-3.5){
\epsfig{file=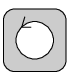}}
\end{picture}}}

\def\bnd{\mathop{\setlength{\unitlength}{1pt}  
\begin{picture}(24,18)(0,0)
\put(-4,-3.5){
\epsfig{file=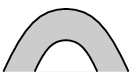}}
\end{picture}}}

\def\caleft{\mathop{\setlength{\unitlength}{1pt}
\begin{picture}(12,14)(0,0)
\put(-4,-3.5){
\epsfig{file=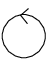}}
\end{picture}}}

\def\caright{\mathop{\setlength{\unitlength}{1pt}
\begin{picture}(12,14)(0,0)
\put(-4,-3.5){
\epsfig{file=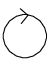}}
\end{picture}}}

\def\tangle#1{
        \setlength{\unitlength}{1pt}
        \begin{picture}(10,14)(0,0)
        \put(-4,-8.0){
                        \epsfig{file=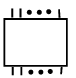}
                     }
        \end{picture}
        \begin{picture}(10,14)(0,0)
        \put(-5,0.8){
                        \makebox(0,-2){\smash{$\textstyle#1$}}               
                 }
        \end{picture}
             }
\def\tangler#1{
        \setlength{\unitlength}{1pt}
        \begin{picture}(10,14)(0,0)
        \put(-4,-8.0){
                        \epsfig{file=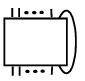}
                     }
        \end{picture}
        \begin{picture}(10,14)(0,0)
        \put(-5,0.8){
                        \makebox(0,-2){\smash{$\textstyle#1$}}
                 }
        \end{picture}
             }
\def\tanglel#1{
        \setlength{\unitlength}{1pt}
        \begin{picture}(10,14)(0,0)
        \put(-4,-8.0){
                        \epsfig{file=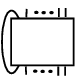}
                     }
        \end{picture} 
        \begin{picture}(10,14)(0,0)
        \put(-2,0.8){
                        \makebox(0,-2){\smash{$\textstyle#1$}}
                 }
        \end{picture}
             }

\def\tanglearrow#1{
        \setlength{\unitlength}{1pt}
        \begin{picture}(14,14)(0,0)
        \put(-4,-7.5){
                        \epsfig{file=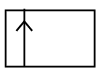}
                     }
        \end{picture}
        \begin{picture}(10,14)(0,0)
        \put(0,0.8){
                        \makebox(0,-2){\smash{$\textstyle#1$}}
                 }
        \end{picture}
             }
\def\tangleudarrow#1{
        \setlength{\unitlength}{1pt}
        \begin{picture}(24,24)(0,0)
        \put(-4,-8.0){
                        \epsfig{file=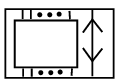}
                     }
        \end{picture}
        \begin{picture}(10,14)(0,0)
        \put(-14, 1.5){
                        \makebox(0,-2){\smash{$\textstyle#1$}}
                 }
        \end{picture}
             }

\def\boxarrow#1{
        \setlength{\unitlength}{1pt}
        \begin{picture}(30,24)(0,0)
        \put(-4,-8.0){
                        \epsfig{file=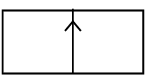}
                     }
        \end{picture}
        \begin{picture}(10,14)(0,0)
        \put(-14, 1.5){
                        \makebox(0,-2){\smash{$\textstyle#1$}}
                 }
        \end{picture}
             }
\def\dtangle#1{ 
        \setlength{\unitlength}{1pt}
        \begin{picture}(30,26)(0,0)
        \put(-4,-11.0){
                        \epsfig{file=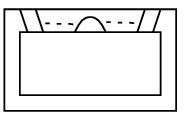}
                     }
        \end{picture}
        \begin{picture}(10,14)(0,0)
        \put(-7, 1.5){ 
                        \makebox(0,-2){\smash{$\textstyle#1$}}
                 }
        \end{picture}
             }

\def\blankbox#1{
        \setlength{\unitlength}{1pt}
        \begin{picture}(10,14)(0,0)
        \put(-4,-6.0){
                        \epsfig{file=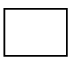}
                     }
        \end{picture}
        \begin{picture}(10,14)(0,0)
        \put(-5,0.8){
                        \makebox(0,-2){\smash{$\textstyle#1$}}
                 }
        \end{picture} 
             }
\def\boxlloop#1{
        \setlength{\unitlength}{1pt}
        \begin{picture}(16,28)(0,0)
        \put(-4,-10.5){
                        \epsfig{file=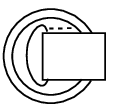}
                     }
        \end{picture}
        \begin{picture}(16,16)(0,0)
        \put(0, 1.5){
                        \makebox(0,-2){\smash{$\textstyle#1$}}
                 }
        \end{picture}
             }
\def\boxrloop#1{    
        \setlength{\unitlength}{1pt}
        \begin{picture}(16,28)(0,0)
        \put(-4,-10.5){
                        \epsfig{file=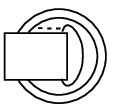}
                     }
        \end{picture}
        \begin{picture}(16,16)(0,0)
        \put(-10, 1.5){ 
                        \makebox(0,-2){\smash{$\textstyle#1$}}
                 }
        \end{picture}
             }
\def\blright#1{ 
        \setlength{\unitlength}{1pt}
        \begin{picture}(16,28)(0,0)
        \put(-4,-10.5){
                        \epsfig{file=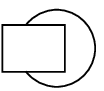} 
                     }  
        \end{picture}   
        \begin{picture}(16,16)(0,0)
        \put(-10, 1.5){
                        \makebox(0,-2){\smash{$\textstyle#1$}}
                 }
        \end{picture}
             }
\def\blrighto#1{
        \setlength{\unitlength}{1pt}
        \begin{picture}(16,18)(0,0)
        \put(-4,-10.5){
                        \epsfig{file=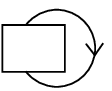}
                     }
        \end{picture}
        \begin{picture}(16,16)(0,0)
        \put(-10, 1.5){
                        \makebox(0,-2){\smash{$\textstyle#1$}}
                 }
        \end{picture}
             }
\def\bllefto#1{  
        \setlength{\unitlength}{1pt}
        \begin{picture}(16,18)(0,0)
        \put(-4,-10.5){
                        \epsfig{file=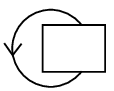}
                     }  
        \end{picture}   
        \begin{picture}(16,16)(0,0)
        \put(0, 1.5){
                        \makebox(0,-2){\smash{$\textstyle#1$}}
                 }
        \end{picture}
             }

\def\boxas#1{
        \setlength{\unitlength}{1pt}
        \begin{picture}(10,24)(0,0)
        \put(-4,-8.0){
                        \epsfig{file=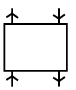}
                     }
        \end{picture}
        \begin{picture}(10,14)(0,0)
        \put(-5,0.8){
                        \makebox(0,-2){\smash{$\textstyle#1$}}
                 }
        \end{picture}
             }
\def\boxasR{  
        \setlength{\unitlength}{1pt}
        \begin{picture}(10,24)(0,0)
        \put(-4,-8.0){
                        \epsfig{file=ps/bar.ps}
                     }
        \end{picture}
        \begin{picture}(10,14)(0,0)
        \put(-5,10){ 
                        \makebox(0,-2){\smash{$\textstyle\rotatebox{180}{R}$}}
                 }
        \end{picture}
             }

\def\boxa#1{
        \setlength{\unitlength}{1pt}
        \begin{picture}(10,15)(0,0)  
        \put(-4,-9.0){
                        \epsfig{file=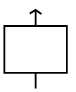}
                     }
        \end{picture}
        \begin{picture}(10,14)(0,0)
        \put(-5,0.8){
                        \makebox(0,-2){\smash{$\textstyle#1$}}
                 }
        \end{picture}
             }
\def\boxb#1{
        \setlength{\unitlength}{1pt}
        \begin{picture}(10,15)(0,0)
        \put(-4,-9.0){
                        \epsfig{file=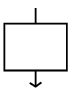}
                     }
        \end{picture}
        \begin{picture}(10,14)(0,0)
        \put(-5,0.8){
                        \makebox(0,-2){\smash{$\textstyle#1$}}
                 }
        \end{picture}
             }

\def\boxe#1{
        \setlength{\unitlength}{1pt}
        \begin{picture}(10,15)(0,0)
        \put(-4,-9.0){
                        \epsfig{file=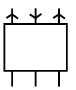}
                     }
        \end{picture}
        \begin{picture}(10,14)(0,0)
        \put(-5,0.8){
                        \makebox(0,-2){\smash{$\textstyle#1$}}
                 }
        \end{picture}
             }
\def\boxf#1{
        \setlength{\unitlength}{1pt}
        \begin{picture}(10,15)(0,0)
        \put(-4,-9.0){
                        \epsfig{file=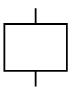}  
                     }
        \end{picture} 
        \begin{picture}(10,14)(0,0)
        \put(-5,0.8){
                        \makebox(0,-2){\smash{$\textstyle#1$}}
                 }
        \end{picture}
             }
\def\boxg#1{
        \setlength{\unitlength}{1pt}
        \begin{picture}(10,24)(0,0)
        \put(-4,-9.0){
                        \epsfig{file=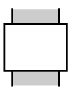}
                     }
        \end{picture}
        \begin{picture}(10,14)(0,0)
        \put(-5,0.8){   
                        \makebox(0,-2){\smash{$\textstyle#1$}}
                 }
        \end{picture}
             }
\def\boxh#1{
        \setlength{\unitlength}{1pt}
        \begin{picture}(10,24)(0,0)
        \put(-4,-11.0){
                        \epsfig{file=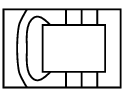} 
                     }
        \end{picture}
        \begin{picture}(10,14)(0,0)
        \put(7,0.8){
                        \makebox(0,-2){\smash{$\textstyle#1$}}
                 }
        \end{picture}\hskip15pt
             }
\def\boxi#1{
        \setlength{\unitlength}{1pt}
        \begin{picture}(10,24)(0,0)
        \put(-4,-11.0){
                        \epsfig{file=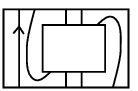}
                     }
        \end{picture}
        \begin{picture}(10,14)(0,0)
        \put(7,0.8){
                        \makebox(0,-2){\smash{$\textstyle#1$}}
                 }
        \end{picture}\hskip16pt
             }

\def\boxj{
        \setlength{\unitlength}{1pt}
        \begin{picture}(10,18)(0,0)
        \put(-4,-5){   
                        \epsfig{file=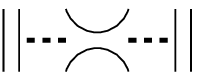}
                     }
        \end{picture}
        \begin{picture}(10,18)(0,0)
        \put(12,-8){ 
                        \makebox(0,-4){\smash{$\scriptstyle{i\ \ \  i+1}$}}
                 }
        \end{picture}
             \ \hskip33pt  }
\def\boxk#1{
        \setlength{\unitlength}{1pt}
        \begin{picture}(10,24)(0,0)
        \put(-4,-11){
                        \epsfig{file=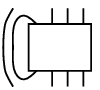}
                     }
        \end{picture}
        \begin{picture}(10,14)(0,0)
        \put(5,0.8){ 
                        \makebox(0,-2){\smash{$\textstyle#1$}}
                 }
        \end{picture}\hskip6pt
             }
\def\boxl#1{
        \setlength{\unitlength}{1pt}
        \begin{picture}(10,24)(0,0)
        \put(-4,-11){
                        \epsfig{file=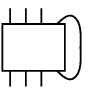}
                     }
        \end{picture}
        \begin{picture}(10,14)(0,0)
        \put(-5,0.8){
                        \makebox(0,-2){\smash{$\textstyle#1$}}
                 }
        \end{picture}\hskip3pt
             }
\def\boxlr#1{
        \setlength{\unitlength}{1pt}
        \begin{picture}(10,24)(0,0)
        \put(-4,-11){
                        \epsfig{file=ps/bbzkkrr.ps}
                     }
        \end{picture}
        \begin{picture}(10,14)(0,0)
        \put(3,0.8){
                        \makebox(0,-2){\smash{$\textstyle#1$}}
                 }
        \end{picture}\hskip3pt
             }
\def\bigtan#1{
        \setlength{\unitlength}{1pt}
        \begin{picture}(20,24)(0,0)
        \put(-8,-10){
                        \epsfig{file=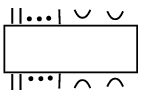}
                     }
        \end{picture}
        \begin{picture}(14,14)(0,0)
        \put(-5,0.8){
                        \makebox(0,-2){\smash{$\textstyle#1$}}
                 }
        \end{picture}\hskip3pt
             }
\def\bigtank{\hskip6pt
        \setlength{\unitlength}{1pt}
        \begin{picture}(20,24)(0,0)
        \put(-8,-5){   
                        \epsfig{file=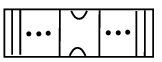}
                     }
        \end{picture}\hskip20pt 
             }

\def\butau{\hskip6pt
        \setlength{\unitlength}{1pt}
        \begin{picture}(24,24)(0,0)
        \put(-8,-10.0){
                        \epsfig{file=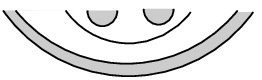}
                     }
        \end{picture}\hskip45pt
             }
\def\butad{\hskip6pt
        \setlength{\unitlength}{1pt}
        \begin{picture}(24,24)(0,0)
        \put(-8,-10.0){
                        \epsfig{file=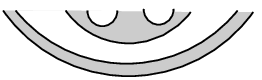}
                     }
        \end{picture}\hskip45pt
             }

\def\vln{
        \setlength{\unitlength}{1pt}
        \begin{picture}(7,16)(0,0)
        \put(0,-8.0){
                        \epsfig{file=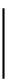}
                     }
        \end{picture}
        }

\def\fedgs{
        \setlength{\unitlength}{1pt}
        \begin{picture}(14,24)(0,0)
        \put(-8,-6.0){
                        \epsfig{file=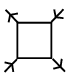}
                     }
        \end{picture}
             }
\def\crss{
        \setlength{\unitlength}{1pt}
        \begin{picture}(14,18)(0,0)
        \put(-8,-6.0){
                        \epsfig{file=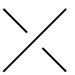}
                     }
        \end{picture}
             }
\def\smth{
        \setlength{\unitlength}{1pt}
        \begin{picture}(14,18)(0,0)
        \put(-8,-6.0){
                        \epsfig{file=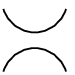}
                     }
        \end{picture}
             }
\def\smtho{
        \setlength{\unitlength}{1pt}
        \begin{picture}(14,18)(0,0)
        \put(-8,-6.0){  
                        \epsfig{file=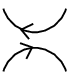}
                     }  
        \end{picture}
             }

\def\smthv{
        \setlength{\unitlength}{1pt}
        \begin{picture}(14,18)(0,0)
        \put(-8,-6.0){
                        \epsfig{file=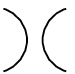}
                     }
        \end{picture}
             }
\def\crsso{
        \setlength{\unitlength}{1pt}
        \begin{picture}(14,18)(0,0)
        \put(-8,-6.0){
                        \epsfig{file=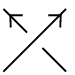}
                     }
        \end{picture}
             }
\def\crssor{
        \setlength{\unitlength}{1pt}
        \begin{picture}(14,18)(0,0)
        \put(-8,-6.0){
                        \epsfig{file=ps/crosso.ps}
                     }
        \end{picture}
             }
\def\crssr{ 
        \setlength{\unitlength}{1pt}
        \begin{picture}(14,18)(0,0)
        \put(-8,-6.0){
                        \epsfig{file=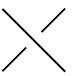} 
                     }  
        \end{picture}   
             }       

\def\crsst{
        \setlength{\unitlength}{1pt}
        \begin{picture}(14,24)(0,0)
        \put(-8,-6.0){
                        \epsfig{file=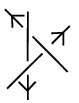}
                     }
        \end{picture}
             }

\def\smthvo{
        \setlength{\unitlength}{1pt}
        \begin{picture}(14,18)(0,0)
        \put(-8,-6.0){
                        \epsfig{file=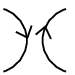}
                     }
        \end{picture}
             }
\def\smthvor{
        \setlength{\unitlength}{1pt} 
        \begin{picture}(14,18)(0,0)
        \put(-8,-6.0){
                        \epsfig{file=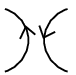}
                     }  
        \end{picture}   
             }       

\def\paro{  
        \setlength{\unitlength}{1pt}
        \begin{picture}(14,24)(0,0)
        \put(-8,-6.0){ 
                        \epsfig{file=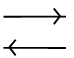} 
                     }
        \end{picture}
             }

\def\blong{\
        \setlength{\unitlength}{1pt}
        \begin{picture}(16,40)(0,0)
        \put(0,-12){
                        \epsfig{file=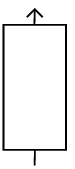}
                     }
        \end{picture}
        \begin{picture}(10,40)(0,0)
        \put(-4,14){
                        \makebox(0,-12){\smash
                {$
                \textstyle{
                        \begin{array}{c}
                        \gamma\vspace{6pt}\\
                                \gamma
                        \end{array}
                          }
                 $}
                                        }
                     }
        \end{picture}
             }
\def\bsp{
        \setlength{\unitlength}{1pt}
        \begin{picture}(40,20)(0,0)
        \put(0,-5){
                        \epsfig{file=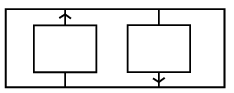}
                     }
        \end{picture}
        \begin{picture}(40,20)(0,0)
        \put(-6,10){
                        \makebox(0,-10){\smash
                {$
                \textstyle{
                        \gamma\ \ \ \ \ \delta  
                          }
                 $}
                                        }
                     }
        \end{picture}
             }

\def\op#1{\bf\rm{#1}} 

\def\qed{{\unskip\nobreak\hfil\penalty50\hskip2em\vadjust{}\nobreak\hfil
    $\square$\parfillskip=0pt\finalhyphendemerits=0\par}}

\def\normalbaselines{\baselineskip20pt
  \lineskip3pt \lineskiplimit3pt }
\def\eqdef#1{\smash{\mathop{=}\limits^{#1}}} 
             
        \newcommand\x{\times}
        \def\<{\langle}
        \def\>{\rangle}
        \renewcommand\a{\alpha}
        \renewcommand\b{{\cal B}}
        \newcommand\bs{\bigskip}
        \newcommand\dt{\cdot}
        \renewcommand\d{\delta} 
        \newcommand\D{\Delta}
        \newcommand\e{\varepsilon}
        \newcommand\g{\gamma} 
        \newcommand\G{\Gamma}
        \newcommand\J{{\cal J}} 
        \renewcommand\l{\lambda} 
        \renewcommand\L{\Lambda}
        \renewcommand\ni{\noindent}
        \newcommand\var{\varphi}
        \renewcommand\r{\mathbb R^2}
        \newcommand\s{\sigma}
        \renewcommand\th{\theta}
        \renewcommand\O{\Omega}
        \renewcommand\o{\omega}  
        \newcommand\z{\mathbb Z}
        \renewcommand\i{\infty}
        \newcommand\p{\partial}
        \def\text#1{\scriptstyle{\rm #1}}       
        \renewcommand\P{{\cal P}}

\baselineskip=16 pt
\begin{abstract}
We introduce a notion of planar
algebra, the simplest example of which is a vector
space of tensors, closed under planar
contractions. A planar algebra with suitable
positivity properties produces a finite index
subfactor of a II${}_1$ factor, and vice versa.
\end{abstract}

\baselineskip=16 pt

\noindent
{\large\bf 0. Introduction}

\medskip At first glance there is nothing planar
about a subfactor. A factor $M$ is a unital
$*$-algebra of bounded linear operators on a
Hilbert space, with trivial centre and closed in
the topology of pointwise convergence. The factor
$M$ is of type II${}_1$ if it admits a
(normalized) trace, a linear function $tr$:
$M\to\Bbb C$ with $tr(ab)=tr(ba)$ and
tr$(1)=1$. In [J1] we defined the notion of {\it
index} $[M:N]$ for II${}_1$ factors $N\subset M$.
The most surprising result of [J1] was that
$[M:N]$ is ``quantized" --- to be precise, if
$[M:N] <4$ there is an integer $n\geq 3$ with
$[M:N]=4\cos^2\pi/n$. This led to a surge of
interest in subfactors and the major theorems of
Pimsner, Popa and Ocneanu ([PP],[Po1],[O1]). These
results turn around a ``standard invariant" for
finite index subfactors, also known variously as
the ``tower of relative commutants", the
``paragroup", or the ``$\l$-lattice". In favorable
cases the standard invariant allows one to
reconstruct the subfactor, and both the paragroup
and $\l$-lattice approaches give complete
axiomatizations of the standard invariant. In this
paper we give, among other things, yet another
axiomatization which has the advantage of
revealing an underlying {\it planar} structure not
apparent in other approaches.  It also places the standard
invariant in a larger mathematical context. In
particular we give a rigorous justification for
pictorial proofs of subfactor theorems.
Non-trivial results have already been obtained
from such arguments in [BJ1]. The standard
invariant is sufficiently rich to justify several
axiomatizations --- it has led to the discovery of
invariants in knot theory ([J2]), 3-manifolds
([TV]) and combinatorics ([NJ]), and is of
considerable interest in conformal and algebraic
quantum field theory ([Wa],[FRS],[Lo]).

Let us now say exactly what we mean by a planar
algebra. The best language to use is that of
operads ([Ma]). We define the {\it planar operad},
each element of which determines a multilinear
operation on the standard invariant.

A planar $k$-tangle will consist of the unit disc
$D$ $(=D_0)$ in $\Bbb C$ together with a finite
(possibly empty) set of disjoint subdiscs
$D_1,D_2,\dots ,D_n$ in the interior of $D$. Each
disc $D_i$, $i\geq 0$, will have an even number
$2k_i\geq 0$ of marked points on its boundary
(with $k=k_0$). Inside $D$ there is also a finite
set of disjoint smoothly embedded curves called
{\it strings} which are either closed curves or
whose boundaries are marked points of the $D_i$'s.
Each marked point is the boundary point of some
string, which meets the boundary of the
corresponding disc transversally. The strings all
lie in the complement of the interiors $\stackrel{\circ}{D_i}$
of the $D_i$, $i\geq 0$. The connected components of
the complement of the strings in 
$\dsize{\stackrel{\circ}{D}\backslash \bigcup^n_{i=1}D_i}$ 
are called regions
and are shaded black and white so that regions
whose closures meet have different shadings. The
shading is part of the data of the tangle, as is
the choice, at every $D_i$, $i\geq 0$, of a white
region whose closure meets that disc. 
The case $k=0$ is exceptional - 
there are two kinds of 0-tangle, according to
whether the region near the boundary is shaded
black or white. An example of a planar 4-tangle,
where the  chosen white regions are marked with a
$*$ close to their respective discs, is given
below.

\[
\input{xfig/pic1}
\]
The {\it planar operad} $\Bbb P$ is the set of all
orientation-preserving diffeomorphism classes of
planar $k$ tangles, $k$ being arbitrary. The
diffeomorphisms preserve the boundary of $D$ but
may move the $D_i$'s, $i>1$.

Given a planar $k$ tangle $T$, a $k'$-tangle $S$,
and a disk $D_i$ of $T$ with $k_i=k'$ we define
the  $k$ tangle $T\circ_i S$ by isotoping $S$ so
that its boundary, together with the marked
points, coincides with that of $D_i$, and the
chosen white regions for $D_i$ (in $T$) and $S$ share a
boundary segment.  The strings may then be joined
at the boundary of $D_i$ and smoothed. The
boundary of $D_i$ is then removed to obtain the
tangle $T\circ_i S$ whose diffeomorphism class
clearly depends only on those of $T$ and $S$.
This gives $\Bbb P$ the structure of a coloured
operad, where each $D_i$ for $i>0$ is assigned
the colour $k_i$ and composition is only allowed
when the colours match. There are two distinct
colours for $k=0$ according to the shading near
the boundary. The $D_i$'s for $i\geq 1$ are to be
thought of as inputs and $D_0$ is the output. (In
the usual definition of an operad the inputs are
labelled and the symmetric group $S_n$ acts on
them. Because of the colours, $S_n$ is here replaced by
$S_{n_1}\x S_{n_2}\x\dots\x S_{n_p}$ where  $n_j$ is the
number of internal discs coloured $j$.
Axioms for such a coloured operad could be given
along the lines of [Ma] but we do not need them
since we have a concrete example.) The picture
below exhibits the composition

\[
	\input{xfig/pic2}
\]

The most general notion of a planar algebra that
we will contemplate is that of an algebra over
$\Bbb P$ in the sense of [Ma].  That is to say, first
of all, a disjoint union $V_k$ of vector spaces
for $k>0$ and two vector spaces
$V_0^{\text{white}}$ and $V_0^{\text{black}}$
(which we will call $P_0$ and $P_{1,1}$ later on).
Linear maps between tensor powers of these vector
spaces form a coloured operad \ Hom \ in the
obvious way under composition of maps and the
planar algebra structure  on the $V$'s is given by
a morphism of coloured operads from $\Bbb P$ to Hom. 
In
practice this means that, to a $k$-tangle $T$ in
$\Bbb P$ there is a linear map
$Z(T):\bigotimes^n_{i=1}V_{k_i}\to V_k$ such that
$Z(T\circ_i S)=Z(T)\circ_i Z(S)$ where the
$\circ_i$ on the right-hand side is composition of
linear maps in \ Hom \ .

Note that the vector spaces $V_0^{\text{white}}$ and $V_0^{\text{black}}$
may be different. This is the case for the ``spin
models" of $\S$3. Both these $V_0$'s become
commutative associative algebras using the tangles

\[
	\begin{picture}(0,0)%
\epsfig{file=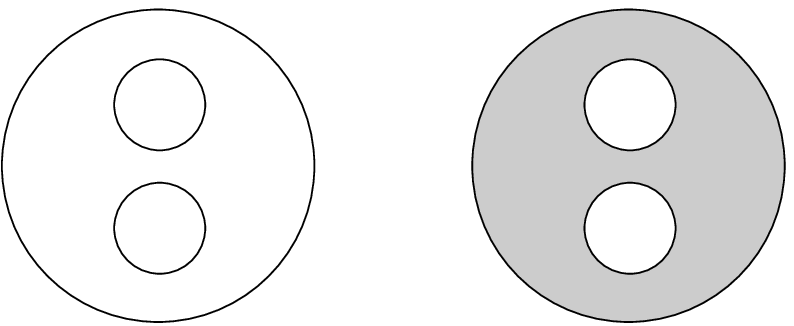}%
\end{picture}%
\setlength{\unitlength}{0.00041700in}%
\begingroup\makeatletter\ifx\SetFigFont\undefined
\def\x#1#2#3#4#5#6#7\relax{\def\x{#1#2#3#4#5#6}}%
\expandafter\x\fmtname xxxxxx\relax \def\y{splain}%
\ifx\x\y   
\gdef\SetFigFont#1#2#3{%
  \ifnum #1<17\tiny\else \ifnum #1<20\small\else
  \ifnum #1<24\normalsize\else \ifnum #1<29\large\else
  \ifnum #1<34\Large\else \ifnum #1<41\LARGE\else
     \huge\fi\fi\fi\fi\fi\fi
  \csname #3\endcsname}%
\else
\gdef\SetFigFont#1#2#3{\begingroup
  \count@#1\relax \ifnum 25<\count@\count@25\fi
  \def\x{\endgroup\@setsize\SetFigFont{#2pt}}%
  \expandafter\x
    \csname \romannumeral\the\count@ pt\expandafter\endcsname
    \csname @\romannumeral\the\count@ pt\endcsname
  \csname #3\endcsname}%
\fi
\fi\endgroup
\begin{picture}(7545,3028)(871,-4570)
\put(2161,-3821){\makebox(0,0)[lb]{\smash{\SetFigFont{12}{14.4}{rm}$D_2$}}}
\put(2146,-2621){\makebox(0,0)[lb]{\smash{\SetFigFont{12}{14.4}{rm}$D_1$}}}
\put(4396,-3161){\makebox(0,0)[lb]{\smash{\SetFigFont{12}{14.4}{rm}and}}}
\put(6676,-3821){\makebox(0,0)[lb]{\smash{\SetFigFont{12}{14.4}{rm}$D_2$}}}
\put(6661,-2621){\makebox(0,0)[lb]{\smash{\SetFigFont{12}{14.4}{rm}$D_1$}}}
\end{picture}

\]

To handle tangles with no internal discs we decree
that the tensor product over the empty set be the
field $K$ and identify Hom$(K,V_k)$ with $V_k$ so
that each $V_k$ will contain a privileged subset
which is $Z(\{ k$-tangles with no internal
discs$\})$. This is the ``unital"
 structure (see [Ma]).

One may want to impose various conditions such as
dim$(V_k)<\i$ for all $k$. The condition
dim$(V_0^{\text{white}})=1=$dim$(V_0^{\text{black}})$
is significant and we impose it in our formal
definition of planar algebra (as opposed to
general planar algebra) later on.  It implies that
there is a unique way to identify each $V_0$ with
$K$ as algebras, and $Z(\bigcirc)=1=Z(\deyea )$.
There are thus also two scalars associated to a
planar algebra, $\d_1=Z(\weye )$ and
$\d_2=Z(\beye )$ (the inner circles are strings,
not discs!). It follows that  $Z$ is
multiplicative on connected components, i.e., if a
part of a tangle $T$ can be surrounded by a disc
so that $T=T'\circ_i S$ for a tangle $T'$ and a
0-tangle $S$, then $Z(T)=Z(S)Z(T')$ where $Z(S)$
is a multilinear map into the field $K$.

Two simple examples serve as the keys to
understanding the notion of a planar algebra. The
first is the Temperley-Lieb algebra $TL$, some
vestige of which is present in every planar
algebra. The vector spaces $TL_k$ are:
$$
TL_0^{\text{black}}\simeq
TL_0^{\text{white}}\simeq K
$$
and $TL_k$ is the vector space  whose basis is the
set of diffeomorphism classes of connected planar
$k$-tangles with no internal discs, of which there
are $\frac{1}{k+1}{2k\choose k}$. The
action of a planar $k$-tangle on $TL$ is almost
obvious --- when one fills the internal discs of a
tangle with basis elements of $TL$ one obtains
another basis element, except for some simple
closed curves. Each closed curve counts a
multiplicative factor of $\d$ and then is removed.
It is easily verified that this defines an action
of $\Bbb P$ on $TL$.
As we have observed, any planar algebra contains
elements corresponding to the $TL$ basis. They are
not necessarily linearly independent. See [GHJ]
and $\S$2.1.

 The second key example of a planar algebra is
given by tensors. We think of a tensor as an
object which yields a number each time its indices
are specified. Let $V_k$ be the vector space of
tensors with $2k$ indices. An element of $\Bbb P$
gives a scheme for contracting tensors, once a
tensor is assigned to each internal disc. The
indices lie on the strings and are locally
constant thereon. The boundary indices are fixed
and are the indices of the output tensor. All
indices  on strings not touching $D$ are summed
over and one contracts by taking, for a given set
of indices, the product of the values of the
tensors in the internal discs. One recognizes the
partition function of a statistical mechanical
model ([Ba]), the boundary index values being the
boundary conditions and the tensor values being
the Boltzmann weights. This diagrammatic
contraction calculus for tensors is well known
([Pe]) but here we are only considering {\it
planar} contraction systems. If the whole planar
algebra of all tensors were the only example this
subject would be of no interest, but in fact there
is a huge family of planar subalgebras --- vector
spaces of tensors closed under planar contractions
--- rich enough to contain the theory of finitely
generated groups and their Cayley graphs. See
$\S$2.7.

The definition of planar algebra we give in $\S$1
is not the operadic one. When the planar algebra
structure first revealed itself, the $V_k$'s
already had an associative algebra structure
coming from the von Neumann algebra context. Thus
our definition will be in terms of a universal planar
algebra on some set of generators (labels) which
can be combined in arbitrary planar fashion. The
discs we have used above become boxes in section one, reflecting
the  specific algebra structure we began with. The
equivalence of the two definitions is completed in
Proposition 1.20. The main ingredient of the equivalence is that the
planar operad $\Bbb P$ is generated by the Temperley-Lieb algebra
and tangles of two kinds:
\begin{enumerate}
\item Multiplication, which is the following
tangle (illustrated for k =5)

\[
	\begin{picture}(0,0)%
\epsfig{file=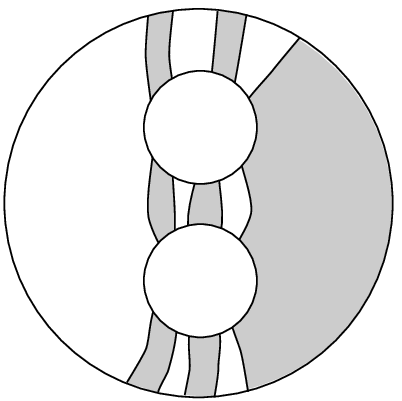}%
\end{picture}%
\setlength{\unitlength}{0.00041700in}%
\begingroup\makeatletter\ifx\SetFigFont\undefined
\def\x#1#2#3#4#5#6#7\relax{\def\x{#1#2#3#4#5#6}}%
\expandafter\x\fmtname xxxxxx\relax \def\y{splain}%
\ifx\x\y   
\gdef\SetFigFont#1#2#3{%
  \ifnum #1<17\tiny\else \ifnum #1<20\small\else
  \ifnum #1<24\normalsize\else \ifnum #1<29\large\else
  \ifnum #1<34\Large\else \ifnum #1<41\LARGE\else
     \huge\fi\fi\fi\fi\fi\fi
  \csname #3\endcsname}%
\else
\gdef\SetFigFont#1#2#3{\begingroup
  \count@#1\relax \ifnum 25<\count@\count@25\fi
  \def\x{\endgroup\@setsize\SetFigFont{#2pt}}%
  \expandafter\x
    \csname \romannumeral\the\count@ pt\expandafter\endcsname
    \csname @\romannumeral\the\count@ pt\endcsname
  \csname #3\endcsname}%
\fi
\fi\endgroup
\begin{picture}(3759,3756)(661,-4406)
\put(2319,-1916){\makebox(0,0)[lb]{\smash{\SetFigFont{12}{14.4}{rm}$D_1$}}}
\put(2318,-3417){\makebox(0,0)[lb]{\smash{\SetFigFont{12}{14.4}{rm}$D_2$}}}
\put(1718,-1039){\makebox(0,0)[lb]{\smash{\SetFigFont{12}{14.4}{rm}$*$}}}
\put(1823,-1624){\makebox(0,0)[lb]{\smash{\SetFigFont{12}{14.4}{rm}$*$}}}
\put(1823,-3048){\makebox(0,0)[lb]{\smash{\SetFigFont{12}{14.4}{rm}$*$}}}
\end{picture}

\]

\item Annular tangles: ones with only one internal
disc, e.g.,
\end{enumerate}

\[
	\begin{picture}(0,0)%
\epsfig{file=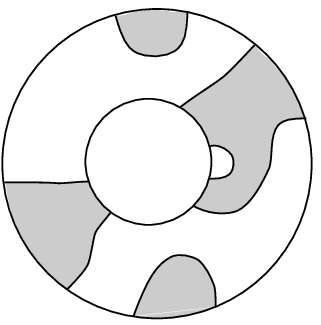}%
\end{picture}%
\setlength{\unitlength}{0.00041700in}%
\begingroup\makeatletter\ifx\SetFigFont\undefined
\def\x#1#2#3#4#5#6#7\relax{\def\x{#1#2#3#4#5#6}}%
\expandafter\x\fmtname xxxxxx\relax \def\y{splain}%
\ifx\x\y   
\gdef\SetFigFont#1#2#3{%
  \ifnum #1<17\tiny\else \ifnum #1<20\small\else
  \ifnum #1<24\normalsize\else \ifnum #1<29\large\else
  \ifnum #1<34\Large\else \ifnum #1<41\LARGE\else
     \huge\fi\fi\fi\fi\fi\fi
  \csname #3\endcsname}%
\else
\gdef\SetFigFont#1#2#3{\begingroup
  \count@#1\relax \ifnum 25<\count@\count@25\fi
  \def\x{\endgroup\@setsize\SetFigFont{#2pt}}%
  \expandafter\x
    \csname \romannumeral\the\count@ pt\expandafter\endcsname
    \csname @\romannumeral\the\count@ pt\endcsname
  \csname #3\endcsname}%
\fi
\fi\endgroup
\begin{picture}(3000,3018)(308,-3695)
\put(1193,-1076){\makebox(0,0)[lb]{\smash{\SetFigFont{12}{14.4}{rm}$*$}}}
\put(1351,-3061){\makebox(0,0)[lb]{\smash{\SetFigFont{12}{14.4}{rm}$*$}}}
\end{picture}

\]

The universal planar algebra is useful for
constructing planar algebras and restriction to
the two generating tangles sometimes makes it
shorter to check that a given structure is a
planar  algebra.
The algebra structure we begin with in $\S$1 corresponds of
course to the multiplication tangle given above.

The original algebra structure has been studied in
some detail (see $\S$3.1) but it should be quite
clear that the operad provides a vast family of
algebra structures on a planar algebra which we
have only just begun to appreciate. For instance,
the annular tangles above form an algebra over
which all the $V_k$'s in a planar algebra are
modules. This structure alone seems quite rich
([GL]) and we exploit it just a little to get
information on principal graphs of subfactors in
4.2.11. We have obtained more sophisticated
results in terms of generating functions  which we
will present in a future paper.

We present several examples of planar algebras in
$\S$2, but it is the connection with subfactors
that has been our main motivation and guide for
this work. The two leading theorems occur in
$\S$4. The first one shows how to obtain a planar
algebra from  a finite index subfactor $N\subset
M$. The vector space $V_k$ is the set of
$N$-central vectors in the $N-N$ bimodule
$M_{k-1}=M\otimes_N M\otimes_N\dots\otimes_N M$
($k$ copies of $M$), which, unlike $M_{k-1}$
itself, is finite dimensional. The planar algebra
structure on these $V_k$'s is obtained by a method
reminiscent of topological quantum field theory.
Given a planar  $k$-tangle $T$ whose internal
discs are labelled by elements of the $V_j$'s, we
have to show how to construct an element of $V_k$,
associated with the boundary of $T$, in a natural
way. One starts with a very small circle (the
``bubble") in the distinguished white region of
$T$, tangent to the boundary of $D$. We then allow
this circle to bubble out until it gets to the
boundary. On its way the bubble will have to cross
strings of the tangle and envelop internal discs.
As it does so it acquires shaded intervals which
are its intersections  with the shaded regions of
$T$. Each time the bubble envelops an internal
disc $D_i$, it acquires $k_i$ such shaded
intervals and, since an element of $V_{k_i}$ is a
tensor in $\otimes^{k_i}_N M$, we assign elements
of $M$ to the shaded region according to this
tensor.  There are also rules for assigning and
contracting tensors as the bubble crosses strings
of the tangle. At the end we have an element of
$\otimes^k_N M$ assigned to the boundary. This is
the action of the operad element on the vectors in
$V_{k_i}$. Once the element of $\otimes^k_N M$ has
been constructed and shown to be invariant under
diffeomorphisms, the formal operadic properties
are immediate.

One could try to carry out this procedure for an
arbitrary inclusion $A\subset B$ of rings, but
there are a few obstructions involved in showing
that our bubbling process is well defined. Finite
index (extremal) subfactors have all the special
properties required, though there are surely other
families of subrings for which the procedure is
possible.

The following tangle:

\[
	\begin{picture}(0,0)%
\epsfig{file=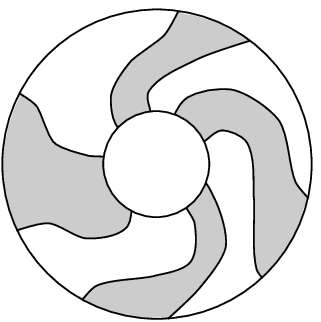}%
\end{picture}%
\setlength{\unitlength}{0.00041700in}%
\begingroup\makeatletter\ifx\SetFigFont\undefined
\def\x#1#2#3#4#5#6#7\relax{\def\x{#1#2#3#4#5#6}}%
\expandafter\x\fmtname xxxxxx\relax \def\y{splain}%
\ifx\x\y   
\gdef\SetFigFont#1#2#3{%
  \ifnum #1<17\tiny\else \ifnum #1<20\small\else
  \ifnum #1<24\normalsize\else \ifnum #1<29\large\else
  \ifnum #1<34\Large\else \ifnum #1<41\LARGE\else
     \huge\fi\fi\fi\fi\fi\fi
  \csname #3\endcsname}%
\else
\gdef\SetFigFont#1#2#3{\begingroup
  \count@#1\relax \ifnum 25<\count@\count@25\fi
  \def\x{\endgroup\@setsize\SetFigFont{#2pt}}%
  \expandafter\x
    \csname \romannumeral\the\count@ pt\expandafter\endcsname
    \csname @\romannumeral\the\count@ pt\endcsname
  \csname #3\endcsname}%
\fi
\fi\endgroup
\begin{picture}(3004,3012)(304,-3679)
\put(1846,-1609){\makebox(0,0)[lb]{\smash{\SetFigFont{12}{14.4}{rm}$*$}}}
\put(1441,-919){\makebox(0,0)[lb]{\smash{\SetFigFont{12}{14.4}{rm}$*$}}}
\end{picture}

\]

\noindent defines a rotation of period $k$ ($2k$
boundary points) so it is a consequence of the
planar algebra structure that the rotation
$x_1\otimes x_2\otimes\dots\otimes x_k\mapsto
x_2\otimes x_3\otimes\dots\otimes x_k\otimes x_1$,
which makes no sense on $M_{k-1}$, is well defined
on $N$-central vectors and has period $k$. This
result is in fact an essential technical
ingredient of the proof of Theorem 4.2.1.

Note that we seem to have avoided the use of
correspondences in the sense of Connes ([Co1]) by
working in the {\it purely algebraic} tensor
product. But the avoidance of $L^2$-analysis,
though extremely convenient, is a little illusory
since the proof of the existence and periodicity
of the rotation uses $L^2$ methods. The Ocneanu
approach ([EK]) uses the $L^2$ definition and the vector spaces $V_k$ are defined
as $Hom_{N,N}(\otimes^j_N M)$ and $Hom_{N,M}(\otimes^j_N M)$
depending on the parity of k. It is no doubt possible
to give a direct proof of Theorem 4.2.1 using this
definition - this would be the ``hom" 
version, our method being the ``$\otimes$" 
method.

To identify the operad structure with the usual
algebra structure on $M_{k-1}$ coming from the
``basic construction" of [J1], we show that the
multiplication tangle above does indeed define the
right formula. This, and a few similar details, is
suprisingly involved and accounts for some
unpleasant looking formulae in $\S$4. Several other
subfactor notions, e.g.~tensor product, are shown
to correspond to their planar algebra
counterparts, already abstractly defined in $\S$3.
Planar algebras also inspired, in joint work with
D.~Bisch, a notion of free product. We give the
definition here and will explore this notion in a
forthcoming paper with Bisch.

The second theorem of $\S$4 shows that one can
construct a subfactor from a planar algebra with
$*$-structure and suitable ``reflection"
positivity. It is truly remarkable that the axioms
needed by Popa for his construction of subfactors
in [Po2] follow so closely the axioms of planar
algebra, at least as formulated using boxes and
the universal planar algebra. For Popa's construction is
quite different from the ``usual" one of
[J1], [F+], [We1], [We2]. Popa uses an
amalgamated free product construction which
introduces an unsatisfactory element in the
correspondence between planar algebras and
subfactors. For although it is true that the
standard invariant of Popa's subfactor  is indeed
the planar algebra from which the subfactor was
constructed, it is {\it not} true that, if one
begins with a subfactor $N\subset M$, even
hyperfinite, and applies Popa's procedure to the
standard invariant, one obtains $N\subset M$ as a
result. There are many difficult questions here,
the main one of which is to decide when a given
planar algebra arises from a subfactor of the
Murray-von Neumann hyperfinite type II${}_1$
factor ([MvN]).

There is a criticism that has and should be made of 
our definition of a planar algebra - that it is
too restrictive. By enlarging the class of tangles in the   
planar operad, say so as to include oriented edges and 
boundary points, or discs with an odd number of 
boundary points, one would obtain a notion of planar
algebra applicable to more examples. For instance,
if the context were the study of group representations
our definition would have us studying say, SU(n) by 
looking at tensor powers of the form ${V\otimes\bar V\otimes V\otimes \bar V\dots}$
(where V is the defining representation on $\Bbb C^n$)

whereas a full categorical treatment would insist on
arbitrary tensor products. In fact, more general notions
already exist in the literature. Our planar algebras could
be formulated as a rather special kind of ``spider" in
the sense of Kuperberg in [Ku], or one could place them in
the context of pivotal and spherical categories ([FY],[BW]), 
and the theory of $C^*-$tensor categories even has the 
ever desirable positivity ([LR],[We3]). Also, in the semisimple
case at least, the work in section 3.3 on cabling and reduction 
shows how to extend our planar diagrams to ones with labelled
edges.

But it is the very restrictive nature of our definition of planar
algebras that should be its great virtue. We have good reasons
for limiting the generality. The most compelling is the 
equivalence with subfactors, which has been our guiding light.
We have tried to introduce as little formalism as possible
compatible with exhibiting quite clearly the planar nature of
subfactor theory. Thus our intention has been to give pride
of place to the pictures. But subfactors are not the only
reason for our procedure. By restricting the scope of the theory
one hopes to get to the most vital examples as quickly as
possible. And we believe that we will see, in some form,
{\it all} the examples in our restricted theory anyway. 
Thus the Fuss-Catalan algebras of [BJ2] (surely among the 
most basic planar algebras, whatever one's definition)
first appeared with our strict axioms. Yet at the same time,
as we show in section 2.5, the {\sc homfly} polynomial, for
which one might have thought oriented strings essential, can
be completely captured within our unoriented framework.

It is unlikely that any other restriction of some more general
operad is as rich as the one we use here. To see why, note that 
in the operadic picture, the role of the identity is played by 
tangles without internal discs -see [Ma]. In our case we get the
whole Temperley-Lieb algebra corresponding to the identity 
whereas any orientation restriction will reduce the size of
this ``identity". The beautiful structure of the Temperley-Lieb
algebra is thus always at our disposal. This leads to the 
following rather telling reason for looking carefully 
at our special planar algebras among more general ones:
if we introduce the generating function for the dimensions
of a planar algebra,  $\sum^\i_{n=0}\dim (V_n)z^n$,
we shall see that if the planar algebra satisfies reflection 
positivity, then this power series has non-zero radius of
convergence. By contrast, if we take the natural oriented planar
algebra structure given by the {\sc homfly} skein, it is a 
result of Ocneanu and Wenzl ([We1],[F+]) that there is a positive
definite Markov trace on the whole algebra, even though the
generating function has zero radius of convergence.

In spite of the previous polemic, it would be foolish to 
neglect the fact that our planar algebra formalism fits
into a more general one. Subfactors can be constructed
with arbitrary orientations by the procedure of [We1],[F+]
and it should be possible to calculate their planar algebras
by planar means.

We end this introduction by discussing three of our motivations
for the introduction of planar algebras as we have defined them.

{\it Motivation 1 }  Kauffman gave his now well-known pictures
for the Temperley-Lieb algebra in [Ka1]. In the mid 1980's he 
asked the author if it was possible to give a pictorial
representation of all elements in the tower of algebras of [J1].
We have only developed the planar algebra formalism for the
sub-tower of relative commutants, as the all-important rotation is
not defined on the whole tower. Otherwise this paper constitutes
an answer to Kauffman's question.

{\it Motivation 2 } One of the most extraordinary developments in 
subfactors was the discovery by Haagerup in [Ha] of a subfactor
of index $(5 + \sqrt{13})/2$, along with the proof that this is 
the smallest index value, greater than $4$, of a finite depth
subfactor. As far as we know there is no way to obtain Haagerup's
``sporadic" subfactor from the conformal field theory/quantum group
methods of [Wa],[We3],[X],[EK]. It is our hope that the planar algebra
context will put Haagerup's subfactor in at least one natural
family, besides yielding tools for its study that are more
general than those of [Ha]. For instance it follows from 
Haagerup's results that the planar algebra of his subfactor is
generated by a single element in $V_4$ (a ``4-box"). The small
dimensionality of the planar algebra forces extremely strong 
conditions on this 4-box. The only two simpler such planar algebras
(with reflection positivity) are those of the $D_6$ subfactor
of index $4\cos^2\pi/10$ and the $\tilde E_7$ subfactor of
index $4$. There are analogous planar algebras generated by
2-boxes and 3-boxes. The simplest 2-box case comes from the
$D_4$ subfactor (index $3$) and the two simplest 3-box cases from
$E_6$ and $\tilde E_6$ (indices $4\cos^2\pi/12$ and $4$).
Thus we believe there are a handful of planar algebras for
each k, generated by a single k-box, satisfying extremely 
strong relations. Common features among these relations
should yield a unified calculus for constructing and manipulating
these planar algebras. In this direction we have classified 
with Bisch in [BJ1] {\it all} planar algebras generated by
a 2-box and tightly restricted in dimension. A result of
D.Thurston shows an analogous result should exist for 3-boxes - see section 2.5.
The 4-box case has yet to be attempted.

In general one would like to understand all systems of relations
on planar algebras that cause the free planar algebra to collapse
to finite dimensions. This is out of sight at the moment. 
Indeed it is know from [BH] that subfactors of index $6$ are 
``wild" in some technical sense, but up to $3+\sqrt{3}$ they
appear to be ``tame". It would be significant to know for what   
index value subfactors first become wild.

{\it Motivation 3 } Since the earliest days of subfactors it
has been known that they can be constructed from certain 
finite data known as a commuting square (see [GHJ]). A theorem 
of Ocneanu (see [JS] or [EK]) reduced the problem of  
calculating the planar algebra component $V_k$ of such a subfactor
to the solution of a finite system of linear equations in finitely
many unknowns. Unfortunately the number of equations 
grows exponentially with $k$ and it is unknown at present whether 
the most simple questions concerning these $V_k$ are solvable
in polynomial time or not. On the other hand the planar algebra
gives interesting invariants of the original combinatorial 
data and it was a desire to exploit this information that led
us to consider planar algebras. First it was noticed that there
is a suggestive planar notation for the linear equations themselves.
Then the invariance of the solution space under the action
of planar tangles was observed. It then became clear that 
one should consider other ways of constructing planar 
algebras from combinatorial data, such as the planar
algebra generated by a tensor in the tensor planar algebra.

These ideas were the original motivation for introducing planar
algebras. We discuss these matters in more detail in section
2.11 which is no doubt the most important part of this work.
The significance of Popa's result on $\lambda -$ lattices became
apparent as the definition evolved. Unfortunately we have
not yet been able to use planar algebras in a convincing
way as a tool in the calculation of the planar algebra for
specific commuting squares.

This paper has been written over a period of several years and
many people have contributed. In particular I would like to
thank Dietmar Bisch, Pierre de la Harpe, Roland Bacher, Sorin Popa,
Dylan Thurston, Bina Bhattacharya, Zeph Landau, Adrian Ocneanu, Gib Bogle and 
Richard Borcherds. Deborah Craig for her 
patience and first-rate typing, and Tsukasa Yashiro for the pictures.

\bs\bs
\ni{\large\bf 1. The Formalism}

\medskip
{\bf Definition 1.1.} If $k$ is a non-negative integer, the
{\it standard} $k${\it -box}, $\b _k$, is
$\{ (x,y)\in\r\mid 0\leq x\leq k\!+\!1, \
0\leq y\leq 1\}$, together with the $2k$ marked points,
$1\!=\!(1,1)$, $2\!=\!(2,1)$, $3=(3,1),\dots ,k=(k,1)$,
$k+1=(k,0)$, \ $k+2=(k-1,0),\dots ,2k=(1,0)$.

\bs{\bf Definition 1.2.}
A {\it planar network} $\cal N$ will be a subset of $\r$ consisting
of the union of a finite set of disjoint images of $\b_k$'s
(with $k$ varying) under smooth orientation-preserving
diffeomorphisms of $\r$, and a finite number of oriented disjoint
curves, smoothly embedded, which may be closed (i.e. isotopic to
circles), but if not their endpoints coincide with marked points of
the boxes. Otherwise the curves are disjoint from the boxes.
All the marked points are endpoints of curves, which meet the boxes
transversally. The orientations of the curves must satisfy the
following two conditions.
\begin{itemize}
\item[a)] A curve meeting a box at an odd marked point
must exit the box at that point.
\item[b)] The connected components of $\r\backslash\cal N$
may be oriented in such a way that the orientation of a curve
coincides with the orientation induced as part of the boundary of a
connected component.
\end{itemize}

{\bf Remark.} Planar networks are of two kinds according to the
orientation of the unbounded region.

\bs
Let $L_i$, $i=1,2,\dots$ be sets and $L=\coprod_i L_i$ be their
disjoint union. $L$ will be called the set of ``labels".

\bs{\bf Definition 1.3.} A {\it labelled planar network}
(on $L$) will be a planar network together with a function from its
$k$-boxes to $L_k$, for all $k$ with $L_k\neq\emptyset$.

\bs If the labelling set consists of asymmetric letters, we may
represent the labelling function diagrammatically by placing the
corresponding letter in its box, with the understanding that the
first marked point is at the top left. This allows us to ignore the
orientations on the edges and the specification of the marked
points. In Fig.~1.4 we give an example of a labelled planar network
with $L_1=\{ P\}$, $L_2=\{ R\}$, $L_3=\{ Q\}$. Here the unbounded
region is positively oriented and, in order to make the conventions
quite clear, we have explicitly oriented the edges and numbered the
marked points of the one 3-boxes labelled $Q$.
\[
        \input{xfig/pic7}
\]
\begin{center}
Figure 1.4
\end{center}

\bs\ni
Note that the same picture as in Fig.~1.4, but with an $R$ upside
down, would be a different labelled planar tangle since the marked
points would be different.
With or without labels, it is only  necessary to say which
distinguished boundary point is first. 

\bs{\bf Remark 1.5.} By shrinking each $k$-box to a point as in
Fig.~1.6 one obtains from a planar network a system of immersed
curves with transversal multiple point singularities.
\[
        \begin{picture}(0,0)%
\epsfig{file=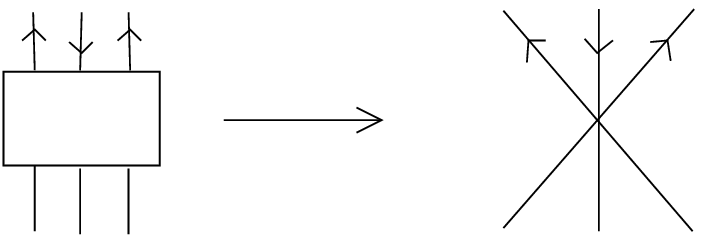}%
\end{picture}%
\setlength{\unitlength}{0.00041700in}%
\begingroup\makeatletter\ifx\SetFigFont\undefined
\def\x#1#2#3#4#5#6#7\relax{\def\x{#1#2#3#4#5#6}}%
\expandafter\x\fmtname xxxxxx\relax \def\y{splain}%
\ifx\x\y   
\gdef\SetFigFont#1#2#3{%
  \ifnum #1<17\tiny\else \ifnum #1<20\small\else
  \ifnum #1<24\normalsize\else \ifnum #1<29\large\else
  \ifnum #1<34\Large\else \ifnum #1<41\LARGE\else
     \huge\fi\fi\fi\fi\fi\fi
  \csname #3\endcsname}%
\else
\gdef\SetFigFont#1#2#3{\begingroup
  \count@#1\relax \ifnum 25<\count@\count@25\fi
  \def\x{\endgroup\@setsize\SetFigFont{#2pt}}%
  \expandafter\x
    \csname \romannumeral\the\count@ pt\expandafter\endcsname
    \csname @\romannumeral\the\count@ pt\endcsname
  \csname #3\endcsname}%
\fi
\fi\endgroup
\begin{picture}(6674,2204)(1179,-4023)
\end{picture}

\]
\begin{center}
Figure 1.6
\end{center}
\ni
Cusps can also be handled by labelled 1-boxes. To reverse the
procedure requires a choice of incoming curve at each multiple
point but we see that our object is similar to that of Arnold in
[A]. In particular, in what follows we will construct a huge
supply of invariants for systems of immersed curves. It remains to
be seen whether these invariants are of interest in singularity
theory, and whether Arnold's invariants may be used to construct
planar algebras with the special properties we shall describe.

\bs{\bf Definition 1.7.} A planar $k$-tangle $T$ (for
$k=0,1,2,\dots$) is the intersection of a planar network $\cal N$
with the standard $k$-box $\b_k$, with the condition that the
boundary of $\b_k$ meets $\cal N$ transversally precisely in the
set of marked points of $\b_k$, which are points on the curves of
$\cal N$ other than endpoints. The orientation induced by $\cal N$
on a neighborhood of (0,0) is required to be positive. A labelled
planar $k$-tangle is defined in the obvious way.

The connected curves in a tangle $T$ will be called the {\it
strings} of $T$. 

The set of smooth isotopy classes of labelled planar $k$-tangles,
with isotopies being the identity on the boundary of $\b_k$, is
denoted $T_k(L)$.

{\bf Note.} $T_0(L)$ is naturally identified with the set of planar
isotopy classes of labelled networks with unbounded region
positively oriented.

\bs{\bf Definition 1.8.} The associative algebra
$\P_k(L)$ over the field $K$ is the vector space having $T_k(L)$ as
basis, with multiplication defined as follows. If $T_1,T_2\in
T_k(L)$, let
$\tilde T_2$ be $T_2$ translated in the negative $y$ direction by
one unit. After isotopy if necessary we may suppose that the union
of the curves in $T_1$ and $\tilde T_2$ define smooth curves.
Remove  $\{(x,0)\mid 0\leq x\leq k+1, \ x\not\in\z\}$ from
$T_1\cup\tilde T_2$ and finally rescale by multiplying the
$y$-coordinates $\frac 12$, then adding $\frac 12$.
The resulting isotopy class of labelled planar $k$-tangles is
$T_1T_2$. See Figure 1.9 for an example.
\[
        \begin{picture}(0,0)%
\epsfig{file=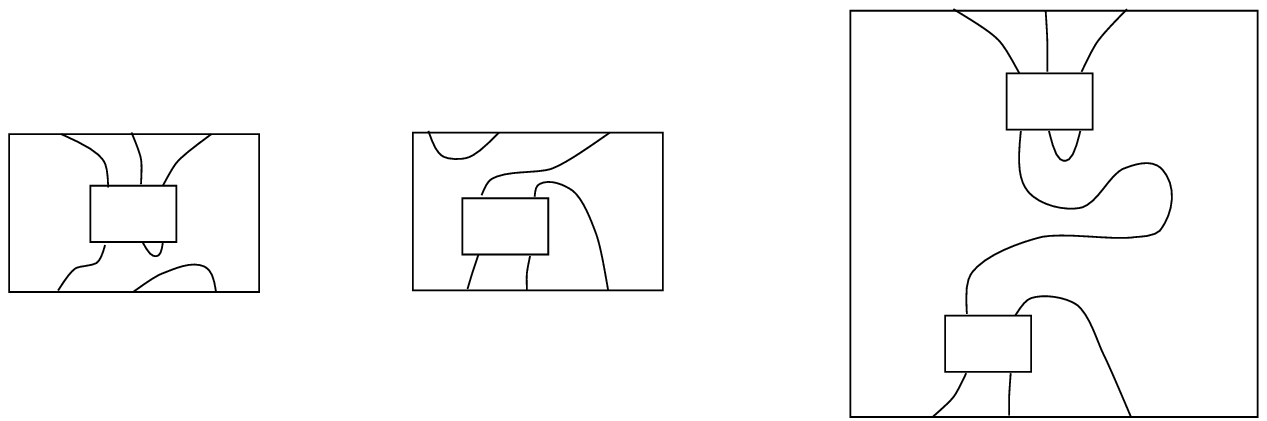}%
\end{picture}%
\setlength{\unitlength}{0.00041700in}%
\begingroup\makeatletter\ifx\SetFigFont\undefined
\def\x#1#2#3#4#5#6#7\relax{\def\x{#1#2#3#4#5#6}}%
\expandafter\x\fmtname xxxxxx\relax \def\y{splain}%
\ifx\x\y   
\gdef\SetFigFont#1#2#3{%
  \ifnum #1<17\tiny\else \ifnum #1<20\small\else
  \ifnum #1<24\normalsize\else \ifnum #1<29\large\else
  \ifnum #1<34\Large\else \ifnum #1<41\LARGE\else
     \huge\fi\fi\fi\fi\fi\fi
  \csname #3\endcsname}%
\else
\gdef\SetFigFont#1#2#3{\begingroup
  \count@#1\relax \ifnum 25<\count@\count@25\fi
  \def\x{\endgroup\@setsize\SetFigFont{#2pt}}%
  \expandafter\x
    \csname \romannumeral\the\count@ pt\expandafter\endcsname
    \csname @\romannumeral\the\count@ pt\endcsname
  \csname #3\endcsname}%
\fi
\fi\endgroup
\begin{picture}(13006,3952)(511,-5178)
\put(511,-3356){\makebox(0,0)[lb]{\smash{\SetFigFont{12}{14.4}{rm}$T_1=$}}}
\put(4036,-3431){\makebox(0,0)[lb]{\smash{\SetFigFont{12}{14.4}{rm},}}}
\put(4426,-3401){\makebox(0,0)[lb]{\smash{\SetFigFont{12}{14.4}{rm}$T_2=$}}}
\put(7966,-3446){\makebox(0,0)[lb]{\smash{\SetFigFont{12}{14.4}{rm},}}}
\put(8341,-3476){\makebox(0,0)[lb]{\smash{\SetFigFont{12}{14.4}{rm}$T_1T_2=$}}}
\put(2531,-3266){\makebox(0,0)[lb]{\smash{\SetFigFont{12}{14.4}{rm}$Q$}}}
\put(11387,-2187){\makebox(0,0)[lb]{\smash{\SetFigFont{12}{14.4}{rm}$Q$}}}
\put(10826,-4527){\makebox(0,0)[lb]{\smash{\SetFigFont{12}{14.4}{rm}$R$}}}
\put(6176,-3401){\makebox(0,0)[lb]{\smash{\SetFigFont{12}{14.4}{rm}$R$}}}
\end{picture}

\]
\begin{center}
        Figure 1.9
\end{center}
\bs{\bf Remark.} The algebra $\P_k(L)$ has an obvious unit and
embeds unitally in $\P_{k+1}(L)$ by adding the line
$\{(k+1,t)\mid 0\leq t\leq 1\}$ to an element of $\P_k(L)$.
Since isotopies are the identity on the boundary this gives an
injection from the basis of $\P_{k}(L)$ to that of
$\P_{k+1}(L)$.

If there is no source of confusion we will suppress the explicit
dependence on $L$.

\bs{\bf Definition 1.10.} The {\it universal planar algebra}
$\P(L)$ on $L$ is the filtered algebra given by the union of all the
$\P_k$'s ($k=0,1,2,\dots$) with $\P_k$ included in $\P_{k+1}$ as in
the preceding remark.

\bs
A planar algebra will be basically a filtered quotient of $\P(L)$
for some $L$, but in order to reflect the planar structure we need
to impose a condition of {\it annular}
invariance.

\bs{\bf Definition 1.11.} The $j\!-\!k$ {\it annulus} $A_{j,k}$ will
be the complement of the interior of $\b_j$ in
$(j+2)\b_k-(\frac 12,\frac 12)$. So there are $2j$ marked points
on the inner boundary of ${\cal A}_{j,k}$ and $2k$ marked points on
the outer one. An {\it annular} $j-k$ {\it tangle} is the
intersection of a planar network $\cal N$ with $A_{j,k}$ such that
the boundary of $A_{j,k}$ meets $\cal N$ transversally precisely
in the set of marked points of $A_{j,k}$, which are
points on the curves of $\cal N$ other than
endpoints.  The orientation induced by $\cal N$ in
neighborhoods of
$(-\frac 12,-\frac 12)$ and $(0,0)$ are
 required to be positive. Labeling is as
usual.

\bs{\bf Warning.} The diagram in Fig.~1.11(a) is {\it not} an
annular 2--1 tangle, whereas the diagram in Fig.~1.11 (b) is.
\[
        \begin{picture}(0,0)%
\epsfig{file=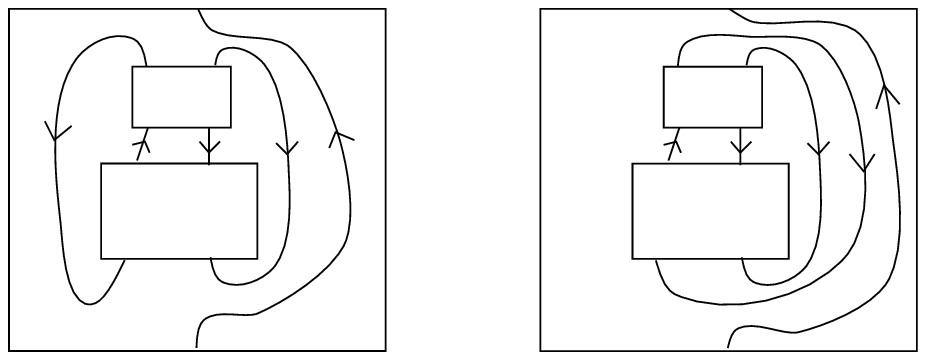}%
\end{picture}%
\setlength{\unitlength}{0.00041700in}%
\begingroup\makeatletter\ifx\SetFigFont\undefined
\def\x#1#2#3#4#5#6#7\relax{\def\x{#1#2#3#4#5#6}}%
\expandafter\x\fmtname xxxxxx\relax \def\y{splain}%
\ifx\x\y   
\gdef\SetFigFont#1#2#3{%
  \ifnum #1<17\tiny\else \ifnum #1<20\small\else
  \ifnum #1<24\normalsize\else \ifnum #1<29\large\else
  \ifnum #1<34\Large\else \ifnum #1<41\LARGE\else
     \huge\fi\fi\fi\fi\fi\fi
  \csname #3\endcsname}%
\else
\gdef\SetFigFont#1#2#3{\begingroup
  \count@#1\relax \ifnum 25<\count@\count@25\fi
  \def\x{\endgroup\@setsize\SetFigFont{#2pt}}%
  \expandafter\x
    \csname \romannumeral\the\count@ pt\expandafter\endcsname
    \csname @\romannumeral\the\count@ pt\endcsname
  \csname #3\endcsname}%
\fi
\fi\endgroup
\begin{picture}(9368,3330)(871,-5479)
\put(871,-3821){\makebox(0,0)[lb]{\smash{\SetFigFont{12}{14.4}{rm}$(a)$}}}
\put(3001,-3041){\makebox(0,0)[lb]{\smash{\SetFigFont{12}{14.4}{rm}$Q$}}}
\put(8131,-3026){\makebox(0,0)[lb]{\smash{\SetFigFont{12}{14.4}{rm}$Q$}}}
\put(5956,-3806){\makebox(0,0)[lb]{\smash{\SetFigFont{12}{14.4}{rm}$(b)$}}}
\end{picture}

\]
\begin{center}
Figure 1.11
\end{center}

The set of all isotopy classes (isotopies being the identity on the
boundary) of labelled annular $j-k$ tangles, ${\cal A}(L)=\bigcup_{j,k}
{\cal A}_{j,k}(L)$ forms a category whose objects are the sets of
$2j$-marked points of $\b_j$. To compose $A_1\in{\cal A}_{j,k}$ and
$A_2\in{\cal A}_{k,\ell}$, rescale and move $A_1$ so that its outside
boundary coincides with the inside boundary of $A_2$, and the $2k$ 
boundary points match up. Join the strings of $A_1$ to those of $A_2$
at their common boundary and smooth them. Remove that part of the common
boundary that is not strings. Finally rescale the whole annulus so that it
is the standard one. 
The result will depend only on the
isotopy classes of $A_1$ and $A_2$ and defines an element
$A_2A_1$ in ${\cal A}(L)$.

Similarly, an $A\in{\cal A}_{j,k}(L)$ determines a map
$\pi_A:T_j(L)\to T_k(L)$ by surrounding $T\in T_j(L)$ with $A$ and
rescaling. Obviously $\pi_A\pi_B=\pi_{AB}$, and the action of
${\cal A}(L)$ extends to $\P(L)$ by linearity.

\bs
{\bf Definition 1.12} A {\it general planar algebra} will be a
filtered algebra $P=\cup_k P_k$, together with a surjective
homomorphism of filtered algebras, $\Phi:\P(L)\to P$, for some label
set $L$, $\Phi(\P_k)=P_k$, with ker $\Phi$ invariant under ${\cal A}(L)$ 
in the sense that, if $\Phi(x)=0$ for $x\in \P_j$, and
$A\in{\cal A}_{j,k}$ then $\Phi(\pi_A(x))=0$, we say $\Phi$
presents $P$ on $L$.

\bs{\bf Note.} Definition 1.12 ensures that ${\cal A}(L)$ acts on $P$
via $\pi_A(\Phi(x))\eqdef{\rm def}\Phi(\pi_A(x))$.
In particular ${\cal A}(\emptyset)$ ($\emptyset=$ emptyset) acts on
any planar algebra.

\bs
The next results show that this action
extends multilinearly to planar surfaces with several boundary
components.

If $T$ is a planar $k$-tangle (unlabelled), number its boxes
$b_1,b_2,\dots ,b_n$. Then given labelled tangles $T_1,\dots ,T_n$
with $T_i$ having the same number of boundary points as $b_i$, we
may form a labelled planar $k$-tangle $\pi_T(T_1,T_2,\dots ,T_n)$ by
filling each $b_i$ with $T_i$ --- by definition $b_i$ is the image
under a planar diffeomorphism $\th$ of $\b_j$ (for some $j$), and
$T_i$ is in $\b_j$, so replace $b_i$ with $\th(T_i)$ and remove the
boundary (apart from marked points, smoothing the curves at the
marked points). None of this depends on isotopy so the isotopy
class of $T$ defines a multilinear map
$\pi_T: \P_{j_1}\x \P_{j_2}\x\dots\x \P_{j_n}\to \P_k$.
Though easy, the following result is fundamental and its conclusion
is the definition given in the introduction of planar algebras
based on the operad defined by unlabeled planar tangles.

\bigskip
\ni{\bf Proposition 1.13} If $P$ is a general planar algebra
presented on $L$, by $\Phi$, $\pi_T$ defines a multilinear map
$P_{j_1}\x P_{j_2}\x\dots\x P_{j_n}\to P_k$.

{\sc Proof}. It suffices to show that, if all the $T_i$'s but
one, say $i_0$, are fixed in $\P_{j_i}$, then the linear map
$\a:\P_{j_{i_0}}\to \P_k$, induced by $\pi_T$, is zero on ker
$\Phi$. By multilinearity the $T_i$'s can be supposed to be isotopy classes
of labelled tangles. So fill all the boxes other than the $i_0$'th
box with the
$T_i$.  Then we may isotope the resulting picture so that
$b_{i_0}$ is the inside box of a $j_{i_0}-k$ annulus. The map $\a$
is then the map $\pi_A$ for some annular tangle $A$ so
ker $\Phi\subseteq {\op{ker}} \ \pi_A$ by Definition 1.12. \qed

\bigskip
\ni{\bf Proposition 1.14} Let $P$ be a general planar
algebra presented on $L$ by $\Phi$. For each $k$ let $S_k$ be a set
and
$\a:S_k\to P_k$ be a function. Put $S=\coprod_k S_k$. Then there is
a unique filtered algebra homomorphism
$\Theta_S:\P(S)\to P$ with ${\op{ker}} \ \Theta_S$ invariant under
${\cal A}(S)$, intertwining the ${\cal A}(\emptyset)$ actions and such that
$\Theta_S(\tangle{R})=\a(R)$ for $R\in S$.

{\sc Proof}. Let $T$ be a tangle in $P(S)$ with boxes
$b_1,\dots ,b_n$ and let $f(b_i)$ be the label of $b_i$.
We set $\Theta_S(T)=\pi_T(\a(f(b_1)),\a(f(b_2)),\dots ,\a(f(b_n))$
with $\pi_T$ as in 1.13. For the homomorphism property, observe
that $\pi_{T_1T_2}$ and $\pi_{T_1}\dt\pi_{T_2}$ are both
multilinear maps agreeing on a basis. For the annular invariance of
ker $\Theta_S$, note that $\Theta_S$ factors through $\P(L)$, say
$\Theta_S=\Phi\circ \theta$, so that
$\Theta_S(x)=0\Longleftrightarrow
\th(x)\in{\op{ker}} \ \Phi$. Moreover, if $A\in{\cal A}(S)$,
$\theta\!\circ\!\pi_A$ is a linear combination of $\pi_{A'}$'s for
$A'$ in ${\cal A}(L)$. Hence $\Phi(\theta(\pi_A(x))=0$ if
$\theta(x)\in {\op{ker}} \ \Phi$.

Finally we must show that $\Theta_S$ is unique. Suppose we are
given a tangle $T\in \P(S)$. Then we may isotope $T$ so that all its
boxes occur in a vertical stack, as in Figure~1.15.
\[
        \begin{picture}(0,0)%
\epsfig{file=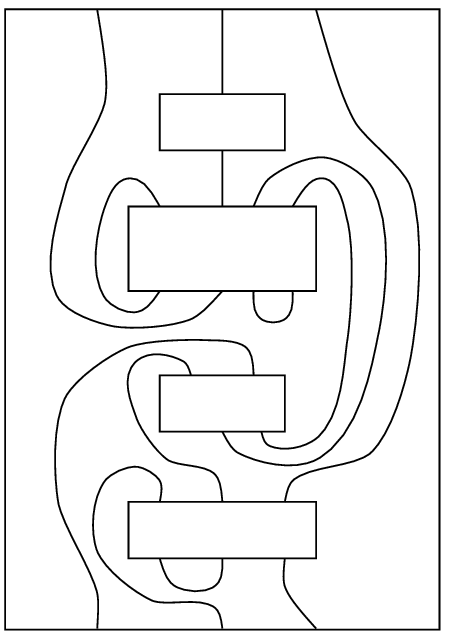}%
\end{picture}%
\setlength{\unitlength}{0.00041700in}%
\begingroup\makeatletter\ifx\SetFigFont\undefined
\def\x#1#2#3#4#5#6#7\relax{\def\x{#1#2#3#4#5#6}}%
\expandafter\x\fmtname xxxxxx\relax \def\y{splain}%
\ifx\x\y   
\gdef\SetFigFont#1#2#3{%
  \ifnum #1<17\tiny\else \ifnum #1<20\small\else
  \ifnum #1<24\normalsize\else \ifnum #1<29\large\else
  \ifnum #1<34\Large\else \ifnum #1<41\LARGE\else
     \huge\fi\fi\fi\fi\fi\fi
  \csname #3\endcsname}%
\else
\gdef\SetFigFont#1#2#3{\begingroup
  \count@#1\relax \ifnum 25<\count@\count@25\fi
  \def\x{\endgroup\@setsize\SetFigFont{#2pt}}%
  \expandafter\x
    \csname \romannumeral\the\count@ pt\expandafter\endcsname
    \csname @\romannumeral\the\count@ pt\endcsname
  \csname #3\endcsname}%
\fi
\fi\endgroup
\begin{picture}(4215,5999)(1794,-6626)
\put(3601,-5761){\makebox(0,0)[lb]{\smash{\SetFigFont{12}{14.4}{rm}$R_1$}}}
\put(3601,-4561){\makebox(0,0)[lb]{\smash{\SetFigFont{12}{14.4}{rm}$R_2$}}}
\put(3526,-2986){\makebox(0,0)[lb]{\smash{\SetFigFont{12}{14.4}{rm}$R_3$}}}
\put(3601,-1786){\makebox(0,0)[lb]{\smash{\SetFigFont{12}{14.4}{rm}$R_4$}}}
\end{picture}

\]
\begin{center}
Figure 1.15
\end{center}

\bs In between each box cut horizontally along a level for which
there are no critical points for the height function along the
curves. Then the tangle becomes a product of single labelled boxes
surrounded by ${\cal A}(\emptyset)$ elements. By introducing kinks as
necessary, as depicted in Figure~1.16, all the surrounded boxes may
be taken in $P_k(S)$ for some large fixed $k$.
\[
        \begin{picture}(0,0)%
\epsfig{file=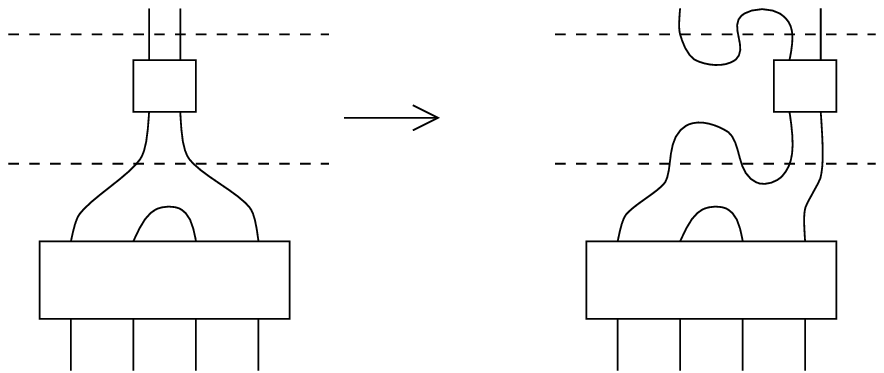}%
\end{picture}%
\setlength{\unitlength}{0.00041700in}%
\begingroup\makeatletter\ifx\SetFigFont\undefined
\def\x#1#2#3#4#5#6#7\relax{\def\x{#1#2#3#4#5#6}}%
\expandafter\x\fmtname xxxxxx\relax \def\y{splain}%
\ifx\x\y   
\gdef\SetFigFont#1#2#3{%
  \ifnum #1<17\tiny\else \ifnum #1<20\small\else
  \ifnum #1<24\normalsize\else \ifnum #1<29\large\else
  \ifnum #1<34\Large\else \ifnum #1<41\LARGE\else
     \huge\fi\fi\fi\fi\fi\fi
  \csname #3\endcsname}%
\else
\gdef\SetFigFont#1#2#3{\begingroup
  \count@#1\relax \ifnum 25<\count@\count@25\fi
  \def\x{\endgroup\@setsize\SetFigFont{#2pt}}%
  \expandafter\x
    \csname \romannumeral\the\count@ pt\expandafter\endcsname
    \csname @\romannumeral\the\count@ pt\endcsname
  \csname #3\endcsname}%
\fi
\fi\endgroup
\begin{picture}(8948,3521)(1351,-4460)
\put(1351,-1186){\makebox(0,0)[lb]{\smash{\SetFigFont{12}{14.4}{rm}cut}}}
\put(1351,-2461){\makebox(0,0)[lb]{\smash{\SetFigFont{12}{14.4}{rm}cut}}}
\put(6601,-1261){\makebox(0,0)[lb]{\smash{\SetFigFont{12}{14.4}{rm}cut}}}
\put(6601,-2461){\makebox(0,0)[lb]{\smash{\SetFigFont{12}{14.4}{rm}cut}}}
\end{picture}

\]
\begin{center}
Figure 1.16
\end{center}

\bs\ni Since $\Theta_S$ is required to intertwine the $\cal
A(\emptyset)$ action and is an algebra homomorphism, it is determined on
all the surrounded boxes by its value on $\{\tangle{R} : R\in
S\}$, and their products. The beginning and end of $T$ may involve a
change in the value of $k$, but they are represented by an element
of ${\cal A}(\emptyset)$ applied to the product of the surrounded boxes.
So $\Theta_S$ is completely determined on $T$. \qed

\bs{\bf Definition 1.17.} Let $P^1,P^2$  be general planar algebras
presented by $\Phi_1,\Phi_2$ on $L^1,L^2$ respectively.

If $\a: L^1_k\to P^2_k$, as in 1.14, is such that \
ker $\Theta_\a\supseteq {\op{ker}} \ \Phi_1$, then the resulting
homomorphism of filtered algebras $\Gamma_\a: P^1\to P^2$ is called
a {\it planar algebra homomorphism.} A {\it planar subalgebra}
of a general planar algebra is the image of a planar algebra
homomorphism. A planar algebra homomorphism that is bijective is
called a {\it planar algebra isomorphism.} Two presentations
$\Phi_1$ and $\Phi_2$ of a planar algebra will be considered to
define the same planar algebra structure if the identity map is a
planar algebra homomorphism.

\bs{\bf Remarks.} (i) It is obvious that
planar algebra homomorphisms intertwine the ${\cal A}(\emptyset)$ actions.

(ii) By 1.14, any presentation of a general planar algebra $P$ can
be altered to one whose labelling set is the whole algebra itself,
defining the same planar algebra structure and such that
$\Phi(\tangle{R})=R$ for all $R\in P$.  Thus
there is a canonical, if somewhat unexciting, labelling set.
We will abuse notation by using the same letter $\Phi$ for the
extension of a labelling set to all of $P$. Two presentations
defining the same planar algebra structure will define the same
extensions to all of $P$ as labelling set.

\bigskip
\ni{\bf Proposition 1.18} Let $P$ be a  general planar algebra,
and let $C_n\subseteq P_n$ be unital subalgebras invariant under
${\cal A}(\emptyset)$ (i.e., $\pi_A(C_j)\subseteq C_k$ for $A\in{\cal
A}_{j,k}(\emptyset)$). Then $C=\cup C_n$ is a planar subalgebra of $P$.

\bigskip
{\sc Proof.} As a labelling set for $C$ we choose $C$ itself.
We have to show that $\Theta_C(\cal P(C))\subseteq C$. But this follows
immediately from the argument for the uniqueness of $\Theta_S$ in
1.14. (Note that $C_n\subseteq C_{n+1}$ as subalgebras of
$P_{n+1}$ is automatic from invariance under ${\cal A}(\emptyset)$.) \qed

The definition of isomorphism
was asymmetric. The next result shows that the notion {\bf is}
symmetric.

\bs
\ni{\bf Proposition 1.19} If $\Gamma_{\a}:P^1\to P^2$ is an
isomorphism of planar algebras, so is $(\Gamma_{\a})^{-1}$.

\bs
{\sc Proof.} Define \ $\a^{-1}:L^2_k\to P^1_k$ \ by \
$\a^{-1}(R)=(\Gamma_{\a})^{-1}(\Phi_2(\tangle{R})$.
Then $\G_\a\circ\Theta_{\a^{-1}}$ is a filtered algebra
homomorphism intertwining the ${\cal A}(\emptyset)$ actions so it equals
$\Phi_1$ by 1.14. Thus ker $\Phi_2\subseteq {\op{ker}} \
\Theta_{\a^{-1}}$ \ and \ $\G_{\a^{-1}}=(\G_\a)^{-1}$. \qed

\bs The definitions of planar algebra homomorphisms, etc., as above
are a little clumsy.
The meaning of the following result is that this operadic
definition of the introduction would give the same notion as the one we have defined.

\bigskip
\ni{\bf Proposition 1.20} If $P_i,\Phi_i,L_i$ for $i=1,2$ are as
in Definition {\rm{1.17}}, then linear maps $\G: P^1_k\to P^2_k$
define a planar algebra homomorphism iff
$$
\pi_T(\G(x_1),\G(x_2),\dots ,\G(x_n)) =\G(\pi_T(x_1,x_2,\dots ,x_n))
$$
for every unlabelled tangle $T$ as in {\rm{1.13}}.

\bigskip
{\sc Proof.} Given  $\G$, define $\a: L_1\to P^2$ by
$\a(R)=\G(\Phi_1(\tangle{R} ))$.
Then $\Theta_\a=\G\circ\Phi_1$ by the uniqueness criterion of 1.14
(by choosing $T$ appropriately it is clear that $\G$ is a
homomorphism of filtered algebras intertwining $\cal
A(\emptyset)$-actions). On the other hand, a planar algebra isomorphism
provides linear maps $\G$ which satisfy the intertwining condition
with $\pi_T$. \qed

\bs{\bf Definition 1.21.} For each $j,k=0,1,2,\dots$ with $j\leq
k$, $\P_{j,k}(L)$ will be the subalgebra of $\P_k(L)$ spanned by
tangles for which all marked points are connected by vertical
straight lines  except those having $x$ coordinates $j+1$ through
$k$. (Thus $\P_{0,k}=\P_k$.) If $B$ is a general planar algebra, put
$P_{j,k}=\Phi(\P_{j,k})$ for some, hence any, presenting map $\Phi$.

\bs{\bf Definition 1.22.} A {\it planar algebra} will be a
general planar algebra $P$ with
 \ dim $P_0=1={\op{dim}} \ P_{1,1}$ and
$\Phi(\caright\ ),\Phi(\caleft\ )$ both nonzero.

\bs A  planar algebra $P$, with presenting map
$\Phi: \P(L)\to P$, defines a planar isotopy invariant of
labelled planar networks, ${\cal N}\mapsto Z_{\Phi}(\cal N)$ by
$Z_{\Phi}(\cal N){\op{id}}=\Phi (\blankbox{\cal N})\in P_0$ if
the unbounded region of ${\cal N}$ is positively oriented 
(and ${\cal N}$ is moved inside $\b_0$ by an isotopy), and 
$Z_{\Phi}({\cal N}){\op{id}}=\Phi(\tanglearrow{\cal N}\  )\in P_{1,1}$
in the other case ($\cal N$ has been isotoped into the right
half of $\b_1$). The invariant $Z$ is called the {\it partition
function.} It is multiplicative in the following sense.

\bs\ni{\bf Proposition 1.23} Let $P$  be a planar algebra with
partition function $Z$. If $T$ is a labelled tangle containing
a planar network $\cal N$ as a connected component, then
$$
\Phi(T)=Z({\cal N})\Phi(T')
$$
where $T'$ is the tangle $T$ from which $\cal N$ has been
removed.

\bs
{\sc Proof.} If we surround $\cal N$ by a 0-box (after isotopy
if necessary) we see that $T$ is just $\blankbox{\cal N}$ to which
a $0-k$ annular tangle has been applied. But
$\Phi( \blankbox{\cal N})=Z({\cal N})\Phi(\blankbox\ )$,
so by annular invariance, $\Phi(T)=Z({\cal N})T$ \qed

\bs
A planar algebra has two scalar parameters,
$\d_1=Z(\caleft\ )$ and $\d_2=Z(\caright\ )$
which we have supposed to be non-zero.

\bs We present two useful procedures to construct  planar algebras.
The first is from an invariant and is analogous to the GNS method
in operator algebras.

Let $Z'$ be a planar isotopy invariant of labelled planar networks
for some labelling set $L$. Extend $Z'$ to $\P_0(L)$ by linearity.
Assume $Z'$ is multiplicative on connected components and that
$Z'(\caright\ )\neq 0$, $Z'(\caleft\ )\neq 0$.
For each $k$ let ${\cal J}_k=\{x\in \P_k(L)\mid Z'(A(T))\!=\!0
\ \forall A\in{\cal A}_{k,0}\}$. Note that $Z'$ (empty network) =1.

\bs
\ni{\bf Proposition 1.24} {\rm{(i)}} ${\cal J}_k$ is a 2-sided
ideal of $\P_k(L)$ and ${\cal J}_{k+1}\cap  \P_k(L)={\cal J}_k$.

{\rm{(ii)}} Let $P_k=\P_k(L)/ {\cal J}_k$ and let $\Phi$ be the
quotient map. Then $P=\cup P_k$ becomes a planar algebra presented
by $\Phi$ with partition function $Z_{\Phi}=Z'$.

{\rm{(iii)}} If $x\in P_k$ then $x=0$ iff $Z_{\Phi}(A(x))\!=\!0
\quad\forall \ A\in{\cal A}_{k,0}$.

\bs
{\sc Proof}. (i) If $T_1$ and $T_2$ are tangles in $\P_k(L)$, the
map $x\mapsto T_1x T_2$ is given by an element $T$ of
${\cal A}_{k,k}(L)$, and if $A\in{\cal A}_{k,0}$ then
$Z'(A(Tx))=Z'(AT)(x))=0$ if $x\in{\cal J}_k$. Hence ${\cal J}_k$ is an
ideal.

It is obvious that ${\cal J}_k\subset{\cal J}_{k+1}\cap{\P}_k(L)$.
So suppose  $x\in{\cal J}_{k+1}\cap \P_k$. Then for some $y\in \P_k$,
$x=\tangleudarrow{y}$, the orientation of the
last straight line depending on the parity of $k$. We want to show
that
$y\in{\cal J}_k$. Take an $A\in{\cal A}_{k,0}$ and form the element
$\tilde A$ in ${\cal A}_{k+1,0}$ which joins the rightmost two
points, inside the annulus, close to the inner boundary.
Then $\tilde A(x)$ will be $A(y)$ with a circle inserted close to
the right extremity of  $\tangle{y}$. So by multiplicativity,
$Z'(\tilde A(x))=Z'(\caleft\ )Z'(A(y))$.
Since $Z'(\caleft\ )\neq 0$, $Z'(A(y))=0$ and
$y\in{\cal J}_k$.

(ii) By (i) we have a natural inclusion of $P_k$ in $P_{k+1}$.
Invariance of the ${\cal J}_k$'s under ${\cal A}$ is immediate.
To show that \ dim $P_0=1={\op{dim}} \ P_{1,1}$, define maps 
$U:\P_0\to K$ ($K=$ the field) and $V:\P_{1,1}\to K$ by linear extensions of
$U(\blankbox{\cal N})=Z'(\cal N)$ and 
$V(\boxarrow{{\cal N}\ \ {\cal M}}\ )= 
Z'({\cal N})Z'({\cal M})$. Observe that $U({\cal J}_0)=0$ and if
${\cal N}_i,{\cal M}_i,\lambda_i (\in K)$ satisfy
$Z'(\sum_i\lambda_i A( \boxarrow{{\cal N}_i\ \ {\cal M}_i}\ )=0$ for
all
$A\in{\cal A}_{1,0}$, then by multiplicativity,
$Z'(\caright\ )(\sum\lambda_i Z({\cal N}_i)Z({\cal M}_i))=0$,
so that $U$ and $V$ define maps from $P_0$ and $P_{1,1}$ to $K$,
respectively. In particular both $U$ and $V$ are surjective since
$U(\blankbox{})=1$, $V(\blankbox{\uparrow})=1$. We need only show
injectivity. So take a linear combination  $\sum\lambda_i\blankbox{{\cal N}_i}$ 
with $\sum\lambda_i Z'({\cal N}_i)=0$. Then if
$A\in{\cal A}_{0,0}$, $Z'(\sum\lambda_i\dt A (\blankbox{{\cal N}_i})
 )\!=\!0$ by multiplicativity so $\sum\lambda_i\dt \blankbox{{\cal N}_i}
\in\cal J_0$. Similarly for 
$\sum_i\lambda_i A(\boxarrow{{\cal N}_i\ \ {\cal M}_i}\ ) \in P_{1,1}$.

Thus \ dim $P_0=1={\op{dim}} \ P_{1,1}$ and by construction, $Z=Z'$.

(iii) This is the definition of ${\cal J}_k$ (and $Z_{\Phi}=Z'$). \qed

\bs{\bf Remark.} If one tried to make the construction of 1.24 for
an invariant that was not multiplicative, one
would rapidly conclude that the resulting algebras all have
dimension zero.

\bs{\bf Definition 1.25.} A planar algebra satisfying condition
(iii) of 1.24 will be called non-degenerate.

\bs The second construction procedure is by generators and
relations. Given a label set $L$ and a subset $R\subseteq \P(L)$,
let
$\J_j(R)$ be the linear span of
$\dsize{\bigcup_{T\in R \atop T\in \P_k(L)}}{\cal A}_{k,j}(L)(T)$. 
It is immediate that
$\J_{j+1}(R)\cap \P_j(L)=\J_j(R)$ (just apply an element of
${\cal A}(\emptyset)$ to kill off the last string), and $\J_j(R)$ is
invariant under ${\cal A}(L)$ by construction.

\bs{\bf Definition 1.26.} With notation as above, set
$P_n(L,R)=\dsize{\frac{\P_n(L)}{\J_n(R)}}$. Then
$P(L,R)=\cup_n P_n(L,R)$ will be called the planar algebra with
generators $L$ and relations $R$.

This method of constructing planar algebras suffers the same
drawbacks as constructing groups by generators and relations.
It is not clear how big $\J_n(R)$ is inside $\P_n(L)$. It is a very
interesting problem to find relation sets $R$ for which
$0 < {\op{dim}} \ P_n(L,R) < \infty$ for each $n$. Knot theory
provides some examples as we shall see.

\bs{\bf Definition 1.27.} A  planar algebra is called {\it spherical} 
if its partition function $Z$ is an invariant of
networks on the two-sphere $S^2$ (obtained from $\mathbb R^2$ by
adding a point at infinity).

\bs The definition of non-degeneracy of a planar algebra involves
all ways of closing a tangle. For a spherical algebra these
closures can be arranged in a more familiar way as follows.

\bs{\bf Definition 1.28.} Let $P$ be a  planar algebra
with partition function $Z$.\linebreak Define two traces tr${}_L$
and tr${}_R$ on $P_k$ by
\[
 {\op{tr}}_L(\tangle{R})=Z(\boxlloop{R})\hskip 15pt  {\rm and}\hskip15pt 
{\op{tr}}_R(\tangle{ R})=Z(\boxrloop{R}).
\]

\bs{\bf Note.} For a spherical planar algebra $P$, $\d_1\!=\!\d_2$
and we shall use $\d$ for this quantity. Similarly
Tr${}_L={\op{Tr}}_R$ and we shall use Tr. If we define
tr$(x)=\frac{1}{\d^n} {\op{Tr}}(x)$ for $x\in P_n$ then tr is
compatible with the inclusions $P_n\subseteq P_{n+1}$
(and tr$(1)\!=\!1$), so defines a trace on $P$ itself.

\bs\ni
{\bf Proposition 1.29} A spherical planar algebra is
nondegenerate iff ${\op{Tr}}$ defines a nondegenerate bilinear form
on $P_k$ for each $k$.

\bs
{\sc Proof.} $(\Leftarrow)$ The picture defining Tr is the
application of a particular element $A$ of ${\cal A}_{k,0}$ to
$x\in P_k$.

$(\Rightarrow)$ It suffices to show that, for any $A\in{\cal A}_{k,0}(L)$
there is a $y\in P_k$ such that Tr$(xy)=Z(A(x))$.
By  spherical invariance one may arrange $A(x)$ so that the box containing
$x$ has no strings to its left. The part of $A(x)$ outside
that box can then be isotoped into a $k$-box which contains the element $y$.
\qed

\bs{\bf Remark 1.30.} One of the significant consequences of
1.29 is that, for nondegenerate $P$, if one can find a finite
set of tangles which linearly span $P_k$, the calculation of
dim $P_k$ is reduced to the {\it finite} problem of calculating
the rank of the bilinear form defined on $P_k$ by Tr. Of course
this may not be easy!

\bs
\ni{\bf Corollary 1.31} A
nondegenerate planar algebra is semisimple.

\bs\ni{\bf Positivity}

 For the rest of this section suppose the field is $\Bbb R$ or
$\Bbb C$.

Suppose we are given an involution $R\to R^*$ on the set of labels
$L$. Then $\P(L)$ becomes a $*$-algebra as follows. If $T$ is a
tangle in $T_k(L)$ we reflect the underlying unlabelled tangle in
the line $y=\frac 12$ and reverse all the orientations of the
strings. The first boundary point for a box in the reflected 
unlabelled tangle is the one that was the last boundary point for that box
in the original unlabelled tangle.The new tangle $T^*$ is then obtained  by assigning the
label $R^*$ to a box that was labelled $R$. This operation is
extended sesquilinearly to all of $\P_k(L)$. If $\Phi$ presents a
general planar algebra, $*$ preserves ker $\Phi$ and defines a
$*$-algebra structure on $\Phi(\P_k(L))$. The operation $*$ on
$T_0(L)$ also gives a well-defined map on isotopy classes of planar
networks and we say an invariant $Z$ is sesquilinear if
$Z(\cal N^*)=\overline{Z(\cal N)}$.

\bs{\bf Definition 1.32.} A $*$-algebra $P$ is called a (general)
planar $*$-algebra if it is presented by a $\Phi$ on  $\P(L)$, $L$
with involution $*$, such that $\Phi$ is a $*$-homomorphism.

Note that if $P$ is planar, $Z$ is sesquilinear. Moreover if $Z$
is a sesquilinear multiplicative invariant, the construction of
1.24 yields a planar $*$-algebra.
The partition function on a planar algebra will be called {\it
positive} if tr${}_L(x^*x) \geq 0$ for $x\in P_k$, $k$ arbitrary.

\bs
\ni{\bf Proposition 1.33} Let $P$ be a  planar
$*$-algebra with positive partition function $Z$. The following
are equivalent:
\begin{itemize}
\item[(i)] $P$ is non-degenerate {\rm{(Def.~1.24).}}

\item[(ii)] ${\op{tr}}_R(x^*x) > 0$ \ for \ $x\neq 0$.

\item[(iii)] ${\op{tr}}_L(x^*x) > 0$ \ for \ $x\neq 0$.
\end{itemize}

{\sc Proof}. For (ii)$\Leftrightarrow$(iii), argue first that
$\d_1={\op{tr}}_R(  \blankbox{\uparrow}) > 0$ and
$\d_2=Z(\caleft\ )=\frac{1}{\d_1}
{\op{tr}}_R( \blankbox{\ \uparrow\downarrow} \ ) > 0$
and then define antiautomorphisms $j$ of $P_{2n}$ by
$j(\boxas{R})=\boxasR$, so that
$j(x^*)=j(x)^*$ and ${\op{tr}}_L(j(x))={\op{tr}}_R(x)$.

(ii)$\Rightarrow$(i) is immediate since
${\op{tr}}_R(x^*x)=\sum\lambda_i A_i(x^*)$ where
$A_i\in{\cal A}_{k,0}$ is the annular tangle of figure $1.34$,
\[
        \begin{picture}(0,0)%
\epsfig{file=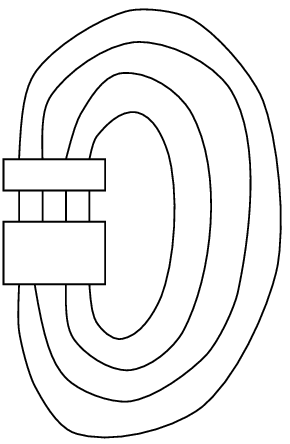}%
\end{picture}%
\setlength{\unitlength}{0.00041700in}%
\begingroup\makeatletter\ifx\SetFigFont\undefined
\def\x#1#2#3#4#5#6#7\relax{\def\x{#1#2#3#4#5#6}}%
\expandafter\x\fmtname xxxxxx\relax \def\y{splain}%
\ifx\x\y   
\gdef\SetFigFont#1#2#3{%
  \ifnum #1<17\tiny\else \ifnum #1<20\small\else
  \ifnum #1<24\normalsize\else \ifnum #1<29\large\else
  \ifnum #1<34\Large\else \ifnum #1<41\LARGE\else
     \huge\fi\fi\fi\fi\fi\fi
  \csname #3\endcsname}%
\else
\gdef\SetFigFont#1#2#3{\begingroup
  \count@#1\relax \ifnum 25<\count@\count@25\fi
  \def\x{\endgroup\@setsize\SetFigFont{#2pt}}%
  \expandafter\x
    \csname \romannumeral\the\count@ pt\expandafter\endcsname
    \csname @\romannumeral\the\count@ pt\endcsname
  \csname #3\endcsname}%
\fi
\fi\endgroup
\begin{picture}(2710,4161)(2079,-5764)
\put(2476,-3961){\makebox(0,0)[lb]{\smash{\SetFigFont{12}{14.4}{rm}$R_i$}}}
\end{picture}

\]
\begin{center}
        Figure 1.34
\end{center}
\ni writing $x=\Phi(\,\sum_i\lambda_i\tangle{ R_i})\in B_k$.

(i)$\Rightarrow$(ii) Suppose $x\in P_k$ satisfies
${\op{tr}}_R(x^*x)=0$. Then if $A\in{\cal A}_{k,0}$, we may isotope
$A(x)$ so it looks like Figure 1.35
\[
        \begin{picture}(0,0)%
\epsfig{file=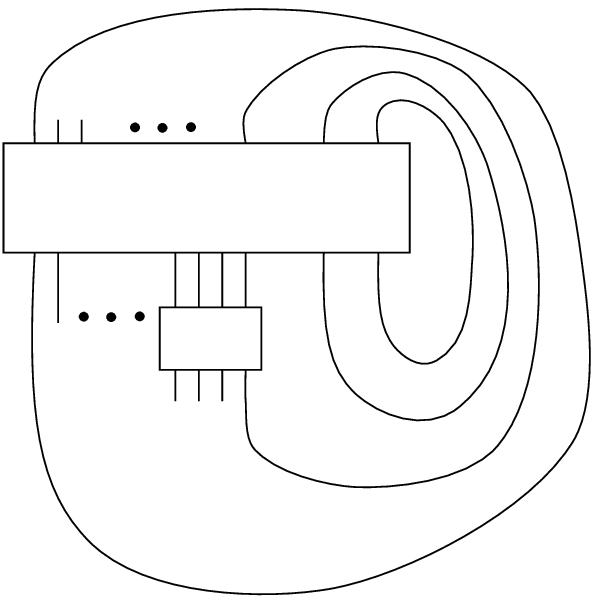}%
\end{picture}%
\setlength{\unitlength}{0.00041700in}%
\begingroup\makeatletter\ifx\SetFigFont\undefined
\def\x#1#2#3#4#5#6#7\relax{\def\x{#1#2#3#4#5#6}}%
\expandafter\x\fmtname xxxxxx\relax \def\y{splain}%
\ifx\x\y   
\gdef\SetFigFont#1#2#3{%
  \ifnum #1<17\tiny\else \ifnum #1<20\small\else
  \ifnum #1<24\normalsize\else \ifnum #1<29\large\else
  \ifnum #1<34\Large\else \ifnum #1<41\LARGE\else
     \huge\fi\fi\fi\fi\fi\fi
  \csname #3\endcsname}%
\else
\gdef\SetFigFont#1#2#3{\begingroup
  \count@#1\relax \ifnum 25<\count@\count@25\fi
  \def\x{\endgroup\@setsize\SetFigFont{#2pt}}%
  \expandafter\x
    \csname \romannumeral\the\count@ pt\expandafter\endcsname
    \csname @\romannumeral\the\count@ pt\endcsname
  \csname #3\endcsname}%
\fi
\fi\endgroup
\begin{picture}(5672,5685)(1779,-6455)
\put(3676,-4036){\makebox(0,0)[lb]{\smash{\SetFigFont{12}{14.4}{rm}$x$}}}
\put(3601,-2536){\makebox(0,0)[lb]{\smash{\SetFigFont{12}{14.4}{rm}$y$}}}
\end{picture}

\]
\begin{center}
        Figure1.35
\end{center}
\bs\ni where $y\in P_n$, $n\geq k$. Thus
$Z(A(x))={\op{tr}}_R(\tilde x y)$ where $\tilde x$ denotes $x$ with
$n-k$ vertical straight lines to the right and left of it.
By the Cauchy-Schwartz inequality,
$|{\op{tr}}_R(\tilde x y)|\leq \sqrt{{\op{tr}}_R(\tilde x^*\tilde
x)} \sqrt{{\op{tr}}_R(y^*y)}$, so if tr${}_R(x^*x)=0$, \
$Z(A(x))=0$. \qed

\bs We will call a general planar algebra $P$ {\it
finite-dimensional} if dim $P_k < \infty$ for all $k$.

\bs
\ni{\bf Corollary 1.36} If $P$ is a non-degenerate
finite-dimensional  planar $*$-algebra with positive partition
function then $P_k$ is semisimple for all $k$, so there is a unique
norm $\| \ \|$ on $P_k$ making it into a $C^*$-algebra.

{\sc Proof}. Each $P_k$ is semisimple since tr$(x^*x) > 0$
means there are no nilpotent ideals. The rest is standard. 
\qed

\bs{\bf Definition 1.37.} We call a planar algebra (over $\Bbb R$
or $\Bbb C$) a {\it $C^*$-planar algebra}
if it satisfies the conditions of corollary 1.34.

\bs\bs
\ni{\large\bf{2. Examples}}

{\bf Example 2.1: Temperley-Lieb algebra.}
If $\d_1$ and $\d_2$ are two non-zero scalars, one defines
$TL(n,\d_1,\d_2)$ as being the subspace of ${\cal P}_n(\emptyset)$
spanned by the tangles with no closed loops. Defining
multiplication on $TL(n,\d_1,\d_2)$ by multiplication as in
${\cal P}_n(\emptyset)$ except that one multiplies by a factor
$\d_1$ for each loop $\caright$ and $\d_2$
for each loop $\caleft$, then discarding the loop.
Clearly the map from ${\cal P}_n(\emptyset)$
to $TL(n,\d_1,\d_2)$ given by multiplying by $\d_1$'s or $\d_2$'s
then discarding loops, gives a $\Phi$ exhibiting
 $TL(n,\d_1,\d_2)$ as a planar algebra. For general values of
$\d_1$ and $\d_2$, $TL(n)$ is not non-degenerate. An extreme case is
$\d_1=\d_2=1$ where $Z(c(T_1-T_2))=0$ for all relevant tangles
$c,T_1,T_2$. In fact the structure of the algebras
$TL(n,\d_1,\d_2)$ (forgetting $\Phi$), depends only on $\d_1\d_2$.
To see this, show as in [GHJ] that $TL(n,\d_1,\d_2)$ is presented as
an algebra by $E_i$, $i\leq 1,\dots ,n-1$ with
$E^2_i=\d_1E_i$ for $i$ odd, $E^2_i=\d_2E_i$ for $i$ even, and
$E_iE_{i\pm 1}E_i=E_i$ and
$E_iE_j=E_jE_i$ for $|i-j|\geq 2$. Then setting
$e_i=\frac{1}{\d_1}E_i$ ($i$ odd), $e_i=\frac{1}{\d_2}E_i$ ($i$
even), the relations become $e^2_i=e_i$, \
$e_ie_{i\pm 1}e_i=\frac{1}{\d_1\d_2} e_i$, \
$e_ie_j=e_je_i$ for $|i-j|\geq 2$. If $\d_1=\d_2=\d$, we write
$TL(\d_1,\d_2)=TL(\d)$.

One may also obtain $TL$ via invariants, as a planar algebra on
one box, in several ways.

\medskip
$\underline{\mbox{(i) The chromatic polynomial.}}$ A planar network
$\cal N$ on
$L=L_2$ with \#$(L_2)\!=\!1$ determines a planar graph $G(\cal N)$
by choosing as vertices the positively oriented regions of $\Bbb
R^2\backslash\cal N$ and replacing the 2-boxes by edges joining the
corresponding vertices (thus \ $\fedgs \rightsquigarrow
\bullet$---$\bullet$). Fix $Q\in\Bbb C-\{0\}$ and let $Z(\cal N)=$
(chromatic polynomial of $G({\cal N})$ as evaluated at $Q)\times f$,
where $f\!=\!1$ if the outside region is negatively oriented and
$f=Q^{-1}$ if the outside region is positively oriented. To see that
${P}_Z$ is Temperley-Lieb, define the map $\alpha:{\cal P}(L)\to
TL(1,Q)$ by
$\alpha( \fedgs )=\smthvor  - \smtho $ (extended by multilinearity to
${\cal P}(L))$. It is easy to check that $\alpha$ makes $TL(1,Q)$ a
planar algebra on $L$ and the corresponding partition function is
$Z$ as above. Thus ${P}_Z$ is the non-degenerate quotient of
$TL(1,Q)$.

\medskip
$\underline{\mbox{(ii) The knot polynomial of [J2].}}$
Given a planar network $\cal N$ on one 2-box, replace the 2-box
$\fedgs$ by $\crss$ to get an unoriented link diagram.
Define $Z(\cal N)$ to be the Kauffman bracket ([Ka1]) of this
diagram. Sending $\crss$ to $A\  \smthv +A^{-1}\  \smth$
we see that this defines a map
from ${\cal P}(L)$ to $TL(-A^2-A^2)$ with the Temperley-Lieb
partition function.

\medskip
Both (i) and (ii) are generalized by the dichromatic polynomial
(see [Tut]).

\bigskip
{\bf Example 2.2: Planar algebras on 1-boxes.}
If $A$ is an associative algebra with identity and a trace
functional tr: $A\to K$, tr$(ab)={\op{tr}}(ba)$, tr$(1)=\d$, we may
form a kind of ``wreath product" of $A$ with $TL(n,\d)$.
In terms of generators and relations, we put
$L=L_1=A$ and
\[
        \input{xfig/pic15}
\]
One may give a direct construction of this planar algebra using a
basis as follows. Choose a basis $\{a_i\mid i\in I\}$ of $A$ with
$a_ia_j=\sum c^k_{ij}a_k$ for scalars $c^k_{ij}$.  (Assume
$1\in\{a_i\}$ for convenience.) Let $P^A_n$ be the vector space
whose basis is the set of all Temperley-Lieb basis $n$-tangles together
with  a function from the strings of the tangle to $\{a_i\}$.
Multiply these basis elements as for Temperley-Lieb except that,
when a string labelled $a_i$ is joined with one labelled $a_j$, the
result gives a sum over $j$ of $c^k_{ij}$ times the same underlying
Temperley-Lieb tangle with the joined string labelled $a_k$. In the
resulting sum of at most \#$(I)^n$ terms, if a closed loop is
labelled $a_k$, remove it and multiply by a factor of tr$(a_k)$.
This gives an associative algebra structure on each $P^A_n$. It
becomes a planar algebra on $A$ in the obvious way with $\Phi$
mapping $\boxa{a}$ to a linear combination of strings labelled
$a_j$, the coefficients being those of $a$ in the basis $\{a_i\}$.

If a string in $\P(A)$ has no 1-box on it, it is sent to the same
string labelled with 1. One may check that the kernel of $\Phi$ is
precisely the ideal generated by our relations $R$, so
$P^A=\P(A)/\cal J(R)$.

Observe how $P^A_n$ is a sum, over Temperley-Lieb basis tangles, of
tensor powers of $A$. When $n\!=\!2$ this gives an associative
algebra structure on $A\otimes A\oplus A\otimes A$. Explicitly,
write $(a\otimes b)\oplus 0$ as $a\otimes b$ and $0\oplus (x\otimes y)$
as $x\otimes y$. Multiplication is then determined by the rules:
\begin{eqnarray*}
(a_1\otimes b_1)(a_2\otimes b_2) & = a_1a_2\otimes b_2b_1 \\
(x_1\otimes y_1)(x_2\otimes y_2) & = {\op{tr}}(y_1x_2)x_1\otimes
y_2 \\
(a\otimes b)(x\otimes y) & = 0\oplus axb\otimes y \\
(x\otimes y)(a\otimes b) & = 0\oplus x\otimes bya
\end{eqnarray*}

The planar algebra $P^A$ may be degenerate,
even when tr on $A$ is non-degenerate and $\d$ is such that
$TL(\d)$ is non-degenerate. We will give more details on the
structure of $P^A$ in $\S$3.1.

\bigskip
{\bf Example 2.3: The Fuss-Catalan algebras} (see [BJ2]).
If $a_1,a_2,\dots ,a_k\in K-\{0\}$, $FC(n,a_1,\dots ,a_k)$ is the
algebra having as basis the Temperley-Lieb diagrams in $TL(nk)$ for
which, for each $p=1,2,\dots ,k$, the set of all boundary points
(counting from the left) indexed by
$\{jk+(-1)^jp+(\sin^2\frac{j\pi}{2})(k+j)\mid j=01,2,\dots (n-1)\}$
are connected among themselves. Assign a colour to each $p=1,2,\dots
,k$ so we think of the Temperley-Lieb strings as being coloured.
Then multiplication preserves colours so that closed loops will have
colours. Removing a closed loop coloured $m$ contributes a
multiplicative factor $a_m$. To see that $FC(n,a_1,\dots ,a_k)$ is
a planar algebra, begin with the case $k=2$. We claim $FC(n,a,b)$
is planar on one 2-box. We draw the box symbolically as
\[
        \begin{picture}(0,0)%
\epsfig{file=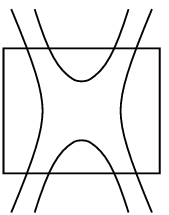}%
\end{picture}%
\setlength{\unitlength}{0.00041700in}%
\begingroup\makeatletter\ifx\SetFigFont\undefined
\def\x#1#2#3#4#5#6#7\relax{\def\x{#1#2#3#4#5#6}}%
\expandafter\x\fmtname xxxxxx\relax \def\y{splain}%
\ifx\x\y   
\gdef\SetFigFont#1#2#3{%
  \ifnum #1<17\tiny\else \ifnum #1<20\small\else
  \ifnum #1<24\normalsize\else \ifnum #1<29\large\else
  \ifnum #1<34\Large\else \ifnum #1<41\LARGE\else
     \huge\fi\fi\fi\fi\fi\fi
  \csname #3\endcsname}%
\else
\gdef\SetFigFont#1#2#3{\begingroup
  \count@#1\relax \ifnum 25<\count@\count@25\fi
  \def\x{\endgroup\@setsize\SetFigFont{#2pt}}%
  \expandafter\x
    \csname \romannumeral\the\count@ pt\expandafter\endcsname
    \csname @\romannumeral\the\count@ pt\endcsname
  \csname #3\endcsname}%
\fi
\fi\endgroup
\begin{picture}(1544,1980)(2379,-2551)
\end{picture}

\]
This shows in fact how to define the corresponding 
$\Phi: {\cal P}_n\to FC(n,a,b)$: double all the strings and replace all the
2-boxes according to the diagram. Thus for instance the tangle 
$\cal N$ below (with boxes shrunk to points, there being only one
2-box),
\[
        \begin{picture}(0,0)%
\epsfig{file=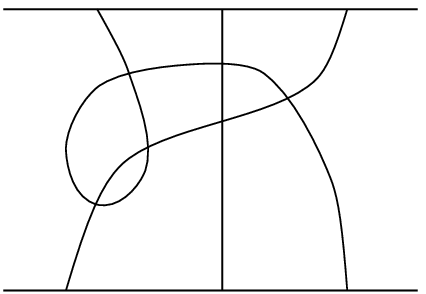}%
\end{picture}%
\setlength{\unitlength}{0.00041700in}%
\begingroup\makeatletter\ifx\SetFigFont\undefined
\def\x#1#2#3#4#5#6#7\relax{\def\x{#1#2#3#4#5#6}}%
\expandafter\x\fmtname xxxxxx\relax \def\y{splain}%
\ifx\x\y   
\gdef\SetFigFont#1#2#3{%
  \ifnum #1<17\tiny\else \ifnum #1<20\small\else
  \ifnum #1<24\normalsize\else \ifnum #1<29\large\else
  \ifnum #1<34\Large\else \ifnum #1<41\LARGE\else
     \huge\fi\fi\fi\fi\fi\fi
  \csname #3\endcsname}%
\else
\gdef\SetFigFont#1#2#3{\begingroup
  \count@#1\relax \ifnum 25<\count@\count@25\fi
  \def\x{\endgroup\@setsize\SetFigFont{#2pt}}%
  \expandafter\x
    \csname \romannumeral\the\count@ pt\expandafter\endcsname
    \csname @\romannumeral\the\count@ pt\endcsname
  \csname #3\endcsname}%
\fi
\fi\endgroup
\begin{picture}(4019,2744)(1779,-2783)
\end{picture}

\]
is sent to $\Phi(\cal N)$ below
\[
        \input{xfig/pic18}
\]
\noindent It is clear that $\Phi$ defines an algebra homomorphism
and surjectivity follows from [BJ2]. That ker
$\Phi$ is annular invariant is straightforward.
The general case of $FC(n,a_1,\dots ,a_k)$ is similar.
One considers the $k-1$ 2-boxes drawn symbolically as
\[
        \begin{picture}(0,0)%
\epsfig{file=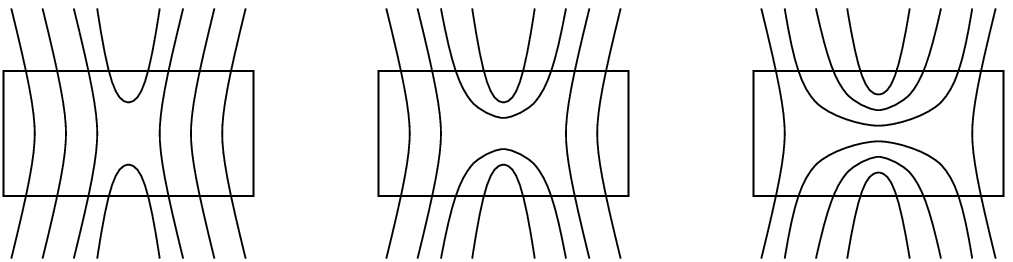}%
\end{picture}%
\setlength{\unitlength}{0.00041700in}%
\begingroup\makeatletter\ifx\SetFigFont\undefined
\def\x#1#2#3#4#5#6#7\relax{\def\x{#1#2#3#4#5#6}}%
\expandafter\x\fmtname xxxxxx\relax \def\y{splain}%
\ifx\x\y   
\gdef\SetFigFont#1#2#3{%
  \ifnum #1<17\tiny\else \ifnum #1<20\small\else
  \ifnum #1<24\normalsize\else \ifnum #1<29\large\else
  \ifnum #1<34\Large\else \ifnum #1<41\LARGE\else
     \huge\fi\fi\fi\fi\fi\fi
  \csname #3\endcsname}%
\else
\gdef\SetFigFont#1#2#3{\begingroup
  \count@#1\relax \ifnum 25<\count@\count@25\fi
  \def\x{\endgroup\@setsize\SetFigFont{#2pt}}%
  \expandafter\x
    \csname \romannumeral\the\count@ pt\expandafter\endcsname
    \csname @\romannumeral\the\count@ pt\endcsname
  \csname #3\endcsname}%
\fi
\fi\endgroup
\begin{picture}(9644,2430)(1779,-2776)
\end{picture}

\]

\noindent One proceeds as above, replacing the single strings in an
$\cal N$ by $k$ coloured strings. Surjectivity follows from
[La].

Note that these planar algebras give
invariants of systems of immersed curves with generic
singularities, and/or planar graphs. The most general such
invariant may be obtained by introducing a single 2-box which is a
linear combination of the $k-1$\quad 2-boxes described above. This will
generalize the dichromatic polynomial.

\bigskip{\bf Example 2.4: The BMW algebra.}
Let $L=L_2=\{R,Q\}$ and define the planar algebra $\cal B\cal M\cal
W$ on $L$ by the relations

        \input{xfig/pic201}

        \begin{picture}(0,0)%
\epsfig{file=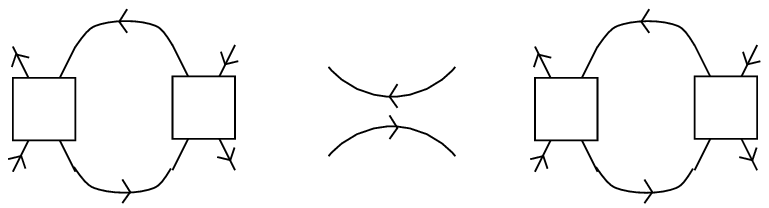}%
\end{picture}%
\setlength{\unitlength}{0.00041700in}%
\begingroup\makeatletter\ifx\SetFigFont\undefined
\def\x#1#2#3#4#5#6#7\relax{\def\x{#1#2#3#4#5#6}}%
\expandafter\x\fmtname xxxxxx\relax \def\y{splain}%
\ifx\x\y   
\gdef\SetFigFont#1#2#3{%
  \ifnum #1<17\tiny\else \ifnum #1<20\small\else
  \ifnum #1<24\normalsize\else \ifnum #1<29\large\else
  \ifnum #1<34\Large\else \ifnum #1<41\LARGE\else
     \huge\fi\fi\fi\fi\fi\fi
  \csname #3\endcsname}%
\else
\gdef\SetFigFont#1#2#3{\begingroup
  \count@#1\relax \ifnum 25<\count@\count@25\fi
  \def\x{\endgroup\@setsize\SetFigFont{#2pt}}%
  \expandafter\x
    \csname \romannumeral\the\count@ pt\expandafter\endcsname
    \csname @\romannumeral\the\count@ pt\endcsname
  \csname #3\endcsname}%
\fi
\fi\endgroup
\begin{picture}(8785,1905)(2281,-3969)
\put(2281,-3106){\makebox(0,0)[lb]{\smash{\SetFigFont{12}{14.4}{rm}$(ii)$}}}
\put(6250,-3165){\makebox(0,0)[lb]{\smash{\SetFigFont{12}{14.4}{rm}$=$}}}
\put(8290,-3165){\makebox(0,0)[lb]{\smash{\SetFigFont{12}{14.4}{rm}$=$}}}
\put(4067,-3090){\makebox(0,0)[lb]{\smash{\SetFigFont{12}{14.4}{rm}$R$}}}
\put(5626,-3091){\makebox(0,0)[lb]{\smash{\SetFigFont{12}{14.4}{rm}$Q$}}}
\put(9106,-3061){\makebox(0,0)[lb]{\smash{\SetFigFont{12}{14.4}{rm}$Q$}}}
\put(10699,-3076){\makebox(0,0)[lb]{\smash{\SetFigFont{12}{14.4}{rm}$R$}}}
\end{picture}

        \begin{picture}(0,0)%
\epsfig{file=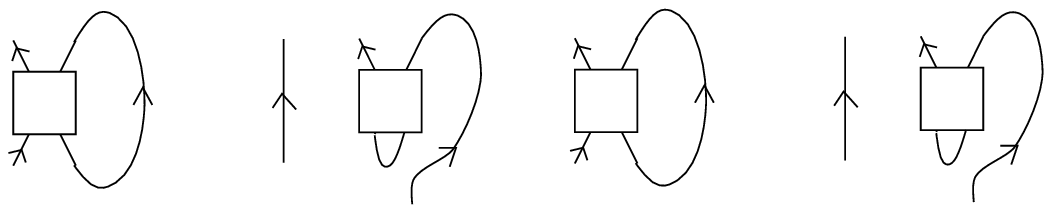}%
\end{picture}%
\setlength{\unitlength}{0.00041700in}%
\begingroup\makeatletter\ifx\SetFigFont\undefined
\def\x#1#2#3#4#5#6#7\relax{\def\x{#1#2#3#4#5#6}}%
\expandafter\x\fmtname xxxxxx\relax \def\y{splain}%
\ifx\x\y   
\gdef\SetFigFont#1#2#3{%
  \ifnum #1<17\tiny\else \ifnum #1<20\small\else
  \ifnum #1<24\normalsize\else \ifnum #1<29\large\else
  \ifnum #1<34\Large\else \ifnum #1<41\LARGE\else
     \huge\fi\fi\fi\fi\fi\fi
  \csname #3\endcsname}%
\else
\gdef\SetFigFont#1#2#3{\begingroup
  \count@#1\relax \ifnum 25<\count@\count@25\fi
  \def\x{\endgroup\@setsize\SetFigFont{#2pt}}%
  \expandafter\x
    \csname \romannumeral\the\count@ pt\expandafter\endcsname
    \csname @\romannumeral\the\count@ pt\endcsname
  \csname #3\endcsname}%
\fi
\fi\endgroup
\begin{picture}(12138,1904)(2266,-7063)
\put(2266,-6121){\makebox(0,0)[lb]{\smash{\SetFigFont{12}{14.4}{rm}$(iii)$}}}
\put(10816,-6129){\makebox(0,0)[lb]{\smash{\SetFigFont{12}{14.4}{rm}$=$}}}
\put(8494,-6253){\makebox(0,0)[lb]{\smash{\SetFigFont{12}{14.4}{rm},}}}
\put(7471,-6121){\makebox(0,0)[lb]{\smash{\SetFigFont{12}{14.4}{rm}$Q$}}}
\put(9556,-6121){\makebox(0,0)[lb]{\smash{\SetFigFont{12}{14.4}{rm}$Q$}}}
\put(12904,-6106){\makebox(0,0)[lb]{\smash{\SetFigFont{12}{14.4}{rm}$R$}}}
\put(12211,-6159){\makebox(0,0)[lb]{\smash{\SetFigFont{12}{14.4}{rm}$=$}}}
\put(11167,-6148){\makebox(0,0)[lb]{\smash{\SetFigFont{12}{14.4}{rm}$a^{-1}$}}}
\put(5425,-6150){\makebox(0,0)[lb]{\smash{\SetFigFont{12}{14.4}{rm}$=$}}}
\put(6076,-6184){\makebox(0,0)[lb]{\smash{\SetFigFont{12}{14.4}{rm}$a$}}}
\put(6775,-6165){\makebox(0,0)[lb]{\smash{\SetFigFont{12}{14.4}{rm}$=$}}}
\put(4129,-6151){\makebox(0,0)[lb]{\smash{\SetFigFont{12}{14.4}{rm}$R$}}}
\put(13984,-6223){\makebox(0,0)[lb]{\smash{\SetFigFont{12}{14.4}{rm},}}}
\put(14404,-6193){\makebox(0,0)[lb]{\smash{\SetFigFont{12}{14.4}{rm}$(a\in C-\{0\})$}}}
\end{picture}

        \begin{picture}(0,0)%
\epsfig{file=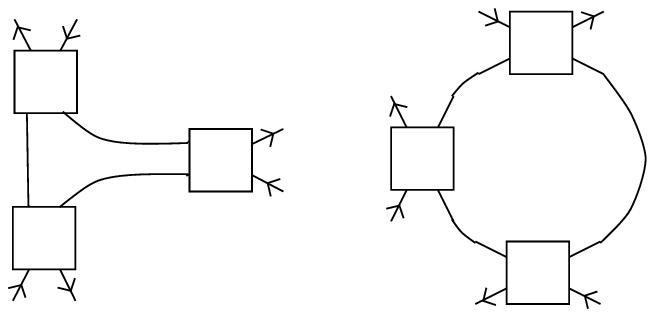}%
\end{picture}%
\setlength{\unitlength}{0.00041700in}%
\begingroup\makeatletter\ifx\SetFigFont\undefined
\def\x#1#2#3#4#5#6#7\relax{\def\x{#1#2#3#4#5#6}}%
\expandafter\x\fmtname xxxxxx\relax \def\y{splain}%
\ifx\x\y   
\gdef\SetFigFont#1#2#3{%
  \ifnum #1<17\tiny\else \ifnum #1<20\small\else
  \ifnum #1<24\normalsize\else \ifnum #1<29\large\else
  \ifnum #1<34\Large\else \ifnum #1<41\LARGE\else
     \huge\fi\fi\fi\fi\fi\fi
  \csname #3\endcsname}%
\else
\gdef\SetFigFont#1#2#3{\begingroup
  \count@#1\relax \ifnum 25<\count@\count@25\fi
  \def\x{\endgroup\@setsize\SetFigFont{#2pt}}%
  \expandafter\x
    \csname \romannumeral\the\count@ pt\expandafter\endcsname
    \csname @\romannumeral\the\count@ pt\endcsname
  \csname #3\endcsname}%
\fi
\fi\endgroup
\begin{picture}(7756,3100)(2267,-2575)
\put(2267,-1142){\makebox(0,0)[lb]{\smash{\SetFigFont{12}{14.4}{rm}$(iv)$}}}
\put(6827,-1187){\makebox(0,0)[lb]{\smash{\SetFigFont{12}{14.4}{rm}$=$}}}
\put(4187,-467){\makebox(0,0)[lb]{\smash{\SetFigFont{12}{14.4}{rm}$Q$}}}
\put(6062,-1202){\makebox(0,0)[lb]{\smash{
\SetFigFont{12}{14.4}{rm}$Q$
}}}
\put(9152,-122){\makebox(0,0)[lb]{\smash{
\SetFigFont{12}{14.4}{rm}$Q$
}}}
\put(9152,-2342){\makebox(0,0)[lb]{\smash{
\SetFigFont{12}{14.4}{rm}$R$
}}}
\put(4202,-2012){\makebox(0,0)[lb]{\smash{\SetFigFont{12}{14.4}{rm}$R$}}}
\put(7757,-1187){\makebox(0,0)[lb]{\smash{\SetFigFont{12}{14.4}{rm}$R$}}}
\end{picture}

        \begin{picture}(0,0)%
\epsfig{file=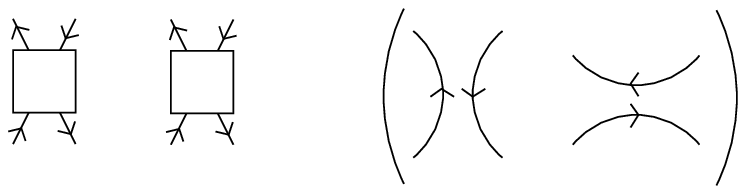}%
\end{picture}%
\setlength{\unitlength}{0.00041700in}%
\begingroup\makeatletter\ifx\SetFigFont\undefined
\def\x#1#2#3#4#5#6#7\relax{\def\x{#1#2#3#4#5#6}}%
\expandafter\x\fmtname xxxxxx\relax \def\y{splain}%
\ifx\x\y   
\gdef\SetFigFont#1#2#3{%
  \ifnum #1<17\tiny\else \ifnum #1<20\small\else
  \ifnum #1<24\normalsize\else \ifnum #1<29\large\else
  \ifnum #1<34\Large\else \ifnum #1<41\LARGE\else
     \huge\fi\fi\fi\fi\fi\fi
  \csname #3\endcsname}%
\else
\gdef\SetFigFont#1#2#3{\begingroup
  \count@#1\relax \ifnum 25<\count@\count@25\fi
  \def\x{\endgroup\@setsize\SetFigFont{#2pt}}%
  \expandafter\x
    \csname \romannumeral\the\count@ pt\expandafter\endcsname
    \csname @\romannumeral\the\count@ pt\endcsname
  \csname #3\endcsname}%
\fi
\fi\endgroup
\begin{picture}(8615,1725)(2297,-1374)
\put(4847,-528){\makebox(0,0)[lb]{\smash{\SetFigFont{12}{14.4}{rm}$+$}}}
\put(8826,-714){\makebox(0,0)[lb]{\smash{\SetFigFont{12}{14.4}{rm}$+$}}}
\put(2297,-543){\makebox(0,0)[lb]{\smash{\SetFigFont{12}{14.4}{rm}$(v)$}}}
\put(4157,-452){\makebox(0,0)[lb]{\smash{\SetFigFont{12}{14.4}{rm}$R$}}}
\put(5702,-452){\makebox(0,0)[lb]{\smash{\SetFigFont{12}{14.4}{rm}$Q$}}}
\put(6407,-558){\makebox(0,0)[lb]{\smash{\SetFigFont{12}{14.4}{rm}$=$}}}
\put(7235,-584){\makebox(0,0)[lb]{\smash{\SetFigFont{12}{14.4}{rm}$x$}}}
\end{picture}

Note that we could use relation (v) to express $\cal B\cal M\cal
W$ using only the one label $R$, but the relations would then be
more complicated. At this stage $\cal B\cal M\cal
W$ could be zero or infinite dimensional, but we may define a homomorphism from $\cal B\cal M\cal
W$ to the algebra $BMW$ of [BiW],[Mu] by sending \
$\blankbox{R}$ to $\crss$ and  $\blankbox{Q}$ to $\crssr$ .
This homomorphism is obviously surjective and one may use the
dimension count of [BiW] to show also that
dim ${\cal BMW}(n)\leq 1.3.5.\dots(2n-1)$ so that
${\cal BMW}\cong BMW$ as algebras. Thus $BMW$ is planar.
It is also connected and the invariant of planar networks is the
Kauffman regular-isotopy two-variable polynomial of [Ka2].

\bigskip{\bf Remark.} Had we presented $BMW$ on the single 2-box
$R$, the Reidemeister type III move (number (iv) above) would have
been
\[
        \input{xfig/pic21}
\]
\noindent
This leads us to consider the general planar algebra $B_n$ with the
following three conditions:
\begin{verse}
(1) $B_n$ is planar on one 2-box.

(2) dim $B_2=3$.

(3) dim $B_3\leq 15$.
\end{verse}
If one lists 16 tangles in $B_3$ then generically any one of them
will have to be a linear combination of the other $15$. Looking at
the 15th and 16th tangles in a listing according to the number of
2-boxes occuring in the tangle, we will generically obtain a type III Reidemeister 
move, or Yang-Baxter equation, modulo terms with less 2-boxes,as above.
It is not hard to show that these conditions force
dim $B_n\leq 1.3.5\cdot\dots \cdot (2n-1)$ since there are
necessarily Reidemeister-like moves of types I and II.
Note that $FC(n,a,b)$ satisfies these conditions as well as $BMW$!
For $C^*$-planar algebras, we have shown with Bisch
([BJ1]) that the only $B$'s with (1) and (2) as above,
and dim $B_3\leq 12$ are the Fuss-Catalan algebras  (with one
exception, when dim $B_3=9$).

\bigskip{\bf Example 2.5: a Hecke-algebra related example.}
The {\sc homfly} polynomial of [F+] is highly sensitive to the
orientation of a link and we may not proceed to use it to define a
planar algebra as in Example 2.4. In particular, a crossing in the
{\sc homfly} theory is necessarily oriented as $\crsso$.
Thus it does {\it not} yield a 2-box in our planar algebra context.
Nevertheless it is possible to use the {\sc homfly} skein theory to
define a planar algebra. We let $P^H_k$ be the usual {\sc homfly}
skein algebra of linear combinations of (3-dimensional) isotopy
classes of oriented tangles in the product of the $k$-box with an
interval, with orientations alternating out-in, modulo the
{\sc homfly} skein relation\newline
 $t\crsso -t^{-1}\crssor =x\smthv$ 
where $t\neq 0$ and $x$ are scalars. Projected onto the $k$-box,
such a tangle could look as in Figure 2.5.1,

\[
	\begin{picture}(0,0)%
\epsfig{file=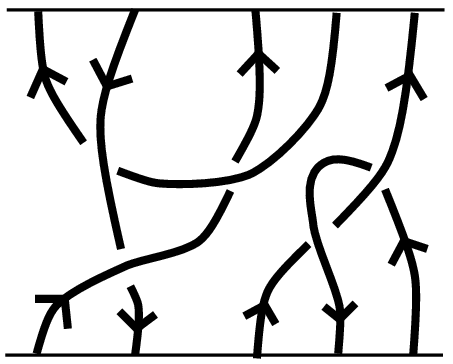}%
\end{picture}%
\setlength{\unitlength}{0.00041700in}%
\begingroup\makeatletter\ifx\SetFigFont\undefined
\def\x#1#2#3#4#5#6#7\relax{\def\x{#1#2#3#4#5#6}}%
\expandafter\x\fmtname xxxxxx\relax \def\y{splain}%
\ifx\x\y   
\gdef\SetFigFont#1#2#3{%
  \ifnum #1<17\tiny\else \ifnum #1<20\small\else
  \ifnum #1<24\normalsize\else \ifnum #1<29\large\else
  \ifnum #1<34\Large\else \ifnum #1<41\LARGE\else
     \huge\fi\fi\fi\fi\fi\fi
  \csname #3\endcsname}%
\else
\gdef\SetFigFont#1#2#3{\begingroup
  \count@#1\relax \ifnum 25<\count@\count@25\fi
  \def\x{\endgroup\@setsize\SetFigFont{#2pt}}%
  \expandafter\x
    \csname \romannumeral\the\count@ pt\expandafter\endcsname
    \csname @\romannumeral\the\count@ pt\endcsname
  \csname #3\endcsname}%
\fi
\fi\endgroup
\begin{picture}(4281,3435)(553,-3751)
\end{picture}

\]
\begin{center}
Figure 2.5.1
\end{center}

If we take the labeling set $L_k=P^H_k$, then $P^H$ is a general planar
algebra since planar isotopy of projections implies 3-dimensional
isotopy (invariance under the annular category is easy).
Standard {\sc homfly} arguments show dim $(P^H_k)\leq k!$
and a specialization could be used to obtain equality. Thus the
algebra is planar and the invariant is clearly the
{\sc homfly} polynomial of the
oriented link diagram given by a labeled network in $P^H_0$.
If we used this invariant to define the algebra as in $\S$1, we
would only obtain the same algebra for generic values of $t$ and
$x$. Note that $P^H_k$ is {\it not} isomorphic to the Hecke algebra
for $k\geq 4$, e.g. $P^H_4$ has an irreducible 4-dimensional
representation. In fact $P^H_k$ is, for generic $(t,x)$ and large
$n$, isomorphic to End${}_{SU(n)}
\underbrace{(V\otimes\bar V\otimes V\otimes \bar V\dots)}_{k\text{ \ vector\  spaces}}$,
where $V=\Bbb C^n$, the obvious $SU(n)$-module. This isomorphism is
only an algebra isomorphism, not a planar algebra isomorphism.

It is clear that the labeling set for $P^H$ could be reduced to a
set of $k!$ isotopy classes of tangles for $P^H_k$. But in fact a
single 3-label suffices as we now show

\bs
\ni{\bf Theorem 2.5.2} Any tangle {\it in the knot-theoretic sense}
with alternating in and out boundary orientations is
isotopic to a tangle with a diagram where all crossings occur in
disjoint discs which contain the pattern 
\[
	\begin{picture}(0,0)%
\epsfig{file=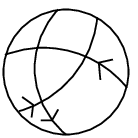}%
\end{picture}%
\setlength{\unitlength}{0.00041700in}%
\begingroup\makeatletter\ifx\SetFigFont\undefined
\def\x#1#2#3#4#5#6#7\relax{\def\x{#1#2#3#4#5#6}}%
\expandafter\x\fmtname xxxxxx\relax \def\y{splain}%
\ifx\x\y   
\gdef\SetFigFont#1#2#3{%
  \ifnum #1<17\tiny\else \ifnum #1<20\small\else
  \ifnum #1<24\normalsize\else \ifnum #1<29\large\else
  \ifnum #1<34\Large\else \ifnum #1<41\LARGE\else
     \huge\fi\fi\fi\fi\fi\fi
  \csname #3\endcsname}%
\else
\gdef\SetFigFont#1#2#3{\begingroup
  \count@#1\relax \ifnum 25<\count@\count@25\fi
  \def\x{\endgroup\@setsize\SetFigFont{#2pt}}%
  \expandafter\x
    \csname \romannumeral\the\count@ pt\expandafter\endcsname
    \csname @\romannumeral\the\count@ pt\endcsname
  \csname #3\endcsname}%
\fi
\fi\endgroup
\begin{picture}(1230,1230)(886,-976)
\end{picture}

\]
with some non-alternating choice of crossings.

\bs
{\sc Proof.} We begin with a tangle without boundary, i.e.~an
oriented link $L$. Choose a diagram for $L$ and add a parallel
double $L'$ of $L$ to the left of $L$ and oppositely oriented, with
crossings chosen so that

\begin{itemize}
\item[a)] $L'$ is always under $L$

\item[b)] $L'$  itself is an unlink.
\end{itemize}
An example of the resulting diagram (for the Whitehead Link) is
given in Figure 2.5.3
\[
	\begin{picture}(0,0)%
\epsfig{file=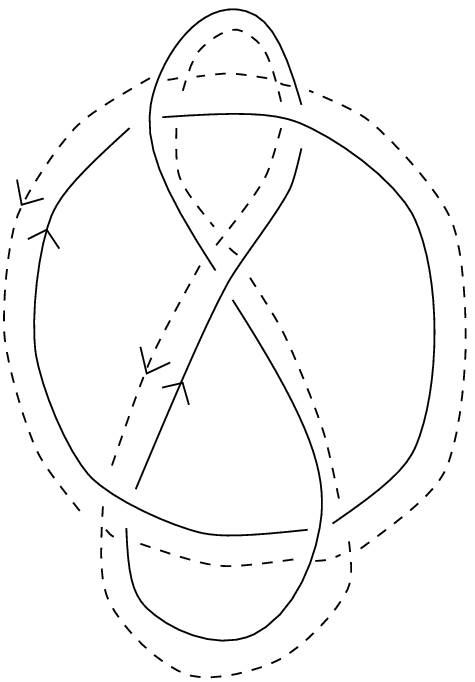}%
\end{picture}%
\setlength{\unitlength}{0.00041700in}%
\begingroup\makeatletter\ifx\SetFigFont\undefined
\def\x#1#2#3#4#5#6#7\relax{\def\x{#1#2#3#4#5#6}}%
\expandafter\x\fmtname xxxxxx\relax \def\y{splain}%
\ifx\x\y   
\gdef\SetFigFont#1#2#3{%
  \ifnum #1<17\tiny\else \ifnum #1<20\small\else
  \ifnum #1<24\normalsize\else \ifnum #1<29\large\else
  \ifnum #1<34\Large\else \ifnum #1<41\LARGE\else
     \huge\fi\fi\fi\fi\fi\fi
  \csname #3\endcsname}%
\else
\gdef\SetFigFont#1#2#3{\begingroup
  \count@#1\relax \ifnum 25<\count@\count@25\fi
  \def\x{\endgroup\@setsize\SetFigFont{#2pt}}%
  \expandafter\x
    \csname \romannumeral\the\count@ pt\expandafter\endcsname
    \csname @\romannumeral\the\count@ pt\endcsname
  \csname #3\endcsname}%
\fi
\fi\endgroup
\begin{picture}(4473,6470)(1684,-6412)
\end{picture}

\]
\begin{center}
	Figure 2.5.3
\end{center}
\ni Now join $L$ to $L'$, component by component, by replacing
$\paro$  by $\smthvo$, at some point well away from any crossings. 
Since $L'$ is
an unlink below $L$, the resulting link is isotopic to $L$. 
All the crossings in $L\cup L'$ occur in disjoint discs containing the pattern
\begin{picture}(0,0)%
\epsfig{file=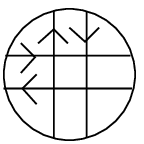}%
\end{picture}%
\setlength{\unitlength}{0.00041700in}%
\begingroup\makeatletter\ifx\SetFigFont\undefined
\def\x#1#2#3#4#5#6#7\relax{\def\x{#1#2#3#4#5#6}}%
\expandafter\x\fmtname xxxxxx\relax \def\y{splain}%
\ifx\x\y   
\gdef\SetFigFont#1#2#3{%
  \ifnum #1<17\tiny\else \ifnum #1<20\small\else
  \ifnum #1<24\normalsize\else \ifnum #1<29\large\else
  \ifnum #1<34\Large\else \ifnum #1<41\LARGE\else
     \huge\fi\fi\fi\fi\fi\fi
  \csname #3\endcsname}%
\else
\gdef\SetFigFont#1#2#3{\begingroup
  \count@#1\relax \ifnum 25<\count@\count@25\fi
  \def\x{\endgroup\@setsize\SetFigFont{#2pt}}%
  \expandafter\x
    \csname \romannumeral\the\count@ pt\expandafter\endcsname
    \csname @\romannumeral\the\count@ pt\endcsname
  \csname #3\endcsname}%
\fi
\fi\endgroup
\begin{picture}(1298,1288)(548,-1305)
\end{picture}
 
which can be isotoped to the pattern
\begin{picture}(0,0)%
\epsfig{file=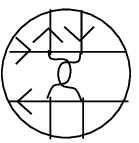}%
\end{picture}%
\setlength{\unitlength}{0.00041700in}%
\begingroup\makeatletter\ifx\SetFigFont\undefined
\def\x#1#2#3#4#5#6#7\relax{\def\x{#1#2#3#4#5#6}}%
\expandafter\x\fmtname xxxxxx\relax \def\y{splain}%
\ifx\x\y   
\gdef\SetFigFont#1#2#3{%
  \ifnum #1<17\tiny\else \ifnum #1<20\small\else
  \ifnum #1<24\normalsize\else \ifnum #1<29\large\else
  \ifnum #1<34\Large\else \ifnum #1<41\LARGE\else
     \huge\fi\fi\fi\fi\fi\fi
  \csname #3\endcsname}%
\else
\gdef\SetFigFont#1#2#3{\begingroup
  \count@#1\relax \ifnum 25<\count@\count@25\fi
  \def\x{\endgroup\@setsize\SetFigFont{#2pt}}%
  \expandafter\x
    \csname \romannumeral\the\count@ pt\expandafter\endcsname
    \csname @\romannumeral\the\count@ pt\endcsname
  \csname #3\endcsname}%
\fi
\fi\endgroup
\begin{picture}(1260,1289)(571,-1313)
\end{picture}
 
containing two discs of the required form. 
Alternating patterns can be
avoided by keeping the top string on top in this isotopy.

For a tangle with boundary we make the doubling curve follow the
boundary, turning right just before it would hit it and right
again as it nears the point where the next string exits the
tangle, as in Figure 2.5.4
\[
	\begin{picture}(0,0)%
\epsfig{file=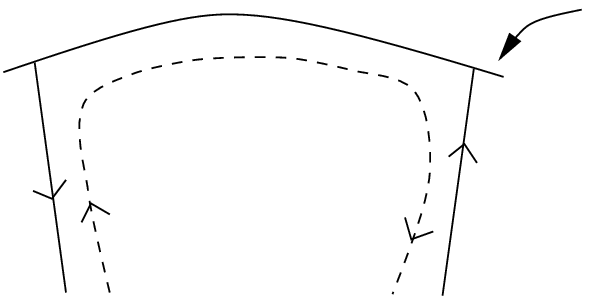}%
\end{picture}%
\setlength{\unitlength}{0.00041700in}%
\begingroup\makeatletter\ifx\SetFigFont\undefined
\def\x#1#2#3#4#5#6#7\relax{\def\x{#1#2#3#4#5#6}}%
\expandafter\x\fmtname xxxxxx\relax \def\y{splain}%
\ifx\x\y   
\gdef\SetFigFont#1#2#3{%
  \ifnum #1<17\tiny\else \ifnum #1<20\small\else
  \ifnum #1<24\normalsize\else \ifnum #1<29\large\else
  \ifnum #1<34\Large\else \ifnum #1<41\LARGE\else
     \huge\fi\fi\fi\fi\fi\fi
  \csname #3\endcsname}%
\else
\gdef\SetFigFont#1#2#3{\begingroup
  \count@#1\relax \ifnum 25<\count@\count@25\fi
  \def\x{\endgroup\@setsize\SetFigFont{#2pt}}%
  \expandafter\x
    \csname \romannumeral\the\count@ pt\expandafter\endcsname
    \csname @\romannumeral\the\count@ pt\endcsname
  \csname #3\endcsname}%
\fi
\fi\endgroup
\begin{picture}(5745,2926)(586,-3998)
\put(6331,-1336){\makebox(0,0)[lb]{\smash{\SetFigFont{12}{14.4}{rm}boundary of tangle}}}
\end{picture}

\]
\begin{center}
	Figure 2.5.4
\end{center}
\ni Join the original tangle to $L'$ one string at a time and
proceed as before. \qed

\bs\ni{\bf Corollary 2.5.5} The planar algebra $P^H$ is generated by
the single 3-box \begin{picture}(0,0)%
\epsfig{file=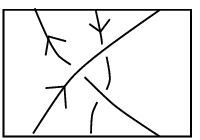}%
\end{picture}%
\setlength{\unitlength}{0.00041700in}%
\begingroup\makeatletter\ifx\SetFigFont\undefined
\def\x#1#2#3#4#5#6#7\relax{\def\x{#1#2#3#4#5#6}}%
\expandafter\x\fmtname xxxxxx\relax \def\y{splain}%
\ifx\x\y   
\gdef\SetFigFont#1#2#3{%
  \ifnum #1<17\tiny\else \ifnum #1<20\small\else
  \ifnum #1<24\normalsize\else \ifnum #1<29\large\else
  \ifnum #1<34\Large\else \ifnum #1<41\LARGE\else
     \huge\fi\fi\fi\fi\fi\fi
  \csname #3\endcsname}%
\else
\gdef\SetFigFont#1#2#3{\begingroup
  \count@#1\relax \ifnum 25<\count@\count@25\fi
  \def\x{\endgroup\@setsize\SetFigFont{#2pt}}%
  \expandafter\x
    \csname \romannumeral\the\count@ pt\expandafter\endcsname
    \csname @\romannumeral\the\count@ pt\endcsname
  \csname #3\endcsname}%
\fi
\fi\endgroup
\begin{picture}(1844,1267)(879,-1883)
\end{picture}
.

\bs
{\sc Proof.} The {\sc homfly} relations can be used to go between
the various possible choices of crossings in the 3-box of Theorem
2.5.2. \qed

\bs This corollary was first proved  by W.~B.~R.~Lickorish using an
argument adapted to the {\sc homfly} skein. His argument is much
more efficient in producing a skein element involving only the
above 3-box. The tangles may be chosen
{\it alternating} in Theorem 2.5.2.

\bs{\bf Remarks.} 1) Another way of stating Theorem 2.5.2 is to say
that any tangle can be projected with only simple tiple point
singularities.  One may ask if there are a set of ``Reidemeister
moves" for such non-generic projections.

2) A related question would be to find a presentation of $P^H$ on
the above 3-box.

\bs{\bf Discussion 2.5.6}  In the remark of Example 2.4 we
introduced relations on the planar algebra generated by a 2-box,
which force finite dimensionality of all the $P_n$'s. One should
explore the possibilities for the planar algebra generated by a
single 3-box. The dimension restrictions analogous to the 1,3 $\leq
15$ values of Example 2.4 are 1,2,6,$\leq 24$ and we conjecture,
somewhat weakly, the following

\bs
\ni
{\bf Conjecture 2.5.7} Let $(P,\Phi)$ be a planar algebra with
labelling set $L=L_3$, $\#(L_3)=1$. Suppose $\dim P_n\leq n!$ for
$n\leq 4$. Let $V$ be the subspace of ${\cal P}_4(L)$ spanned by
tangles with at most two labeled 3-boxes, and let
$R=V\cap{\op{ker}} \ \Phi$ be relations. Then
\[
\dim\left( \frac{{\cal P}_n(L)}{{\cal J}_n(R)}\right)\leq n! \quad
{\text{for all}} \ n.
\]

\bigskip
There is some evidence for the conjecture. It would imply in
particular the $n=0$ case which implies the following
result, proved by D. Thurston about hexagons:

``Consider all graphs with hexagonal faces that may be drawn on
$S^2$ with non-intersecting edges. Let $M$ be the move of Figure
2.5.8 on the set of all such graphs (where the 8 external vertices
are connected in an arbitrary way to the rest of the graph).
\[
	\begin{picture}(0,0)%
\epsfig{file=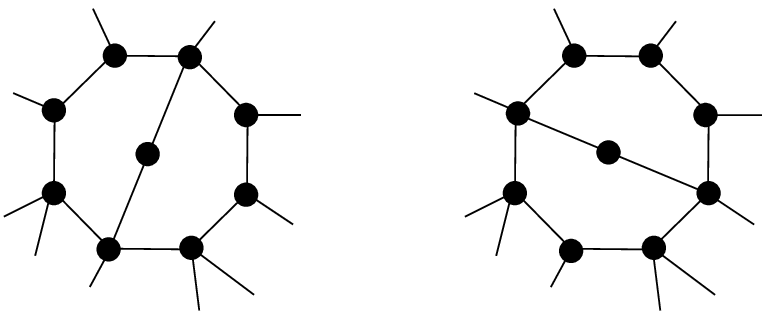}%
\end{picture}%
\setlength{\unitlength}{0.00041700in}%
\begingroup\makeatletter\ifx\SetFigFont\undefined
\def\x#1#2#3#4#5#6#7\relax{\def\x{#1#2#3#4#5#6}}%
\expandafter\x\fmtname xxxxxx\relax \def\y{splain}%
\ifx\x\y   
\gdef\SetFigFont#1#2#3{%
  \ifnum #1<17\tiny\else \ifnum #1<20\small\else
  \ifnum #1<24\normalsize\else \ifnum #1<29\large\else
  \ifnum #1<34\Large\else \ifnum #1<41\LARGE\else
     \huge\fi\fi\fi\fi\fi\fi
  \csname #3\endcsname}%
\else
\gdef\SetFigFont#1#2#3{\begingroup
  \count@#1\relax \ifnum 25<\count@\count@25\fi
  \def\x{\endgroup\@setsize\SetFigFont{#2pt}}%
  \expandafter\x
    \csname \romannumeral\the\count@ pt\expandafter\endcsname
    \csname @\romannumeral\the\count@ pt\endcsname
  \csname #3\endcsname}%
\fi
\fi\endgroup
\begin{picture}(7320,2940)(850,-2725)
\end{picture}

\]
\begin{center}
	Figure 2.5.8: The move $M$
\end{center}

\medskip\ni Then one may find a finite number of applications of
the move $M$ leading to a graph with two adjacent 2-valent
vertices." 

\bs
\bs
{\bf Example 2.6: Tensors.}
Let $V$ be a finite dimensional vector space with dual $\hat V$.
We will define a planar algebra $P^{\otimes}\!\!=\!\cup_k
P^{\otimes}_k$ with dim $P_0^{\otimes}\!=\!0$ and
$P_k^{\otimes}\!\!={\op{End}}(V\!\otimes\!\hat V\!\otimes
\!V\!\otimes\!
\hat V\!\otimes\!\dots)$ where there are $k$ vector spaces in the
tensor product. The planar structure on $P^{\otimes}$ can be defined
invariantly using the canonical maps $V\otimes\hat V\to K$ and
$\hat V\otimes  V\to K$, which are applied to any pair of vector
spaces connected by an internal edge in a planar tangle, where $V$
and $\hat V$ are associated with the marked points of a $k$-box in
an alternating fashion with $V$ associated to $*$. One could
also think of $\Bbb C\subset$End$(V)$ as a finite factor
and use the method of Theorem 4.2.1.
 It is perhaps easiest to understand this structure using
a basis
$(v_1,v_2,\dots ,v_n)$ of $V$, with corresponding dual basis.
An element of
$P^\otimes_k$ is then the same as a tensor
$X^{j_1j_2\dots j_k}_{i_1i_2\dots i_k}$. The labelling set $L_k$ is
$P^\otimes_k$ itself and the presenting map $\Phi:{\cal P}_k(L)\to\
P^\otimes_k$ is defined by summing (``contracting") over all the
internal indices  in a labelled planar tangle. The first marked
point in a box corresponds to the ``$j_1$" above. To be more
precise, one defines a {\it state} $\sigma$ of a planar tangle $T$
to be a function from the connected components of the set $S(T)$
of strings in $T$, \ $\sigma: S(T)\to\{1,2,\dots ,n\}$, to the basis
elements of $V$. A state defines a set of indices around every box
$B$ in $T$, and since the label associated to $B$ is a tensor, with
the appropriate number of indices, to each labelled box, the state
$\sigma$ associates a number, $\sigma(B)$. A state also induces a
function $\p\sigma$ from the marked points on the boundary of $T$
to $\{1,2,\dots , n\}$. Now we associate a tensor $\Phi(T)$ with
$T$ as follows: let $f:$ \{marked points $(T)\}\to\{1,2,\dots ,n\}$
denote the indices 
$\dsize{{j_1\ldots j_k}\choose {i_1 \ldots  i_k}}$
of a tensor in $P^\otimes_k\quad (f(p,0)=i_p,f(p,1)=j_p,$ for $ 1\leq p \leq k)$. Then
$$
\Phi(T)^{j_1\dots j_k}_{i_1\dots i_k}=\sum_{\sigma:\p\sigma =f}
\left( \prod_{B\in\{ {\text{labelled}}\atop {\text{boxes of}} \ T\}
} \sigma(B)\right) \ .
$$
An empty sum is zero and an empty product is 1. One easily checks
that $\Phi$ defines an algebra homomorphism
and that ker $\Phi$ is invariant under the annular category.

This planar algebra has an obvious
$*$-structure. The invariant $Z$ is recognizable as the partition
function for the ``vertex model" defined by the labelled network,
the labels supplying the Boltzmann weights (see [Ba]).

For further discussion we introduce the following notation ---
consider the indices as a (finite) set $\D$. Given a function
$\dsize{{\g_1\ldots \g_k}\choose {\d_1 \ldots  \d_k}}$  
from the marked
points of a $k$-box to $\D$, we define the corresponding
{\it basic tensor} to be
\[ 
T^{j_1\dots j_k}_{i_1\dots i_k}= \left\{
		\begin{array}{ll}
		1 & {\mbox{if}}\ \ i_1=\d_1, \ i_2=\d_2 {\mbox{ \ etc.}}\\ 
		{} & {\mbox{and}} \ \ j_1=\g_1, \ j_2=\g_2 {\mbox{ \ etc.}}\\
		0 & {\mbox{otherwise.}}
		\end{array}
	\right.
\]

If $\frak S(\D)$ is the free semigroup on $\D$, $\p T$ is then the
word \newline $\g_1\g_2\g_3\dots \g_k\d_k\d_{k-1}\dots\d_1$ in $\frak S(\D)$, and we will use
the notation
\fbox{\parbox{2.3cm}{
$\dsize{
\begin{array}{c}
\g_1\ \g_2\ldots \g_k \\
\d_1\ \d_2 \ldots  \d_k
\end{array}
	}$
}} 
for this basic tensor.

The planar algebra $P^\otimes$ is not terribly interesting by
itself (and there seems  to be no reason to limit the contractions
allowed to planar ones). But one may look for planar subalgebras.
One way is to take a set $\{A_i\in P^\otimes_{k_i}\}$ and look at
the planar subalgebra $P_k(A_i)$ they generate. The calculation of
$P_k(A_i)$ as a function of the $A_i$'s can be extremely difficult.
While it is easy enough to decide if the $A_i$'s are in the TL
subalgebra, we will see in the next example that the question of
whether $P_k(A_i)\neq P_k^{\otimes}$ is undecidable, even for $k=1$.

Note that if the tensors $A_i$ have only 0--1 entries, the
partition function will be simply the number of ``edge colourings"
of the network by $n$ colours with colourings allowed only if they
correspond to a non-zero entry of the tensor label at each box.

The next example gives a situation where we {\it can} say
$P_k(A_i)\neq P^{\otimes}_k$.

\bigskip{\bf Example 2.7: Finitely generated groups.}
As in 2.6, if $\D$ is a set, $\frak S(\D)$ will denote the free
semigroup on $\D$, and $F(\D)$ will denote the free group on $\D$.
We define the map
alt:$\frak S(\D)\to F(\D)$ by alt$(\g_1\dots\g_m)=
\g_1\g_2^{-1}\g_3\g_4^{-1}\dots \g_m^{\pm 1}$,
(where the $+$ sign occurs if $m$ is odd, -- if $m$ is
even). Note that alt is only a homomorphism from the subsemigroup
of words of even length.

Now let $\G$ be a discrete group and $\D$ a finite set, together
with a function $\d\mapsto\tilde\d$ from $\D$ to $\G$. There is
then a natural map $\phi: F(\D)\to\G$ defined by
$\phi(\d)=\tilde\d$. Let
$V$ be the vector space with basis $\D$. Use $V$ and the basis $\D$
to form the planar algebra $P^{\otimes}$ of $\S$2.6. Recall that,
for a basic tensor
$T\in P^{\otimes}_k$, \ $\p T$ is the element of $\frak S(\D)$
obtained by reading around the boundary of $T$. Let $*$ denote the
involution on $\frak S(\D)$ given by writing words backwards.

\bigskip{\bf Definition.} $P^{\G,\D}=\cup_k P_k^{\G,\D}$ is the
linear span of all basic tensors $T$ such that
$\phi({\op{alt}}(\p T))=1$ in $\G$.

\ni\bs{\bf Proposition 2.7.1} $P^{\G,\D}$ is a planar $*$-subalgebra
of $P^{\otimes}$.

\bs
{\bf Proof}. By 1.18, it suffices to show $P^{\G,\D}$ is a unital
subalgebra invariant under the annular category ${\cal A} (\emptyset)$. If $T$ is a basic
tensor in $P_k^{\G,\D}$, let $\p^+T$ (resp.~$\p^-T$) be the element
of $\frak S(\D)$ obtained by reading along the top of $T$
(resp.~the bottom), so $\p T=\p^+ T(\p^-T)^*$. Then for the
product $T_1T_2$ to be non-zero, $\p^-T_1=\p^+T_2$.
In the product alt$(T_1)$alt$(T_2)$, the last letter of
$(\p^-T_1)^*$ is then the same as the first letter of $\p T_2$, but
with opposite sign. Thus the contribution of $\p^-T_1$ cancels with
that of $\p^+T_2$ and $\phi({\op{alt}}(T_1){\op{alt}}(T_2))=1$.
So $P^{\G,\D}$ is a subalgebra, clearly unital and self-adjoint.

Now consider a typical ${\cal A} (\emptyset)$ element $C$ applied to a basic tensor
$T$ as in Figure~2.7.2. 
\[
	\begin{picture}(0,0)%
\epsfig{file=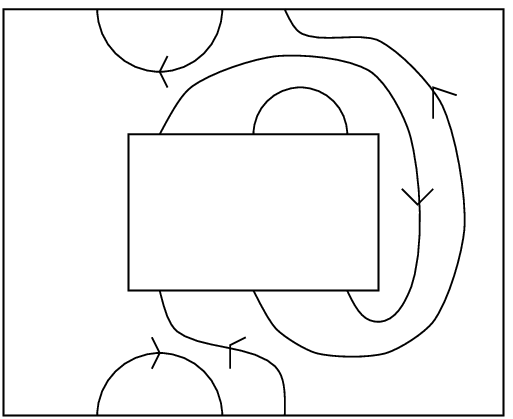}%
\end{picture}%
\setlength{\unitlength}{0.00041700in}%
\begingroup\makeatletter\ifx\SetFigFont\undefined
\def\x#1#2#3#4#5#6#7\relax{\def\x{#1#2#3#4#5#6}}%
\expandafter\x\fmtname xxxxxx\relax \def\y{splain}%
\ifx\x\y   
\gdef\SetFigFont#1#2#3{%
  \ifnum #1<17\tiny\else \ifnum #1<20\small\else
  \ifnum #1<24\normalsize\else \ifnum #1<29\large\else
  \ifnum #1<34\Large\else \ifnum #1<41\LARGE\else
     \huge\fi\fi\fi\fi\fi\fi
  \csname #3\endcsname}%
\else
\gdef\SetFigFont#1#2#3{\begingroup
  \count@#1\relax \ifnum 25<\count@\count@25\fi
  \def\x{\endgroup\@setsize\SetFigFont{#2pt}}%
  \expandafter\x
    \csname \romannumeral\the\count@ pt\expandafter\endcsname
    \csname @\romannumeral\the\count@ pt\endcsname
  \csname #3\endcsname}%
\fi
\fi\endgroup
\begin{picture}(4844,3944)(879,-3683)
\put(2326,-1411){\makebox(0,0)[lb]{\smash{\SetFigFont{12}{14.4}{rm}$\gamma_1$}}}
\put(3226,-1411){\makebox(0,0)[lb]{\smash{\SetFigFont{12}{14.4}{rm}$\gamma_2$}}}
\put(4051,-1411){\makebox(0,0)[lb]{\smash{\SetFigFont{12}{14.4}{rm}$\gamma_3$}}}
\put(2251,-2161){\makebox(0,0)[lb]{\smash{\SetFigFont{12}{14.4}{rm}$\gamma_6$}}}
\put(3076,-2236){\makebox(0,0)[lb]{\smash{\SetFigFont{12}{14.4}{rm}$\gamma_5$}}}
\put(4051,-2236){\makebox(0,0)[lb]{\smash{\SetFigFont{12}{14.4}{rm}$\gamma_4$}}}
\end{picture}

\]
\begin{center}
 	Figure 2.7.2
\end{center}

\noindent This tensor is a sum of basic tensors, the sum ranging
over all functions from the curves in the diagram to $\D$. If $R$
is a particular basic tensor in the sum, notice that the non through-strings in $C$
contribute to alt$(\p R)$ in two ways --- either they occur in
cancelling pairs or, if their beginning and end are separated by
the left-hand side of the diagram, they change alt$(\p R)$ by
conjugation (eliminate all the cancelling pairs first to see this).
Thus the conjugacy class in $F(\D)$ of alt$(R)$ (and alt$(T)$) is
not changed by removing all non through-strings. Once this is done,
however, the word around the outer boundary is just an even cyclic
permutation of the word $w'$ around the inner boundary, so
$\phi({\op{alt}} \ w)=1\Longleftrightarrow \phi({\op{alt}}
(w'))=1$. Thus $\P^{\G,\D}$ is invariant under ${\cal A} (\emptyset)$. \qed

\bigskip
Note that the basic tensors in $P^{\G,\D}$ are unchanged if we
change $\sim:\D\to\G$ by right multiplication by an element of
$\G$. So we may suppose there is an element $e$ of $\D$ with
$\tilde e=1\in\G$. To denote this situation we will say simply
``$e\in\D$".

Let $G\!\subseteq\G$ (resp.~$G'$) $=\{\phi({\op{alt}}(\p^+ T))\mid T
{\mbox{ \ a\  basic\  tensor in }} P^{\G,\D}_{2k} \
({\text{resp.}} \ P^{\G,\D}_k), \ k\!\in\!\Bbb N\}$.

\bs\ni
{\bf Lemma 2.7.3} $G$ is the subgroup
$\<\tilde\D\tilde\D^{-1}\>$ of $\G$ generated by
$\tilde\D\tilde\D^{-1}$, and if
$e\in\D$, \ $G=G'$.

\bs
{\sc Proof}. The definition of alt implies immediately that
$G\subseteq\<\tilde\D\tilde\D^{-1}\>$. That $G=G^{-1}$ follows
from Figure~2.7.4
\[
        \begin{picture}(0,0)%
\epsfig{file=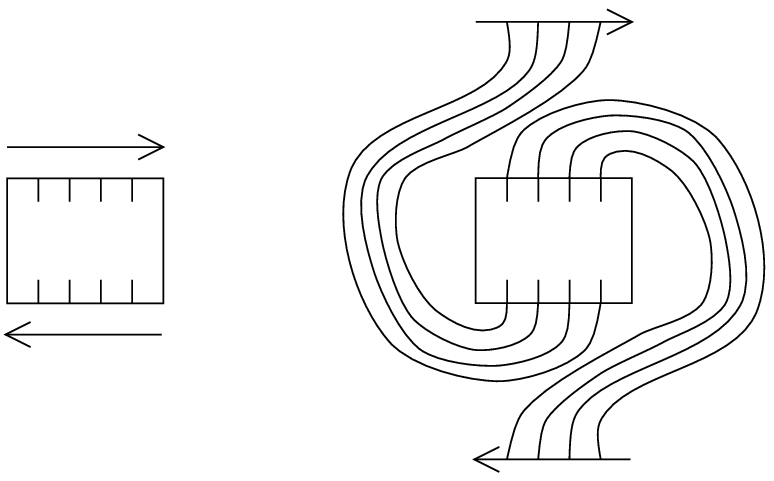}%
\end{picture}%
\setlength{\unitlength}{0.00041700in}%
\begingroup\makeatletter\ifx\SetFigFont\undefined
\def\x#1#2#3#4#5#6#7\relax{\def\x{#1#2#3#4#5#6}}%
\expandafter\x\fmtname xxxxxx\relax \def\y{splain}%
\ifx\x\y   
\gdef\SetFigFont#1#2#3{%
  \ifnum #1<17\tiny\else \ifnum #1<20\small\else
  \ifnum #1<24\normalsize\else \ifnum #1<29\large\else
  \ifnum #1<34\Large\else \ifnum #1<41\LARGE\else
     \huge\fi\fi\fi\fi\fi\fi
  \csname #3\endcsname}%
\else
\gdef\SetFigFont#1#2#3{\begingroup
  \count@#1\relax \ifnum 25<\count@\count@25\fi
  \def\x{\endgroup\@setsize\SetFigFont{#2pt}}%
  \expandafter\x
    \csname \romannumeral\the\count@ pt\expandafter\endcsname
    \csname @\romannumeral\the\count@ pt\endcsname
  \csname #3\endcsname}%
\fi
\fi\endgroup
\begin{picture}(7349,5310)(1165,-5383)
\put(1877,-2762){\makebox(0,0)[lb]{\smash{\SetFigFont{12}{14.4}{rm}$T$}}}
\put(6376,-2761){\makebox(0,0)[lb]{\smash{\SetFigFont{12}{14.4}{rm}$T$}}}
\put(1201,-1561){\makebox(0,0)[lb]{\smash{\SetFigFont{12}{14.4}{rm}$\phi ($alt$(\partial^+T))=w$}}}
\put(1201,-4186){\makebox(0,0)[lb]{\smash{\SetFigFont{12}{14.4}{rm}$\phi ($alt$(\partial^+T))=w^{-1}$}}}
\put(5701,-361){\makebox(0,0)[lb]{\smash{\SetFigFont{12}{14.4}{rm}$\phi ($alt$(\partial^+T))=w^{-1}$}}}
\put(5476,-5311){\makebox(0,0)[lb]{\smash{\SetFigFont{12}{14.4}{rm}$\phi ($alt$(\partial^+T))=w$}}}
\end{picture}

\]
\begin{center}
	Figure 2.7.4
\end{center}
\noindent That $G$ is a group follows
from Figure~2.7.5 
\[
	\begin{picture}(0,0)%
\epsfig{file=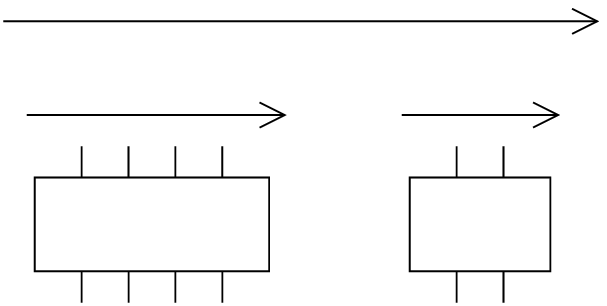}%
\end{picture}%
\setlength{\unitlength}{0.00041700in}%
\begingroup\makeatletter\ifx\SetFigFont\undefined
\def\x#1#2#3#4#5#6#7\relax{\def\x{#1#2#3#4#5#6}}%
\expandafter\x\fmtname xxxxxx\relax \def\y{splain}%
\ifx\x\y   
\gdef\SetFigFont#1#2#3{%
  \ifnum #1<17\tiny\else \ifnum #1<20\small\else
  \ifnum #1<24\normalsize\else \ifnum #1<29\large\else
  \ifnum #1<34\Large\else \ifnum #1<41\LARGE\else
     \huge\fi\fi\fi\fi\fi\fi
  \csname #3\endcsname}%
\else
\gdef\SetFigFont#1#2#3{\begingroup
  \count@#1\relax \ifnum 25<\count@\count@25\fi
  \def\x{\endgroup\@setsize\SetFigFont{#2pt}}%
  \expandafter\x
    \csname \romannumeral\the\count@ pt\expandafter\endcsname
    \csname @\romannumeral\the\count@ pt\endcsname
  \csname #3\endcsname}%
\fi
\fi\endgroup
\begin{picture}(5744,3236)(1179,-3384)
\put(2551,-2536){\makebox(0,0)[lb]{\smash{\SetFigFont{12}{14.4}{rm}$T$}}}
\put(5701,-2536){\makebox(0,0)[lb]{\smash{\SetFigFont{12}{14.4}{rm}$R$}}}
\put(1576,-1261){\makebox(0,0)[lb]{\smash{\SetFigFont{12}{14.4}{rm}$\phi($alt$(\partial^+T))$}}}
\put(5026,-1261){\makebox(0,0)[lb]{\smash{\SetFigFont{12}{14.4}{rm}$\phi$(alt$(\partial^+R))$}}}
\put(2626,-436){\makebox(0,0)[lb]{\smash{\SetFigFont{12}{14.4}{rm}$\phi($alt$(\partial^+T))\phi($alt$(\partial^+R))$}}}
\end{picture}

\]
\begin{center}
	Figure 2.7.5
\end{center}
\noindent Also $G$ contains $\tilde\D\tilde\D^{-1}$ since
$\bsp$ is in $P_2^{\G,\D}$ where
for $\g\in\D$, \ $\boxa{\g}$ is the ``diagonal" tensor
$\blong$. That $G=G'$ if $e\in\D$ is
easily seen by attaching $\boxb{e}$ to the right of basic tensors
in $P_k^{\G,\D}$ when $k$ is even. \qed

\bigskip
We see that, if $e\in\D$, a basis for $P^{\G,\D}$ is formed by all
random walks on $G$, starting and ending at $1\in\G$, where the odd
transitions correspond to multiplying by a $\tilde \d$ for each $\d\in\D$,
and the even ones by $\tilde\d^{-1}$ for $\d\in\D$. If
$\tilde\D=\tilde\D^{-1}$ and $\tilde{}$ is injective, these are
just random walks on the Cayley graph of $G$.

If $e\in\D$, each basic tensor $T\in P^{\G,\D}$ gives the {\it
relation} alt$(\p T)$ in $G$, thinking of $G$ as being presented on
$\D\backslash\{e\}$.

Suppose $G=\<\D\backslash\{e\}\mid
r_1,r_2,\dots\>$ is a presentation of $G$, i.e. the kernel of
the map induced by
$\sim$ from $F(\D\backslash\{e\})$ to $G$ is the normal closure of
the
$r_i$'s. Then each $r_i$ may be represented by a $k$-box, written
$\blankbox{r_i}$, with $\mu({\op{alt}}(\p(\blankbox{r_i})))=r_i$,
for some $k$ with $2k\geq\ell(r)$. (We use $\ell(w)$ to denote the
length of a word $w$.) To do this one may have to use $e\in\D$ so
that the word $r_i$ conforms with the alternating condition. For
instance to represent $\g\d^2\g^{-1}\d$ one might use the basic
tensor 
\fbox{\parbox{1.8cm}{
$\dsize{
\begin{array}{c}
\g\ e\ \d\ e\\
e\ \d\ \g\  \d 
\end{array}
        }$ 
}}.      

Let $\mu:F(\D)\to F(\D\backslash\{e\})$ be the homomorphism defined
by $\mu(e)=1\in F(\D\backslash\{e\})$, \ $\mu|_{\D\backslash\{e\}}=$
id.

\bigskip{\bf Definition.} Let $R=\bigcup^{\infty}_{k=0} R_k$
be the planar subalgebra of $P^{\G,\D}$ generated by
$\{ \boxf{\d}  :\d\in\D\}\cup \{ \blankbox{r}\}
\cup \{\blankbox{\hskip2pt r^{-1}_i}\}$.
Let $H=\{\mu({\op{alt}}(\p T))\mid T$ is a basic tensor in $R\}$.

{\bf Theorem 2.7.6} The set $H$ is a subgroup of
$F(\D\backslash\{e\})$ equal to the normal closure $N$ of $\{r_i\}$
in $F(\D\backslash\{e\})$. Moreover, $P^{G,\D}=R$.

\medskip
{\bf Proof}. That $H$ is multiplicatively closed follows from
Figure 2.7.7.
\[
        \begin{picture}(0,0)%
\epsfig{file=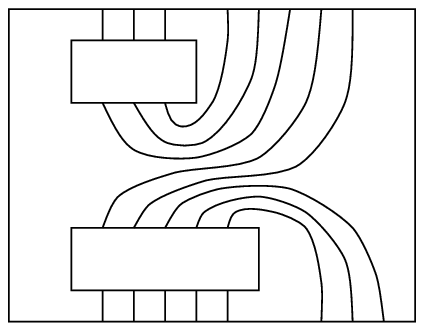}%
\end{picture}%
\setlength{\unitlength}{0.00041700in}%
\begingroup\makeatletter\ifx\SetFigFont\undefined
\def\x#1#2#3#4#5#6#7\relax{\def\x{#1#2#3#4#5#6}}%
\expandafter\x\fmtname xxxxxx\relax \def\y{splain}%
\ifx\x\y   
\gdef\SetFigFont#1#2#3{%
  \ifnum #1<17\tiny\else \ifnum #1<20\small\else
  \ifnum #1<24\normalsize\else \ifnum #1<29\large\else
  \ifnum #1<34\Large\else \ifnum #1<41\LARGE\else
     \huge\fi\fi\fi\fi\fi\fi
  \csname #3\endcsname}%
\else
\gdef\SetFigFont#1#2#3{\begingroup
  \count@#1\relax \ifnum 25<\count@\count@25\fi
  \def\x{\endgroup\@setsize\SetFigFont{#2pt}}%
  \expandafter\x
    \csname \romannumeral\the\count@ pt\expandafter\endcsname
    \csname @\romannumeral\the\count@ pt\endcsname
  \csname #3\endcsname}%
\fi
\fi\endgroup
\begin{picture}(5400,3044)(1201,-4284)
\put(3302,-1937){\makebox(0,0)[lb]{\smash{\SetFigFont{12}{14.4}{rm}$S$}}}
\put(3527,-3737){\makebox(0,0)[lb]{\smash{\SetFigFont{12}{14.4}{rm}$T$}}}
\put(1201,-2836){\makebox(0,0)[lb]{\smash{\SetFigFont{12}{14.4}{rm}$Q=$}}}
\put(6601,-2836){\makebox(0,0)[lb]{\smash{\SetFigFont{12}{14.4}{rm}alt($\partial Q)=$alt$(\partial S)$alt$(\partial T)$}}}
\end{picture}

\]
\begin{center}
Figure 2.7.7
\end{center}
\noindent To see that $H=H^{-1}$, note that the transpose of a
basic tensor gives the inverse boundary word, and a planar algebra
generated by a $*$-closed set of boxes is $*$-closed (note that the
box $\blankbox{\d}$ is self-adjoint). Figure 2.7.8 exhibits
conjugation of $\alpha={\op{alt}}(\p T)$ by $\g\d\g^{-1}$ for
$\g,\d\in\D$, which shows how to prove that $H$ is normal
\[
        \begin{picture}(0,0)%
\epsfig{file=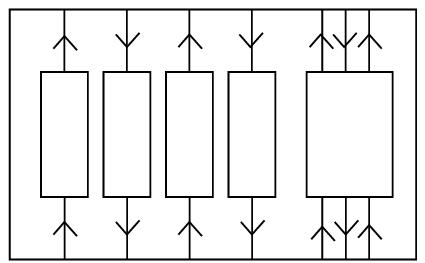}%
\end{picture}%
\setlength{\unitlength}{0.00041700in}%
\begingroup\makeatletter\ifx\SetFigFont\undefined
\def\x#1#2#3#4#5#6#7\relax{\def\x{#1#2#3#4#5#6}}%
\expandafter\x\fmtname xxxxxx\relax \def\y{splain}%
\ifx\x\y   
\gdef\SetFigFont#1#2#3{%
  \ifnum #1<17\tiny\else \ifnum #1<20\small\else
  \ifnum #1<24\normalsize\else \ifnum #1<29\large\else
  \ifnum #1<34\Large\else \ifnum #1<41\LARGE\else
     \huge\fi\fi\fi\fi\fi\fi
  \csname #3\endcsname}%
\else
\gdef\SetFigFont#1#2#3{\begingroup
  \count@#1\relax \ifnum 25<\count@\count@25\fi
  \def\x{\endgroup\@setsize\SetFigFont{#2pt}}%
  \expandafter\x
    \csname \romannumeral\the\count@ pt\expandafter\endcsname
    \csname @\romannumeral\the\count@ pt\endcsname
  \csname #3\endcsname}%
\fi
\fi\endgroup
\begin{picture}(5415,2444)(1201,-3983)
\put(1201,-2836){\makebox(0,0)[lb]{\smash{\SetFigFont{12}{14.4}{rm}$Q=$}}}
\put(2551,-2836){\makebox(0,0)[lb]{\smash{\SetFigFont{12}{14.4}{rm}$\gamma$}}}
\put(3181,-2866){\makebox(0,0)[lb]{\smash{\SetFigFont{12}{14.4}{rm}$e$}}}
\put(3736,-2851){\makebox(0,0)[lb]{\smash{\SetFigFont{12}{14.4}{rm}$\delta$}}}
\put(4321,-2851){\makebox(0,0)[lb]{\smash{\SetFigFont{12}{14.4}{rm}$\gamma$}}}
\put(5281,-2836){\makebox(0,0)[lb]{\smash{\SetFigFont{12}{14.4}{rm}$T$}}}
\put(6616,-3361){\makebox(0,0)[lb]{\smash{\SetFigFont{12}{14.4}{rm}$=\gamma\delta\gamma^{-1}$alt$(\partial Q)\gamma\delta^{-1}\gamma^{-1}$}}}
\put(6601,-2836){\makebox(0,0)[lb]{\smash{\SetFigFont{12}{14.4}{rm}$\mu$(alt($\partial T)$}}}
\end{picture}

\]
\begin{center}
Figure 2.7.8
\end{center}
\noindent Thus $H$ contains the normal closure $N$.

Now the tangle picture of an arbitrary basic tensor in $R$ can be
isotoped so that it is as in Figure 2.7.9.
\[
        \begin{picture}(0,0)%
\epsfig{file=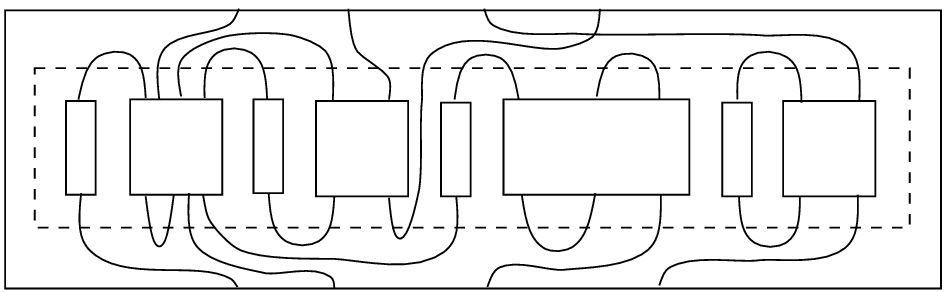}%
\end{picture}%
\setlength{\unitlength}{0.00041700in}%
\begingroup\makeatletter\ifx\SetFigFont\undefined
\def\x#1#2#3#4#5#6#7\relax{\def\x{#1#2#3#4#5#6}}%
\expandafter\x\fmtname xxxxxx\relax \def\y{splain}%
\ifx\x\y   
\gdef\SetFigFont#1#2#3{%
  \ifnum #1<17\tiny\else \ifnum #1<20\small\else
  \ifnum #1<24\normalsize\else \ifnum #1<29\large\else
  \ifnum #1<34\Large\else \ifnum #1<41\LARGE\else
     \huge\fi\fi\fi\fi\fi\fi
  \csname #3\endcsname}%
\else
\gdef\SetFigFont#1#2#3{\begingroup
  \count@#1\relax \ifnum 25<\count@\count@25\fi
  \def\x{\endgroup\@setsize\SetFigFont{#2pt}}%
  \expandafter\x
    \csname \romannumeral\the\count@ pt\expandafter\endcsname
    \csname @\romannumeral\the\count@ pt\endcsname
  \csname #3\endcsname}%
\fi
\fi\endgroup
\begin{picture}(9029,2722)(894,-3383)
\end{picture}

\]
\begin{center}
        Figure 2.7.9
\end{center}

\noindent This is a tangle $T$ all of whose curves are vertical
straight lines surrounded by an element of the category ${\cal A} (\emptyset)$. But it
is easy to see that applying an ${\cal A} (\emptyset)$ element changes alt$(\p T)$ at
most by a conjugation.

The most difficult part of Theorem 2.7.6 is to show that
$P^{\G,\D}=R$. We must show that if $w$ is a word of even length on
$\D$ with $\phi({\op{alt}}(w))=1$ then there is a basic tensor
$T\in R$ with $\p T=w$.

As a first step, observe that if $\mu$ is a homomorphism from
$F(\D)$ to $F(\D-\{e\})$ sending $e$ to the identity and with
$\mu(\d)=\d$ for $\d\neq e$, then if $w_1,w_2\in\frak S(\D)$ are of
even length and $\mu({\op{alt}}(w_1))=\mu({\op{alt}}(w_2))$ then
${\op{alt}}(w_1)={\op{alt}}(w_2)$. This is because
$w_1w^*_2$ ($w^*$ is $w$ written backwards) satisfies
${\op{alt}}(w_1w^*_2)={\op{alt}}(w_1){\op{alt}}(w_2)^{-1}$, thus
${\op{alt}}(w_1w_2^*)\in{\op{ker}} \ \mu$ which is the normal
closure of $e$. The length of $w_1w_2^*$ can be reduced
(if necessary) by eliminating consecutive letters two at a time to
obtain another word $w$, of even length, with
$\ell(w)=\ell({\op{alt}} \ w)$ ($\ell=$ length). By the uniqueness
of reduced words in a free group, $w$ must be a product of words of
the form $x \ e \ y$, which map to conjugates of
$e^{\pm 1}$ in $F(\D)$. But the last letter of $x$ and the first
letter of $y$ must then be the same, and alt will send both these
letters to the same free group element. Thus in the process of
reducing $w_1w^*_2$, all occurrences of $e$ must disappear and
${\op{alt}}(w_1)={\op{alt}}(w_2)$.

A consequence of this observation is that, if $T$ is a basic tensor
with $\phi({\op{alt}}(\p T))=1$ in $\G$ so that
$\mu({\op{alt}}(\p T))$ is in the normal closure of $\{r_i\}$ in
$F(\D\backslash\{e\})$, then alt$(\p T)$ is the normal closure of
$\{{\op{alt}}(\p( \blankbox{r_i}))\}$ in $F(\D)$. Thus it
suffices to show that, if $T_1$ and $T_2$ are basic tensors with
${\op{alt}}(\p T_1)={\op{alt}}(\p T_2)\in F(\D)$, then $T_1=cT_2$
for some $c$ in ${\cal A} (\emptyset)$ (since for any $x\in\frak
S(\D)$ with $\phi |\mu({\op{alt}}(x))=1$ we have shown there is a
$T$ in the planar algebra $R$ with alt$(\p T)={\op{alt}}(x)$).
But this is rather easy --- we may suppose without loss of
generality that no cancellation happens going from $\p T_1$ to
alt$(\p T_1)$ and then use induction on $\ell(\p T_2)$. If
$\ell(\p T_2)=\ell(\p T_1)$ then $\p T_2=\p T_1$. Otherwise there
must be a sequence $\dots \d\d\dots$ in $\p T_2$ for some
$\d\in\D$. Connecting $\d$ to $\d$ in the tangle reduces the length
of $\p T_2$ by 2, and the remaining region is a disc.

\qed
A.Casson has pointed out the connection between
$P^{\G,\D}$ and van Kampen diagrams.

\bigskip
If $e\in\D$, the dimension of $P^{\G,\D}_n$ is the number of ways
of writing $1\in\G$ as a product of elements $\tilde\d$,
$\d\in\D$, with alternating signs. In particular,
dim$(P^{\G,\D}_1)=|\D|^2$ iff $\G$ is trivial. Since the problem of
the triviality of a group with given presentation is undecidable,
we conclude the following.

\bs\ni{\bf Corollary 2.7.10} The calculation of the dimension of a
planar subalgebra of $P^{\otimes}$ is undecidable.

\bs
Since there are
groups which are finitely generated but not finitely presented we
have

\bs\ni{\bf Corollary 2.7.11} There are finite dimensional
planar $*$-algebras which are not finitely generated as planar
algebras.

\bs
{\sc Proof}. If finitely many linear combinations of basic tensors
generated a planar algebra, then certainly the basic tensors
involved would also. But by 2.7.6, the group would then be finitely
presented. \qed

\bs{\bf Example 2.8 \ Spin models.}
We give a general planar algebra that is not a planar algebra, although it
is the planar algebras associated with it that will be of most
interest. In some sense it is a ``square root" of the planar
algebra $P^{\otimes}$ of $\S$2.6.

Let $V$ be a vector space of dimension $Q$ with a basis indexed by
``spin states" $\{1,2,\dots Q\}$. For each odd $n$ let $P^\s_n$ be
the subalgebra of End$(V^{\otimes\frac{n+1}{2}})$ given by
End$(V^{\otimes\frac{n-1}{2}})\otimes\D$ where $\D$ is the
subalgebra of End$(V)$ consisting of linear maps diagonal with
respect to the basis. For $n=0$, $P_0$ is the field $K$ and for $n$
even, $P^\s_n={\op{End}}(V^{\otimes\frac{n}{2}})$. Elements of
$P^\s_n$ will be identified with functions from $\{1,2,\dots ,Q\}^n$
to $K$, the value of the function on $(i_1,i_2,\dots ,i_n)$ being
the coefficient of the basic tensor
\fbox{\parbox{4.3cm}{
$\dsize{
\begin{array}{clcr}
i_1 & i_2  \ldots  \ldots i_{m-1} & i_m \\
i_n & i_{n-1}  \ldots .i_{m+2} & i_{m+1} 
\end{array}
        }$ 
}}      

for
$n=2m$, and the coefficient of
\fbox{\parbox{4.3cm}{
$\dsize{
\begin{array}{clcr}
i_1  & i_2  \ldots \ldots i_{m-1} & i_m \\
i_n &  i_{n-1}  \ldots .i_{m+2} & i_{m}      
\end{array}
        }$
}}

for
$n=2m-1$. (See $\S$2.6 for notation.) We shall make $P^\s$ into a
general planar algebra in two slightly different ways. In both cases the
labelling set will be $P^\s$ itself.

\bs{\bf First planar structure on $P^\s$.}
Take a tangle $T$ in ${\cal P}_k(L)$. We will define
$\Phi_0(T)\in{\cal P}^\s_k$ as follows.

First, shade the connected components of ${\cal B}_k\backslash T$
(called regions) black and white so that the region containing a
neighborhood of (0,0) is white, and so that regions whose closures
intersect (i.e.~which share an edge) have different colours. In
other words, regions whose boundary induces the positive
orientation of $\Bbb R^2$ are coloured white and negatively oriented
ones are black. Observe that the top and bottom of ${\cal B}_k$
consists of segments of length one forming parts of the boundaries
of regions alternately coloured white and black. If $k$ is odd, the
right-hand boundary of ${\cal B}_k$ can be joined with the rightmost
top and bottom segments to form part of the boundary of a black
region.  This way the boundary of ${\cal B}_k$ always has $k$
segments attached to black regions whose closure meets the
boundary. Number these segments cyclically $1,2,\dots ,k$ starting
form the top left and going clockwise.  To define an element of
$P^\s_k$ from $T$ we must give a function 
$\Phi_0(T):\{1,2,\dots Q\}^k\to K$. It is
\[
\Phi_0(T)(i_1,i_2,\dots i_k)=\sum_\s \
\prod_{B\in\left\{
{\text{labelled\  boxes}}\atop  {\text{of }} T\right\}} \s (B)
\]
where $\s$ runs over all functions from the black regions of $T$ to
$\{1,2,\dots ,Q\}$ which take the value $i_j$ on the black region
whose closure contains the $j^{\text{th}}$ boundary segment, for
all $j=1,2,\dots ,k$. Given a labeled box $B$ of $T$, and such a
$\s$, the boundary segments of $B$ which meet closures of black
regions are numbered 1 to $k_B$ so $\s$ defines an element of
$\{1,2,\dots ,Q\}^{k_B}$, and thus the label of $B$ gives a scalar
$\s(B)$ in $K$. As usual empty sums are zero and empty products are
1. This completes the definition of $\Phi_0$ and it is easily
checked that $\Phi_0$ presents $P^\s$ as a planar algebra. The
induced representation of ${\cal P}(\phi)$ gives a representation of
TL with $\d_1\!=\!Q$, $\d_2\!=\!1$. It is precisely the
representation associated with the Potts model used by Temperley
and Lieb in [TL].

 To be sure of relevance to subfactors, we
now show how to adjust these parameters so that
$\d_1=\d_2=\sqrt{Q}$.

\bs{\bf Second planar structure on $P^\s$.}
If $T$ is a labeled tangle in ${\cal P}_k(L)$, we define
a tangle $\tilde T$ in ${\cal P}_k(\phi)$ by ``smoothing" all the
boxes of $\tilde T$, i.e.~replacing $\boxe{R}$ by
$\uparrow \ \downarrow \ \uparrow$, and shrinking all
non-through-strings to semicircles near the top or bottom of ${\cal B}_k$. 
Put $f(T)=Q^{\frac 12(n_+-n_-)+\frac 14(n^\p_+-n^\p_-)}$
where $n_+$ and $n_-$ are the numbers of positively and negatively
oriented circles in $\tilde T$ respectively, and $n_{\pm}^\p$ are
similarly the numbers of positively and negatively oriented
semicircles near the top and bottom.  Thus defined, $f(T)$ is
clearly an isotopy invariant, so we could redefine it by assuming
all the boxes are parallel to the $x$-axis. Assuming all maxima and
minima of the $y$-coordinate restricted to the strings of $T$ are
nondegenerate, $2(n_+-n_-)+(n^\p_+-n^\p_-)$ is just
$p_++q_+-p_--q_-$ where $p_+,p_-$ are the numbers of local maxima
of $y$ oriented to the left and right respectively and similarly
$q_+$ and $q_-$ count minima to the right and left respectively.
It follows that $T\to f(T)$ is multiplicative and indeed that if $A$
is in the annular category one may define $f(A)$ so that
$f(AT)=f(A)f(T)$.

The normalisation constant
$n_+-n_-+\frac{1}{2}(n_+^{\partial}-n_-^{\partial})$
may seem mysteious. What is actually being calculated is
the isotopy invariant $\int d\theta$
where the intgral is taken over the strings of the tangle
and $d\theta$ is the change of angle or curvature 1-form, normalised so that
integrating over a positively oriented circle counts one.
The above factor is then this integral when all strings meet
all boxes at right angles. Thus, by following shaded regions at every internal
box, another formula for this normalisation factor for a $k-$tangle is

\[
        [\frac{k+1}{2}]-b-\sum_{i\geq 1}[\frac{n_i}{2}].
\]

If all boxes are 2- or 3-boxes we get
\[
        \frac{1}{2}(\#(\mbox{black regions})-\#(\mbox{boxes}))+
        \frac{1}{4}(n_+^{\partial}-n_-^{\partial}),
\]
where now $n_+^{\partial}$ and $n_-^{\partial}$ are calculated
by eliminating the boxes by following the black regions rather
than going straight through the box.

\bs
\ni{\bf Proposition 2.8.1} The map $\Phi^\s:{\cal P}(L)\to P^\s$,
$\Phi(T)=f(T)\Phi_0(T)$ (linearly extended) presents $P^\s$ as a
general planar algebra, $\Phi |_{{\cal P}(\phi)}$  presents $TL$ with
$\d_1=\d_2=\sqrt{Q}$.

\bs
{\sc Proof.} Annular invariance follows from the relation
$f(AT)=f(A)f(T)$ and the annular invariance of $\Phi_0$.
If $\deye$ is a part of a tangle $T$ then
$\Phi(T)=QQ^{-\frac 12} \Phi(\tilde T)$ where $\deye$
has been removed from $T$. If $\deyer$\ \ \  is part of $T$,
$\Phi(T)=Q^{\frac 12} \Phi(\tilde T)$. \qed

\bs
When we refer to $P^\s$, we will mean $P^\s$,
together with $\Phi$.

Although dim $P^\s_0=1$, so that $P^\s$ gives an invariant of
labeled planar networks with unbounded region positively oriented,
dim$(P^\s_{1,1})=Q$ so $P^\s$ is not planar. However, $P^\s$
does have the obvious $*$ structure and tr${}_R$ is defined and
positive definite, so that any self-adjoint planar
subalgebra of
$P^\s$ will be a $C^*$-planar algebra.

\bs\ni
{\bf Proposition 2.8.2} A planar subalgebra $P$ of
$P^\s$ is spherical.

\bs
{\sc Proof.} Given a planar network $N$ in $P$ with positively
oriented unbounded region, we need only show that
\[
        \begin{picture}(0,0)%
\epsfig{file=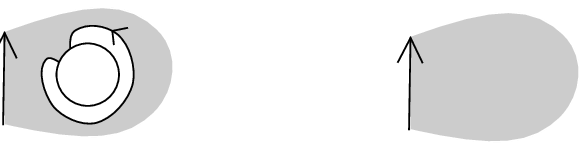}%
\end{picture}%
\setlength{\unitlength}{0.00041700in}%
\begingroup\makeatletter\ifx\SetFigFont\undefined
\def\x#1#2#3#4#5#6#7\relax{\def\x{#1#2#3#4#5#6}}%
\expandafter\x\fmtname xxxxxx\relax \def\y{splain}%
\ifx\x\y   
\gdef\SetFigFont#1#2#3{%
  \ifnum #1<17\tiny\else \ifnum #1<20\small\else
  \ifnum #1<24\normalsize\else \ifnum #1<29\large\else
  \ifnum #1<34\Large\else \ifnum #1<41\LARGE\else
     \huge\fi\fi\fi\fi\fi\fi
  \csname #3\endcsname}%
\else
\gdef\SetFigFont#1#2#3{\begingroup
  \count@#1\relax \ifnum 25<\count@\count@25\fi
  \def\x{\endgroup\@setsize\SetFigFont{#2pt}}%
  \expandafter\x
    \csname \romannumeral\the\count@ pt\expandafter\endcsname
    \csname @\romannumeral\the\count@ pt\endcsname
  \csname #3\endcsname}%
\fi
\fi\endgroup
\begin{picture}(5578,1341)(1175,-1455)
\put(1917,-901){\makebox(0,0)[lb]{\smash{\SetFigFont{12}{14.4}{rm}$N$}}}
\put(3207,-901){\makebox(0,0)[lb]{\smash{\SetFigFont{12}{14.4}{rm}$=$}}}
\put(4032,-931){\makebox(0,0)[lb]{\smash{\SetFigFont{12}{14.4}{rm}$Z(N)$}}}
\end{picture}

\]
\ni But since $P$ is planar, the sum over all internal spins in
\begin{picture}(0,0)%
\epsfig{file=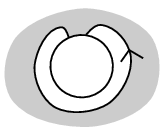}%
\end{picture}%
\setlength{\unitlength}{0.00041700in}%
\begingroup\makeatletter\ifx\SetFigFont\undefined
\def\x#1#2#3#4#5#6#7\relax{\def\x{#1#2#3#4#5#6}}%
\expandafter\x\fmtname xxxxxx\relax \def\y{splain}%
\ifx\x\y   
\gdef\SetFigFont#1#2#3{%
  \ifnum #1<17\tiny\else \ifnum #1<20\small\else
  \ifnum #1<24\normalsize\else \ifnum #1<29\large\else
  \ifnum #1<34\Large\else \ifnum #1<41\LARGE\else
     \huge\fi\fi\fi\fi\fi\fi
  \csname #3\endcsname}%
\else
\gdef\SetFigFont#1#2#3{\begingroup
  \count@#1\relax \ifnum 25<\count@\count@25\fi
  \def\x{\endgroup\@setsize\SetFigFont{#2pt}}%
  \expandafter\x
    \csname \romannumeral\the\count@ pt\expandafter\endcsname
    \csname @\romannumeral\the\count@ pt\endcsname
  \csname #3\endcsname}%
\fi
\fi\endgroup
\begin{picture}(1543,1191)(1214,-1327)
\put(1842,-811){\makebox(0,0)[lb]{\smash{\SetFigFont{12}{14.4}{rm}$N$}}}
\end{picture}

is independent of the spin value in the unbounded
region, and each term in the sum for 
\begin{picture}(0,0)%
\epsfig{file=ps/pic37.ps}%
\end{picture}%
\setlength{\unitlength}{0.00041700in}%
\begingroup\makeatletter\ifx\SetFigFont\undefined
\def\x#1#2#3#4#5#6#7\relax{\def\x{#1#2#3#4#5#6}}%
\expandafter\x\fmtname xxxxxx\relax \def\y{splain}%
\ifx\x\y   
\gdef\SetFigFont#1#2#3{%
  \ifnum #1<17\tiny\else \ifnum #1<20\small\else
  \ifnum #1<24\normalsize\else \ifnum #1<29\large\else
  \ifnum #1<34\Large\else \ifnum #1<41\LARGE\else
     \huge\fi\fi\fi\fi\fi\fi
  \csname #3\endcsname}%
\else
\gdef\SetFigFont#1#2#3{\begingroup
  \count@#1\relax \ifnum 25<\count@\count@25\fi
  \def\x{\endgroup\@setsize\SetFigFont{#2pt}}%
  \expandafter\x
    \csname \romannumeral\the\count@ pt\expandafter\endcsname
    \csname @\romannumeral\the\count@ pt\endcsname
  \csname #3\endcsname}%
\fi
\fi\endgroup
\begin{picture}(1543,1191)(1214,-1327)
\put(1842,-811){\makebox(0,0)[lb]{\smash{\SetFigFont{12}{14.4}{rm}$N$}}}
\end{picture}

is $Q$ times the
corresponding term for
\begin{picture}(0,0)%
\epsfig{file=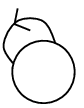}%
\end{picture}%
\setlength{\unitlength}{0.00041700in}%
\begingroup\makeatletter\ifx\SetFigFont\undefined
\def\x#1#2#3#4#5#6#7\relax{\def\x{#1#2#3#4#5#6}}%
\expandafter\x\fmtname xxxxxx\relax \def\y{splain}%
\ifx\x\y   
\gdef\SetFigFont#1#2#3{%
  \ifnum #1<17\tiny\else \ifnum #1<20\small\else
  \ifnum #1<24\normalsize\else \ifnum #1<29\large\else
  \ifnum #1<34\Large\else \ifnum #1<41\LARGE\else
     \huge\fi\fi\fi\fi\fi\fi
  \csname #3\endcsname}%
\else
\gdef\SetFigFont#1#2#3{\begingroup
  \count@#1\relax \ifnum 25<\count@\count@25\fi
  \def\x{\endgroup\@setsize\SetFigFont{#2pt}}%
  \expandafter\x
    \csname \romannumeral\the\count@ pt\expandafter\endcsname
    \csname @\romannumeral\the\count@ pt\endcsname
  \csname #3\endcsname}%
\fi
\fi\endgroup
\begin{picture}(713,945)(1587,-1029)
\put(1842,-811){\makebox(0,0)[lb]{\smash{\SetFigFont{12}{14.4}{rm}$N$}}}
\end{picture}
. Taking the sum over all $Q$
spin states in the shaded region we are done. \qed

\bs There are ways to obtain planar subalgebras of
$P^\s$. An obvious place to look is association schemes 
where one is given a family of (0,1) $Q\times Q$ matrices $A_i$,
$i=1,\dots ,d$, whose linear span is closed under the operations
given by the tangles
\medskip
\begin{picture}(0,0)%
\epsfig{file=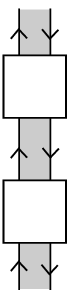}%
\end{picture}%
\setlength{\unitlength}{0.00041700in}%
\begingroup\makeatletter\ifx\SetFigFont\undefined
\def\x#1#2#3#4#5#6#7\relax{\def\x{#1#2#3#4#5#6}}%
\expandafter\x\fmtname xxxxxx\relax \def\y{splain}%
\ifx\x\y   
\gdef\SetFigFont#1#2#3{%
  \ifnum #1<17\tiny\else \ifnum #1<20\small\else
  \ifnum #1<24\normalsize\else \ifnum #1<29\large\else
  \ifnum #1<34\Large\else \ifnum #1<41\LARGE\else
     \huge\fi\fi\fi\fi\fi\fi
  \csname #3\endcsname}%
\else
\gdef\SetFigFont#1#2#3{\begingroup
  \count@#1\relax \ifnum 25<\count@\count@25\fi
  \def\x{\endgroup\@setsize\SetFigFont{#2pt}}%
  \expandafter\x
    \csname \romannumeral\the\count@ pt\expandafter\endcsname
    \csname @\romannumeral\the\count@ pt\endcsname
  \csname #3\endcsname}%
\fi
\fi\endgroup
\begin{picture}(644,2744)(1179,-2933)
\end{picture}

(matrix multiplication) and
\begin{picture}(0,0)%
\epsfig{file=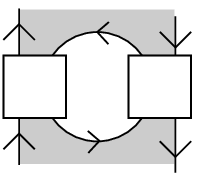}%
\end{picture}%
\setlength{\unitlength}{0.00041700in}%
\begingroup\makeatletter\ifx\SetFigFont\undefined
\def\x#1#2#3#4#5#6#7\relax{\def\x{#1#2#3#4#5#6}}%
\expandafter\x\fmtname xxxxxx\relax \def\y{splain}%
\ifx\x\y   
\gdef\SetFigFont#1#2#3{%
  \ifnum #1<17\tiny\else \ifnum #1<20\small\else
  \ifnum #1<24\normalsize\else \ifnum #1<29\large\else
  \ifnum #1<34\Large\else \ifnum #1<41\LARGE\else
     \huge\fi\fi\fi\fi\fi\fi
  \csname #3\endcsname}%
\else
\gdef\SetFigFont#1#2#3{\begingroup
  \count@#1\relax \ifnum 25<\count@\count@25\fi
  \def\x{\endgroup\@setsize\SetFigFont{#2pt}}%
  \expandafter\x
    \csname \romannumeral\the\count@ pt\expandafter\endcsname
    \csname @\romannumeral\the\count@ pt\endcsname
  \csname #3\endcsname}%
\fi
\fi\endgroup
\begin{picture}(1845,1618)(1179,-2107)
\end{picture}

(Hadamard product).
\medskip\ni
The requirement that this linear span (the
``Bose-Mesner algebra") be closed under all planar contractions is
presumably much more stringent. If the requirement is satisfied, and the row and
column sums of each $A_i$ do not depend on the row or column, we
will have a planar subalgebra of
$P^\s$. A particularly simple example of this comes from transitive
actions of a finite group $G$ on a set $S$. Then the orbits of $G$
on $S\times S$ define an association scheme whose Bose-Mesner
algebra is the fixed points for the action on $M_{|S|}(\Bbb C)$ by
conjugation. We get a planar subalgebra of $P^\s$ either
by taking the fixed points for the $G$-action on $P^\s$ or the
planar subalgebra generated by the assocation scheme. They are
different in general. A case where they are the same is for the
dihedral group on a set with five elements (see [J4]). They
are different for Jaeger's Higman-Sims model ([Ja],[dlH])
 --- although the dimensions of the two planar algebras agree for
a while, they have different asymptotic growth rates, one being that
of the commutant of Sp(4) on $(\Bbb C^4)^{\otimes k}$ and the other
being $100^k$. 

Here is an interesting example for a doubly transitive group.
It connects with Example 2.5 and gives a new kind of ``spin model"
for link invariants from links projected with only triple point
singularities.

The alternating group $A_4$ is doubly transitive on the set
$\{1,2,3,4\}$ but there are two orbits on the set of ordered
triples $(a,b,c)$ of distinct elements according to whether
$1\mapsto a$, $2\mapsto b$, $3\mapsto c$, $4\mapsto d$
(with $\{a,b,c,d\}=\{1,2,3,4\}$) is an even or odd permutation.
Let \ $\boxe{e}\in{\cal P}_3(L)$ be such that
$\Phi^\s(\boxe{e})$ is the characteristic
function  of the even orbit. Define a mapping from
${\cal P}(\ \crsst\ )$ (the universal planar algebra on a single 3-box)
to $P^\s$ by sending $\crsst$ to
$\boxe{e}-\frac 12 \vln\ \  \smth$.
It is possible to prove that this map passes to the quotient
$P^H$ (the planar algebra of 2.5) with parameters  $t=i=x$.
This is equivalent to showing that twice the value of the {\sc
homfly} polynomial of a link obtained by connecting 3-boxes
$\crsst$ (at 1,--1 in $\ell-m$ variables) in an oriented way is
the partition function in $P^\s$ (with $Q=4$) given by filling the
same three boxes with $\boxe{e}-\frac 12 \vln\ \ \smth$. 
We give a sample calculation below
which illustrates all the considerations. Note that, for $t=i=x$,
the value of a single circle in the {\sc homfly} skein is 2.
\[
        \begin{picture}(0,0)%
\epsfig{file=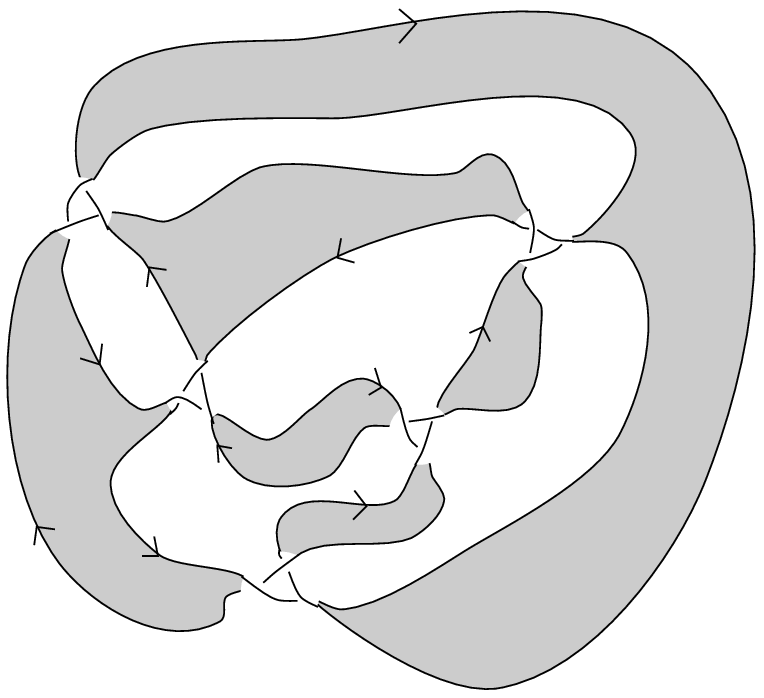}%
\end{picture}%
\setlength{\unitlength}{0.00041700in}%
\begingroup\makeatletter\ifx\SetFigFont\undefined
\def\x#1#2#3#4#5#6#7\relax{\def\x{#1#2#3#4#5#6}}%
\expandafter\x\fmtname xxxxxx\relax \def\y{splain}%
\ifx\x\y   
\gdef\SetFigFont#1#2#3{%
  \ifnum #1<17\tiny\else \ifnum #1<20\small\else
  \ifnum #1<24\normalsize\else \ifnum #1<29\large\else
  \ifnum #1<34\Large\else \ifnum #1<41\LARGE\else
     \huge\fi\fi\fi\fi\fi\fi
  \csname #3\endcsname}%
\else
\gdef\SetFigFont#1#2#3{\begingroup
  \count@#1\relax \ifnum 25<\count@\count@25\fi
  \def\x{\endgroup\@setsize\SetFigFont{#2pt}}%
  \expandafter\x
    \csname \romannumeral\the\count@ pt\expandafter\endcsname
    \csname @\romannumeral\the\count@ pt\endcsname
  \csname #3\endcsname}%
\fi
\fi\endgroup
\begin{picture}(7277,6618)(88,-5970)
\put(676,-3811){\makebox(0,0)[lb]{\smash{\SetFigFont{12}{14.4}{rm}$b$}}}
\put(2881,-3676){\makebox(0,0)[lb]{\smash{\SetFigFont{12}{14.4}{rm}$c$}}}
\put(3166,-4411){\makebox(0,0)[lb]{\smash{\SetFigFont{12}{14.4}{rm}$c$}}}
\put(5086,-5266){\makebox(0,0)[lb]{\smash{\SetFigFont{12}{14.4}{rm}$a$}}}
\put(4846,-2701){\makebox(0,0)[lb]{\smash{\SetFigFont{12}{14.4}{rm}$c$}}}
\put(2461,-1696){\makebox(0,0)[lb]{\smash{\SetFigFont{12}{14.4}{rm}$d$}}}
\end{picture}

\]
\begin{center}
        Figure 2.8.3
\end{center}

Smoothing all the 3-boxes leads to a single negatively oriented
circle so we must divide the final partition function by 2.
Replacing the 3-boxes by $\boxe{e}-\frac12 \vln\ \  \smth$ we look for spin states, i.e. functions
from the shaded regions to $\{1,2,3,4\}$ for which each 3-box
yields a non-zero contribution to the partition function. Around
each 3-box this means that either the three spin values are in the
even orbit under
$A_4$, or they are all the same. The first case contributes +1  to
the product over boxes, the second case contributes --1 ({\bf not}
$-\frac 12$ because of the maxima and minima in the box). If the box labeled
($\dagger$) is surrounded by the same spin value, all the spin
states must be the same for a nonzero contribution to $Z$. This
gives a factor $4\times (-1)^5$. On the other hand, if the spins at
($\dagger$) are as in Figure 2.8.3 with $(a,b,c)$ in the even
orbit, the other spin
choices are forced (where $\{a,b,c,d\}=\{1,2,3,4\}$), for a
contribution of --1. The orbit is of size 12 so the partition
function is $\frac 12(-12-4)=-8$. For this link the value of the
{\sc homfly} polynomial $P_L(1,-1)$ is --4. The factor of 2 is
accounted for by the fact that our partition function is 2 on the
unknot. Thus our answer is correct. Note how few spin patterns
actually contributed to $Z$!

If we wanted to use non-alternating 3-boxes we could simply use the
{\sc homfly} skein relation to modify the 3-box. For instance
\[
        \begin{picture}(0,0)%
\epsfig{file=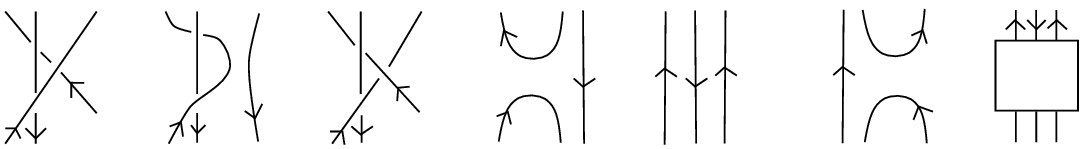}%
\end{picture}%
\setlength{\unitlength}{0.00041700in}%
\begingroup\makeatletter\ifx\SetFigFont\undefined
\def\x#1#2#3#4#5#6#7\relax{\def\x{#1#2#3#4#5#6}}%
\expandafter\x\fmtname xxxxxx\relax \def\y{splain}%
\ifx\x\y   
\gdef\SetFigFont#1#2#3{%
  \ifnum #1<17\tiny\else \ifnum #1<20\small\else
  \ifnum #1<24\normalsize\else \ifnum #1<29\large\else
  \ifnum #1<34\Large\else \ifnum #1<41\LARGE\else
     \huge\fi\fi\fi\fi\fi\fi
  \csname #3\endcsname}%
\else
\gdef\SetFigFont#1#2#3{\begingroup
  \count@#1\relax \ifnum 25<\count@\count@25\fi
  \def\x{\endgroup\@setsize\SetFigFont{#2pt}}%
  \expandafter\x
    \csname \romannumeral\the\count@ pt\expandafter\endcsname
    \csname @\romannumeral\the\count@ pt\endcsname
  \csname #3\endcsname}%
\fi
\fi\endgroup
\begin{picture}(10342,1339)(130,-1341)
\put(1186,-811){\makebox(0,0)[lb]{\smash{\SetFigFont{12}{14.4}{rm}$=$}}}
\put(4306,-781){\makebox(0,0)[lb]{\smash{\SetFigFont{12}{14.4}{rm}$=$}}}
\put(2791,-781){\makebox(0,0)[lb]{\smash{\SetFigFont{12}{14.4}{rm}$-$}}}
\put(5956,-826){\makebox(0,0)[lb]{\smash{\SetFigFont{12}{14.4}{rm}$-$}}}
\put(9196,-811){\makebox(0,0)[lb]{\smash{\SetFigFont{12}{14.4}{rm}$-$}}}
\put(9991,-751){\makebox(0,0)[lb]{\smash{\SetFigFont{12}{14.4}{rm}$e$}}}
\put(7291,-811){\makebox(0,0)[lb]{\smash{\SetFigFont{12}{14.4}{rm}$+$}}}
\put(7876,-811){\makebox(0,0)[lb]{\smash{\SetFigFont{12}{14.4}{rm}$\frac{1}{2}$}}}
\end{picture}

\]
\ni In general by [LM], $P_L(1,-1)$ is $(-1)^{c-1}(-2)^{\frac
12 d}$ where $c$ is the number of components of $L$ and $d$ is the
dimension of the first homology group (with $\Bbb Z/2\Bbb Z$
coefficients) of the triple branched cover of $S^3$, branched over
$L$. It would be reassuring to be able to see directly why our
formula gives this value. This would also prove directly that the
map $\crsst\mapsto \boxe{e}-\frac{1}{2} \vln\ \ \smth$
passes to the {\sc homfly} quotient.
Our proof of this is a little indirect --- one shows that the
planar subalgebras $P^H$ and $(P^\s)^{A_4}$ are the same by showing
they arise as centralizer towers from the same subfactor
(constructed in [GHJ]). Thus there must be a 3-box
corresponding to: $\crsst$  and we obtained the explicit
expression for it by solving an obvious set of equations.

As far as we know, this is the first genuine ``3-spin interaction"
statistical mechanical model for a link invariant.
Of course one may produce 3-spin interaction models by taking a
2-spin one and summing over the internal spin $\s$ in the picture
\[
        \begin{picture}(0,0)%
\epsfig{file=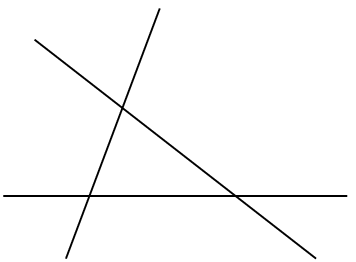}%
\end{picture}%
\setlength{\unitlength}{0.00041700in}%
\begingroup\makeatletter\ifx\SetFigFont\undefined
\def\x#1#2#3#4#5#6#7\relax{\def\x{#1#2#3#4#5#6}}%
\expandafter\x\fmtname xxxxxx\relax \def\y{splain}%
\ifx\x\y   
\gdef\SetFigFont#1#2#3{%
  \ifnum #1<17\tiny\else \ifnum #1<20\small\else
  \ifnum #1<24\normalsize\else \ifnum #1<29\large\else
  \ifnum #1<34\Large\else \ifnum #1<41\LARGE\else
     \huge\fi\fi\fi\fi\fi\fi
  \csname #3\endcsname}%
\else
\gdef\SetFigFont#1#2#3{\begingroup
  \count@#1\relax \ifnum 25<\count@\count@25\fi
  \def\x{\endgroup\@setsize\SetFigFont{#2pt}}%
  \expandafter\x
    \csname \romannumeral\the\count@ pt\expandafter\endcsname
    \csname @\romannumeral\the\count@ pt\endcsname
  \csname #3\endcsname}%
\fi
\fi\endgroup
\begin{picture}(3344,2444)(1479,-3683)
\put(2326,-1711){\makebox(0,0)[lb]{\smash{\SetFigFont{12}{14.4}{rm}$\sigma_1$}}}
\put(1651,-3436){\makebox(0,0)[lb]{\smash{\SetFigFont{12}{14.4}{rm}$\sigma3$}}}
\put(2626,-2761){\makebox(0,0)[lb]{\smash{\SetFigFont{12}{14.4}{rm}$\sigma$}}}
\put(4426,-3436){\makebox(0,0)[lb]{\smash{\SetFigFont{12}{14.4}{rm}$\sigma_2$}}}
\end{picture}

\]
\ni but that is of little interest. One may check quite easily that
the above model does not factorize in this way.

\bs\ni{\bf Example 2.9. Finite groups}

A special case of the subalgebra of $P^\s$ of Example 2.8 where a
group acts on a set $S$ is where $S=G$ itself, the action being
left multiplication. The algebra $P_2$ is then the group algebra
$\Bbb C G$, linearly spanned by elements $\boxg{g}$ for $g$
in $G$, which are by definition the matrices
$g^a_b=\d_{ag,b}$. They have the properties
\[
        \input{xfig/pic43}
\]
\ni which can be used to present $P$ as a planar algebra (see
[La]). In particular $P_{1,3}=\ell^{\i}(G)$ and $P_2=M_{|G|}(\Bbb
C)$. It is obvious that $G$ can be recovered from the abstract
planar algebra $P$ by using the minimal idempotents of $P_{1,3}$.

A noncommutative, finite-dimensional Hopf algebra  gives a planar
algebra but it cannot be a subalgebra of $P^\s$.

\bs{\bf Example 2.10 \ Invariant planar algebras.} Given an
invertible element $u\in P_1$ in a general planar algebra $P$ we
define $u^{\otimes k}$ to be the element of $P_k$ defined by the
following $k$-tangle
\[
        \begin{picture}(0,0)%
\epsfig{file=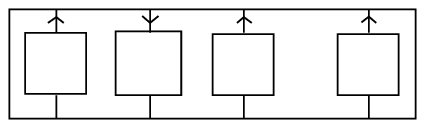}%
\end{picture}%
\setlength{\unitlength}{0.00041700in}%
\begingroup\makeatletter\ifx\SetFigFont\undefined
\def\x#1#2#3#4#5#6#7\relax{\def\x{#1#2#3#4#5#6}}%
\expandafter\x\fmtname xxxxxx\relax \def\y{splain}%
\ifx\x\y   
\gdef\SetFigFont#1#2#3{%
  \ifnum #1<17\tiny\else \ifnum #1<20\small\else
  \ifnum #1<24\normalsize\else \ifnum #1<29\large\else
  \ifnum #1<34\Large\else \ifnum #1<41\LARGE\else
     \huge\fi\fi\fi\fi\fi\fi
  \csname #3\endcsname}%
\else
\gdef\SetFigFont#1#2#3{\begingroup
  \count@#1\relax \ifnum 25<\count@\count@25\fi
  \def\x{\endgroup\@setsize\SetFigFont{#2pt}}%
  \expandafter\x
    \csname \romannumeral\the\count@ pt\expandafter\endcsname
    \csname @\romannumeral\the\count@ pt\endcsname
  \csname #3\endcsname}%
\fi
\fi\endgroup
\begin{picture}(5326,1094)(1951,-2932)
\put(3708,-2446){\makebox(0,0)[lb]{\smash{\SetFigFont{12}{14.4}{rm}$u$}}}
\put(5508,-2459){\makebox(0,0)[lb]{\smash{\SetFigFont{12}{14.4}{rm}$u$}}}
\put(4951,-2228){\makebox(0,0)[lb]{\smash{
\SetFigFont{12}{14.4}{rm}$u$
}}}
\put(6708,-2459){\makebox(0,0)[lb]{\smash{\SetFigFont{12}{14.4}{rm}$u$}}}
\put(5957,-2505){\makebox(0,0)[lb]{\smash{\SetFigFont{12}{14.4}{rm}$\cdots$}}}
\put(1951,-2536){\makebox(0,0)[lb]{\smash{\SetFigFont{12}{14.4}{rm}$u^{\otimes k}=$}}}
\end{picture}

\]
($k$ is odd in the picture). The $u$'s and upside down $u^{-1}$'s alternate.

\medskip
{\bf Proposition 2.10.1} If $P$ is a general planar algebra
and $S$ is a set of invertible elements of $P_1$, set
$$
P^S_k=\{x\in P_k\mid u^{\otimes k}x=xu^{\otimes k} \ \forall \
u\in S\} \ .
$$
Then $P^S$ is a general planar subalgebra of $P$ (planar if $P$ is)
and a $*$-planar algebra if $P$ is, and $S=S^*$.

\bigskip
{\sc Proof.} That $P^S$ is a unital filtered subalgebra is
obvious. In the $*$ case note that
$(u^{\otimes k})^*=(u^*)^{\otimes k}$. So by Lemma 1.18 we only
have to check invariance under ${\cal A}(\phi)$. Given an
$A\in{\cal A}_{k,n}(\phi)$ consider the tangle representing
$(\otimes^nu)\pi_A(x)(\otimes ^nu)^{-1}$ in Figure~\ref{pic66}

\begin{figure}[h]
\[
        \input{xfig/pic66}
\]
\caption[2.10.2]{}\label{pic66}
\end{figure}
\ni Each string of $\otimes^ku$, at the top and bottom, either
connects to another external boundary point, in which case the $u$
cancels with $u^{-1}$, or it connects to an internal boundary
point of the annulus. Isotoping each such $u$ and $u^{-1}$ close to
the internal boundary and inserting cancelling pairs of $u$ and
$u^{-1}$ on strings connecting internal boundary points, we see
that the tangle of Figure 2.10.2 gives the same element of $P$ as
the one where the only instances of $u$ and $u^{-1}$ surround the
$\blankbox{x}$ in an alternating fashion.
Since $x\in P^S$, these $u$'s may be eliminated and we are left
with $\pi_A(x)$. \qed

\bs
This gives a useful way of constructing planar subalgebras. In
particular $P^S_1$ is the commutant of $S$ in $P_1$ which may be
much smaller than $P_1$. Of special interest is the case where
$P=P^{\otimes}$ and $S$ is a subgroup $G$ of the unitary group.
In this case $P^G_k={\op{End}}_G(V\otimes V^*\otimes V\otimes
V^*\otimes \dots)$ ($k$ copies of $V$ or $V^*$) where $V^*$
represents the contragredient representation of $G$. Other cases of
interest can be constructed by cabling as in $\S$3 and then picking
some set of invertible elements in the original  $P_n$.

To obtain a more general construction one may replace
$\otimes^nu$ with tangles of the form
\[
        \begin{picture}(0,0)%
\epsfig{file=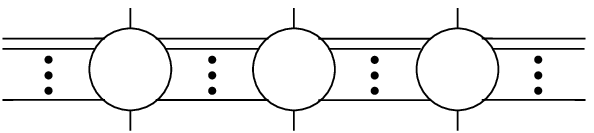}%
\end{picture}%
\setlength{\unitlength}{0.00041700in}%
\begingroup\makeatletter\ifx\SetFigFont\undefined
\def\x#1#2#3#4#5#6#7\relax{\def\x{#1#2#3#4#5#6}}%
\expandafter\x\fmtname xxxxxx\relax \def\y{splain}%
\ifx\x\y   
\gdef\SetFigFont#1#2#3{%
  \ifnum #1<17\tiny\else \ifnum #1<20\small\else
  \ifnum #1<24\normalsize\else \ifnum #1<29\large\else
  \ifnum #1<34\Large\else \ifnum #1<41\LARGE\else
     \huge\fi\fi\fi\fi\fi\fi
  \csname #3\endcsname}%
\else
\gdef\SetFigFont#1#2#3{\begingroup
  \count@#1\relax \ifnum 25<\count@\count@25\fi
  \def\x{\endgroup\@setsize\SetFigFont{#2pt}}%
  \expandafter\x
    \csname \romannumeral\the\count@ pt\expandafter\endcsname
    \csname @\romannumeral\the\count@ pt\endcsname
  \csname #3\endcsname}%
\fi
\fi\endgroup
\begin{picture}(5638,1222)(2004,-1261)
\put(3076,-736){\makebox(0,0)[lb]{\smash{\SetFigFont{12}{14.4}{rm}$u_1$}}}
\put(4651,-736){\makebox(0,0)[lb]{\smash{\SetFigFont{12}{14.4}{rm}$u_2$}}}
\put(6226,-736){\makebox(0,0)[lb]{\smash{\SetFigFont{12}{14.4}{rm}$u_3$}}}
\end{picture}

\]
\ni for $u$'s satisfying appropriate equations. We will use this
approach with $n=1$ to pick up some important cases of commuting
squares in Example 2.11.

\newpage
\ni{\bf 2.11. Binunitaries}

Commuting squares have been used to construct subfactors (see
[GHJ],[Ha]) and the general theory of calculating the
subfactor planar algebra of $\S$4.2 from a commuting square will be
dealt with in a future paper. The treatment uses the language of
statistical mechanical models with some attention paid to critical
points as in chapter 4 of [JS]. Here we give a different
approach which seems more natural from a planar point of view and
will allow us to capture, as special cases, spin model commuting
squares and some vertex model ones with no extension of the planar
algebra formalism. The main concept is that of a bi-invertible
element in a planar algebra which is the next step in the
hierarchy discussed at the end of Example 2.10.

\bs{\bf Definition 2.11.1.} Let $P$ be a general planar algebra. An
invertible element $u\in P_2$ will be called {\it bi-invertible} if
\[
        \begin{picture}(0,0)%
\epsfig{file=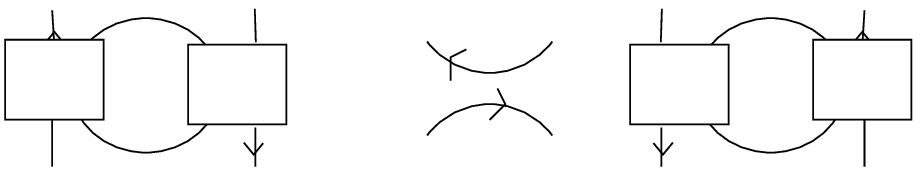}%
\end{picture}%
\setlength{\unitlength}{0.00041700in}%
\begingroup\makeatletter\ifx\SetFigFont\undefined
\def\x#1#2#3#4#5#6#7\relax{\def\x{#1#2#3#4#5#6}}%
\expandafter\x\fmtname xxxxxx\relax \def\y{splain}%
\ifx\x\y   
\gdef\SetFigFont#1#2#3{%
  \ifnum #1<17\tiny\else \ifnum #1<20\small\else
  \ifnum #1<24\normalsize\else \ifnum #1<29\large\else
  \ifnum #1<34\Large\else \ifnum #1<41\LARGE\else
     \huge\fi\fi\fi\fi\fi\fi
  \csname #3\endcsname}%
\else
\gdef\SetFigFont#1#2#3{\begingroup
  \count@#1\relax \ifnum 25<\count@\count@25\fi
  \def\x{\endgroup\@setsize\SetFigFont{#2pt}}%
  \expandafter\x
    \csname \romannumeral\the\count@ pt\expandafter\endcsname
    \csname @\romannumeral\the\count@ pt\endcsname
  \csname #3\endcsname}%
\fi
\fi\endgroup
\begin{picture}(8744,1558)(429,-1584)
\put(7276,-586){\makebox(0,0)[lb]{\smash{
\SetFigFont{12}{14.4}{rm}$u^{-1}$
}}}
\put(901,-811){\makebox(0,0)[lb]{\smash{\SetFigFont{12}{14.4}{rm}$u$}}}
\put(3076,-586){\makebox(0,0)[lb]{\smash{
\SetFigFont{12}{14.4}{rm}$u^{-1}$
}}}
\put(3451,-886){\makebox(0,0)[lb]{\smash{\SetFigFont{12}{14.4}{rm}$=$}}}
\put(5776,-886){\makebox(0,0)[lb]{\smash{\SetFigFont{12}{14.4}{rm}$=$}}}
\put(8626,-811){\makebox(0,0)[lb]{\smash{\SetFigFont{12}{14.4}{rm}$u$}}}
\put(4051,-886){\makebox(0,0)[lb]{\smash{\SetFigFont{12}{14.4}{rm}$\Delta$}}}
\end{picture}

\]
\ni for some non-zero scalar $\D$. (If $P$ is planar $\D$ is
necessarily $\d_1 \ (\d_2)$. If $P$ is a planar $*$-algebra, a
 bi-invertible $u$ is called biunitary  if $u^*=u^{-1}$.
Bi-invertible elements define planar subalgebras as we now
describe. It will be convenient to consider labelled tangles
containing certain distinguished curves joining boundary points,
which intersect with only the other strings in a tangle and do not
meet any internal boxes. From such a tangle, and a
bi-invertible element $u$, we construct an honest labelled tangle,
in the sense of $\S$1, in two steps.

\medskip (i) Orient the distinguished curves. The global
orientation will be denoted $\twoheadrightarrow\!\!-$ and the tangle
orientation by $\rightarrow\!\!-$.

\smallskip (ii) At a point of intersection between the
distinguished curves and the strings of the tangle, insert 2-boxes
containing $u$ or $u^{-1}$ according to the following conventions
\[
	\begin{picture}(0,0)%
\epsfig{file=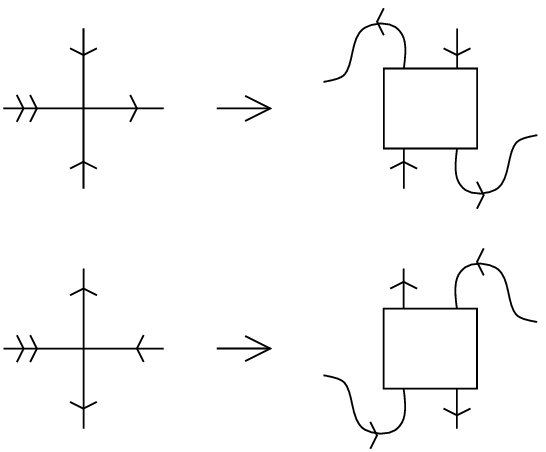}%
\end{picture}%
\setlength{\unitlength}{0.00041700in}%
\begingroup\makeatletter\ifx\SetFigFont\undefined
\def\x#1#2#3#4#5#6#7\relax{\def\x{#1#2#3#4#5#6}}%
\expandafter\x\fmtname xxxxxx\relax \def\y{splain}%
\ifx\x\y   
\gdef\SetFigFont#1#2#3{%
  \ifnum #1<17\tiny\else \ifnum #1<20\small\else
  \ifnum #1<24\normalsize\else \ifnum #1<29\large\else
  \ifnum #1<34\Large\else \ifnum #1<41\LARGE\else
     \huge\fi\fi\fi\fi\fi\fi
  \csname #3\endcsname}%
\else
\gdef\SetFigFont#1#2#3{\begingroup
  \count@#1\relax \ifnum 25<\count@\count@25\fi
  \def\x{\endgroup\@setsize\SetFigFont{#2pt}}%
  \expandafter\x
    \csname \romannumeral\the\count@ pt\expandafter\endcsname
    \csname @\romannumeral\the\count@ pt\endcsname
  \csname #3\endcsname}%
\fi
\fi\endgroup
\begin{picture}(5497,4270)(579,-4684)
\put(4651,-1486){\makebox(0,0)[lb]{\smash{\SetFigFont{12}{14.4}{rm}$u$}}}
\put(5026,-3436){\makebox(0,0)[lb]{\smash{
\SetFigFont{12}{14.4}{rm}$u^{-1}$
}}}
\put(6076,-1036){\makebox(0,0)[lb]{\smash{\SetFigFont{12}{14.4}{rm}(global orientation)}}}
\put(6076,-2215){\makebox(0,0)[lb]{\smash{\SetFigFont{12}{14.4}{rm}going out)}}}
\put(6076,-3436){\makebox(0,0)[lb]{\smash{\SetFigFont{12}{14.4}{rm}(global orientation)}}}
\put(6076,-1711){\makebox(0,0)[lb]{\smash{\SetFigFont{12}{14.4}{rm}(=tangle orientation}}}
\put(6076,-4111){\makebox(0,0)[lb]{\smash{\SetFigFont{12}{14.4}{rm}(=tangel orientation}}}
\put(6076,-4615){\makebox(0,0)[lb]{\smash{\SetFigFont{12}{14.4}{rm}going in)}}}
\put(4876,-4036){\makebox(0,0)[lb]{\smash{\SetFigFont{12}{14.4}{rm}$*$}}}
\put(4351,-1261){\makebox(0,0)[lb]{\smash{\SetFigFont{12}{14.4}{rm}$*$}}}
\end{picture}

\]
\bs
\ni Thus along a distinguished curve one alternately meets $u$ or
$u^{-1}$.

From now on we suppose for convenience that $\D = 1$.

\bs\ni
{\bf Lemma 2.11.2} The Reidemeister type  II moves are
satisfied, i.e.,
\[
	\begin{picture}(0,0)%
\epsfig{file=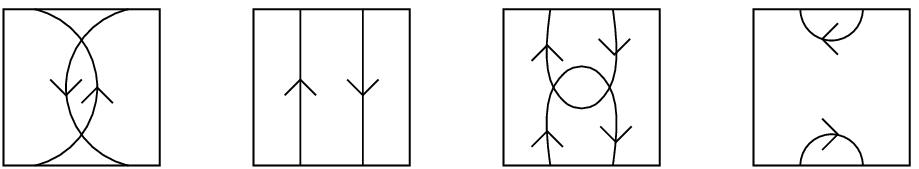}%
\end{picture}%
\setlength{\unitlength}{0.00041700in}%
\begingroup\makeatletter\ifx\SetFigFont\undefined
\def\x#1#2#3#4#5#6#7\relax{\def\x{#1#2#3#4#5#6}}%
\expandafter\x\fmtname xxxxxx\relax \def\y{splain}%
\ifx\x\y   
\gdef\SetFigFont#1#2#3{%
  \ifnum #1<17\tiny\else \ifnum #1<20\small\else
  \ifnum #1<24\normalsize\else \ifnum #1<29\large\else
  \ifnum #1<34\Large\else \ifnum #1<41\LARGE\else
     \huge\fi\fi\fi\fi\fi\fi
  \csname #3\endcsname}%
\else
\gdef\SetFigFont#1#2#3{\begingroup
  \count@#1\relax \ifnum 25<\count@\count@25\fi
  \def\x{\endgroup\@setsize\SetFigFont{#2pt}}%
  \expandafter\x
    \csname \romannumeral\the\count@ pt\expandafter\endcsname
    \csname @\romannumeral\the\count@ pt\endcsname
  \csname #3\endcsname}%
\fi
\fi\endgroup
\begin{picture}(8744,1544)(879,-2183)
\put(5026,-1411){\makebox(0,0)[lb]{\smash{\SetFigFont{12}{14.4}{rm},}}}
\put(2776,-1486){\makebox(0,0)[lb]{\smash{\SetFigFont{12}{14.4}{rm}$=$}}}
\put(7561,-1456){\makebox(0,0)[lb]{\smash{\SetFigFont{12}{14.4}{rm}$=$}}}
\end{picture}

\]
\ni where either of the two curves is distinguished and its global
orientation is arbitrary.

\bs
{\sc Proof.} This is just a re-expression of bi-invertibility. \qed

\bs\ni
{\bf Theorem 2.11.3} If $P$ is a general planar algebra and
$u$ is bi-invertible, let $P^u_k$ be
\bs
\[
\{x\in P_k\ \text{s.t.}\begin{picture}(0,0)%
\epsfig{file=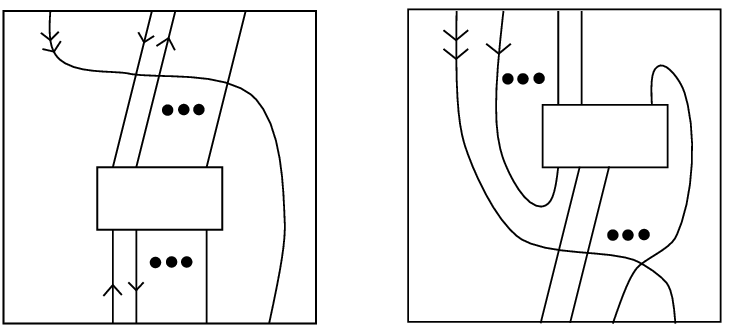}%
\end{picture}%
\setlength{\unitlength}{0.00041700in}%
\begingroup\makeatletter\ifx\SetFigFont\undefined
\def\x#1#2#3#4#5#6#7\relax{\def\x{#1#2#3#4#5#6}}%
\expandafter\x\fmtname xxxxxx\relax \def\y{splain}%
\ifx\x\y   
\gdef\SetFigFont#1#2#3{%
  \ifnum #1<17\tiny\else \ifnum #1<20\small\else
  \ifnum #1<24\normalsize\else \ifnum #1<29\large\else
  \ifnum #1<34\Large\else \ifnum #1<41\LARGE\else
     \huge\fi\fi\fi\fi\fi\fi
  \csname #3\endcsname}%
\else
\gdef\SetFigFont#1#2#3{\begingroup
  \count@#1\relax \ifnum 25<\count@\count@25\fi
  \def\x{\endgroup\@setsize\SetFigFont{#2pt}}%
  \expandafter\x
    \csname \romannumeral\the\count@ pt\expandafter\endcsname
    \csname @\romannumeral\the\count@ pt\endcsname
  \csname #3\endcsname}%
\fi
\fi\endgroup
\begin{picture}(6929,3059)(879,-3083)
\put(4216,-1591){\makebox(0,0)[lb]{\smash{\SetFigFont{12}{14.4}{rm}$=$}}}
\put(2356,-1876){\makebox(0,0)[lb]{\smash{\SetFigFont{12}{14.4}{rm}$x$}}}
\put(6571,-1246){\makebox(0,0)[lb]{\smash{\SetFigFont{12}{14.4}{rm}$y$}}}
\end{picture}
 
\mbox{\ \ for some}\ \ y\in  P_k\}.
\]
\bs
Then $P^u$ is a  general planar subalgebra, $P^u$ is a (general)
planar $*$-subalgebra if $P$ is a planar $*$-algebra and $u$ is
bi-unitary. The properties of being planar, $C^*$ and spherical
are inherited from $P$.

\bs
{\sc Proof.} That $P^u$ is a subalgebra is obvious. The
$*$-property is more interesting. Applying $*$ to the pictures we
obtain (using $u^*=u^{-1}$)
\[
	\begin{picture}(0,0)%
\epsfig{file=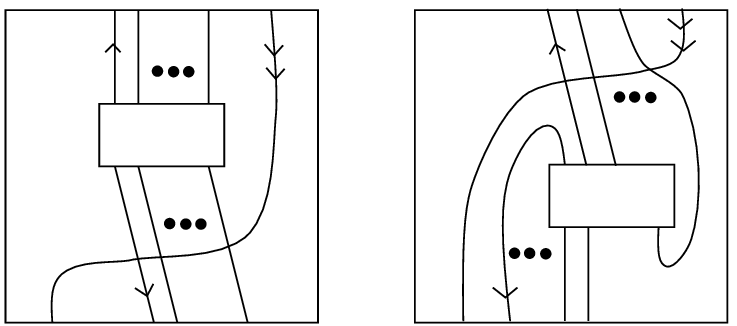}%
\end{picture}%
\setlength{\unitlength}{0.00041700in}%
\begingroup\makeatletter\ifx\SetFigFont\undefined
\def\x#1#2#3#4#5#6#7\relax{\def\x{#1#2#3#4#5#6}}%
\expandafter\x\fmtname xxxxxx\relax \def\y{splain}%
\ifx\x\y   
\gdef\SetFigFont#1#2#3{%
  \ifnum #1<17\tiny\else \ifnum #1<20\small\else
  \ifnum #1<24\normalsize\else \ifnum #1<29\large\else
  \ifnum #1<34\Large\else \ifnum #1<41\LARGE\else
     \huge\fi\fi\fi\fi\fi\fi
  \csname #3\endcsname}%
\else
\gdef\SetFigFont#1#2#3{\begingroup
  \count@#1\relax \ifnum 25<\count@\count@25\fi
  \def\x{\endgroup\@setsize\SetFigFont{#2pt}}%
  \expandafter\x
    \csname \romannumeral\the\count@ pt\expandafter\endcsname
    \csname @\romannumeral\the\count@ pt\endcsname
  \csname #3\endcsname}%
\fi
\fi\endgroup
\begin{picture}(6974,3055)(865,-3084)
\put(4216,-1591){\makebox(0,0)[lb]{\smash{\SetFigFont{12}{14.4}{rm}$=$}}}
\put(2251,-1336){\makebox(0,0)[lb]{\smash{\SetFigFont{12}{14.4}{rm}$x^*$}}}
\put(6556,-1891){\makebox(0,0)[lb]{\smash{\SetFigFont{12}{14.4}{rm}$y^*$}}}
\end{picture}

\]
\ni Now we surround these pictures with the annular tangle
\[
	\begin{picture}(0,0)%
\epsfig{file=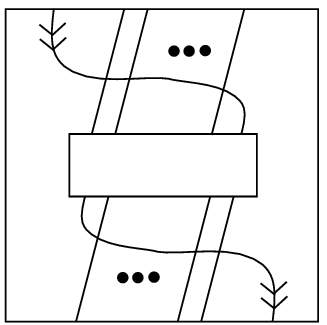}%
\end{picture}%
\setlength{\unitlength}{0.00041700in}%
\begingroup\makeatletter\ifx\SetFigFont\undefined
\def\x#1#2#3#4#5#6#7\relax{\def\x{#1#2#3#4#5#6}}%
\expandafter\x\fmtname xxxxxx\relax \def\y{splain}%
\ifx\x\y   
\gdef\SetFigFont#1#2#3{%
  \ifnum #1<17\tiny\else \ifnum #1<20\small\else
  \ifnum #1<24\normalsize\else \ifnum #1<29\large\else
  \ifnum #1<34\Large\else \ifnum #1<41\LARGE\else
     \huge\fi\fi\fi\fi\fi\fi
  \csname #3\endcsname}%
\else
\gdef\SetFigFont#1#2#3{\begingroup
  \count@#1\relax \ifnum 25<\count@\count@25\fi
  \def\x{\endgroup\@setsize\SetFigFont{#2pt}}%
  \expandafter\x
    \csname \romannumeral\the\count@ pt\expandafter\endcsname
    \csname @\romannumeral\the\count@ pt\endcsname
  \csname #3\endcsname}%
\fi
\fi\endgroup
\begin{picture}(3044,3045)(3867,-3085)
\end{picture}

\]
\ni and then apply type II Reidemeister moves. We see that $x^*\in
P^u$ if $x$ does (though note that the ``$y$'' for $x^*$ is
$y^*$, but rotated.

To show that $P^u$ is a general planar subalgebra, we only have to
show by 1.18 that it is invariant under $\cal A(\phi)$. But if we
arrange the annular tangle $A$ so that all critical points of the
height function on strings are local maxima and minima, the
distinguished line
\[
	\begin{picture}(0,0)%
\epsfig{file=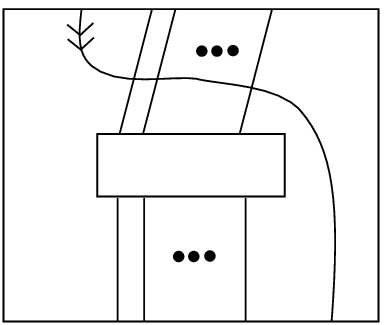}%
\end{picture}%
\setlength{\unitlength}{0.00041700in}%
\begingroup\makeatletter\ifx\SetFigFont\undefined
\def\x#1#2#3#4#5#6#7\relax{\def\x{#1#2#3#4#5#6}}%
\expandafter\x\fmtname xxxxxx\relax \def\y{splain}%
\ifx\x\y   
\gdef\SetFigFont#1#2#3{%
  \ifnum #1<17\tiny\else \ifnum #1<20\small\else
  \ifnum #1<24\normalsize\else \ifnum #1<29\large\else
  \ifnum #1<34\Large\else \ifnum #1<41\LARGE\else
     \huge\fi\fi\fi\fi\fi\fi
  \csname #3\endcsname}%
\else
\gdef\SetFigFont#1#2#3{\begingroup
  \count@#1\relax \ifnum 25<\count@\count@25\fi
  \def\x{\endgroup\@setsize\SetFigFont{#2pt}}%
  \expandafter\x
    \csname \romannumeral\the\count@ pt\expandafter\endcsname
    \csname @\romannumeral\the\count@ pt\endcsname
  \csname #3\endcsname}%
\fi
\fi\endgroup
\begin{picture}(3644,3043)(3579,-3083)
\put(5026,-1636){\makebox(0,0)[lb]{\smash{\SetFigFont{12}{14.4}{rm}$\pi_A(x)$}}}
\end{picture}

\]
\ni can be moved through $\pi_A(x)$, close to \ $\tangle{x}$ 
using only planar isotopy and type II Reidemeister moves.
It can then go past $x$, producing a $y$, and down to the bottom
via type II Reidemeister moves. Thus  $\pi_A(x)$ is in $P^u$.

Planarity, positivity and sphericity are all
inherited. \qed

\bs {\bf Notes.} (i) We will see that $P^u$ may be planar even
when $P$ is only general planar.

(ii) When $k=2$ the equation of Theorem 2.11.3 is an abstract
version of the Yang-Baxter equation ([Ba]).

\bs We now recast the equation of Theorem 2.11.3 in some equivalent
forms which reveal some of its structure. The idea of the equation,
as seen clearly in the proof of 2.11.3 is just that the
distinguished lines can move freely past the boxes. As stated this
requires a special configuration as the distinguished line
approaches a box, but by Reidemeister type II invariance any
approach will do. We record this below, keeping the notation of
2.11.3. (Note that it is somewhat cumbersome to force all the
pictures to fit appropriately into the standard $k$-box. We use a
disk with $2k$ boundary points-as in the introduction.)

\medskip
{\bf Proposition 2.11.4} An $x$ in $P_k$ is in $P^u_k$ iff
there is a $y$ in $P^u_k$ with
\[
	\begin{picture}(0,0)%
\epsfig{file=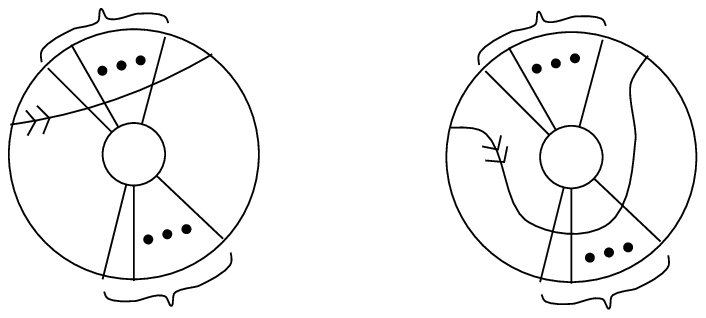}%
\end{picture}%
\setlength{\unitlength}{0.00041700in}%
\begingroup\makeatletter\ifx\SetFigFont\undefined
\def\x#1#2#3#4#5#6#7\relax{\def\x{#1#2#3#4#5#6}}%
\expandafter\x\fmtname xxxxxx\relax \def\y{splain}%
\ifx\x\y   
\gdef\SetFigFont#1#2#3{%
  \ifnum #1<17\tiny\else \ifnum #1<20\small\else
  \ifnum #1<24\normalsize\else \ifnum #1<29\large\else
  \ifnum #1<34\Large\else \ifnum #1<41\LARGE\else
     \huge\fi\fi\fi\fi\fi\fi
  \csname #3\endcsname}%
\else
\gdef\SetFigFont#1#2#3{\begingroup
  \count@#1\relax \ifnum 25<\count@\count@25\fi
  \def\x{\endgroup\@setsize\SetFigFont{#2pt}}%
  \expandafter\x
    \csname \romannumeral\the\count@ pt\expandafter\endcsname
    \csname @\romannumeral\the\count@ pt\endcsname
  \csname #3\endcsname}%
\fi
\fi\endgroup
\begin{picture}(7380,3750)(136,-3958)
\put(4111,-1921){\makebox(0,0)[lb]{\smash{\SetFigFont{12}{14.4}{rm}$=$}}}
\put(2521,-3886){\makebox(0,0)[lb]{\smash{\SetFigFont{12}{14.4}{rm}$q$ points}}}
\put(6766,-3856){\makebox(0,0)[lb]{\smash{\SetFigFont{12}{14.4}{rm}$q$ points}}}
\put(4381,-511){\makebox(0,0)[lb]{\smash{\SetFigFont{12}{14.4}{rm}$p$ points}}}
\put(136,-496){\makebox(0,0)[lb]{\smash{\SetFigFont{12}{14.4}{rm}$p$ points}}}
\put(2041,-1906){\makebox(0,0)[lb]{\smash{\SetFigFont{12}{14.4}{rm}$x$}}}
\put(6256,-1936){\makebox(0,0)[lb]{\smash{\SetFigFont{12}{14.4}{rm}$y$}}}
\end{picture}

\]
\ni where $p+q=2k$, where the pictures may be isotoped in any way
so that the annular region becomes the standard annular region
(with $*$ anywhere allowed by the orientations).

\medskip
{\sc Proof.} If $p>q$, surround the two pictures with the annular
tangle
\[
	\begin{picture}(0,0)%
\epsfig{file=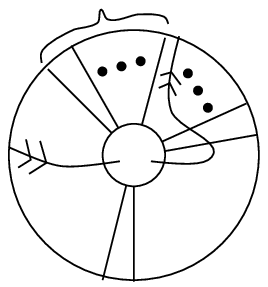}%
\end{picture}%
\setlength{\unitlength}{0.00041700in}%
\begingroup\makeatletter\ifx\SetFigFont\undefined
\def\x#1#2#3#4#5#6#7\relax{\def\x{#1#2#3#4#5#6}}%
\expandafter\x\fmtname xxxxxx\relax \def\y{splain}%
\ifx\x\y   
\gdef\SetFigFont#1#2#3{%
  \ifnum #1<17\tiny\else \ifnum #1<20\small\else
  \ifnum #1<24\normalsize\else \ifnum #1<29\large\else
  \ifnum #1<34\Large\else \ifnum #1<41\LARGE\else
     \huge\fi\fi\fi\fi\fi\fi
  \csname #3\endcsname}%
\else
\gdef\SetFigFont#1#2#3{\begingroup
  \count@#1\relax \ifnum 25<\count@\count@25\fi
  \def\x{\endgroup\@setsize\SetFigFont{#2pt}}%
  \expandafter\x
    \csname \romannumeral\the\count@ pt\expandafter\endcsname
    \csname @\romannumeral\the\count@ pt\endcsname
  \csname #3\endcsname}%
\fi
\fi\endgroup
\begin{picture}(3180,2890)(136,-3098)
\put(136,-496){\makebox(0,0)[lb]{\smash{\SetFigFont{12}{14.4}{rm}$k$ points}}}
\end{picture}

\]
\ni and  use isotopy, rotation (if necessary to get $*$ in the
right place), and type II Reidemeister moves to obtain the same
picture as in 2.11.3. \qed

\bs{\bf Definition 2.11.5.} Given $P$ and $u$ as above, define
$\s_u:P_k\to P_{k+1}$ by the tangle below:
\[
	\begin{picture}(0,0)%
\epsfig{file=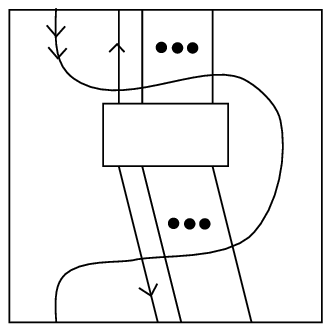}%
\end{picture}%
\setlength{\unitlength}{0.00041700in}%
\begingroup\makeatletter\ifx\SetFigFont\undefined
\def\x#1#2#3#4#5#6#7\relax{\def\x{#1#2#3#4#5#6}}%
\expandafter\x\fmtname xxxxxx\relax \def\y{splain}%
\ifx\x\y   
\gdef\SetFigFont#1#2#3{%
  \ifnum #1<17\tiny\else \ifnum #1<20\small\else
  \ifnum #1<24\normalsize\else \ifnum #1<29\large\else
  \ifnum #1<34\Large\else \ifnum #1<41\LARGE\else
     \huge\fi\fi\fi\fi\fi\fi
  \csname #3\endcsname}%
\else
\gdef\SetFigFont#1#2#3{\begingroup
  \count@#1\relax \ifnum 25<\count@\count@25\fi
  \def\x{\endgroup\@setsize\SetFigFont{#2pt}}%
  \expandafter\x
    \csname \romannumeral\the\count@ pt\expandafter\endcsname
    \csname @\romannumeral\the\count@ pt\endcsname
  \csname #3\endcsname}%
\fi
\fi\endgroup
\begin{picture}(4448,3053)(1576,-3082)
\put(4366,-1334){\makebox(0,0)[lb]{\smash{\SetFigFont{12}{14.4}{rm}$x$}}}
\put(1576,-1561){\makebox(0,0)[lb]{\smash{\SetFigFont{12}{14.4}{rm}$\sigma_u(x)=$}}}
\end{picture}

\]

\bs\ni
{\bf Proposition 2.11.6} The map $\s_u$ is a unital
endomorphism of the filtered algebra $P$ (a $*$-endomorphism if $u$
is unitary) and $P^u_k=\{x\mid \s_u(x)\in P_{1,k+1}\}$.

\bs
{\sc Proof.}  That $\s_u$ preserves multiplication follows from
type II Reidemeister moves. The alternative definition of $P^u$ is
just the case $q=0$, with $*$ appropriately placed, in 2.11.4. \qed

\bs Note that the endomorphism $\s_u$ is the obvious ``shift de un"
when restricted to the Temperley-Lieb subalgebra.

The condition of $2.11.5$ involved a pair $(x,y)$. In fact $x$ is
determined by $y$ and vice versa as we now record.

\bs\ni
{\bf Proposition 2.11.7} Suppose $P$ is a general planar
algebra. If $x\in P^u$ with \ $\s_u(x)=$ $\boxi{y}$.
Then

\medskip
\[
	\begin{picture}(0,0)%
\epsfig{file=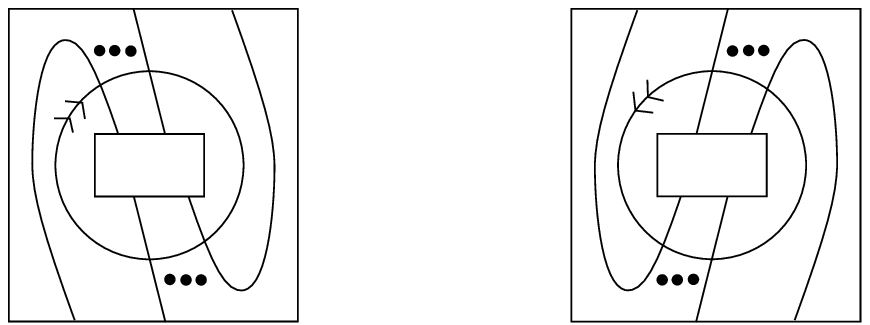}%
\end{picture}%
\setlength{\unitlength}{0.00041700in}%
\begingroup\makeatletter\ifx\SetFigFont\undefined
\def\x#1#2#3#4#5#6#7\relax{\def\x{#1#2#3#4#5#6}}%
\expandafter\x\fmtname xxxxxx\relax \def\y{splain}%
\ifx\x\y   
\gdef\SetFigFont#1#2#3{%
  \ifnum #1<17\tiny\else \ifnum #1<20\small\else
  \ifnum #1<24\normalsize\else \ifnum #1<29\large\else
  \ifnum #1<34\Large\else \ifnum #1<41\LARGE\else
     \huge\fi\fi\fi\fi\fi\fi
  \csname #3\endcsname}%
\else
\gdef\SetFigFont#1#2#3{\begingroup
  \count@#1\relax \ifnum 25<\count@\count@25\fi
  \def\x{\endgroup\@setsize\SetFigFont{#2pt}}%
  \expandafter\x
    \csname \romannumeral\the\count@ pt\expandafter\endcsname
    \csname @\romannumeral\the\count@ pt\endcsname
  \csname #3\endcsname}%
\fi
\fi\endgroup
\begin{picture}(10148,3055)(451,-3088)
\put(5401,-1711){\makebox(0,0)[lb]{\smash{\SetFigFont{12}{14.4}{rm}and }}}
\put(451,-1711){\makebox(0,0)[lb]{\smash{\SetFigFont{12}{14.4}{rm}$y$}}}
\put(976,-1711){\makebox(0,0)[lb]{\smash{\SetFigFont{12}{14.4}{rm}$=$}}}
\put(1726,-1711){\makebox(0,0)[lb]{\smash{\SetFigFont{12}{14.4}{rm}$\frac{1}{\delta_2}$}}}
\put(3676,-1561){\makebox(0,0)[lb]{\smash{\SetFigFont{12}{14.4}{rm}$x$}}}
\put(9076,-1561){\makebox(0,0)[lb]{\smash{\SetFigFont{12}{14.4}{rm}$y$}}}
\put(7276,-1711){\makebox(0,0)[lb]{\smash{\SetFigFont{12}{14.4}{rm}$=$}}}
\put(6376,-1711){\makebox(0,0)[lb]{\smash{\SetFigFont{12}{14.4}{rm}$x\frac{1}{\delta_1}$}}}
\end{picture}

\]

\bs
{\sc Proof.}  Just apply the appropriate annular tangles
and use Reidemeister moves. \qed

\bs Thus we could rewrite equations for $x\in P^u$ entirely in
terms of $x$. The least obvious reformulation of these equations
involves less boundary points than above and requires positivity.

\bs\ni
{\bf Theorem 2.11.8} Let $P$ be a spherical finite-dimensional
$C^*$-planar algebra and $u\in P_2$ be bi-unitary. Then $x\in P^u$
iff
\[
	\begin{picture}(0,0)%
\epsfig{file=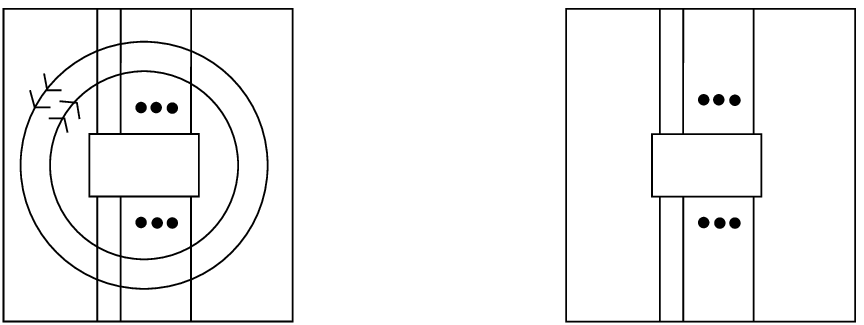}%
\end{picture}%
\setlength{\unitlength}{0.00041700in}%
\begingroup\makeatletter\ifx\SetFigFont\undefined
\def\x#1#2#3#4#5#6#7\relax{\def\x{#1#2#3#4#5#6}}%
\expandafter\x\fmtname xxxxxx\relax \def\y{splain}%
\ifx\x\y   
\gdef\SetFigFont#1#2#3{%
  \ifnum #1<17\tiny\else \ifnum #1<20\small\else
  \ifnum #1<24\normalsize\else \ifnum #1<29\large\else
  \ifnum #1<34\Large\else \ifnum #1<41\LARGE\else
     \huge\fi\fi\fi\fi\fi\fi
  \csname #3\endcsname}%
\else
\gdef\SetFigFont#1#2#3{\begingroup
  \count@#1\relax \ifnum 25<\count@\count@25\fi
  \def\x{\endgroup\@setsize\SetFigFont{#2pt}}%
  \expandafter\x
    \csname \romannumeral\the\count@ pt\expandafter\endcsname
    \csname @\romannumeral\the\count@ pt\endcsname
  \csname #3\endcsname}%
\fi
\fi\endgroup
\begin{picture}(8220,3047)(2379,-3084)
\put(3676,-1561){\makebox(0,0)[lb]{\smash{\SetFigFont{12}{14.4}{rm}$x$}}}
\put(9077,-1562){\makebox(0,0)[lb]{\smash{\SetFigFont{12}{14.4}{rm}$x$}}}
\put(5986,-1711){\makebox(0,0)[lb]{\smash{\SetFigFont{12}{14.4}{rm}$=$}}}
\put(7051,-1681){\makebox(0,0)[lb]{\smash{\SetFigFont{12}{14.4}{rm}$\delta^2$}}}
\end{picture}
.
\]

\bs
{\sc Proof.}  $(\Rightarrow)$ This is easy and requires no
positivity.

$(\Leftarrow)$ We begin by observing that
$\d\s^*_u(z)$ is given by the tangle below
\[
	\begin{picture}(0,0)%
\epsfig{file=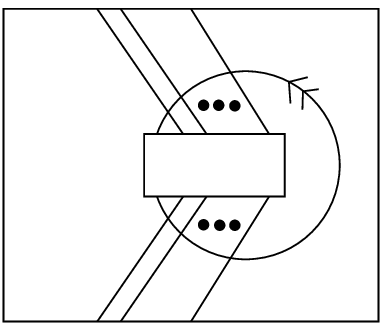}%
\end{picture}%
\setlength{\unitlength}{0.00041700in}%
\begingroup\makeatletter\ifx\SetFigFont\undefined
\def\x#1#2#3#4#5#6#7\relax{\def\x{#1#2#3#4#5#6}}%
\expandafter\x\fmtname xxxxxx\relax \def\y{splain}%
\ifx\x\y   
\gdef\SetFigFont#1#2#3{%
  \ifnum #1<17\tiny\else \ifnum #1<20\small\else
  \ifnum #1<24\normalsize\else \ifnum #1<29\large\else
  \ifnum #1<34\Large\else \ifnum #1<41\LARGE\else
     \huge\fi\fi\fi\fi\fi\fi
  \csname #3\endcsname}%
\else
\gdef\SetFigFont#1#2#3{\begingroup
  \count@#1\relax \ifnum 25<\count@\count@25\fi
  \def\x{\endgroup\@setsize\SetFigFont{#2pt}}%
  \expandafter\x
    \csname \romannumeral\the\count@ pt\expandafter\endcsname
    \csname @\romannumeral\the\count@ pt\endcsname
  \csname #3\endcsname}%
\fi
\fi\endgroup
\begin{picture}(3644,3046)(2379,-3083)
\put(4426,-1636){\makebox(0,0)[lb]{\smash{\SetFigFont{12}{14.4}{rm}$z$}}}
\end{picture}

\]
\ni Here the adjoint of $\s_u$ is as a map between the
finite-dimensional Hilbert spaces $P_k$ and $P_{k+1}$ with inner
products given by the normalized traces.  This formula for $\s^*_u$
thus follows from the equality by isotopy of the following two
networks, the first of which is, up to a power of $\d$,  \
$\<\s_u(x),z^*\>$
\[
	\begin{picture}(0,0)%
\epsfig{file=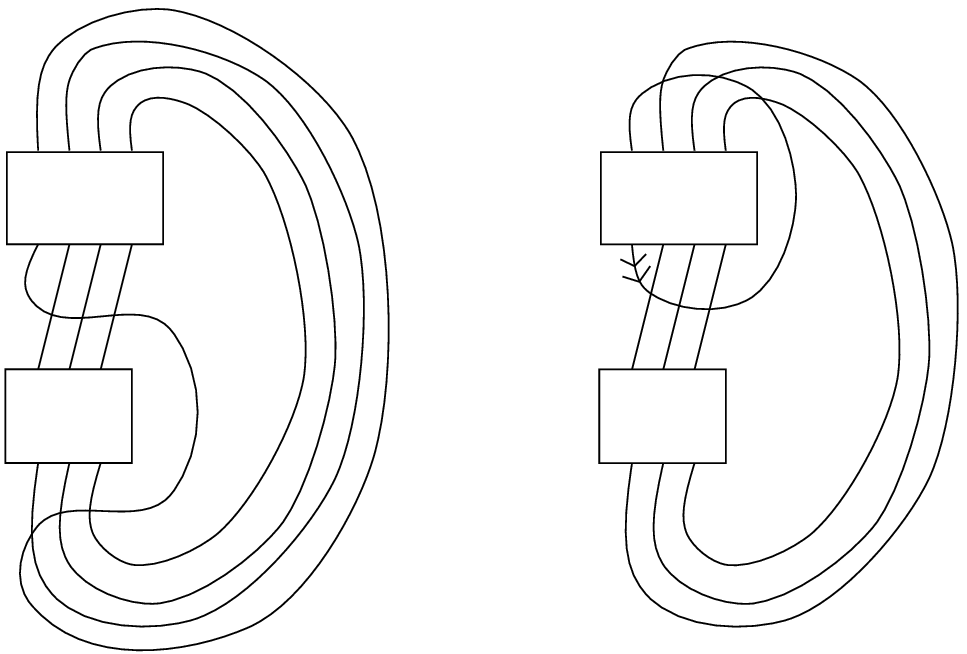}%
\end{picture}%
\setlength{\unitlength}{0.00041700in}%
\begingroup\makeatletter\ifx\SetFigFont\undefined
\def\x#1#2#3#4#5#6#7\relax{\def\x{#1#2#3#4#5#6}}%
\expandafter\x\fmtname xxxxxx\relax \def\y{splain}%
\ifx\x\y   
\gdef\SetFigFont#1#2#3{%
  \ifnum #1<17\tiny\else \ifnum #1<20\small\else
  \ifnum #1<24\normalsize\else \ifnum #1<29\large\else
  \ifnum #1<34\Large\else \ifnum #1<41\LARGE\else
     \huge\fi\fi\fi\fi\fi\fi
  \csname #3\endcsname}%
\else
\gdef\SetFigFont#1#2#3{\begingroup
  \count@#1\relax \ifnum 25<\count@\count@25\fi
  \def\x{\endgroup\@setsize\SetFigFont{#2pt}}%
  \expandafter\x
    \csname \romannumeral\the\count@ pt\expandafter\endcsname
    \csname @\romannumeral\the\count@ pt\endcsname
  \csname #3\endcsname}%
\fi
\fi\endgroup
\begin{picture}(9194,6211)(1764,-5783)
\put(2461,-1366){\makebox(0,0)[lb]{\smash{\SetFigFont{12}{14.4}{rm}$z$}}}
\put(2266,-3451){\makebox(0,0)[lb]{\smash{\SetFigFont{12}{14.4}{rm}$x$}}}
\put(8221,-1351){\makebox(0,0)[lb]{\smash{\SetFigFont{12}{14.4}{rm}$z$}}}
\put(8041,-3496){\makebox(0,0)[lb]{\smash{\SetFigFont{12}{14.4}{rm}$x$}}}
\put(6226,-2776){\makebox(0,0)[lb]{\smash{\SetFigFont{12}{14.4}{rm}$=$}}}
\end{picture}

\]
\ni Thus orthogonal projection $E$ onto $\s_u(P_k)$ is
$\s_u\s^*_u$ which is given on $z\in P_k$ by $\frac 1\d$ times the
following picture
\[
	\begin{picture}(0,0)%
\epsfig{file=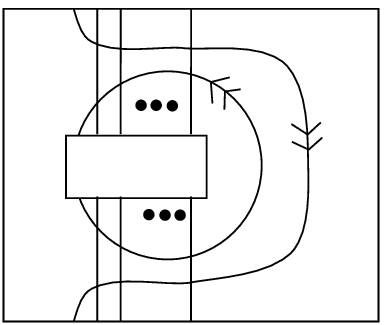}%
\end{picture}%
\setlength{\unitlength}{0.00041700in}%
\begingroup\makeatletter\ifx\SetFigFont\undefined
\def\x#1#2#3#4#5#6#7\relax{\def\x{#1#2#3#4#5#6}}%
\expandafter\x\fmtname xxxxxx\relax \def\y{splain}%
\ifx\x\y   
\gdef\SetFigFont#1#2#3{%
  \ifnum #1<17\tiny\else \ifnum #1<20\small\else
  \ifnum #1<24\normalsize\else \ifnum #1<29\large\else
  \ifnum #1<34\Large\else \ifnum #1<41\LARGE\else
     \huge\fi\fi\fi\fi\fi\fi
  \csname #3\endcsname}%
\else
\gdef\SetFigFont#1#2#3{\begingroup
  \count@#1\relax \ifnum 25<\count@\count@25\fi
  \def\x{\endgroup\@setsize\SetFigFont{#2pt}}%
  \expandafter\x
    \csname \romannumeral\the\count@ pt\expandafter\endcsname
    \csname @\romannumeral\the\count@ pt\endcsname
  \csname #3\endcsname}%
\fi
\fi\endgroup
\begin{picture}(3644,3046)(2379,-3083)
\put(3676,-1561){\makebox(0,0)[lb]{\smash{\SetFigFont{12}{14.4}{rm}$z$}}}
\end{picture}

\]
\ni Orthogonal projection $F$ onto $P_{1,k+1}$ is given by
$\d F(z)= \quad \boxh{z}$.
An element $w$ of $\s_u(P_k)$ is in $P_{1,k+1}$  iff
$EF(w)=w$. But if $x$ satisfies the condition of the theorem we have
\bs
\[
	\begin{picture}(0,0)%
\epsfig{file=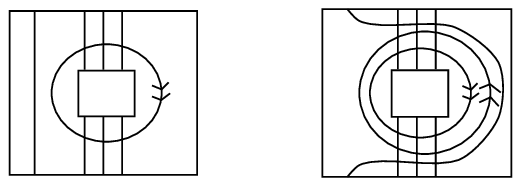}%
\end{picture}%
\setlength{\unitlength}{0.00041700in}%
\begingroup\makeatletter\ifx\SetFigFont\undefined
\def\x#1#2#3#4#5#6#7\relax{\def\x{#1#2#3#4#5#6}}%
\expandafter\x\fmtname xxxxxx\relax \def\y{splain}%
\ifx\x\y   
\gdef\SetFigFont#1#2#3{%
  \ifnum #1<17\tiny\else \ifnum #1<20\small\else
  \ifnum #1<24\normalsize\else \ifnum #1<29\large\else
  \ifnum #1<34\Large\else \ifnum #1<41\LARGE\else
     \huge\fi\fi\fi\fi\fi\fi
  \csname #3\endcsname}%
\else
\gdef\SetFigFont#1#2#3{\begingroup
  \count@#1\relax \ifnum 25<\count@\count@25\fi
  \def\x{\endgroup\@setsize\SetFigFont{#2pt}}%
  \expandafter\x
    \csname \romannumeral\the\count@ pt\expandafter\endcsname
    \csname @\romannumeral\the\count@ pt\endcsname
  \csname #3\endcsname}%
\fi
\fi\endgroup
\begin{picture}(7875,1654)(151,-2800)
\put(151,-2086){\makebox(0,0)[lb]{\smash{\SetFigFont{12}{14.4}{rm}$EF(\sigma_u(x))=\frac{1}{\delta}($}}}
\put(3826,-2011){\makebox(0,0)[lb]{\smash{\SetFigFont{12}{14.4}{rm}$x$}}}
\put(4951,-2086){\makebox(0,0)[lb]{\smash{\SetFigFont{12}{14.4}{rm}$) = \frac{1}{\delta^2}$}}}
\put(6826,-2011){\makebox(0,0)[lb]{\smash{\SetFigFont{12}{14.4}{rm}$x$}}}
\put(8026,-2086){\makebox(0,0)[lb]{\smash{\SetFigFont{12}{14.4}{rm}$=\sigma_u(x)$.}}}
\end{picture}

\]
Hence $x\in P^u$. \qed

\bs {\bf Remark.} We did not use the full force of the hypotheses.
The result will hold in a finite-dimensional general $C^*$-planar
algebra provided isolated circles can be removed with a
multiplicative factor $\d$, that the tangle formula for orthogonal
projection onto $P_{1,k+1}$ is correct, and dim $P_1=1$.

We will see in the case of spin models that $P^u$ may be planar
although $P$ is not. (But the conditions of the above remark are
satisfied by $P^\s$.)

It is easy to check that a bi-invertible $u\in P_2$ may be altered
by four invertible elements $A,B,C,D$ in $P_1$ as in Figure 2.11.9
\begin{figure}
\[
	\input{xfig/pic60}
\]
\caption[2.11.9]{}\label{pic60}
\end{figure}
\bs

{\bf Definition 2.11.10.} Two bi-invertibles are said to differ by
a {\it gauge transformation} if one is obtained from the other as
in Figure 2.11.9.

\bs Gauge transformations have an inessential effect on $P^u$; $A$
and $C$ change absolutely nothing, $B$ and $D$ change $P^u$ by a
planar algebra isomorphism (induced by one on $P$-conjugation by
\ $\boxa{d} \ \boxb{d} \ \boxa{d}\dots$).

In the $*$ case, gauge transformations on bi-unitary matrices are
ones with $A,B,C,D$ unitary.

A significant observation about the equations defining $P^u$ above
is that they are {\it linear} so the calculation of $P^u_k$, given
$P_k$ and $u$, is a {\it finite} problem, unlike the calculation of
the planar subalgebra generated by some set, which requires
consideration of infinitely many tangles. In practice, however, the
brute force calculation, even just of dim $P^u_k$, runs into a
serious problem. For the dimension of $P_k$ grows exponentially
with $k$. For $k=2$ the calculation is usually easy enough (indeed
we give an entirely satisfactory general solution for $k=2$ when
$P=P^\s$, below) and somewhat harder for $k=3$. For $k=4$ it tends
to be very demanding even for relatively ``small" $P$'s.
On the other hand, we are dealing with objects  with a lot of
structure. For instance, once we have calculated $P^u_2$ by brute
force or otherwise, the fact that $P^u$ is a planar algebra means
that every unlabelled 2-tangle gives a {\it nonlinear constraint.}
For if a 2-tangle is given labels with elements in  $P$, in order
for the corresponding element of $P_2$ to be in $P^u$, it must lie
in the linear subspace of $P_2$ already calculated. It was the
desire to systematically exploit these highly interesting nonlinear
constraints that led to the theory of planar algebras --- their
generality was only appreciated afterwards.

There are good reasons for wanting to calculate $P^u$. In general
the calculation is greatly facilitated by the presence of group
symmetry but there are many cases of bi-invertible elements with no
apparent symmetry. We hope that the planar algebra $P^u$ plays the
role of ``higher", non-group-like symmetries which reveal
structural properties of the combinatorial object $u$.

We turn now to a special case where this program has had some
partial success, namely in Hadamard matrices.  The theory is no
different for generalized Hadamard matrices, which occur as
biunitaries for spin models. Consider the spin model  $P^\s$ with
$Q$ spins and its spherical planar algebra structure ($\S$2.8). A
bi-invertible element $u$ of $P^\s_2$ is an invertible $Q\x Q$
matrix $u^b_a$ such that
\[
	\begin{picture}(0,0)%
\epsfig{file=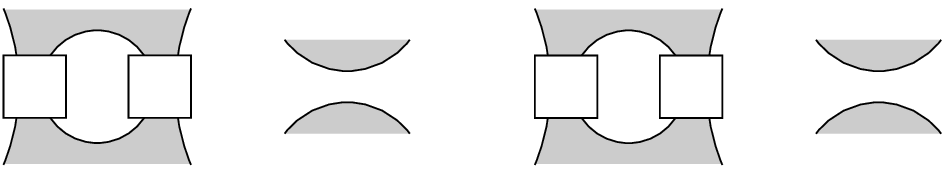}%
\end{picture}%
\setlength{\unitlength}{0.00041700in}%
\begingroup\makeatletter\ifx\SetFigFont\undefined
\def\x#1#2#3#4#5#6#7\relax{\def\x{#1#2#3#4#5#6}}%
\expandafter\x\fmtname xxxxxx\relax \def\y{splain}%
\ifx\x\y   
\gdef\SetFigFont#1#2#3{%
  \ifnum #1<17\tiny\else \ifnum #1<20\small\else
  \ifnum #1<24\normalsize\else \ifnum #1<29\large\else
  \ifnum #1<34\Large\else \ifnum #1<41\LARGE\else
     \huge\fi\fi\fi\fi\fi\fi
  \csname #3\endcsname}%
\else
\gdef\SetFigFont#1#2#3{\begingroup
  \count@#1\relax \ifnum 25<\count@\count@25\fi
  \def\x{\endgroup\@setsize\SetFigFont{#2pt}}%
  \expandafter\x
    \csname \romannumeral\the\count@ pt\expandafter\endcsname
    \csname @\romannumeral\the\count@ pt\endcsname
  \csname #3\endcsname}%
\fi
\fi\endgroup
\begin{picture}(9038,1545)(2679,-1734)
\put(4876,-1036){\makebox(0,0)[lb]{\smash{\SetFigFont{12}{14.4}{rm}$=$}}}
\put(9977,-1037){\makebox(0,0)[lb]{\smash{\SetFigFont{12}{14.4}{rm}$=$}}}
\put(2776,-886){\makebox(0,0)[lb]{\smash{\SetFigFont{12}{14.4}{rm}$*$}}}
\put(4276,-1261){\makebox(0,0)[lb]{\smash{\SetFigFont{12}{14.4}{rm}$*$}}}
\put(3001,-1186){\makebox(0,0)[lb]{\smash{\SetFigFont{12}{14.4}{rm}$u$}}}
\put(8176,-1186){\makebox(0,0)[lb]{\smash{\SetFigFont{12}{14.4}{rm}$*$}}}
\put(9076,-886){\makebox(0,0)[lb]{\smash{\SetFigFont{12}{14.4}{rm}$*$}}}
\put(9226,-1111){\makebox(0,0)[lb]{\smash{\SetFigFont{12}{14.4}{rm}$u$}}}
\put(6976,-1036){\makebox(0,0)[lb]{\smash{\SetFigFont{12}{14.4}{rm}and}}}
\put(4531,-736){\makebox(0,0)[lb]{\smash{
\SetFigFont{12}{14.4}{rm}$u^{-1}$
}}}
\put(8461,-736){\makebox(0,0)[lb]{\smash{
\SetFigFont{12}{14.4}{rm}$u^{-1}$
}}}
\end{picture}

\]
\ni So if $(u^{-1})^b_a=v^b_a$ we have $u_{a,b}v_{b,a}=1/Q$ where
the factor $1/Q$ comes from counting oriented circles after
smoothing. If $u$ is biunitary, $v_{a,b}=\overline{u_{b,a}}$, so
the condition is precisely
\[
 |u_{a,b}|=\frac{1}{\sqrt{Q}} 
\]
We call a unitary matrix satisfying 2.11.11 a {\it
generalized Hadamard matrix.} A Hadamard matrix is just $\sqrt{Q}$
times a real generalized Hadamard matrix.

Gauge transformations alter a generalized Hadamard matrix by
multiplying rows and columns by scalars of modulus one ($\pm 1$ in
the Hadamard case). This, together with permutations of the rows
and columns, gives what is called {\it Hadamard equivalence} of
(generalized) Hadamard matrices. Row and column permutations
are easily seen by
2.11.6 to produce equivalent $P^u$'s so any information about $u$
obtained from $P^u$ alone will be {\it invariant under  Hadamard
equivalence.} (The endomorphism $\s_u$ of 2.11.5 itself is more
information than just $P^u$.)

\bs\ni
{\bf Proposition 2.11.12} If $u$ is a generalized Hadamard
matrix, $P^u$ is planar, hence a spherical $C^*$-planar algebra.
Moreover, ${\op{dim}} \ P^u_1=1$, and $P^u_2$ and $P^u_{1,3}$ are
abelian.

\bs
{\sc Proof.} Obviously dim $P^u_1=1$ implies planarity, so
consider a tangle $T$ representing an element of $P_1$. It consists
of a vertical straight line with a 1-box on it, and networks to the
left and right. The networks to the left have exterior shaded white
so only contribute scalars. The picture below is the condition for
such an element to be in  $P^u_1$ (for some element $S$)
\[
	\begin{picture}(0,0)%
\epsfig{file=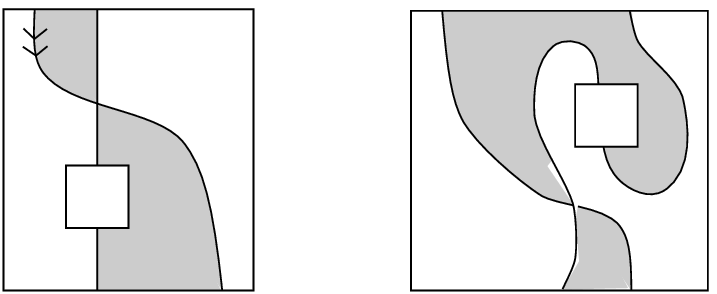}%
\end{picture}%
\setlength{\unitlength}{0.00041700in}%
\begingroup\makeatletter\ifx\SetFigFont\undefined
\def\x#1#2#3#4#5#6#7\relax{\def\x{#1#2#3#4#5#6}}%
\expandafter\x\fmtname xxxxxx\relax \def\y{splain}%
\ifx\x\y   
\gdef\SetFigFont#1#2#3{%
  \ifnum #1<17\tiny\else \ifnum #1<20\small\else
  \ifnum #1<24\normalsize\else \ifnum #1<29\large\else
  \ifnum #1<34\Large\else \ifnum #1<41\LARGE\else
     \huge\fi\fi\fi\fi\fi\fi
  \csname #3\endcsname}%
\else
\gdef\SetFigFont#1#2#3{\begingroup
  \count@#1\relax \ifnum 25<\count@\count@25\fi
  \def\x{\endgroup\@setsize\SetFigFont{#2pt}}%
  \expandafter\x
    \csname \romannumeral\the\count@ pt\expandafter\endcsname
    \csname @\romannumeral\the\count@ pt\endcsname
  \csname #3\endcsname}%
\fi
\fi\endgroup
\begin{picture}(6806,2744)(2079,-3683)
\put(3001,-2761){\makebox(0,0)[lb]{\smash{\SetFigFont{12}{14.4}{rm}$T$}}}
\put(7858,-1981){\makebox(0,0)[lb]{\smash{\SetFigFont{12}{14.4}{rm}$S$}}}
\put(5158,-2431){\makebox(0,0)[lb]{\smash{\SetFigFont{12}{14.4}{rm}$=$}}}
\end{picture}

\]
\ni If the bottom shaded region is assigned a spin $a$, and the top
region a spin $b$, the left-hand side gives $u^b_aT_a$ and the
right-hand side gives $u^b_aS_b$, so $T_a$ is independent of $a$,
and dim  $P^u_1=1$.  $P^u_{1,3}$ is abelian because $P_{1,3}$ is
and $P^u_2$ is abelian since it is $\s^{-1}_u(\s_u(P_2)\cap P_{1,3})$
by 2.11.6. \qed

\bs So by $\S$4.3, a generalized Hadamard matrix $u$ yields a
subfactor whose planar algebra invariant is $P^u$. In fact such a
subfactor was the starting point of the theory of planar algebras, as the equations for $P^u$ are
those for the relative commutants of a spin model commuting square
given in [JS]. Note that the original subfactor is
hyperfinite whereas the one obtained from 4.3 is not! We now
determine $P^u_2$ for a generalized Hadamard matrix $u$.

\bs {\bf Definition 2.11.13.} Given a $Q\x Q$
generalized Hadamard matrix $u^b_a$ we define the $Q^2\x Q^2$ {\it
profile matrix} \ Prof$(u)$ by
$$
{\op{Prof}}(u)^{c,d}_{a,b} =\sum_x
u^x_a \overline{u^x_b} \ \overline{u^x_c} \ u^x_d \ .
$$
The profile matrix is used in the theory of Hadamard matrices.
We will see that it determines $P^u$.

\bs {\bf Definition 2.11.14.} Given the $Q^2\x Q^2$
matrix \ Prof$(u)$, define the directed graph $\G_u$ on $Q^2$
vertices by $(a,b)\rightarrow\!\!\!- (c,d)$ iff \
Prof$(u)^{c,d}_{a,b}\neq 0$.

\bs The isomorphism class of $\G_u$ is an invariant of Hadamard
equivalence.

\bs\ni
{\bf Theorem 2.11.15} If $u$ is a $Q\x Q$
generalized Hadamard matrix thought of as a biunitary for the spin
model $P^\s$, then the minimal projections of the abelian
$C^*$-algebra $P^u_2$ are in bijection with the connected
components of the graph $\G_u$. Moreover the (normalized) trace of
such a projection is $n/Q^2$ where $n$ is the size of the connected
component, which is necessarily a multiple of $Q$.

\bs
{\sc Proof.} For matrices $x^b_a, \ y^b_a$, the equations of
Theorem 2.11.3 are the ``star-triangle" equations
\[
\sum_d u^d_a \ \overline{u^d_b} \ x^c_d =
u^c_a  \ \overline{u^c_b} \ y^b_a 
\]
which amount to saying that, for each $(a,b)$, the vector
$v_{(a,b)}$ whose $d^{\text{th}}$ component is the
$u^d_a\overline{u^d_b}$ is an eigenvector of the matrix $x^c_d$
with eigenvalue $y^b_a$. The profile matrix is just the matrix of
inner products
$\<v_{(a,b)},v_{(c,d)}\>$ so the orthogonal projection onto the
linear span of $v_{(a,b)}$'s in a connected component is in $P^u_2$
and is necessarily minimal since eigenvectors for distinct minimal
projections are orthogonal.

If the matrix $x$ is an orthogonal projection, $y^b_a$ is either 1
or 0 depending on whether $v_{(a,b)}$ is in the connected
component or not. Consider the picture
\[
	\begin{picture}(0,0)%
\epsfig{file=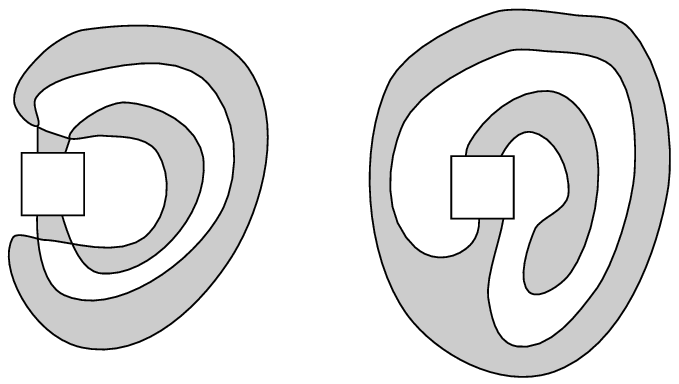}%
\end{picture}%
\setlength{\unitlength}{0.00041700in}%
\begingroup\makeatletter\ifx\SetFigFont\undefined
\def\x#1#2#3#4#5#6#7\relax{\def\x{#1#2#3#4#5#6}}%
\expandafter\x\fmtname xxxxxx\relax \def\y{splain}%
\ifx\x\y   
\gdef\SetFigFont#1#2#3{%
  \ifnum #1<17\tiny\else \ifnum #1<20\small\else
  \ifnum #1<24\normalsize\else \ifnum #1<29\large\else
  \ifnum #1<34\Large\else \ifnum #1<41\LARGE\else
     \huge\fi\fi\fi\fi\fi\fi
  \csname #3\endcsname}%
\else
\gdef\SetFigFont#1#2#3{\begingroup
  \count@#1\relax \ifnum 25<\count@\count@25\fi
  \def\x{\endgroup\@setsize\SetFigFont{#2pt}}%
  \expandafter\x
    \csname \romannumeral\the\count@ pt\expandafter\endcsname
    \csname @\romannumeral\the\count@ pt\endcsname
  \csname #3\endcsname}%
\fi
\fi\endgroup
\begin{picture}(6432,3604)(1830,-3782)
\put(2218,-1936){\makebox(0,0)[lb]{\smash{\SetFigFont{12}{14.4}{rm}$x$}}}
\put(6343,-1906){\makebox(0,0)[lb]{\smash{\SetFigFont{12}{14.4}{rm}$y$}}}
\put(4738,-2041){\makebox(0,0)[lb]{\smash{\SetFigFont{12}{14.4}{rm}$=$}}}
\end{picture}

\]
\ni
Applying Reidemeister type II moves and summing we obtain the
assertion about the trace. (It is a multiple of $1/Q$ since $x$ is
a $Q\x Q$ matrix.) \qed

\bs If $G$ is a finite abelian group and $g\mapsto\hat g$ is an
isomorphism of $G$ with its dual $\hat G$ (=Hom$(G,\Bbb C^*)$), we
obtain a generalized Hadamard matrix $u$, with $Q=|G|$, by setting
$u^h_g=\frac{1}{\sqrt{Q}} \ \hat h(g)$. We call this a {\it
standard} generalized Hadamard matrix.  It is Hadamard if
$G=(\z/2\z)^n$ for some $n$. We leave it to the reader to check
that if $u$ is standard $P^u$ is exactly the planar
algebra of
$\S$2.9 for the group $G$. In particular, dim$(P^u_k)=Q^k$. It is
well known in subfactor theory that any subfactor with
$N'\cap M_1=\Bbb C^{[M:N]}$ comes from a  group. It can
also be seen directly from association schemes that if
dim$(P^u_2)=Q$ then $u$ is standard up to gauge equivalence
(recall that $P^u$ is always an assocation scheme as remarked in
$\S$2.8).

We have, together with R.~Bacher,
 P.~de la Harpe, and M.G.V.~Bogle performed many computer
calculations. So far we have not found a
generalized Hadamard matrix $u$ for which dim$(P^u_2)=2$ but
dim$(P^u_3) >5$.  Such an example would be a confirmation of our
non-group symmetry program as group-like symmetries tend to show up
in $P_2$. In particular the five $16\x 16$ Hadamard matrices have
dim$ P^u_2=$ 16,8,5,3 and 3, and are completely distinguished by
the trace. There are group-like symmetries in all cases
corresponding to the presence of normalizer in the subfactor
picture.

Haagerup has shown how to construct many interesting examples and
given a complete classification for $Q=5$. In the circulant case he
has shown there are only finitely many examples for fixed prime
$Q$  (see [\qquad ]).

Perhaps somewhat surprisingly, the  presence of a lot of symmetry
in $u$ can cause $P^u_2$ to be small! The kind of biunitary
described in the following result is quite common --- the Paley
type Hadamard matrices give an example.

\bs\ni
{\bf Proposition 2.11.17} Suppose $Q-1$ is prime and let $u$
be a $Q\x Q$ generalized Hadamard matrix with the following two
properties (the first of which is always true up to gauge
equivalence):
\begin{verse}
(i) There is an index $*$ with $u^a_*=u^*_a=1$ for all $a$.

(ii) The group $\z/(Q-1)\z$ acts transitively on the
spins other than $*$, and $u^{ga}_{gb}=u^a_b$ for all
$g\in\z/(Q-1)\z$.
\end{verse}
Then dim$(P^u_2)=2$ or $u$ is gauge equivalent to a standard
matrix.

\bs
{\sc Proof.} The nature of the equations 2.11.15 makes it clear
that $\z/(Q-1)\z$ acts by automorphisms on $P^u_2$, obviously
fixing the projection $e_1$ which is the matrix
$x^b_a=1/Q$. Thus the action preserves $(1-e_1)P^u_2(1-e_1)$.
Since $(Q-1)$ is prime there are  only two possibilities:
either the action is non-trivial and dim$(P^u_2)=Q$ so $P^u$ is
standard, or every solution of 2.11.15 is fixed by
$\z/(Q-1)\z$. In the latter case let $x^b_a,y^b_a$ be a solution
of 2.11.15. Then putting $c=*$ we obtain $\sum_d u^d_a \
\overline{u^d_b} x^*_d=y^b_a$, so $y^b_a$ is determined by the two
numbers
$x^*_*$ and $x^*_d$, $d\neq *$. So by 2.11.7 we are done. \qed

\bs Note that the standard case in the above result can occur. The
$8\x 8$ Hadamard matrix is of the required form, but it is Hadamard
equivalent to a standard matrix. For $Q=$ 12 and 24 this cannot be the
case and dim $P_2^u=2$.

We have very few general results on $P^u_k$ for $k >2$.
We only record the observation that $P^u_k$ is the $\d^2$
eigenspace for the $Q^k\x Q^k$ matrix  given by the ``transfer
matrix with periodic horizontal boundary conditions" for the
$Q$-spin vertex model having the profile matrix as Boltzmann
weights. The transfer matrix is given by the picture:
\[
	\begin{picture}(0,0)%
\epsfig{file=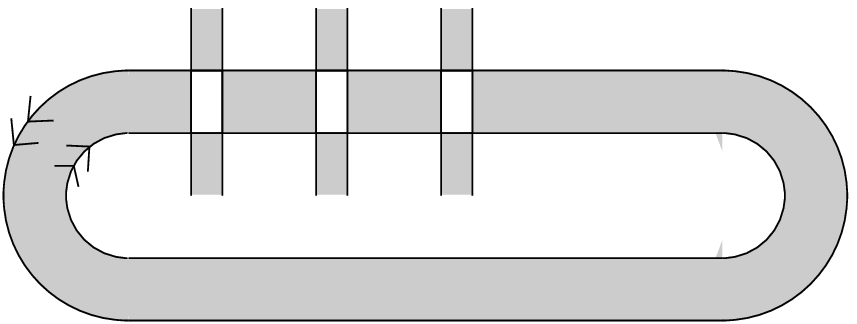}%
\end{picture}%
\setlength{\unitlength}{0.00041700in}%
\begingroup\makeatletter\ifx\SetFigFont\undefined
\def\x#1#2#3#4#5#6#7\relax{\def\x{#1#2#3#4#5#6}}%
\expandafter\x\fmtname xxxxxx\relax \def\y{splain}%
\ifx\x\y   
\gdef\SetFigFont#1#2#3{%
  \ifnum #1<17\tiny\else \ifnum #1<20\small\else
  \ifnum #1<24\normalsize\else \ifnum #1<29\large\else
  \ifnum #1<34\Large\else \ifnum #1<41\LARGE\else
     \huge\fi\fi\fi\fi\fi\fi
  \csname #3\endcsname}%
\else
\gdef\SetFigFont#1#2#3{\begingroup
  \count@#1\relax \ifnum 25<\count@\count@25\fi
  \def\x{\endgroup\@setsize\SetFigFont{#2pt}}%
  \expandafter\x
    \csname \romannumeral\the\count@ pt\expandafter\endcsname
    \csname @\romannumeral\the\count@ pt\endcsname
  \csname #3\endcsname}%
\fi
\fi\endgroup
\begin{picture}(8132,3043)(585,-2783)
\end{picture}

\]

\ni where of course the internal spins have been summed over.
This is an immediate consequence of 2.11.8.

We would like to make the following two open problems about
matrices quite explicit. Both concern a generalized Hadamard matrix
$u$.
\begin{verse}
(i) Is the calculation of dim $P^u_k$ feasible in the
polynomial time as a function of $k$?

(ii) Is there a $u$ for which dim
$P^u_k=\frac{1}{k+1}{2k\choose k}$? (i.e.,
$P^u_k$ is just the Temperley-Lieb algebra).
\end{verse}

Finally we make some comments on vertex models. There are many
formal connections with Hopf algebras here which is not surprising
since quantum groups arose from vertex models in statistical
mechanical models ([Dr]). Banica has done some interesting
work from this point of view --- see [Ban].

A vertex model, in the above context, is simply a biunitary (or
biinvertible) in the planar algebra $P^{\otimes}$. The equations of
Theorem 2.11.3 are then just the equations for the higher relative
commutants of a subfactor coming from a certain commuting square
(see [JS]). Perhaps the most interesting examples not coming from
the quantum group machinery are the Krishan-Sunder ``bipermutation
matrices" where $u$ is a permutation matrix with respect to some
basis of the underlying vector space (see [KS]).
B.Bhattacharya has exhibited a planar algebra which is bigger than
that of example 2.3. (Fuss Catalan) and which is necessarily a
planar subalgebra of $P^u$  if $u$ is a bipermutation matrix.


\bs\bs\ni{\large\bf 3. General Structure Theory}

\medskip\ni{\bf 3.1. Algebra structure, Markov trace}

The proof of Theorem 3.1.3 below is routine for those conversant
with [J1] or [GHJ]. We include it since, as stated, it can
be useful in determining principal graphs. Recall that in a planar
algebra $P$, \ $e_k$ denotes the idempotent in $P_k$ equal to
$\frac{1}{\delta}(| |\dots \quad )$.

{\bf Lemma 3.1.1} Let $P$ be a finite-dimensional
spherical nondegenerate planar algebra over an algebraically closed
field. Then for each  $k$, \
$P_{k-1}e_{k-1}P_{k-1}$ is a 2-sided ideal, denoted $I_k$, in $P_k$
and if ${\cal M}_k$ is a set of minimal idempotents in $P_k$
generating all the distinct minimal ideals in $P_k/I_k$, we have
\begin{verse}
(i)   $pP_ke_{k-1}=0$ \ for $p$ in ${\cal M}_k$

(ii)  $pP_kq=0$ for $p\neq q$ in ${\cal M}_k$

(iii) For each $x$ in $P_k\backslash I_k$ there is a
$p\in\cal M_k$ with $xP_kp\neq 0$.

(iv) ${\op{tr}}(p)\neq 0$ for all $p\in {\cal M}_k$.

(v) $I_{k+2}=\bigoplus_{p\in{\cal M}_k}
P_{k+2}pe_{k+1}P_{k+2}$, \ $pe_{k+1}$ being a minimal idempotent in
$P_{k+2}$.
\end{verse}
Moreover, if, for each $k$, $\cal N_k$ is a set of minimal
idempotents of $P_k$ satisfying ${\op{(i)}}\dots{\op{(iv)}}$ (with
${\cal M}_k$ replaced by ${\cal N}_k$), then there is an invertible
$u_k$ in $P_k$ with $u_k\cal N_ku^{-1}_k={\cal M}_k$ (so in
particular (v) is true for ${\cal N}_k$).

\bigskip
{\sc Proof.} To see that $P_{k-1}e_{k-1}P_{k-1}$ is an
ideal, consider the maps $\a,\beta: P_k\to P_{k-1}$ given by the
annular tangles
\[
	\begin{picture}(0,0)%
\epsfig{file=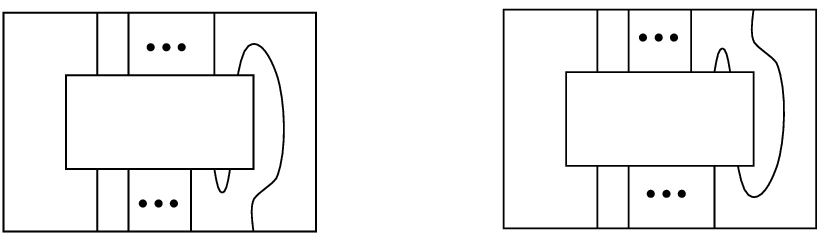}%
\end{picture}%
\setlength{\unitlength}{0.00041700in}%
\begingroup\makeatletter\ifx\SetFigFont\undefined
\def\x#1#2#3#4#5#6#7\relax{\def\x{#1#2#3#4#5#6}}%
\expandafter\x\fmtname xxxxxx\relax \def\y{splain}%
\ifx\x\y   
\gdef\SetFigFont#1#2#3{%
  \ifnum #1<17\tiny\else \ifnum #1<20\small\else
  \ifnum #1<24\normalsize\else \ifnum #1<29\large\else
  \ifnum #1<34\Large\else \ifnum #1<41\LARGE\else
     \huge\fi\fi\fi\fi\fi\fi
  \csname #3\endcsname}%
\else
\gdef\SetFigFont#1#2#3{\begingroup
  \count@#1\relax \ifnum 25<\count@\count@25\fi
  \def\x{\endgroup\@setsize\SetFigFont{#2pt}}%
  \expandafter\x
    \csname \romannumeral\the\count@ pt\expandafter\endcsname
    \csname @\romannumeral\the\count@ pt\endcsname
  \csname #3\endcsname}%
\fi
\fi\endgroup
\begin{picture}(7845,2175)(2079,-2483)
\put(5701,-1561){\makebox(0,0)[lb]{\smash{\SetFigFont{12}{14.4}{rm}and}}}
\end{picture}

\]
\ni respectively. A  diagram shows that
$xe_{k-1}y=\a(x)e_{k-1}\beta(y)$ for $x,y\in P_k$.

By Corollary 1.30, $P_k$ is semisimple and multimatrix since $K$ is
algebraically closed. Thus properties (i), (ii), (iii) and (v) are
obvious for ${\cal M}_k$. If $p\neq 0$ satisfied tr$(p)=0$, then tr
would vanish on the whole matrix algebra containing $p$ which
would then be orthogonal to $P_k$.

Finally, suppose we are given ${\cal N}_k$ satisfying (i)--(iv).
Then since $P_k$ is multimatrix, each $p$ in  ${\cal N}_k$ belongs to
a unique matrix algebra summand in which there is an invertible
$u_p$ with $u_ppu_p^{-1}\in{\cal M}_k$. Putting together the
$u_p$'s, and the identity of $I_k$, we get $u_k$. Property (v) for
${\cal N}_k$ then follows from (iii). \qed

\bs{\bf Definition 3.1.2.} With $P$ and ${\cal N}_k$ as in 3.11, we
define the {\it principal graph} $\G_P$ of $P$ to be the
(bipartite) graph whose vertices are $\bigcup_{k\geq 0}\cal N_k$
with distinguished vertex $*$ so that ${\cal N}_0=\{*\}$, and
dim $(pP_{k+1}q)$ edges between $p\in{\cal N}_k$ and $q\in{\cal N}_{k+1}$. 
Let $d_p$ denote the distance from $p$ to $*$ on $\G_P$.

\bs\ni
{\bf Theorem 3.1.3} As an algebra, $P_k$ is isomorphic to the
algebra whose basis is random walks of length $2k$ on $\G_P$
beginning and ending at $*$ with multiplication rule
$w_1w_2=w_3$ if the first half of the walk $w_2$ is equal to the
second half of $w_1$, and $w_3$ is the first half of $w_1$ followed
by the second half of $w_2$; \ 0 otherwise. Moreover, if $\vec t$
is the function from the vertices of $\G_P$ to $K$, \
$\vec t_p=\d^{d_p}{\op{tr}}(p)$, \ $\vec t$ is an eigenvector for
the adjacency matrix of $\G$, eigenvalue $\d$.

\bigskip
{\sc Proof.} The first assertion is easily equivalent to showing
that the Bratteli diagram (see [GHJ]) of the multimatrix
algebra $P_k$ in $P_{k+1}$ is the bipartite graph consisting of
those vertices $p$ with $d_p\equiv k$ (mod 2) and $d_p\leq k$,
connected to those with
$d_p\equiv (k+1)$ (mod 2) and $d_p\leq k+1$ with appropriate
multiplicities. Observe first that $P_0=\Bbb C$ and ${\cal M}_1$ is a
set of minimal projections, one for each matrix algebra summand of
$P_1$, so the Bratelli diagram is correct for $P_0\subset P_1$. Now
proceed by induction on $k$. The trace on $P_k$ is nondegenerate,
as is its restriction to $P_{k-1}$ so one may perform the abstract
``basic construction" of [J1] to obtain the algebra
$\<P_k,e_{P_{k-1}}\>$ which is multimatrix and isomorphic to
$P_k\otimes_{P_{k-1}} P_k$ as a $P_k-P_k$ bimodule via the map
$x\otimes y\mapsto xe_{P_{k-1}} y$. Moreover the matrix algebra
summands of $\<P_k,e_{P_{k-1}}\>$ are indexed by those of $P_{k-1}$,
which by induction are the vertices of $\G_P$ with $d_p\leq k-1$, \
$d_P\equiv (k+1)$ mod 2.  If one defines the trace tr on
$\<P_k,e_{P_{k-1}}\>$ by $\widetilde{\op{tr}}(xe_{P_{k-1}}y)=
\frac{1}{\d^2} {\op{tr}}(xy)$ then the traces of minimal
projections in $\<P_k,e_{P_{k-1}}\>$ are $\frac{1}{\d^2}$ times
those in $P_{k-1}$. Moreover, setting $\g(xe_{P_{k-1}}y)=xe_ky$
defines an algebra homomorphism from $\<P_k,e_{P_{k-1}}\>$ which is
injective by property (iv) and onto $I_{k+1}$. And
Tr$(xe_ky)=\frac 1\d$ Tr$(xy)$ so tr
$=\widetilde{\op{tr}}\circ\g^{-1}$ on $I_{k+1}$. Properties (i),
(ii) and (iii) ensure that the other  vertices of the Bratteli
diagram for $P_k\subset P_{k+1}$ are labelled by vertices $p$ of
$\G_P$ with $d_P=k+1$. And the number of edges on $\G_P$ connecting
a $p$ in ${\cal M}_k$ to a $q$ in ${\cal M}_{k+1}$ is by definition the
number of edges in the Bratteli diagram. That there are no edges
between ${\cal M}_{k+1}$ and ${\cal M}_j$, $j < k$, follows from
$$
(xe_{k-1}y)p=\frac{1}{\d^2} \ x(e_{k-1}e_ke_{k-1})yp
= \frac{1}{\d^2} \ x(e_{k-1}e_kpe_{k-1})y=0
$$
by (i) for $x,y\in P_{k-1}$ and $p\in {\cal M}_{k+1}$.

Finally, the (normalized) trace of a minimal projection $p$ in
$I_k$ is $\d^{d_P-k}\vec t_p$ so the assertion about the trace
follows as usual (see [J1]). \qed

\bs{\bf Remarks.} (1) Similarly, the algebras $P_{1,k}$ have a
principal graph $\G'_P$ with the trace vector $\vec s$. We call
$\G'_P$ the dual principal graph. Ocneanu has shown, in the $C^*$
case,  how to associate numerical data encoding the ensuing
embedding of the random walk algebra of $\G'_P$ into that of
$\G_P$. This completely captures the planar algebra structure and
is analogous to choosing local coordinates on a manifold. The same
could be done under the hypotheses of Lemma 3.1.1. The principal
graphs alone do not determine the planar algebra --- for instance
the algebras $(P^\s)^{\z/4\z}$ and
$(P^\s)^{\z/2\z\oplus\z/2\z}$ of 2.8 have the same principal graph
but are readily distinguished by counting fixed points under the
rotation.

(2) In fact, the assumption of nondegeneracy on $P$ in 3.1.1 and
3.1.3 could be replaced by the hypothesis $P_k/I_k$ semisimple.
Then conditions (i)--(iv) could be  used to inductively guarantee
nondegeneracy.

(3) If $P$ had been a $C^*$-planar algebra, we would have a
theorem (3.1.3)${}^*$ with all tr$(p)$ positive, all $p$'s
projections, and the obvious adjoint \ ${}^*$ \  on random walks.

\bs
Theorem 3.1.3 can be used to compute the principal graphs for
Temperley Lieb and the Fuss Catalan algebra (it is the ``middle
pattern" method of [BJ2]). We now illustrate its use by
calculating the principal graph of the nondegenerate planar
algebras coming from Example 2.2. We work over $\Bbb C$ for
convenience, and in the
$C^*$ case to simplify life.

Let $(A,{\op{TR}})$ be a finite-dimensional unital $C^*$-algebra
with normalized faithful (positive) trace TR. The labelling set $L$
is $L_1=A$. We choose a number $\d > 0$ and let $\tau_p$ be TR$(p)$
for projections $p\in A$. A labelled network is then a disjoint
union of smoothly embedded circles, each one containing a (possibly
empty) sequence of 1-boxes labelled by elements in $A$. We define
the partition function $Z$ of such a collection of circles to be
$\d^{\# \ {\text{(circles)}}}\prod_{\text{(circles)}}
{\op{TR}}(a_1a_2\dots a_n)$, where
$a_1a_2\dots a_n$ are the labels  on the given circle, numbered in
order around the circle. The partition function $Z$ is obviously
multiplicative so we define $P^{(A,{\op{TR}})}$ to be the
nondegenerate planar algebra, with obvious $*$-structure, defined
by 1.23. It is linearly spanned by Temperley-Lieb diagrams with a
single labelled box on each string. The relations of Example 2.2
hold, noting that
$Z(\blrighto{p})=\d\tau_p$. For certain values of $\d$ and traces
TR we will compute the principal (and dual principal) graphs of
$P^{A,{\op{TR}})}$ and the Markov trace, and show it to be a planar
$C^*$-algebra. Let us first describe the graphs. Let ${\cal M}=\{p\}$
be a set of minimal projections in $A$, one for each matrix algebra
direct summand and let
$n_p={\op{dim}}(pA)$. Let $\frak S({\cal M})$ be the free
semigroup with identity on ${\cal M}$. Let $\l:{\op{Proj}}\to\Bbb
N\cup\{\infty\}$ be a function, and $W_\l$ be the set of words in
 $\frak S({\cal M})$  which contain no consecutive string of $p$'s
longer than $\l(p)$, for each $p$.

\medskip{\bf Definition 3.1.4.} The graph $\G_{A,\l}$ is the rooted
tree having vertices $W_\l$, with $n_p$ edges between $w$ and
$wp$ for every $p\in\cal M$ with $\{w\cup wp\}\subset W_f$. The
root $*$ is the identity of $\frak S({\cal M})$.

\medskip Thus if $\l(p)=1$ for all $p$ and $A$ is abelian,
$\G_{A,\l}$ is the regular tree of valence $|{\cal M}|$. If
$\l(p)=\infty$ for all
$p$ and $A$ is abelian, the root $*$ of the tree $\G_{A,\l}$ has
valence
$|{\cal M}|$ and all other vertices have valence $|{\cal M}|+1$.
If $A=\Bbb Cp+\Bbb Cq$ and $\l(p)=1$, \ $\l(q)=2$, the tree
$\G_{A,\l}$ is as in Figure~\ref{pic69}.
\[
	\begin{picture}(0,0)%
\epsfig{file=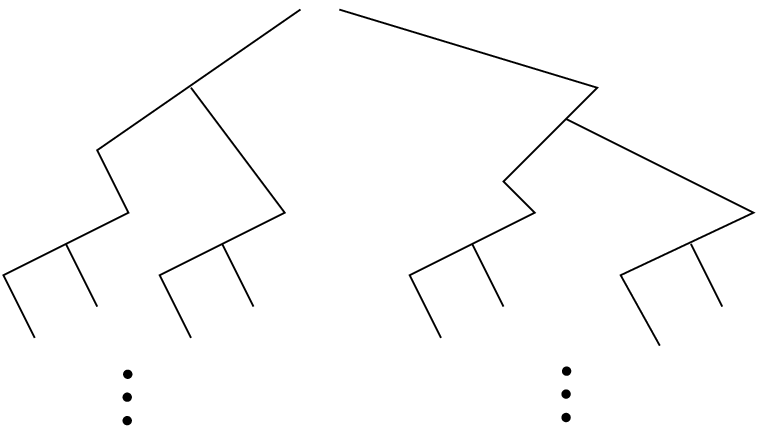}%
\end{picture}%
\setlength{\unitlength}{0.00041700in}%
\begingroup\makeatletter\ifx\SetFigFont\undefined
\def\x#1#2#3#4#5#6#7\relax{\def\x{#1#2#3#4#5#6}}%
\expandafter\x\fmtname xxxxxx\relax \def\y{splain}%
\ifx\x\y   
\gdef\SetFigFont#1#2#3{%
  \ifnum #1<17\tiny\else \ifnum #1<20\small\else
  \ifnum #1<24\normalsize\else \ifnum #1<29\large\else
  \ifnum #1<34\Large\else \ifnum #1<41\LARGE\else
     \huge\fi\fi\fi\fi\fi\fi
  \csname #3\endcsname}%
\else
\gdef\SetFigFont#1#2#3{\begingroup
  \count@#1\relax \ifnum 25<\count@\count@25\fi
  \def\x{\endgroup\@setsize\SetFigFont{#2pt}}%
  \expandafter\x
    \csname \romannumeral\the\count@ pt\expandafter\endcsname
    \csname @\romannumeral\the\count@ pt\endcsname
  \csname #3\endcsname}%
\fi
\fi\endgroup
\begin{picture}(7244,4278)(1479,-4207)
\put(4471,-181){\makebox(0,0)[lb]{\smash{\SetFigFont{12}{14.4}{rm}$*$}}}
\end{picture}

\]
\begin{center}
	Figure 3.1.5
\end{center}
Recall the polynomials $T_n(x)$ of [J1], $T_1\!=\!1$,
$T_2\!=\!1$,
$T_{n+1}=T_n-xT_{n-1}$ and the ``Jones-Wenzl" projections
$f_k\in TL(k)$ with $f^*_k=f^2_k=f_k$, \
$f_ke_i\!=\!0$ for $i\!=\!1,2,\dots ,k-1$ so that $f_k$ are elements
of any (spherical) planar algebra, and
tr$(f_k)=T_{k-2}(\frac{1}{\d^2})$ if $T_j(\frac{1}{\d^2})\neq 0$
for $j < k+2$.

{\bf Theorem 3.1.6} With notation as above, suppose
$\tau_p\d=2\cos \pi/(\l(p)+2)$ ($\tau_p\d\geq 2$ if $\l(p)=\infty$).
Then $P^{(A,{\op{TR}})}$ is a (spherical) $C^*$-planar algebra with
principal and dual principal graphs equal to $\G_{A,\l}$. The
(normalized) trace of the minimal projection in
$P_k^{(A,{\op{TR}})}$ corresponding to the word
$w\!=\!p_1^{m_1}p_2^{m_2}\dots p_r^{m_r}$ (\,$\sum m_i\!=\!k$ and
$p_i\neq p_{i+1}$) is
$\prod^r_{i=1} \tau_{p_i}T_{m_i}(\tau_{p_i}\d)$.

\bigskip
{\sc Proof.} We shall give explicit projections satisfying
conditions (i)--(iv) of 3.1.1. The key observation is that tangles
with a fixed $p$ labelling each string form a Temperley-Lieb
subalgebra $B_p$ with parameter $\d\tau_{p}$ (and identity
$\boxf{p} \ \boxf{p}\dots \boxf{p}$\,). So if $m <\l(p)$ we
consider the projection
\[
	\begin{picture}(0,0)%
\epsfig{file=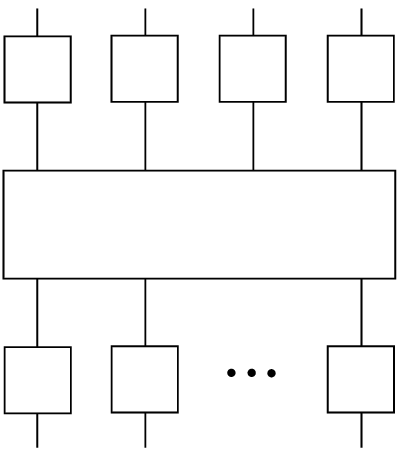}%
\end{picture}%
\setlength{\unitlength}{0.00041700in}%
\begingroup\makeatletter\ifx\SetFigFont\undefined
\def\x#1#2#3#4#5#6#7\relax{\def\x{#1#2#3#4#5#6}}%
\expandafter\x\fmtname xxxxxx\relax \def\y{splain}%
\ifx\x\y   
\gdef\SetFigFont#1#2#3{%
  \ifnum #1<17\tiny\else \ifnum #1<20\small\else
  \ifnum #1<24\normalsize\else \ifnum #1<29\large\else
  \ifnum #1<34\Large\else \ifnum #1<41\LARGE\else
     \huge\fi\fi\fi\fi\fi\fi
  \csname #3\endcsname}%
\else
\gdef\SetFigFont#1#2#3{\begingroup
  \count@#1\relax \ifnum 25<\count@\count@25\fi
  \def\x{\endgroup\@setsize\SetFigFont{#2pt}}%
  \expandafter\x
    \csname \romannumeral\the\count@ pt\expandafter\endcsname
    \csname @\romannumeral\the\count@ pt\endcsname
  \csname #3\endcsname}%
\fi
\fi\endgroup
\begin{picture}(3805,4258)(2379,-4297)
\put(4126,-2236){\makebox(0,0)[lb]{\smash{\SetFigFont{12}{14.4}{rm}$f_m$}}}
\put(3676,-3736){\makebox(0,0)[lb]{\smash{\SetFigFont{12}{14.4}{rm}$p$}}}
\put(2626,-3736){\makebox(0,0)[lb]{\smash{\SetFigFont{12}{14.4}{rm}$p$}}}
\put(5776,-3736){\makebox(0,0)[lb]{\smash{\SetFigFont{12}{14.4}{rm}$p$}}}
\put(3676,-736){\makebox(0,0)[lb]{\smash{\SetFigFont{12}{14.4}{rm}$p$}}}
\put(4726,-736){\makebox(0,0)[lb]{\smash{\SetFigFont{12}{14.4}{rm}$p$}}}
\put(5776,-736){\makebox(0,0)[lb]{\smash{\SetFigFont{12}{14.4}{rm}$p$}}}
\put(2626,-736){\makebox(0,0)[lb]{\smash{\SetFigFont{12}{14.4}{rm}$p$}}}
\end{picture}

\]
\ni where $f_m$ is calculated in Temperley-Lieb with $j$ strings
and $Z(\caright\ )=\tau_p\d$. Now set
\[
	\begin{picture}(0,0)%
\epsfig{file=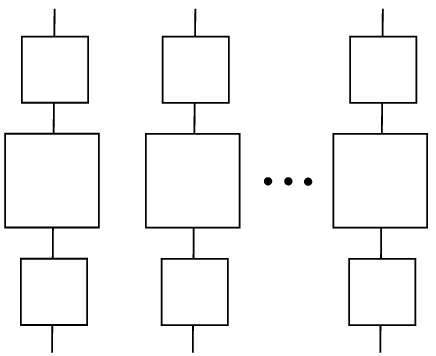}%
\end{picture}%
\setlength{\unitlength}{0.00041700in}%
\begingroup\makeatletter\ifx\SetFigFont\undefined
\def\x#1#2#3#4#5#6#7\relax{\def\x{#1#2#3#4#5#6}}%
\expandafter\x\fmtname xxxxxx\relax \def\y{splain}%
\ifx\x\y   
\gdef\SetFigFont#1#2#3{%
  \ifnum #1<17\tiny\else \ifnum #1<20\small\else
  \ifnum #1<24\normalsize\else \ifnum #1<29\large\else
  \ifnum #1<34\Large\else \ifnum #1<41\LARGE\else
     \huge\fi\fi\fi\fi\fi\fi
  \csname #3\endcsname}%
\else
\gdef\SetFigFont#1#2#3{\begingroup
  \count@#1\relax \ifnum 25<\count@\count@25\fi
  \def\x{\endgroup\@setsize\SetFigFont{#2pt}}%
  \expandafter\x
    \csname \romannumeral\the\count@ pt\expandafter\endcsname
    \csname @\romannumeral\the\count@ pt\endcsname
  \csname #3\endcsname}%
\fi
\fi\endgroup
\begin{picture}(4095,3345)(2229,-3384)
\put(2551,-736){\makebox(0,0)[lb]{\smash{\SetFigFont{12}{14.4}{rm}$p_1$}}}
\put(3901,-736){\makebox(0,0)[lb]{\smash{\SetFigFont{12}{14.4}{rm}$p_2$}}}
\put(5701,-736){\makebox(0,0)[lb]{\smash{\SetFigFont{12}{14.4}{rm}$p_r$}}}
\put(2551,-2911){\makebox(0,0)[lb]{\smash{\SetFigFont{12}{14.4}{rm}$p_1$}}}
\put(3901,-2911){\makebox(0,0)[lb]{\smash{\SetFigFont{12}{14.4}{rm}$p_2$}}}
\put(5701,-2911){\makebox(0,0)[lb]{\smash{\SetFigFont{12}{14.4}{rm}$p_r$}}}
\put(2401,-1786){\makebox(0,0)[lb]{\smash{\SetFigFont{12}{14.4}{rm}$f_{m_1}$}}}
\put(3751,-1786){\makebox(0,0)[lb]{\smash{\SetFigFont{12}{14.4}{rm}$f_{m_2}$}}}
\put(5551,-1861){\makebox(0,0)[lb]{\smash{\SetFigFont{12}{14.4}{rm}$f_{m_r}$}}}
\end{picture}

\]
\ni where we have combined the $m_j$ strings at the top (and
bottom) of $f_{m_j}$ into one. Orientations are completely
forgotten and may be inserted, if required, so as to satisfy
Definition 1.7.

Condition (1) of 3.1.1 is easy: $P_k$ is linearly spanned by
Temperley Lieb diagrams with matrix units  (with the $p$'s among
the diagonal ones) from the simple summands of $A$ in a single box
on each string. If $x$ is such an element, then if
$xe_{k-1}$ is non-zero, the product $p_wxe_{k-1}$ contains
\[
	\begin{picture}(0,0)%
\epsfig{file=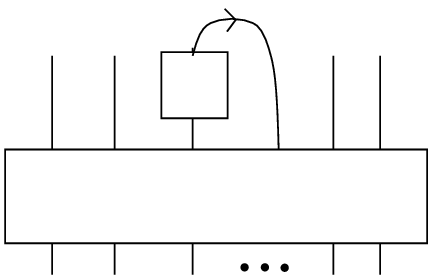}%
\end{picture}%
\setlength{\unitlength}{0.00041700in}%
\begingroup\makeatletter\ifx\SetFigFont\undefined
\def\x#1#2#3#4#5#6#7\relax{\def\x{#1#2#3#4#5#6}}%
\expandafter\x\fmtname xxxxxx\relax \def\y{splain}%
\ifx\x\y   
\gdef\SetFigFont#1#2#3{%
  \ifnum #1<17\tiny\else \ifnum #1<20\small\else
  \ifnum #1<24\normalsize\else \ifnum #1<29\large\else
  \ifnum #1<34\Large\else \ifnum #1<41\LARGE\else
     \huge\fi\fi\fi\fi\fi\fi
  \csname #3\endcsname}%
\else
\gdef\SetFigFont#1#2#3{\begingroup
  \count@#1\relax \ifnum 25<\count@\count@25\fi
  \def\x{\endgroup\@setsize\SetFigFont{#2pt}}%
  \expandafter\x
    \csname \romannumeral\the\count@ pt\expandafter\endcsname
    \csname @\romannumeral\the\count@ pt\endcsname
  \csname #3\endcsname}%
\fi
\fi\endgroup
\begin{picture}(4094,2594)(2229,-2483)
\put(3976,-736){\makebox(0,0)[lb]{\smash{\SetFigFont{12}{14.4}{rm}$e$}}}
\put(4126,-1786){\makebox(0,0)[lb]{\smash{\SetFigFont{12}{14.4}{rm}$P_W$}}}
\end{picture}

\]
\ni with $e$ being a matrix unit. The string containing $e$ enters
the box for $p_w$ either connecting two distinct $f$'s or two
strings of the same one. In the first case the result is zero since
$e$ belongs to precisely one of the direct summands.  In the second
case it is zero because of the properties of the $f$'s.

Condition (ii) follows similarly, noting that a picture like the
above will occur unless all the strings of $x$ are through strings.

For condition (iii) we observe that the ideal $I_k$ is linearly
spanned by $x$'s as above with less than $k$ through-strings, so we
may suppose $x$ is composed of through-strings, each with a matrix
unit in its box. Note that any relation true in $B_p$ is true in
$P_k$. Hence if $j\geq\l(p)$, $f_j=0$ and the identity of $B_p$ is
a linear combination of tangles with less than $\l(p)$
through-strings. So we can suppose that in $x$ there is no sequence
of $\l(p)$ strings in a row whose matrix unit labels are in the
same simple  summand as $p$. Thus by multiplying $x$ to the left
and right by tangles with the appropriate matrix unit labels, we
get $axb=p_w$ for some word $w$ in $\frak S(\cal M)$ of length $k$.
Thus (iii) will follow provided (iv) holds.

We calculate the normalized trace of $p_w$. It is 
\[
	\begin{picture}(0,0)%
\epsfig{file=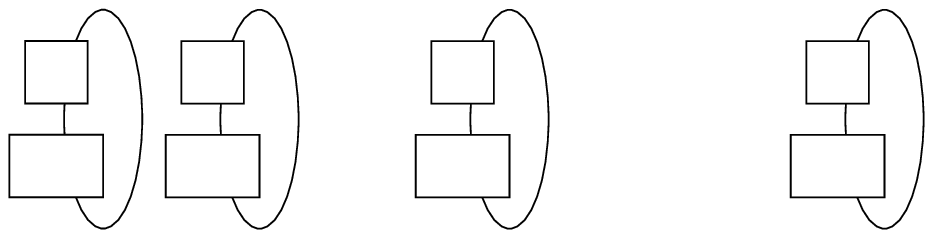}%
\end{picture}%
\setlength{\unitlength}{0.00041700in}%
\begingroup\makeatletter\ifx\SetFigFont\undefined
\def\x#1#2#3#4#5#6#7\relax{\def\x{#1#2#3#4#5#6}}%
\expandafter\x\fmtname xxxxxx\relax \def\y{splain}%
\ifx\x\y   
\gdef\SetFigFont#1#2#3{%
  \ifnum #1<17\tiny\else \ifnum #1<20\small\else
  \ifnum #1<24\normalsize\else \ifnum #1<29\large\else
  \ifnum #1<34\Large\else \ifnum #1<41\LARGE\else
     \huge\fi\fi\fi\fi\fi\fi
  \csname #3\endcsname}%
\else
\gdef\SetFigFont#1#2#3{\begingroup
  \count@#1\relax \ifnum 25<\count@\count@25\fi
  \def\x{\endgroup\@setsize\SetFigFont{#2pt}}%
  \expandafter\x
    \csname \romannumeral\the\count@ pt\expandafter\endcsname
    \csname @\romannumeral\the\count@ pt\endcsname
  \csname #3\endcsname}%
\fi
\fi\endgroup
\begin{picture}(10200,2129)(826,-2176)
\put(5176,-1636){\makebox(0,0)[lb]{\smash{\SetFigFont{12}{14.4}{rm}$\cdots$}}}
\put(7426,-1636){\makebox(0,0)[lb]{\smash{\SetFigFont{12}{14.4}{rm}$)$}}}
\put(7726,-1636){\makebox(0,0)[lb]{\smash{\SetFigFont{12}{14.4}{rm}$= \delta^{-k}\Pi_iZ($}}}
\put(11026,-1636){\makebox(0,0)[lb]{\smash{\SetFigFont{12}{14.4}{rm}$)$}}}
\put(826,-1636){\makebox(0,0)[lb]{\smash{\SetFigFont{12}{14.4}{rm}$\delta^{-k}Z($}}}
\put(2401,-736){\makebox(0,0)[lb]{\smash{\SetFigFont{12}{14.4}{rm}$p_1$}}}
\put(3901,-736){\makebox(0,0)[lb]{\smash{\SetFigFont{12}{14.4}{rm}$p_2$}}}
\put(6301,-736){\makebox(0,0)[lb]{\smash{\SetFigFont{12}{14.4}{rm}$p_r$}}}
\put(6151,-1636){\makebox(0,0)[lb]{\smash{\SetFigFont{12}{14.4}{rm}$f_{m_r}$}}}
\put(3826,-1636){\makebox(0,0)[lb]{\smash{\SetFigFont{12}{14.4}{rm}$f_{m_2}$}}}
\put(2251,-1636){\makebox(0,0)[lb]{\smash{\SetFigFont{12}{14.4}{rm}$f_{m_1}$}}}
\put(9976,-736){\makebox(0,0)[lb]{\smash{\SetFigFont{12}{14.4}{rm}$p_i$}}}
\put(9826,-1636){\makebox(0,0)[lb]{\smash{\SetFigFont{12}{14.4}{rm}$f_{m_i}$}}}
\end{picture}

\]
Now the partition function on $P_m$, restricted to $B_p$,
gives a Markov trace which will be normalized after division by
$Z(\blright{p} )^m=(\d\tau_p)^m$. So
$Z(\blright{f_m})=\d^m\tau_p^m T_m(\tau_{p}\d)$.
Hence tr$(p_w)=\prod^r_{i=1} \tau_{p_i}T_{m_i}(\tau_{p_i}\d)$.

Finally we must calculate the multiplicities dim$(p_vP_kp_w)$ for
$v$ of length $k$ and $w$ of length $k-1$. We must consider
diagrams of the form
\[
	\begin{picture}(0,0)%
\epsfig{file=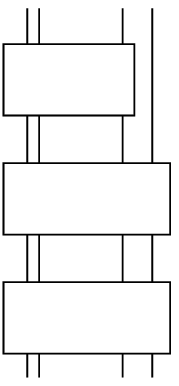}%
\end{picture}%
\setlength{\unitlength}{0.00041700in}%
\begingroup\makeatletter\ifx\SetFigFont\undefined
\def\x#1#2#3#4#5#6#7\relax{\def\x{#1#2#3#4#5#6}}%
\expandafter\x\fmtname xxxxxx\relax \def\y{splain}%
\ifx\x\y   
\gdef\SetFigFont#1#2#3{%
  \ifnum #1<17\tiny\else \ifnum #1<20\small\else
  \ifnum #1<24\normalsize\else \ifnum #1<29\large\else
  \ifnum #1<34\Large\else \ifnum #1<41\LARGE\else
     \huge\fi\fi\fi\fi\fi\fi
  \csname #3\endcsname}%
\else
\gdef\SetFigFont#1#2#3{\begingroup
  \count@#1\relax \ifnum 25<\count@\count@25\fi
  \def\x{\endgroup\@setsize\SetFigFont{#2pt}}%
  \expandafter\x
    \csname \romannumeral\the\count@ pt\expandafter\endcsname
    \csname @\romannumeral\the\count@ pt\endcsname
  \csname #3\endcsname}%
\fi
\fi\endgroup
\begin{picture}(1644,3690)(1479,-3274)
\put(1951,-436){\makebox(0,0)[lb]{\smash{\SetFigFont{12}{14.4}{rm}$P_W$}}}
\put(2101,-1486){\makebox(0,0)[lb]{\smash{\SetFigFont{12}{14.4}{rm}$x$}}}
\put(2026,-2686){\makebox(0,0)[lb]{\smash{\SetFigFont{12}{14.4}{rm}$P_V$}}}
\put(2101,-961){\makebox(0,0)[lb]{\smash{\SetFigFont{12}{14.4}{rm}$\cdots$}}}
\put(2101,-2161){\makebox(0,0)[lb]{\smash{\SetFigFont{12}{14.4}{rm}$\cdots$}}}
\put(2101,-3211){\makebox(0,0)[lb]{\smash{\SetFigFont{12}{14.4}{rm}$\cdots$}}}
\put(2101,164){\makebox(0,0)[lb]{\smash{\SetFigFont{12}{14.4}{rm}$\cdots$}}}
\end{picture}

\]
\ni where $x$ is a Temperley Lieb diagram decorated with matrix
units as before. Arguing on $p_v$, \ $x$ has only through-strings,
$v=w_p$ for some $p\in{\cal M}$, and the first $k-1$ strings of $x$
are labelled by the elements of $\cal M$ in $w$. Thus the diagram
is in fact equal to
\[
	\begin{picture}(0,0)%
\epsfig{file=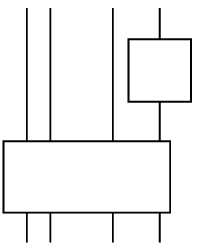}%
\end{picture}%
\setlength{\unitlength}{0.00041700in}%
\begingroup\makeatletter\ifx\SetFigFont\undefined
\def\x#1#2#3#4#5#6#7\relax{\def\x{#1#2#3#4#5#6}}%
\expandafter\x\fmtname xxxxxx\relax \def\y{splain}%
\ifx\x\y   
\gdef\SetFigFont#1#2#3{%
  \ifnum #1<17\tiny\else \ifnum #1<20\small\else
  \ifnum #1<24\normalsize\else \ifnum #1<29\large\else
  \ifnum #1<34\Large\else \ifnum #1<41\LARGE\else
     \huge\fi\fi\fi\fi\fi\fi
  \csname #3\endcsname}%
\else
\gdef\SetFigFont#1#2#3{\begingroup
  \count@#1\relax \ifnum 25<\count@\count@25\fi
  \def\x{\endgroup\@setsize\SetFigFont{#2pt}}%
  \expandafter\x
    \csname \romannumeral\the\count@ pt\expandafter\endcsname
    \csname @\romannumeral\the\count@ pt\endcsname
  \csname #3\endcsname}%
\fi
\fi\endgroup
\begin{picture}(1844,2410)(1479,-3349)
\put(2026,-2686){\makebox(0,0)[lb]{\smash{\SetFigFont{12}{14.4}{rm}$P_V$}}}
\put(2101,-1936){\makebox(0,0)[lb]{\smash{\SetFigFont{12}{14.4}{rm}$\cdots$}}}
\put(2101,-3286){\makebox(0,0)[lb]{\smash{\SetFigFont{12}{14.4}{rm}$\cdots$}}}
\put(2926,-1636){\makebox(0,0)[lb]{\smash{\SetFigFont{12}{14.4}{rm}$e$}}}
\end{picture}

\]
\ni where $e$ is a matrix unit with one subscript fixed.
These diagrams span $p_v P_k p_w$ and the sesquilinear form given
by $(x,y)={\op{tr}}(y^*x)$ is diagonal with non-zero entries.
Hence dim$(p_vP_kp_w)=n_p$. \qed

\bs\ni{\bf 3.2. Duality}

If $P=\cup_k P_k$ is  a planar algebra, the filtered algebras
$\l_n(P)$, where $\l_n(P)_k=P_{n,n+k}$, for fixed $n$, have
natural planar algebra structures. For $n$ even, this is rather
obvious --- just add $n$ straight vertical lines to the left of a
tangle. But if $n$ is odd one must be more careful because of
orientations. In fact $\l_1(P)$ and $P$ are not isomorphic in
general, even as filtered algebras, as one can see from example
2.9. We begin by describing the planar algebra structure on
$\l_1(P)$.

If $T$ is an unlabelled $k$-tangle we define the unlabelled
$(k+1)$-tangle $\tilde T$ to be the tangle consisting of a vertical
straight line from (1,0) to (1,1), and the tangle $T$, with all its
orientations reversed, shifted by 1 in the positive $x$ direction.
Also in $\tilde T$ each internal $p$-box is replaced by a
$(p+1)$-box with the first and last distinguished boundary points
connected by a short curve.  The procedure is illustrated in Figure 3.2.1.
\[
	\input{xfig/pic76}
\]
\begin{center}
Figure 3.2.1
\end{center}
\bs
To each internal box $B$ of $T$ there corresponds in the obvious
way a box $\tilde B$ of $\tilde T$. If $T$ is labelled by
$L=\coprod_{k > 0} L_k$, \
$\tilde T$ will be given the obvious labelling by $\tilde L$, \
$\tilde L_k= L_{k-1}(\tilde L_1=\emptyset)$.

\bs\ni
{\bf Proposition 3.2.2} Let $P=\cup_k P_k$ be a planar algebra
with parameters $\d_1=Z(\caright)$ and
$\d_2=Z(\caleft)$. Assume $P$ is presented on itself by
$\Phi$. Then $\l_1(P)$ is a planar algebra with parameters
$\d_2,\d_1$, presented on itself by $\l_1(\Phi)$ where $\l_1(\Phi)$
is defined by $\l_1(\Phi)(T)=\d^{-p}\Phi(\tilde T)$, \ $p$ being
the number of internal boxes in $T$.

\bs
{\sc Proof.} First note how  the labels in a tangle of
${\cal P}(\l_1(P))$ give valid labels for ${\cal P}(P)$ because of the
inclusion $\l_1(P)_k\subset P_{k+1}$. That $\l_1(\Phi)$ is a
filtered algebra homomorphism is obvious. The annular invariance of
ker $\l_1(\Phi)$ follows immediately from that of $\Phi$, by
representing elements of $P_{1,1+k}$ as linear combinations of
tangles with vertical first string and applying $\sim$ to linear
combinations. Thus $\l_1(P)$ is a general planar algebra.

Now $\l_1(P)_0=P_{1,1}$ and $\l_1(P)_{1,1}=P_{2,2}$ which has
dimension 1 since $P$ is a planar algebra. So
$\l_1(P)$ is planar. The multiplicativity property for $\l_1(P)$
follows immediately from that of $P$, where orientations on
networks without boxes are reversed. \qed

\bs In the next two lemmas, $A$ will be the ${\cal A}(\emptyset)$
element
\[
	\begin{picture}(0,0)%
\epsfig{file=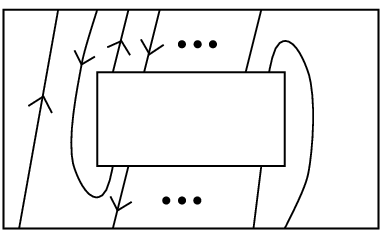}%
\end{picture}%
\setlength{\unitlength}{0.00041700in}%
\begingroup\makeatletter\ifx\SetFigFont\undefined
\def\x#1#2#3#4#5#6#7\relax{\def\x{#1#2#3#4#5#6}}%
\expandafter\x\fmtname xxxxxx\relax \def\y{splain}%
\ifx\x\y   
\gdef\SetFigFont#1#2#3{%
  \ifnum #1<17\tiny\else \ifnum #1<20\small\else
  \ifnum #1<24\normalsize\else \ifnum #1<29\large\else
  \ifnum #1<34\Large\else \ifnum #1<41\LARGE\else
     \huge\fi\fi\fi\fi\fi\fi
  \csname #3\endcsname}%
\else
\gdef\SetFigFont#1#2#3{\begingroup
  \count@#1\relax \ifnum 25<\count@\count@25\fi
  \def\x{\endgroup\@setsize\SetFigFont{#2pt}}%
  \expandafter\x
    \csname \romannumeral\the\count@ pt\expandafter\endcsname
    \csname @\romannumeral\the\count@ pt\endcsname
  \csname #3\endcsname}%
\fi
\fi\endgroup
\begin{picture}(3644,2144)(1779,-2483)
\end{picture}

\]
\ni where the actual number of boundary points is as required by
context.

\bs\ni
{\bf Lemma 3.2.3} If $P$ is a planar algebra, $\pi_A$ defines
a linear isomorphism between $P_k$ and $\l_1(P)_k$, for each $k >0$.

\bigskip
{\sc Proof.} By a little isotopy and the definition of
$P_{1,k+1}$, \ $\pi_A$ is onto. But
\[
	\begin{picture}(0,0)%
\epsfig{file=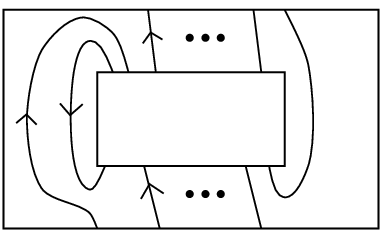}%
\end{picture}%
\setlength{\unitlength}{0.00041700in}%
\begingroup\makeatletter\ifx\SetFigFont\undefined
\def\x#1#2#3#4#5#6#7\relax{\def\x{#1#2#3#4#5#6}}%
\expandafter\x\fmtname xxxxxx\relax \def\y{splain}%
\ifx\x\y   
\gdef\SetFigFont#1#2#3{%
  \ifnum #1<17\tiny\else \ifnum #1<20\small\else
  \ifnum #1<24\normalsize\else \ifnum #1<29\large\else
  \ifnum #1<34\Large\else \ifnum #1<41\LARGE\else
     \huge\fi\fi\fi\fi\fi\fi
  \csname #3\endcsname}%
\else
\gdef\SetFigFont#1#2#3{\begingroup
  \count@#1\relax \ifnum 25<\count@\count@25\fi
  \def\x{\endgroup\@setsize\SetFigFont{#2pt}}%
  \expandafter\x
    \csname \romannumeral\the\count@ pt\expandafter\endcsname
    \csname @\romannumeral\the\count@ pt\endcsname
  \csname #3\endcsname}%
\fi
\fi\endgroup
\begin{picture}(3644,2149)(1779,-2183)
\end{picture}

\]
\ni provides an inverse for $\pi_A$, up to a non-zero scalar.
So $\pi_A$ is an isomorphism. \qed

{\bf Lemma 3.2.4} The subset $S$ of the planar algebra
generates $P$ as a planar algebra iff $\pi_A(S)$ generates
$\l_1(P)$ as a planar algebra.

\bigskip
{\sc Proof.} $(\Rightarrow)$ \  Given a tangle $T$ in
${\cal P}_{1,k+1}(S)$, it suffices to exhibit a tangle $T_A$
in ${\cal P}_k(\pi_A(S))$ with $\Phi(\tilde T_A)$ being a multiple of
$\Phi(T)$. We create $T_A$ from $T$ by eliminating the first
string, reversing all orientations and otherwise changing only in
small neighborhoods of the internal boxes of $G$, sending a box
labelled $R\in S$ in $T$ to the box labelled $\pi_A(R)$ in $T_A$ as
below:
\[
	\begin{picture}(0,0)%
\epsfig{file=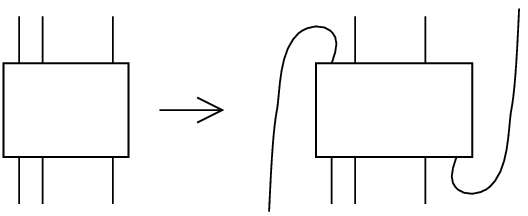}%
\end{picture}%
\setlength{\unitlength}{0.00041700in}%
\begingroup\makeatletter\ifx\SetFigFont\undefined
\def\x#1#2#3#4#5#6#7\relax{\def\x{#1#2#3#4#5#6}}%
\expandafter\x\fmtname xxxxxx\relax \def\y{splain}%
\ifx\x\y   
\gdef\SetFigFont#1#2#3{%
  \ifnum #1<17\tiny\else \ifnum #1<20\small\else
  \ifnum #1<24\normalsize\else \ifnum #1<29\large\else
  \ifnum #1<34\Large\else \ifnum #1<41\LARGE\else
     \huge\fi\fi\fi\fi\fi\fi
  \csname #3\endcsname}%
\else
\gdef\SetFigFont#1#2#3{\begingroup
  \count@#1\relax \ifnum 25<\count@\count@25\fi
  \def\x{\endgroup\@setsize\SetFigFont{#2pt}}%
  \expandafter\x
    \csname \romannumeral\the\count@ pt\expandafter\endcsname
    \csname @\romannumeral\the\count@ pt\endcsname
  \csname #3\endcsname}%
\fi
\fi\endgroup
\begin{picture}(4987,1980)(1179,-1801)
\put(1651,-1636){\makebox(0,0)[lb]{\smash{\SetFigFont{12}{14.4}{rm}$\cdots$}}}
\put(1651,-211){\makebox(0,0)[lb]{\smash{\SetFigFont{12}{14.4}{rm}$\cdots$}}}
\put(4501,-811){\makebox(0,0)[lb]{\smash{\SetFigFont{12}{14.4}{rm}$\pi_A(R)$}}}
\put(1651,-811){\makebox(0,0)[lb]{\smash{\SetFigFont{12}{14.4}{rm}$R$}}}
\put(4696,-211){\makebox(0,0)[lb]{\smash{\SetFigFont{12}{14.4}{rm}$\cdots$}}}
\put(4681,-1621){\makebox(0,0)[lb]{\smash{\SetFigFont{12}{14.4}{rm}$\cdots$}}}
\end{picture}

\]
Then by definition $\tilde T_A$ will be exactly like $T$ except
near its boxes where it will look as below:
\[
	\begin{picture}(0,0)%
\epsfig{file=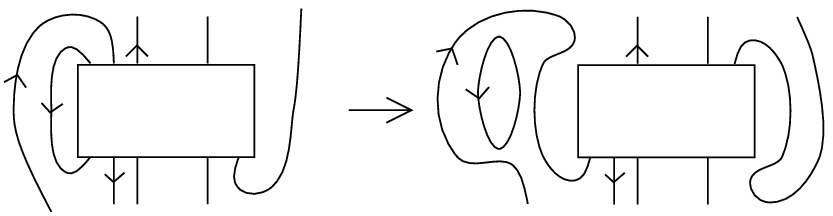}%
\end{picture}%
\setlength{\unitlength}{0.00041700in}%
\begingroup\makeatletter\ifx\SetFigFont\undefined
\def\x#1#2#3#4#5#6#7\relax{\def\x{#1#2#3#4#5#6}}%
\expandafter\x\fmtname xxxxxx\relax \def\y{splain}%
\ifx\x\y   
\gdef\SetFigFont#1#2#3{%
  \ifnum #1<17\tiny\else \ifnum #1<20\small\else
  \ifnum #1<24\normalsize\else \ifnum #1<29\large\else
  \ifnum #1<34\Large\else \ifnum #1<41\LARGE\else
     \huge\fi\fi\fi\fi\fi\fi
  \csname #3\endcsname}%
\else
\gdef\SetFigFont#1#2#3{\begingroup
  \count@#1\relax \ifnum 25<\count@\count@25\fi
  \def\x{\endgroup\@setsize\SetFigFont{#2pt}}%
  \expandafter\x
    \csname \romannumeral\the\count@ pt\expandafter\endcsname
    \csname @\romannumeral\the\count@ pt\endcsname
  \csname #3\endcsname}%
\fi
\fi\endgroup
\begin{picture}(7912,1980)(652,-1802)
\put(1877,-812){\makebox(0,0)[lb]{\smash{\SetFigFont{12}{14.4}{rm}$\pi_A(R)$}}}
\put(2072,-212){\makebox(0,0)[lb]{\smash{\SetFigFont{12}{14.4}{rm}$\cdots$}}}
\put(2057,-1622){\makebox(0,0)[lb]{\smash{\SetFigFont{12}{14.4}{rm}$\cdots$}}}
\put(6874,-214){\makebox(0,0)[lb]{\smash{\SetFigFont{12}{14.4}{rm}$\cdots$}}}
\put(6859,-1624){\makebox(0,0)[lb]{\smash{\SetFigFont{12}{14.4}{rm}$\cdots$}}}
\put(6902,-812){\makebox(0,0)[lb]{\smash{\SetFigFont{12}{14.4}{rm}$R$}}}
\put(4173,-586){\makebox(0,0)[lb]{\smash{\SetFigFont{12}{14.4}{rm}$\Phi$}}}
\end{picture}

\]
\ni Thus $\Phi(\tilde T_A)$ is a multiple of $\Phi(T)$.

$(\Leftarrow)$ \ Given $x\in P$, \ $\pi_A(\tangle{x})$
is in $\l_1(P)$ so by hypothesis it is the image under
$\l_1(\Phi)$ of a linear combination of tangles  labelled by
elements of $\pi_A(S)$ which are in turn images under $\Phi$ of
tangles  labelled by
elements of $S$ (up to nonzero scalars). Using the tangle of Lemma
3.2.3 to invert $\pi_A$, we are done. \qed

\bs By iterating the procedure $P\mapsto\l_1(P)$, we see that all
the $\l_n(P)$ have natural planar algebra structures, but observe
that all the $\l_{2n}(P)$ are isomorphic to $P$ as planar algebras
via the endomorphism (often called ``le shift de deux") defined by
adding two straight vertical strings to the left of a tangle. We
leave the details to the reader.

The planar algebra $\l_1(P)$ is said to be the {\it dual}
of the planar algebra $P$, and we have $\l_1(\l_1(P))\simeq P$
as planar algebras.

In the case of Example 2.9, $P^G_2$ is the group algebra
$\Bbb C G$ and $\l_1(P^G)_2$ is $\ell^\i(G)$. The tangle $\pi_A$
gives a {\it linear} isomorphism between the two. Thus planar
algebra duality extends the duality between a finite group and its
dual object.

\bs\bs\ni{\bf 3.3. Reduction and cabling}

We give two ways to produce new planar algebras from a given one.
The first is a reduction process which makes ``irreducible" planar
algebras --- those with dim $P_1=1$ --- the focus of study.
A planar algebra is not reconstructible in any simple way from its
irreducible reductions, though, as can be seen from example 2
where the irreducible reductions would be trivial.

Given a general planar algebra $P$ presented on itself by $\Phi$,
and an idempotent $p\in  P_1$, we define the {\it reduced} general
planar algebra $pPp$ (by $p$) as follows: for each $k$ we let
$p_k$ be the element
\[
        \begin{picture}(0,0)%
\epsfig{file=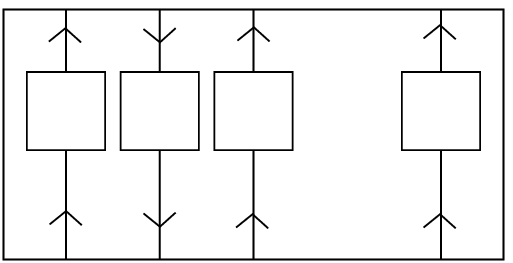}%
\end{picture}%
\setlength{\unitlength}{0.00041700in}%
\begingroup\makeatletter\ifx\SetFigFont\undefined
\def\x#1#2#3#4#5#6#7\relax{\def\x{#1#2#3#4#5#6}}%
\expandafter\x\fmtname xxxxxx\relax \def\y{splain}%
\ifx\x\y   
\gdef\SetFigFont#1#2#3{%
  \ifnum #1<17\tiny\else \ifnum #1<20\small\else
  \ifnum #1<24\normalsize\else \ifnum #1<29\large\else
  \ifnum #1<34\Large\else \ifnum #1<41\LARGE\else
     \huge\fi\fi\fi\fi\fi\fi
  \csname #3\endcsname}%
\else
\gdef\SetFigFont#1#2#3{\begingroup
  \count@#1\relax \ifnum 25<\count@\count@25\fi
  \def\x{\endgroup\@setsize\SetFigFont{#2pt}}%
  \expandafter\x
    \csname \romannumeral\the\count@ pt\expandafter\endcsname
    \csname @\romannumeral\the\count@ pt\endcsname
  \csname #3\endcsname}%
\fi
\fi\endgroup
\begin{picture}(4844,2444)(1179,-2783)
\put(1711,-1351){\makebox(0,0)[lb]{\smash{\SetFigFont{12}{14.4}{rm}$p$}}}
\put(3526,-1351){\makebox(0,0)[lb]{\smash{\SetFigFont{12}{14.4}{rm}$p$}}}
\put(5326,-1351){\makebox(0,0)[lb]{\smash{\SetFigFont{12}{14.4}{rm}$p$}}}
\put(2911,-1141){\makebox(0,0)[lb]{\smash{
\SetFigFont{12}{14.4}{rm}$p$
}}}
\put(4276,-1336){\makebox(0,0)[lb]{\smash{\SetFigFont{12}{14.4}{rm}$\cdots$}}}
\end{picture}

\]
\ni of $P_k$ (illustrated when $k$ is odd) and we set
$(pPp)_k =p_k(P_k)p_k$ with identity $p_k$, and unital
inclusion  $P_k\hookrightarrow P_{k+1}$ given by
$p_kxp_k\hookrightarrow p_{k+1}xp_{k+1}$ (note
$p_{k+1}p_k=p_kp_{k+1}=p_{k+1}$). We make $pPp$ into a planar
algebra on $P$ as follows. Given a tangle $T\in{\cal P(P)}$, define
the tangle ${}_pT_p\in {\cal P}(P)$ by inserting  $\boxa{p}$ in
every string of $TG$. Then
$p\Phi p:{\cal P}(L)\to pPp$ is $p\Phi p(T)=\Phi(pTp)$.
Since $p$ is idempotent, $p\Phi p$ is a filtered algebra
homomorphism with annular invariance, obviously surjective, so
$pPp$ is a general planar algebra. Planarity is inherited from $P$
and $pPp$ has parameters
$Z(\blrighto{p}\ )$ and $Z(\bllefto{p}\ )$, provided these
are nonzero. If $P$ is a $C^*$-planar algebra and $p$ is a
projection
$(p=p^2=p^*)$, $pPp$ is clearly also a $C^*$-planar algebra,
spherical if $P$ is.

\medskip{\bf Note.} We have used the canonical labelling set for $P$
to define $pPp$. If were given another specific
other labelling set $L$, it is not clear that the homomorphism
obtained in the same way from $\cal P(L)$ to $pPp$ is surjective. We
do not have any example of this phenomenon.

\bs To use the reduction process we require dim $P_1 >1$. But even
``irreducible" planar
 algebras can yield this situation by cabling, i.e. grouping
several strings together. The term is borrowed from knot theory.
Given a general planar algebra $P$ we define the $n^{\text{th}}$
cabled (general) planar algebra $\cal C_n(P)$ by ${\cal C}_n(P)_k=P_{nk}$ 
which we endow with a planar algebra structure as
follows. If $\Phi$ presents $P$ on itself, we define
${\cal C}_n(\Phi):{\cal P}({\cal C}_n(P))\to{\cal C}_n(P)$ by taking a
labelled $k$-tangle $T$ in ${\cal P}({\cal C}_n(P))$ and constructing
an $nk$-tangle $\tilde T$ in ${\cal P}(P)$ with the same labels on
boxes, but where every boundary point in $TG$ (both on internal and
external boxes) is replaced by $n$ boundary points, and of course
orientations alternate. Every string in $T$ is then replaced by $n$
parallel strings. The procedure is illustrated in Figure 3.4.1,
where $k=3$ and $n=2$.
\[
        \input{xfig/pic82}
\]
\begin{center}
Figure 3.4.1
\end{center}
\bs\ni Then we define ${\cal C}_n(\Phi)(T)$ to be $\Phi(\tilde T)$.
It is clear that ${\cal C}_n(\Phi)$ is a general planar algebra,
connected and multiplicative if $P$ is, with parameters
$(\d_1\d_2)^{\big[\frac n2\big]}\d_1^{n-\big[\frac n2\big]}$,
$(\d_1\d_2)^{\big[\frac n2\big]}\d_2^{n-\big[\frac n2\big]}$,
(where $\d_1$ and $\d_2$ are the parameters of the planar algebra
$P$ and $\big[\frac n2\big]$ is the  integer part of
$\big[\frac n2\big]$). Also ${\cal C}_n(P)$ is a $C^*$-planar algebra
if $P$ is, spherical if $P$ is.

\bs\bs\ni{\bf 3.4. Tensor product}

Let $P^1=\cup_k P^1_k$ and $P^2=\cup_k P^2_k$ be
general planar algebras. We will endow the filtered algebra
$P^1\otimes P^2=\cup_k P^1_k\otimes P^2_k$ with a
general planar algebra structure on the labelling set
$L=\coprod_{i\geq 1} P^1_i\x P^2_i$. Consider $P^1$ and $P^2$
presented on themselves by $\Phi_1$ and $\Phi_2$ respectively.
First define a linear map ${\cal L}:{\cal P}(L)\to{\cal P}(P_1)\otimes
{\cal P}(P_2)$ by ${\cal L}(T)=T_1\otimes T_2$ where $T$ is a tangle
labelled by $f:{\text{Boxes}}(T)\to P_1\x P_2$, \ $T_i$ have the
same unlabelled tangle as $T$ and they are labelled by $f$ composed
with the projection $P_1\x P_2\to P_i$, $i=1,2$. This ${\cal L}$ is
well defined since the isotopy classes of labelled tangles are a
basis of ${\cal P}(L)$. Now define the presenting map
$\Phi^{P_1\otimes P_2}:{\cal P}(L)\to P_1\otimes P_2$ by
$$
\Phi^{P_1\otimes P_2}=(\Phi_1\otimes\Phi_2)\circ{\cal L} \ .
$$
This obviously gives a homomorphism of filtered algebras.
It is surjective because we may consider the tangle \
$\tangle{}$ labelled by $x\x y$ for an arbitrary pair
$(x,y)$ in $P^1_k\x P^2_k$. This will be sent on to $x\otimes y\in
P^1_k\x P^2_k$. Thus we need only show the annular invariance of
ker $\Phi^{P_1\otimes P_2}$. But if $A$ is an annular tangle in
${\cal A}(L)$ then it is easy to see that 
${\cal L}\circ\pi_A=\pi_{A_1}\otimes \pi_{A_2}\circ{\cal L}$ where $A_1$ and
$A_2$ are the annular tangles having the same unlabelled tangle as
$A$ but labelled by the first and second components of the labels
of $A$ respectively. So if
\begin{eqnarray*}
&x\in{\op{ker}}\in \Phi^{P_1\otimes P_2}\ , \\
&\Phi^{P_1\otimes P_2} \pi_A(x)  = \Phi_1\otimes\Phi_2
(\pi_{A_1}\otimes \pi_{A_2}({\cal L}(x)))  
=\pi_{A_1}\otimes \pi_{A_2}(\Phi_1\otimes\Phi_2({\cal L}(x)))\\ 
&= \pi_{A_1}\otimes \pi_{A_2}(\Phi^{P_1\otimes P_2}(x))=0
\end{eqnarray*}
So $P_1\otimes P_2$ is a general planar algebra.

It is clear that $P_1\otimes P_2$ is connected iff both $P_1$ and
$P_2$ are and that $Z_{P_1\otimes P_2}=Z_{P_1}Z_{P_2}$ in the sense
that a network labelled with $P_1\x P_2$ is the same as two
networks  labelled with $P_1$ and $P_2$ respectively.
Thus $P_1\otimes P_2$ is a planar algebra if $P_1$ and $P_2$ are.
Moreover, nondegeneracy, $*$ structure, positivity and sphericity
are inherited by $P_1\otimes P_2$ from $P_1$ and $P_2$. So the
tensor product of two $C^*$-planar algebras is a $C^*$-planar
algebra.

\bs{\bf Notes.} (i) By representing elements of $P_1$ and $P_2$  by
tangles, one may think of the tensor product planar structure as
being a copy of $P_1$ and one of $P_2$ sitting in boxes on
parallel planes, with no topological interaction between them.
This corresponds to presenting
$P_1\otimes P_2$  on the labelling set ${\cal P}(P_1)\x{\cal P}(P_2)$
in the obvious way.

(ii) It is clear that $P_1\otimes P_2\approx P_2\otimes P_1$ as
(general) planar algebras.

\bs\bs\ni{\bf 3.5. Free product}

The notion of free product of planar algebras was developed in
collaboration with D.~Bisch and will be presented in a future
paper. The free product $P^1\x P^2$ of two planar algebras
$P^1$ and $P^2$ is by definition the subalgebra of the tensor
product linearly spanned by (images of) $T_1$ consisting of a pair
$T_1\in P^1_k$, \ $T_2\in P^2_k$ which can be drawn in a single
$2k-1$-box, with boundary points in pairs, alternately
corresponding to $P^1$ and $P^2$, so that the two tangles
$T_1$ and $T_2$ are disjoint. An example is given by Figure 3.5.1
where we have used the ``colours" to indicate boundary points
belonging to $P_1$ and $P_2$.
\[
        \input{xfig/pic83}
\]
S. Gnerre defined in [Gn] a notion of free product using
detailed connection calculations in the paragroup 
formalism.

\bs The most interesting result so far of the work with Bisch is a
formula, at least for finite dimensional  $C^*$-planar algebras,
for the Poincar\'e series of $P^1\x P^2$ in terms of those of
$P^1$ and $P^2$, using Voiculescu's free multiplicative convolution
([V]).

\bs
\ni{\bf 3.6. Fusion algebra}

The reduced subalgebras of the cables on a planar algebra form a
``fusion algebra" along the lines of [Bi]. This is to be
thought of as part of the {\it graded algebra} structure given by a
planar algebra and will be treated in detail in a future paper with
D.~Bisch.

\bs\bs
\noindent{\large\bf  4. Planar algebras and subfactors}

In this section we show that the centralizer tower for an extremal
finite index type II${}_1$ subfactor admits the structure of a
 spherical $C^*$-planar algebra, and vice versa. We need
several results from subfactors, some of which are well known.

\medskip\ni{\bf 4.1 \ Some facts about subfactors}

Let $N\subset M$ be II${}_1$ factors with $\tau^{-1}=[M:N] < \i$.
We adopt standard notation so that $M_i$, $i=0,1,2,\dots$ is the
tower of [J1] with $M_0=M$, \ $M_1=\< M,e_n\>$, \
$M_{i+1}=\< M_i,e_{i+1}\>$ where $e_i: L^2(M_{i-1})\to
L^2(M_{i-2})$ is orthogonal projection, $e_N=e_1$. Let $B=\{b\}$ be
a finite subset of $M$ (called a basis) with
\[
(4.1.1)\ \ \   \sum_{b\in B}  be_1b^*=1 \ . 
\]

Then by [PP] and [Bi],
$$
(4.1.2)\ \ \ \ \ x(\sum_b be_1b^*)= (\sum_b be_1b^*)x=x \quad {\mbox{for}} \ x\in M.
$$
$$
(4.1.3)\ \ \ \ \ \sum_b bb^*=\tau^{-1} 
$$
$$
(4.1.4)\ \ \ \ \ \sum_b bE_N(b^*x)=x=\sum_b E_N(xb)b^* \quad {\mbox{for}} \ x\in M.
$$
$$
(4.1.5)\ \ \ \  {\mbox{For}} \ x\in N' \ ({\mbox{on}} \ L^2(M)), \ \ \
\tau\sum_b bxb^*=E_{M'}(x) \ . 
$$

\noindent{\bf 4.1.6} \ Recall that the subfactor is called {\it
extremal} if the normalized traces on $N'$ and $M$ coincide on
$N'\cap M$ in which case the traces on $N'\cap M_k$, realized on
$L^2(M)$, coming from $M_k$ and $N'$ coincide for all $k$,
and $E_{M'}(e_1)=\tau$.

\medskip\ni {\bf 4.1.7} \ By standard convex averaging procedures on
$L^2(M_k)$, given a finite subset $X=\{x\}$ of $M_k$ and
$\e >0$, there is a finite set $U=\{u\}$ of unitaries in $M$ and
$\l_u\in\Bbb R^+$, \ $\sum_{u\in U}\l_u=1$ with
$$
\|\sum_u\l_u uxu^*-E_{M'}(x)\|_2 < \e \quad {\mbox{for}} \ x\in X.
$$
Further averaging does not make the estimate worse, so by averaging
again with $\sum_u \l_u {\op{Ad}} \ u^*$ and gathering together
repeated terms if necessary, we may assume
$U^*=U$ and $\l_{u^*}=\l_u$.

It will be convenient to renormalize the $e_i$'s so we set
$E_i=\d e_i$ with $\d^2\tau=1$ $(\d >0)$. (Note that there is a
very slight notational clash with $\S$2, but we will show it to be
consistent.) We then have the formulae
$$
(4.1.8)\ \ \ \ \ E^2_i =\d E_i \ , \quad
E_iE_j=E_jE_i \  \ {\text{if}} \ |i-j|\geq 2, \quad
E_iE_{i\pm 1}E_i=E_i \ , 
$$
$E_1xE_1\!=\d E_N(a)E_1$ and $\sum_b bE_1e_2E_1b^*\!=1$
so that $\{bE_1\mid b\!\in\! B\}$ is a basis for $M_1$ over~$M$.

\bs{\bf Definition 4.1.9.} For $k=1,2,3,\dots$ let
$v_k=E_kE_{k-1}\dots E_1\in N'\cap M_k$.\newline
(Note $v_kxv^*_k=\d E_N(x)E_k$ for $x\in M$, \
$v^*_kv_k=\d E_1$.)

\bs\ni
{\bs Theorem 4.1.10} If $x_i$, $i=1,\dots ,k+1$ are elements
of $M$ then $x_1v_1x_2v_2\dots v_kx_{k+1}=x_1v^*_kx_2v^*_{k-1}\dots
v^*_1x_{k+1}$, and the map $x_1\otimes x_2\otimes\dots\otimes
x_{k+1}\mapsto x_1v_1x_2v_2\dots v_kx_{k+1}$ defines an $M-M$
bimodule isomorphism, written $\th$, from $M\otimes_{\!N}
M\otimes_{\!N}
\dots\otimes_{\!N} M$ (with $k+1$ \ $M$'s), written
$\bigotimes^{k+1}_N M$, onto $M_k$.

\bs
{\sc Proof.} See [J5], Corollary 11.

\bs Recall that if $R$ is a ring and $B$ is an $R-R$ bimodule, an
element $b$ of $B$ is called {\it central} if
$rb=br$ \ $\forall \ r\in R$.

\bs\ni
{\bf Corollary 4.1.11} The centralizer $N'\cap M_k$ is
isomorphic under $\th$ to the vector space $V_{k+1}$
of central vectors in the $N-N$ bimodule $\bigotimes^{k+1}_N M$.

\bs
We now define the most interesting ``new" algebraic ingredient of
subfactors seen from the planar point of view. It is the
``rotation", known to Ocneanu and rediscovered by the author in
specific models. See also [BJ1].

\bs{\bf Definition 4.1.12.} For $x\in M_k$ we define
$$
\rho(x)=\d^2 E_{M_k} (v_{k+1} \ E_{M'}(xv_{k+1}))
$$

\bs\ni
{\bf Proposition 4.1.13} $\rho(M_k)\subseteq N'\cap M_k$
and if $B$ is a basis, $\rho$ coincides on $N'\cap M_k$
with $r:M_k\to M_k$, \ $r(x)=E_{M_k} (v_{k+1} \sum_{b\in B}
bx v_{k+1}b^*)$.

\bigskip
{\sc Proof.} This is immediate from 4.1.5. \qed

\bs\ni
{\bf Lemma 4.1.14} With $\th$ as above,
$$
\th^{-1}r\th(x_1\otimes x_2\otimes\dots\otimes x_{k+1}) =
\sum_{b\in B} E_N(bx_1)x_2\otimes x_3\otimes\dots\otimes
x_{k+1}\otimes b^*
$$

\bigskip
{\sc Proof.}
\begin{eqnarray*}
r(x_1v_1x_2v_2\dots v_kx_{k+1}) &  = &
\sum_b E_{M_k}  (v_{k+1}bx_1 v^*_{k+1}x_2v^*_k\dots v^*_1b^*) \\
& = & \d\sum_b E_{M_k} (E_{k+1}E_N(bx_1)x_2v_k^*\dots v^*_1b^*) \ \
{\text{(by (4.1.9))}} \\
& = & \sum_b E_N(bx_1)x_2v^*_k \dots v^*_1b^* \\
& = & \sum_b\th(E_N(bx_1)x_2\otimes x_3\otimes\dots \otimes
x_{k+1}\otimes b^*)
\end{eqnarray*}
\qed
Note that the rotation on $\bigotimes^{k+1}_{\Bbb C} M$ does {\it
not} pass to the quotient $\bigotimes^{k+1}_N M$, however we have
the following.

{\bf Lemma 4.1.15} Suppose $N\subset M$ is extremal, then
$$
\rho(\th(x_1\otimes x_2\otimes\dots\otimes x_{k+1})) =
E_{N'}(\th(x_2\otimes x_3\otimes\dots\otimes x_{k+1}
\otimes x_1))
$$

\bigskip
{\sc Proof.} If $x=\th(x_1\otimes x_2\otimes\dots\otimes
x_{k+1})$ and $y=\th(x_2\otimes x_3\otimes\dots\otimes
x_1)$, it suffices to show that
tr$(\rho(x)a)={\op{tr}}(ya)$ for all $a$ in $N'\cap M_k$.
Let $\e > 0$ be given and choose by 4.1.7 a finite set $U$
of unitaries in $M$, with $U=U^*$, and $\l_u\in\Bbb R^+$,
$\sum_{u\in U}\l_u=1$, \ $\l_{u^*}=\l_u$ so that
$$
\|\sum_u \l_u v_{k+1} uxv_{k+1}u^* -
v_{k+1}E_{M'} (xv_{k+1})\|_{{}_2} < \e \ ,
$$
and, by extremality, $\|\sum_u \l_u uE_1u^*-\d^{-1}\|_{{}_2} < \e$.
So, if $y\in M$,
$$
(4.1.16)\ \ \ \ |\sum_u \l_u uE_N(u^*y)-\tau y\|_{{}_2} < \e\|y\|
\qquad {\text{(by (4.1.8) and}} \ \|ab\|_2\leq \|a\| \
\|b\|_{{}_2}.  
$$
So
\begin{eqnarray*}
|{\op{tr}}((\rho(x)-y)a| & =  & |{\op{tr}}(\d^2 E_{M_k}(v_{k+1}E_{M'}(xv_{k+1}))-y)a| \\ 
& < & |{\op{tr}}(\d^2 E_{M_k}(v_{k+1}\sum_u \l_u uxv_{k+1}u^*) -y)a| +\d^2\e\|a\|_{{}_2} \\ 
& = & |{\op{tr}}(\d^2 \sum_u \l_u E_N(ux_1)\th(x_2\otimes x_3\otimes\dots\otimes x_{k+1}\otimes u^*)-y)a| +\d^2\e\|a\|_{{}_2}\\ 
&  & \qquad {\mbox{(as in the proof of 4.11)}}\\
& = & |\d^2 \big(\big({\op{tr}} \sum_u \l_u \th
(x_2\otimes
x_3\otimes\dots\otimes x_{k+1}\otimes u^* E_N(u x_1))-\th
(y)\big)a\big)+\d^2\e\|a\|_{{}_2} \\
&  &  \qquad  {\mbox{(since}} \ a\in N')\\
&  = & |\d^2{\op{tr}}(\th(x_2\otimes 
x_3\otimes\dots\otimes x_{k+1}\otimes (\sum_u
\l_u uE_N(u^*x_1)-x_1))a)\|_{{}_2} +\d^2\e\|a\|_{{}_2} \\
&  &  \qquad  {\mbox{(since}} \ U=U^*, \ \ \l_u=\l_{u^*}) \ .
\end{eqnarray*}
For fixed $x_1,x_2,\dots ,x_{k+1}$ and $a$, this can clearly be
made as small as desired by 4.1.16, by choosing $\e$ small. \qed

\medskip
{\bf Corollary 4.1.17} If $x\in N$ and $\xi\in L^2(M_k)$,
\ $\rho(x\xi-\xi x)=0$.

\bigskip
{\sc Proof.} By its definition, $\rho$ extends to a bounded linear
map from $L^2(M_k)$ to itself, so it suffices to show the formula
for $\xi$ of the form $\th(x_1\otimes
x_2\otimes\dots\otimes x_{k+1})$. But if $n\in N$, \
$E_{N'}(\th(x_2\otimes\dots\otimes x_{k+1}\otimes nx_1))=
E_{N'}(\th(x_2\otimes\dots\otimes x_{k+1}n\otimes x_1))$,  so by
4.1.15 we are done. \qed

\medskip
{\bf Theorem 4.1.18}
If $N \subset M$ is extremal,
$\rho^{k+1}={\op{id}}$ \ on \ $N'\cap M_k$.

{\sc Proof.}  Recall that if ${\cal H}$ is an $N-N$ bimodule
(correspondence as in [Co]) then
$\<\eta,x\xi-\xi x\>=\<x^*\eta-\eta x^*,\xi\>$ for $x\in N$
and $\xi,\eta\in{\cal H}$. So $\eta$ is central iff it is orthogonal
to commutators.

Hence Lemma 4.1.15 reads
$$
\rho(\th(x_1\otimes x_2\otimes\dots\otimes x_{k+1}))  =
\th(x_2\otimes x_3\otimes\dots\otimes x_1) +\xi
$$
where $\xi\in \kappa$, which we define to be the closure of the
linear span of commutators in the $N-N$ correspondence $L^2(M_k)$.
Applying $\rho$ to both sides of this equation, and 4.1.17, we
obtain
$$
\rho(\th(x_1\otimes \dots\otimes x_{k+1}))  -
\th(x_3\otimes x_4\otimes\dots\otimes x_1\otimes x_2) \in \kappa
$$
and continuing,
$$
\rho^{k+1}(\th(x_1\otimes\dots\otimes x_{k+1}))
-\th(x_1\otimes x_2\otimes\dots\otimes x_{k+1})\in \kappa \ .
$$
So, by linearity, if $x\in N'\cap M_k$, \
$\rho^{k+1}(x)-x\in \kappa$. But both $\rho^{k+1}(x)$ and $x$ are
central, so orthogonal to $\kappa$ and
$\rho^{k+1}(x)=x$. \qed

\bs We now define five types of maps between centralizers using the
isomorphism $\th$. Choose a basis $B$.

\bs {\bf Definition 4.1.19.} If \ $x=
x_1\otimes x_2\otimes\dots\otimes x_k\in \bigotimes^k_N M$,
\begin{verse}
(1) For $j=2,3,\dots ,k$, \newline
$a_j(x)=\d(x_1\otimes x_2\otimes\dots\otimes x_{j-1}
E(x_j)\otimes\dots\otimes x_k)\in \bigotimes^{k-1}_N M$,\newline
 and \ $a_1(x)=\d E(x_1)x_2 \otimes\dots\otimes x_k$.

(2) For $j=2,3,\dots ,k$, \newline
$\mu_j(x)=x_1\otimes x_2\otimes\dots\otimes x_{j-1}x_j
\otimes\dots\otimes x_k)\in \bigotimes^{k-1}_N M$

(3) For $j=2,3,\dots ,k$, \newline
$\eta_j(x)=x_1\otimes x_2\otimes\dots\otimes x_{j-1}\otimes 1
\otimes x_j
\otimes\dots\otimes x_k)\in \bigotimes^{k-1}_N M$,\newline
and \ $\eta_1(x) = 1\otimes x_1\otimes\dots\otimes x_k$.

(4) For $j=1,2,\dots ,k$,\newline
$\kappa_j(x)=\d^{-1}\sum_{b\in B}
x_1\otimes\dots\otimes x_jb\otimes b^*\otimes
x_{j+1}\otimes\dots\otimes x_k)\in \bigotimes^{k-1}_N M$

(5) If \ $\th(c)\in N'\cap M_{n-1}$ \ define \
$\a_{j,k}\otimes^k_N M\to\otimes^{k+n}_N M$ \ by\newline
$\a_{j,c}(x)=x_1\otimes\dots\otimes x_{j-1}\otimes c
\otimes x_j\otimes\dots\otimes x_k$\newline
$(\a_{1,c}(x)=c\otimes x$, $\a_{k+1,c}(x)=x\otimes c$.
Note also $\a_{j,1}=\eta_j$.)
\end{verse}
Note all these maps are $N$ middle linear (for (5)
this requires $b$ to be central; for (4) we use 4.1.2),
so they are defined on the tensor product over $N$.
They are all $N\!-\!N$ bimodule for maps so they preserve central
vectors and are thus defined on the space $V_k$ of $N$-central
elements of $\otimes^k_N M$ ($=\th^{-1}(N'\cap M_{k-1}))$.
We will use the same notation for the restrictions
of these maps to the $V_k$. Note that $\kappa_j$ does not depend on the
basis, indeed $\kappa_j(x)=\mu_{j+1}\a_{j+1,{\text{id}}}(x)$ where
$\th$(id) is the identity of $N'\cap M$, is a basis-independent
formula for $\kappa_j$.

\bs
\ni{\bf Lemma 4.1.20} If $c\in V_k$ and $d\in V_{k'}$, then
\begin{verse}
(i) \ $\a_{d,j+k}\a_{c,i}=\a_{c,i}\a_{d,j}$ \ for \ 
$i\leq j$.

(ii) For \ $i < j$, \ $a_{j-1}a_i=a_ia_j$,  \
 $\mu_{j-1}\mu_i=\mu_i\mu_j$. \newline
For \ $i\leq j$, \ $\eta_{j+1}\eta_i=\eta_i\eta_j$, \
$\kappa_{j+1}\kappa_i=\kappa_i\kappa_j$.

(iii) For $i\leq j$,
\begin{eqnarray*}
   a_{j+k}\a_{c,i}=\a_{c,i}a_j &
a_{j-1}\a_{c,j}=\a_{c,j-1}a_{i-1} \\
\mu_{j+k+1}\a_{c,i}=\a_{c,i}\mu_{j+1} &
\mu_{i-1}\a_{c,j}=\a_{c,j-1}\mu_{i-1} \\
\eta_{j+k}\a_{c,i}=\a_{c,i}\eta_j &
\eta_{i-1}\a_{c,j-1}=\a_{c,j}\eta_{i-1} \\
\kappa_{j+k}\a_{c,i}=\a_{c,i}\kappa_j & \kappa_{i-1}\a_{c,j}=\a_{c,j+1}
\kappa_{i-1}
\end{eqnarray*}

(iv) For $i < j$,
\begin{eqnarray*}
\mu_j a_i =a_i\mu_{j+1} & \mu_ia_j=a_{j-1}\mu_i\\
\eta_{j-1}a_i=a_i\eta_j &   \eta_{i+1}a_j=a_{j+1}\eta_{i+1}\\
\kappa_{j-1}a_i=a_i\kappa_j &  \kappa_ia_j=a_{j+1}\kappa_i\\
\eta_{j-1}\mu_i=\mu_i\eta_j &  \eta_i\mu_j=\mu_{j+1}\eta_i\\
\kappa_{j-2}\mu_i=\mu_i\kappa_{j-1} &  \kappa_i\mu_j=\mu_{j+1}\kappa_i\\
\kappa_j\eta_i =\eta_i\kappa_{j-1} &  \kappa_i\eta_j=\eta_{j+1}\kappa_i
\end{eqnarray*}

(v) 
\begin{eqnarray*}
a_i\kappa_i={\op{id}} \\
\mu_i\eta_i={\op{id}} \\
a_{i+1}\kappa_i={\op{id}} \\
\mu_{i+1}\eta_i={\op{id}} 
\end{eqnarray*}

(vi) 
\begin{eqnarray*}
a_j\eta_j=\d{\op{id}} \\
\mu_{j+1}\kappa_j=\d{\op{id}}.
\end{eqnarray*}
\end{verse}
(These identities hold when $i, \ j$ and $k$ are such that all the
maps involved are defined by 4.1.19.)

\bigskip
{\sc Proof.} Almost all cases of identities (i)--(iv) are trivial
as they  can be written so as to involve distant tensor product
indices: thus they just amount to a renumbering.
The ones that involve some interaction between the tensor product
components are
$$
\mu_i\mu_i=\mu_i\mu_{i+1} \ , \ \
\kappa_{i+1}\kappa_i=\kappa_i\kappa_i \ , \ \
\kappa_{i-1}\mu_i=\mu_i\kappa_i \ , \ \
\kappa_i\mu_{i+1}=\mu_{i+2}\kappa_i \ .
$$
These all follow easily from associativity of multiplication and\\
$\dsize{\sum_{b\in B}\!b\!\otimes\! b^*x\!=\!\sum_{b\in B}}xb \otimes b^*$
for $x\in M$, which is 4.1.2.

\begin{eqnarray*}
{\mbox{For (v):}}\ a_i\kappa_i& = &{\op{id}}\ {\mbox{follows from}} \sum_{b\in B} E_N(xb)b^*=x \ \ \ (4.1.4)\\
\mu_i\eta_i & = & {\op{id}}\  {\mbox{follows from}}  x1=x\\
a_{i+1}\kappa_i & = & {\op{id}}\  {\mbox{follows from}} 
\sum_{b\in B} bE(b^*x)=x \ \ \ \ (4.1.4)\\
\mu_{i+1}\eta_i & = &{\op{id}}\   {\mbox{follows from}}\  x1=x
\end{eqnarray*}

\begin{eqnarray*}
{\mbox{For (vi):}}\   a_j\eta_j & = & \d{\op{id}}\   {\mbox{follows from}}\ E_N(1)=1\\
\mu_{j+1}\kappa_j&=&\d{\op{id}}\  {\mbox{follows from}}\ \sum_{b\in B} bb^*=\d^2{\op{id}} \ \ \ (4.1.3)
\end{eqnarray*}
\qed

\medskip
\ni{\bf Lemma 4.1.21} \ If \ $x\in M$, \
$2\leq r\leq k$, then \
$v_kxv_r=v_{r-2}xv_k$ (where $v_0=1$), and
$v_kxv_1=\d E_N(x)v_k$.

\bigskip
{\sc Proof.} Simple manipulation of 4.1.8 and 4.1.9. \qed

\medskip
\ni{\bf Lemma 4.1.22} \ If \ $x\in \otimes^k_N M$, \ then \
\begin{verse}
(i)  $\th^{-1}E_{M_{k-2}}\th(x)=\d^{-1}
a_{m+1}(x)$ \ if \ $k$ is odd, $k=2m+1$.

(ii) $\th^{-1}E_{M_{k-2}}\th(x)=\d^{-1}\mu_{m+1}(x)$ \ if \ $k$ is even, $k=2m$.
\end{verse}

\bigskip
{\sc Proof.} Let $x$ be of the form $x_1\otimes
x_2\otimes\dots\otimes x_k$. Then
$$
E_{M_{k-2}}(\th(x))=
E_{M_{k-2}}(x_1v_1x_2v_2\dots v_{k-1}x_k) =
\d^{-1}x_1v_1x_2v_2\dots x_{k-2}v_{k-2}x_{k-1}v_{k-2}x_k
$$

\ni${\underline{\mbox{Case (i)}}}$. If \ $k=2m+1$ we may apply
4.1.21
$m\!-\!1$ times to obtain
\begin{eqnarray*}
E_{M_{k-2}}(\th(x)) & = \d^{-1}x_1v_1\dots x_mv_m x_{m+1}v_1
x_{m+2}v_{m+1}\dots x_{k-1}v_{k-2}x_k\\
&= x_1v_1\dots x_mv_m E(x_{m+1})x_{m+2}v_{m+1}\dots
x_{k-1}v_{k-2}x_k \\
&= \d^{-1}\th(a_{m+1}(x))
\end{eqnarray*}
${\underline{\mbox{Case (ii)}}}$. If \ $k=2m$ we  apply 4.1.21
$m\!-\!1$ times to obtain
\begin{eqnarray*}
E_{M_{k-2}}(\th(x)) & = \d^{-1}x_1v_1\dots x_mv_m x_{m+1}v_0
 x_{m+2}v_{m+1}\dots x_{k-1}v_{k-2}x_k\\
&= \d^{-1}\th(\mu_{m+2}(x)) \quad {\text{(since}} \ v_0\!=\!1).
\end{eqnarray*}
\qed

\bs
\ni{\bf Lemma 4.1.23} \ If \ $y\in \otimes^k_N M$ \ and \
$x,z\in M$, then
$$
xv^*_k\th(y)v_{k+1}z=\th(x\otimes y\otimes z)
$$

\bs
{\sc Proof.} Simple commutation of $E_i$ with $M$ and $E_j$'s
$j\leq i-2$. \qed

\bs
In the next lemma, let $\chi_j=a_j\mu_{j+1}$, \
$j=1,2,\dots k-1$.

\bs
\ni{\bf Lemma 4.1.24} \ If \ $x,y\in \otimes^k_N M$, \ then \
$\th(x)\th(y)=$
\begin{verse}
(i) $\th(\chi_{m+1}\chi_{m+2}\dots\chi_k(x\otimes_N y))$ \
if $k$ is even, $k=2m$.

(ii) $\th(\mu_{m+2}\chi_{m+2}\dots\chi_k(x\otimes_N y))$ \
if $k$ is odd, $k=2m+1$.
\end{verse}

\bs
{\sc Proof.}  By induction on $k$. Let $x=x_1\otimes
x_2\dots\otimes x_{k+1}$, \ $y=y_1\otimes
y_2\dots\otimes y_{k+1}$, \ then

\begin{eqnarray*}
\th(x)\th(y) & = & x_1v^*_k kx_2v^*_{k-1}\dots v^*_1x_{k+1}y_1v_1y_2v_2 \dots y_kv_ky_{k+1} \\
&= & \d x_1v^*_{k-1} E_{M_{k-2}}(\th(x_2\otimes x_3\otimes\dots\otimes x_{k+1})\th(y_1\otimes\dots\otimes y_k))
v_ky_{k+1}
\end{eqnarray*}

(i) If $k$ is even, $k=2m$, by the inductive hypothesis we have
\begin{eqnarray*}
\th(x)\th(y) & = &\d x_1v^*_{k-1} E_{M_{k-2}} (\th(\chi_{m+1}\chi_{m+2}\dots\chi_k(x_2\otimes\dots\otimes x_{k+1}\otimes y_1\otimes\dots\otimes y_k)))v_ky_{k+1} \\
& = & x_1v^*_{k-1} \th(\mu_{m+1}\chi_{m+1}\dots\chi_k (x_2\otimes\dots\otimes y_k))v_ky_{k+1} \ \ \mbox{(by\  Lemma\ 4.1.22)}\\
& = &\th(x_1\otimes (\mu_{m+1}\chi_{m+1}\chi_{m+2}\dots \chi_k(x_2\otimes\dots\otimes y_k)\otimes y_{k+1})) \ \ \mbox{(by Lemma 4.1.23)}\\
& = &\th(\mu_{m+2}\chi_{m+2}\chi_{m+3}\dots \chi_{k+1}(x_1\otimes x_2\otimes x_{k+1}\otimes y_1 \otimes\dots\otimes y_{k+1}))
\end{eqnarray*}

(ii) If $k$ is odd, $k=2m+1$,
\begin{eqnarray*}
\th(x)\th(y) & = & x_1 v^*_{k-1} E_{M_{k-2}} (\th(\mu_{m+2}\chi_{m+2}\chi_{m+3}\dots\chi_k (x_2\otimes\dots x_{k+1}\otimes y_1 \otimes\dots\otimes y_k)))v_ky_{k+1} \\
& & \mbox{(by the induction hypothesis)} \\
& = & x_1v^*_{k-1}\th(a_{m+1}\mu_{m+2}\chi_{m+2}\chi_{m+3}\dots \chi_k(x_2\otimes\dots\otimes y_k))v_ky_{k+1}\\
& & \mbox{(by Lemma 4.1.22)}\\
& = & \th(x_1\otimes \chi_{m+1}\chi_{m+2}\dots \chi_k(x_2\otimes\dots\otimes y_k)\otimes y_{k+1})) \quad \mbox{(by Lemma 4.1.23)}\\
& = &\th(\chi_{m+2}\chi_{m+3}\dots \chi_{k+1}(x_1\otimes x_2 \dots\otimes y_{k+1})) \ .
\end{eqnarray*}
It only remains to check the formula for $k\!=\!1, \ (m\!=\!0)$.
Then $\th(x)=x$, \ $\th(y)=y$ and the formula reads
\[
\th(x)\th(y)=xy=\th(\mu_2(x\otimes y))=\th(xy) \ .
\]
\qed

\bs
\ni{\bf Lemma 4.1.25}
For $m=1,2,\dots$, and $x_{m+1},x_{m+2},\dots x_{2m}\in M$,
$$
v_mv_{m+1}x_{m+1}v_{m+2}x_{m+2}v_{m+3}\dots
v_{2m}x_{2m}=v_mx_{m+1}v_{m+1}x_{m+2}\dots v_{2m-1}x_{2m} \  .
$$

\bigskip
{\sc Proof.} Induction on $m$.

For $m=1$ the formula reads $E_1E_2E_1x_2=E_1x_2$ which is correct.
Now suppose the formula holds for $m$, then
$$
v_{m+1}v_{m+2}x_{m+2} =v_{m+1}x_{m+2}E_{m+2}E_{m+1}\dots E_3 \qquad
\mbox{(by 4.1.8)}
$$
so
$$
v_{m+1}v_{m+2}x_{m+2} v_{m+3}x_{m+3}\dots v_{2m+2}x_{2m+2} =
(v_{m+1}x_{m+2}) V_mV_{m+1} y_{m+1}V_{m+2}y_{m+2}\dots
V_{2m}y_{2m}
$$
where \  $V_n=E_{n+2}E_{n+1}\dots E_3$ \ and \
$y_n=E_2E_1x_{n+2}$.

We may now apply the inductive hypothesis to the subfactor
$M_1\subset M_2$ (for which the $E_i$'s are just those for
$N\subset M$, shifted by 2), to obtain
$$
v_{m+1}x_{m+2}V_my_{m+1}\dots V_{2m}y_{2m}=
v_{m+1}x_{m+2}v_{m+2}x_{m+3}\dots v_{2m+1}x_{2m+2}
$$
\qed

\bs
\ni{\bf Corollary 4.1.26} \ With notation as above, for \
$m=2,3,4,\dots,$
\begin{eqnarray*}
& E_mE_{m-1}\dots E_2(v_{m+1}x_{m+1})(v_{m+2}x_{m+2})\dots
(v_{2m-1}x_{2m-1}) \\
&\quad =
(v_m\dots x_{m+1})(v_{m+1}x_{m+2})\dots (v_{2m-2}x_{2m-1})
\end{eqnarray*}

\bs
{\sc Proof.} Write \ $V_n=E_{n+1}E_{n-1}\dots E_2$, \
$y_n=E_1x_{n+1}$ and apply the 4.1.25 to the subfactor $M\subset
M_1$. \qed

\bs
\ni{\bf Lemma 4.1.27} \ For
$p=1,2,\dots,$
\[
\sum_{b_1,b_2,\dots b_p\in B} (b_1v_1)(b_2v_2)\dots
(b_pv_p)v_{p+1}v_{p+2} b^*_{p}v_{p+3}b^*_{p-1} v_{p+4}\dots
v_{2p+1}b^*_1 = \d^p E_{2p+1}
\]

\bs
{\sc Proof.}
By induction on $p$. For $p=1$ the formula reads
$$
\sum_{b\in B} bE_1E_2E_1E_3E_2E_1b^* = E_3 \sum_{b\in B}
bE_1b^* = \d E_3  \ ,
$$
which is correct. Now observe that
$$
\sum_b bv_p v_{p+1}v_{p+2}b^*v_{p+3}=\d
E_p E_{p-1}\dots E_2 V_{p-1}V_p V_{p+1} E_2E_1
$$
where $V_n=E_{n+2}E_{n+1}\dots E_3$, and that
$$
b_1v_1b_2v_2\dots b_{p-1}v_{p-1}E_pE_{p-1}\dots E_2=
y_1V_1y_2V_2\dots y_{p-2}V_{p-2}y_{p-1}
$$
where $y_i=b_iE_1E_2$, so that
\begin{eqnarray*}
&\dsize{\sum_{b_1,\dots ,b_p\in B} b_1v_1b_2v_2\dots b_pv_p v_{p+1}v_{p+2}b^*_p v_{p+3}b^*_{p-1} v_{p+4}\dots v_{2p+1}b^*_1} \\
& = \dsize{\d \sum_{y_1\dots y_{p-1}\in BE_1E_2} y_1V_1y_2V_2\dots y_{p-2}V_{p-2}y_{p-1}V_{p-1}V_pV_{p+1} y^*_{p-1}V_{p+2}\dots V_{2p-1}y^*_1}
\end{eqnarray*}
so since $\{ bE_1E_2\mid b\in B\}$ is a basis for $M_2$ over $M_1$,
we are through by induction. \qed

\ni{\bf Corollary 4.1.28} \ For
$p=1,2,\dots,$
$$
\sum_{b_1,b_2,\dots ,b_{p+1}\in B}
(b_1v_1)(b_2v_2)\dots (b_pv_p)b^*_p b_{p+1}v_{p+1}
b^*_{p+1}v_{p+2} b^*_{p-1}\dots v_{2p}b^*_1 =\d^{p+1}E_{2p}
$$

\bigskip
{\sc Proof.} Observe that
$$
\sum_{b_p,b_{p+1}\in B} b_pv_p b^*_pb_{p+1} v_{p+1}
b^*_{p+1}v_{p+2} = \d^2 V_{p-1}V_p V_{p+1}E_1 \ ,
$$
where $V_n=E_{n+1}E_n\dots E_2$, so the left-hand side of the
equation becomes
$$
\d^2\sum_{y_1,\dots ,y_{p-1}} y_1V_1y_2V_2\dots y_{p-1}
V_{p-1} V_pV_{p+1} y^*_{p-1}V_{p+2}\dots V_{2p-1}y^*_1
$$
where $y_n=b_nE_1$, and this is $\d^{p+1}E_{2p}$ by 4.1.27 applied
to the subfactor $M\subset M_1$ with basis $\{bE_1\mid b\in B\}$.
\qed

\bs\bs\noindent{\bf 4.2 \ Subfactors give planar algebras}

\smallskip We keep the notation of $\S$4.1.
The next theorem legitimizes the use of pictures to prove
subfactor results.

\ni{\bf Theorem 4.2.1} Let $N\subset M$ be an extremal type
$II_1$ subfactor with $[M:N]^{\frac 12}=\d <\i$. For each $k$ let
$P^{N\subset M}_k=N'\cap M_{k-1}$ (isomorphic via $\th$ to $V_k$,
i.e.
$N$-central vectors in $\bigotimes^k_N M$). Then
$P^{N\subset M}=\bigcup_k P^{N\subset M}_k$ has a
spherical
$C^*$-planar algebra structure (with labelling set $P^{N\subset
M}$) for which
$\Phi(\tangle{x})=x$ and, suppressing the
presenting map~$\Phi$,
\begin{verse}
(i) For $i=1,2,\dots k-1,\boxj = E_i$ \bs

(ii) $\boxh{x} = \d E_{M'}(x)$ , \ \
 $\boxl{x} \ = \d E_{M_{k-2}}(x)$ \bs

(iii) $\tangle{x}\ \  \vln = \tangle{x}$ 
(where on the right, $x$ is considered as an element of $M_{k+1}$)
\bs

(iv) $Z(\boxrloop{x})=\d^k {\op{tr}}(x)$ \ \
$(x\in P^{N\subset M}_k)$
\end{verse}
\bs
Moreover, any other spherical planar algebra structure $\Phi'$ with
$\Phi'(\tangle{x})=x$ and {\rm (i),(ii),(iv)} for
$\Phi'$ is equal to $\Phi$.

\bigskip
{\sc Proof.}  The idea of the proof is fairly simple but it will
involve a lot of details, so we begin with an informal description
of the idea. We must show how to associate an element $\Phi(T)$ in
$N'\cap M_{k-1}$ to a tangle $T$  whose boxes are labelled by
elements of the appropriate $N'\cap M_j$. An example with $k=5$ is
given in Figure 4.2.2.
\[
	\begin{picture}(0,0)%
\epsfig{file=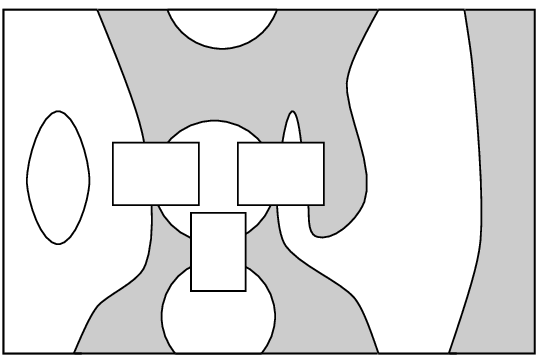}%
\end{picture}%
\setlength{\unitlength}{0.00041700in}%
\begingroup\makeatletter\ifx\SetFigFont\undefined
\def\x#1#2#3#4#5#6#7\relax{\def\x{#1#2#3#4#5#6}}%
\expandafter\x\fmtname xxxxxx\relax \def\y{splain}%
\ifx\x\y   
\gdef\SetFigFont#1#2#3{%
  \ifnum #1<17\tiny\else \ifnum #1<20\small\else
  \ifnum #1<24\normalsize\else \ifnum #1<29\large\else
  \ifnum #1<34\Large\else \ifnum #1<41\LARGE\else
     \huge\fi\fi\fi\fi\fi\fi
  \csname #3\endcsname}%
\else
\gdef\SetFigFont#1#2#3{\begingroup
  \count@#1\relax \ifnum 25<\count@\count@25\fi
  \def\x{\endgroup\@setsize\SetFigFont{#2pt}}%
  \expandafter\x
    \csname \romannumeral\the\count@ pt\expandafter\endcsname
    \csname @\romannumeral\the\count@ pt\endcsname
  \csname #3\endcsname}%
\fi
\fi\endgroup
\begin{picture}(5144,3344)(1779,-3983)
\put(4051,-3286){\makebox(0,0)[lb]{\smash{
\SetFigFont{12}{14.4}{rm}$R_3$
}}}
\put(4276,-2386){\makebox(0,0)[lb]{\smash{\SetFigFont{12}{14.4}{rm}$R_2$}}}
\put(3076,-2386){\makebox(0,0)[lb]{\smash{\SetFigFont{12}{14.4}{rm}$R_1$}}}
\end{picture}

\]
\begin{center}
Figure 4.2.2
\end{center}
\bs Shade the regions black and white and observe that a smooth
oriented curve starting and ending on the left-hand boundary, and
missing the internal boxes, will generically pass through a
certain number of black regions. Number the connected components of
the intersection of the curve with the black
 regions $1,2,\dots
n$ in the order they are crossed.
The regions on the curve will be used to index the tensor product
components in
$\bigotimes^n_N M$.

We will start with our curve close to the boundary so that it
crosses no black regions, and allow it to bubble outwards until it
is very close to the outside boundary at which point it will cross
$k$ black regions. As the curve bubbles out, it will pass through
non-generic  situations with respect to the strings of the tangle,
and it will envelop internal boxes. At generic times we will
associate an $N$-central element of $\bigotimes^n_N M$ with the
curve. As the curve passes through exceptional situations we will
change the element of $\bigotimes^n_N M$ according to certain
rules, the main one of which being that, when the curve envelops a
box labelled by a tensor, we will insert that label into the
tensor on the curve at the appropriate spot, as illustrated in
Figure~\ref{pic85}
\[
        \begin{picture}(0,0)%
\epsfig{file=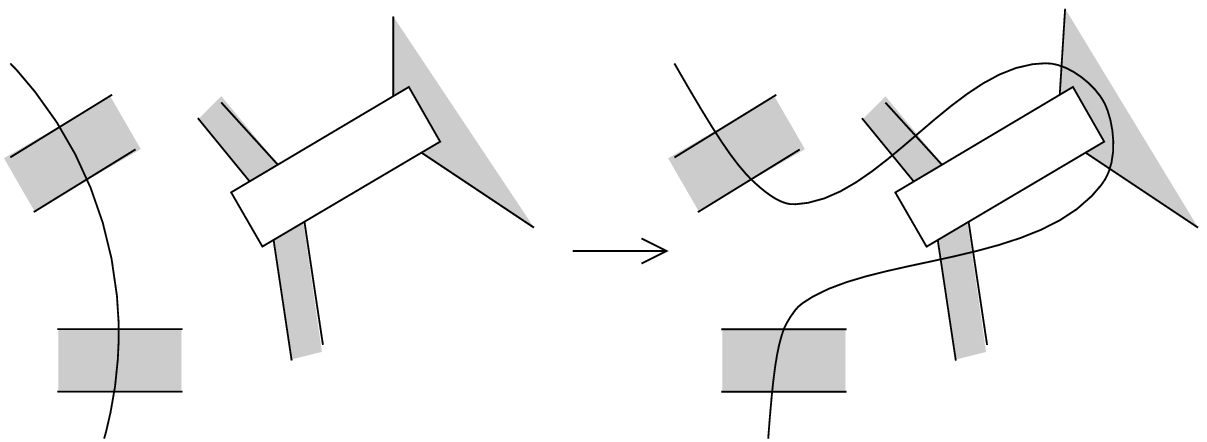}%
\end{picture}%
\setlength{\unitlength}{0.00041700in}%
\begingroup\makeatletter\ifx\SetFigFont\undefined
\def\x#1#2#3#4#5#6#7\relax{\def\x{#1#2#3#4#5#6}}%
\expandafter\x\fmtname xxxxxx\relax \def\y{splain}%
\ifx\x\y   
\gdef\SetFigFont#1#2#3{%
  \ifnum #1<17\tiny\else \ifnum #1<20\small\else
  \ifnum #1<24\normalsize\else \ifnum #1<29\large\else
  \ifnum #1<34\Large\else \ifnum #1<41\LARGE\else
     \huge\fi\fi\fi\fi\fi\fi
  \csname #3\endcsname}%
\else
\gdef\SetFigFont#1#2#3{\begingroup
  \count@#1\relax \ifnum 25<\count@\count@25\fi
  \def\x{\endgroup\@setsize\SetFigFont{#2pt}}%
  \expandafter\x
    \csname \romannumeral\the\count@ pt\expandafter\endcsname
    \csname @\romannumeral\the\count@ pt\endcsname
  \csname #3\endcsname}%
\fi
\fi\endgroup
\begin{picture}(11519,4162)(1704,-3976)
\put(10652,-1862){\makebox(0,0)[lb]{\smash{
\SetFigFont{12}{14.4}{rm}$x\otimes y\otimes z$
}}}
\put(9452,-962){\makebox(0,0)[lb]{\smash{\SetFigFont{12}{14.4}{rm}$a$}}}
\put(9902,-3287){\makebox(0,0)[lb]{\smash{\SetFigFont{12}{14.4}{rm}$b$}}}
\put(4276,-1861){\makebox(0,0)[lb]{\smash{
\SetFigFont{12}{14.4}{rm}$x\otimes y\otimes z$
}}}
\put(3076,-961){\makebox(0,0)[lb]{\smash{\SetFigFont{12}{14.4}{rm}$a$}}}
\put(3526,-3286){\makebox(0,0)[lb]{\smash{\SetFigFont{12}{14.4}{rm}$b$}}}
\end{picture}

\]
\begin{center}
Figure 4.2.3: Tensor on bubbling curve: $\dots\otimes           
a\otimes b\otimes\dots
\ \to \ \dots\otimes a\otimes x \otimes
y\otimes z \otimes b\otimes\dots$
\end{center}

\bs\ni When the curve arrives very close to the boundary it will
have associated to it a central element of $\bigotimes^k_N M$ which
gives an element of $N'\cap M_{k-1}$ via $\th$. This element will
be $\Phi(T)$.

This strategy meets several obstacles.
\begin{verse}
(1) We must show that $\Phi(T)$ is well defined -- note that the
insertions of 4.2.3 are {\bf not} well defined for the tensor
product over $N$.

(2) $\Phi(T)$ must be central. This will require either
enveloping the boxes only starting from the white region touching
the first boundary point, or projecting onto central vectors at
each step. We will adopt the former policy.

(3)$\Phi(T)$ must be independent of isotopy of $T$ and the
choice of the path. Since we must show isotopy invariance
eventually, we might as well suppose that the tangles are in a
convenient standard form, since any two ways of arriving at that
standard form from a given $T$ will only differ by an isotopy.
\end{verse}

We begin the formal proof by describing the standard form.

\medskip A $k${\it -picture} (or just ``picture"
 if the value of $k$ is clear) will be the intersection of the unit
square $[0,1]\x [0,1]$ in the $x\!-\!y$ plane with a system of
smooth curves, called strings, meeting only in finitely many
isolated  singularities, called ``cusps", where $2m$ strings meet
in a single point. The strings must meet the boundary of $[0,1]\x
[0,1]$ transversally in just $2k$ points on the boundary line
$[0,1]\x
\{1\}$.  A cusp $(x,y)$ will be said to be in {\it standard form} \
if, in some neighborhood of $(x,y)$ the $y$-coordinates of points
on the strings are all greater than $y$. A picture may always be
shaded black and white with the $y$-axis being part of the boundary
of a white region.

A picture $\Theta$ will be said to be {\it standard} \ if
\begin{verse}
(i) All its cusps are in standard form, and the region
immediately below the cusp is white.

(ii) The $y$-coordinate, restricted to strings, has only
generic singularities, i.e., isolated maxima and minima.

(iii) The $y$-coordinates of all cusps and  all maxima and
minima are distinct. This set will be written
$\frak S(\Theta)=\{ y_1,y_2,y_3,\dots ,y_c\}$ with
$y_i < y_{i+1}$ for $1\leq i < c$.
\end{verse}

\ni An example of a standard $k$-picture with $k=6$ is in Figure
4.2.4.
\[
	\begin{picture}(0,0)%
\epsfig{file=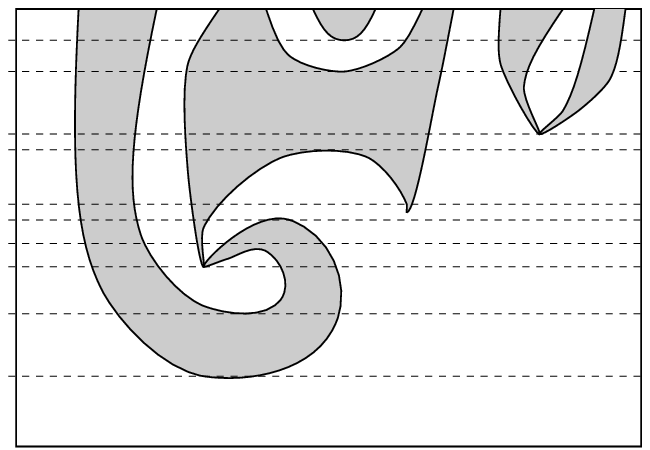}%
\end{picture}%
\setlength{\unitlength}{0.00041700in}%
\begingroup\makeatletter\ifx\SetFigFont\undefined
\def\x#1#2#3#4#5#6#7\relax{\def\x{#1#2#3#4#5#6}}%
\expandafter\x\fmtname xxxxxx\relax \def\y{splain}%
\ifx\x\y   
\gdef\SetFigFont#1#2#3{%
  \ifnum #1<17\tiny\else \ifnum #1<20\small\else
  \ifnum #1<24\normalsize\else \ifnum #1<29\large\else
  \ifnum #1<34\Large\else \ifnum #1<41\LARGE\else
     \huge\fi\fi\fi\fi\fi\fi
  \csname #3\endcsname}%
\else
\gdef\SetFigFont#1#2#3{\begingroup
  \count@#1\relax \ifnum 25<\count@\count@25\fi
  \def\x{\endgroup\@setsize\SetFigFont{#2pt}}%
  \expandafter\x
    \csname \romannumeral\the\count@ pt\expandafter\endcsname
    \csname @\romannumeral\the\count@ pt\endcsname
  \csname #3\endcsname}%
\fi
\fi\endgroup
\begin{picture}(7072,4244)(151,-4283)
\put(676,-3736){\makebox(0,0)[lb]{\smash{\SetFigFont{12}{14.4}{rm}$y_1$}}}
\put(676,-3136){\makebox(0,0)[lb]{\smash{\SetFigFont{12}{14.4}{rm}$y_2$}}}
\put(676,-2761){\makebox(0,0)[lb]{\smash{\SetFigFont{12}{14.4}{rm}$y_3$}}}
\put(151,-2236){\makebox(0,0)[lb]{\smash{\SetFigFont{12}{14.4}{rm}$y_5$}}}
\put(676,-811){\makebox(0,0)[lb]{\smash{\SetFigFont{12}{14.4}{rm}$y_9$}}}
\put(676,-1561){\makebox(0,0)[lb]{\smash{\SetFigFont{12}{14.4}{rm}$y_7$}}}
\put(676,-2386){\makebox(0,0)[lb]{\smash{\SetFigFont{12}{14.4}{rm}$y_4$}}}
\put(676,-1186){\makebox(0,0)[lb]{\smash{\SetFigFont{12}{14.4}{rm}$y_8$}}}
\put(676,-1936){\makebox(0,0)[lb]{\smash{\SetFigFont{12}{14.4}{rm}$y_6$}}}
\put(526,-436){\makebox(0,0)[lb]{\smash{\SetFigFont{12}{14.4}{rm}$y_{10}$}}}
\end{picture}

\]
\begin{center}
	Figure 4.2.4:  A standard 6-picture with
$\frak S(\Theta)=\{ y_1,y_2,y_3,\dots ,y_{10}\}$.
\end{center}
\bs
A standard $k$ picture $\Theta$ will be {\it labelled} \ if there is a
function from the cusps of $\Theta$ to $\coprod_m V_m$ so that a
cusp where $2m$ strings meet is assigned an element of
$V_m$ (or $N'\cap M_{m-1}$, via $\th$). We now describe  how to
associate an element $Z_{\Theta}$ of $N'\cap M_{k-1}$ to a
labelled standard $k$-picture $\Theta$, using the operators of
definition 4.1.19.

Let $\Theta$ be a labelled standard $k$-picture with
$\frak S(\Theta)=\{ y_1,y_2,y_3,\dots ,y_c\}$. We define a locally
constant function $Z:[0,1]\backslash\frak S(\Theta)\to V_{k(y)}$,
where $k(y)$ is the number  of distinct intervals of
$[0,1]\times \{y\}$ which are the connected components of its
intersection with all the black regions of $\Theta$. For instance, if
$\Theta$ is as in Figure 4.2.4 and $y_8 < y < y_9$ then $k(y)=4$.
Obviously $k(y)$ is locally constant and $k(y)=k$ for $y_c < y < 1$.
For $y < y_1$, \ $V_0=\Bbb C$ and we set $Z(y)=1$. There are five
possible ways for $Z$ to change as $y$ goes from a value $y_-$,
just less than $y_i$, to $y_+$, just bigger than $y_i$. We define
$Z(y_+)$ from $Z(y_-)$ in each case:

\medskip${\underline{\mbox{Case (i)}}}$. $y_i$ is the
$y$-coordinate of a cusp, between the $(j-1)^{\text{th}}$ and
$j^{\text{th}}$ connected components of the intersection of
$[0,1]\x y_-$ with the black regions, as below
\[
	\begin{picture}(0,0)%
\epsfig{file=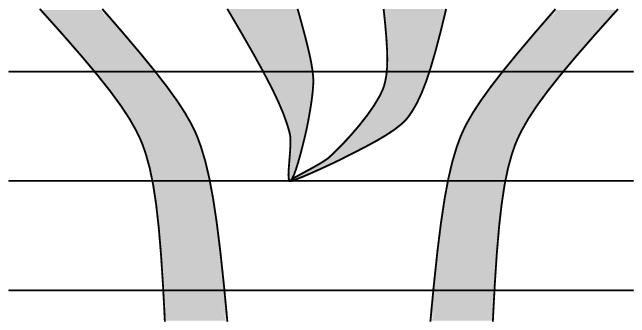}%
\end{picture}%
\setlength{\unitlength}{0.00041700in}%
\begingroup\makeatletter\ifx\SetFigFont\undefined
\def\x#1#2#3#4#5#6#7\relax{\def\x{#1#2#3#4#5#6}}%
\expandafter\x\fmtname xxxxxx\relax \def\y{splain}%
\ifx\x\y   
\gdef\SetFigFont#1#2#3{%
  \ifnum #1<17\tiny\else \ifnum #1<20\small\else
  \ifnum #1<24\normalsize\else \ifnum #1<29\large\else
  \ifnum #1<34\Large\else \ifnum #1<41\LARGE\else
     \huge\fi\fi\fi\fi\fi\fi
  \csname #3\endcsname}%
\else
\gdef\SetFigFont#1#2#3{\begingroup
  \count@#1\relax \ifnum 25<\count@\count@25\fi
  \def\x{\endgroup\@setsize\SetFigFont{#2pt}}%
  \expandafter\x
    \csname \romannumeral\the\count@ pt\expandafter\endcsname
    \csname @\romannumeral\the\count@ pt\endcsname
  \csname #3\endcsname}%
\fi
\fi\endgroup
\begin{picture}(6623,3543)(901,-3883)
\put(901,-1036){\makebox(0,0)[lb]{\smash{\SetFigFont{12}{14.4}{rm}$y_+$}}}
\put(901,-2086){\makebox(0,0)[lb]{\smash{\SetFigFont{12}{14.4}{rm}$y_i$}}}
\put(901,-3136){\makebox(0,0)[lb]{\smash{\SetFigFont{12}{14.4}{rm}$y_-$}}}
\put(3151,-3811){\makebox(0,0)[lb]{\smash{\SetFigFont{12}{14.4}{rm}$j-1$}}}
\put(5701,-3811){\makebox(0,0)[lb]{\smash{\SetFigFont{12}{14.4}{rm}$j$}}}
\end{picture}

\]

\bs

\ni Set \ $Z(y_+)=\a_{j,c}(Z(y_-))$ where $c$ is the label
(in $V_{k(y_-)}$) associated with the cusp.

\bs ${\underline{\mbox{Case (ii)}}}$. $y_i$ is the
$y$-coordinate of a minimum, with numbering and shading as below:
\[
	\begin{picture}(0,0)%
\epsfig{file=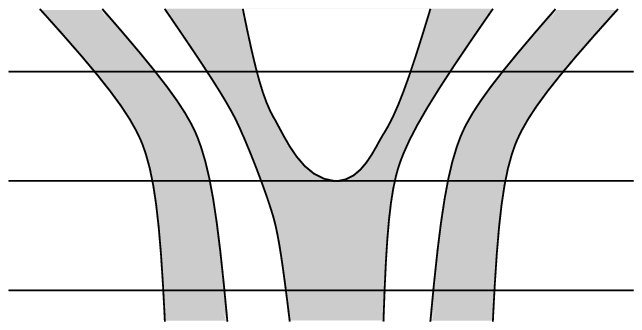}%
\end{picture}%
\setlength{\unitlength}{0.00041700in}%
\begingroup\makeatletter\ifx\SetFigFont\undefined
\def\x#1#2#3#4#5#6#7\relax{\def\x{#1#2#3#4#5#6}}%
\expandafter\x\fmtname xxxxxx\relax \def\y{splain}%
\ifx\x\y   
\gdef\SetFigFont#1#2#3{%
  \ifnum #1<17\tiny\else \ifnum #1<20\small\else
  \ifnum #1<24\normalsize\else \ifnum #1<29\large\else
  \ifnum #1<34\Large\else \ifnum #1<41\LARGE\else
     \huge\fi\fi\fi\fi\fi\fi
  \csname #3\endcsname}%
\else
\gdef\SetFigFont#1#2#3{\begingroup
  \count@#1\relax \ifnum 25<\count@\count@25\fi
  \def\x{\endgroup\@setsize\SetFigFont{#2pt}}%
  \expandafter\x
    \csname \romannumeral\the\count@ pt\expandafter\endcsname
    \csname @\romannumeral\the\count@ pt\endcsname
  \csname #3\endcsname}%
\fi
\fi\endgroup
\begin{picture}(6623,3544)(901,-3883)
\put(901,-1036){\makebox(0,0)[lb]{\smash{\SetFigFont{12}{14.4}{rm}$y_+$}}}
\put(901,-2086){\makebox(0,0)[lb]{\smash{\SetFigFont{12}{14.4}{rm}$y_i$}}}
\put(901,-3136){\makebox(0,0)[lb]{\smash{\SetFigFont{12}{14.4}{rm}$y_-$}}}
\put(3151,-3811){\makebox(0,0)[lb]{\smash{\SetFigFont{12}{14.4}{rm}$j-1$}}}
\put(4501,-3811){\makebox(0,0)[lb]{\smash{\SetFigFont{12}{14.4}{rm}$j$}}}
\put(5701,-3811){\makebox(0,0)[lb]{\smash{\SetFigFont{12}{14.4}{rm}$j+1$}}}
\end{picture}

\]
\bs

\ni Set \ $Z(y_+)=\kappa_j(Z(y_-))$.

\newpage ${\underline{\mbox{Case (iii)}}}$. $y_i$ is the
$y$-coordinate of a minimum, with numbering and shading as below:
\[
	\begin{picture}(0,0)%
\epsfig{file=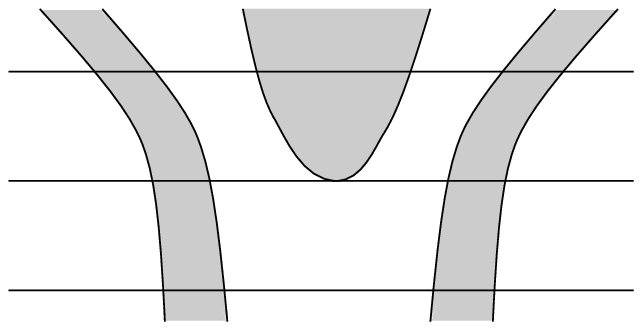}%
\end{picture}%
\setlength{\unitlength}{0.00041700in}%
\begingroup\makeatletter\ifx\SetFigFont\undefined
\def\x#1#2#3#4#5#6#7\relax{\def\x{#1#2#3#4#5#6}}%
\expandafter\x\fmtname xxxxxx\relax \def\y{splain}%
\ifx\x\y   
\gdef\SetFigFont#1#2#3{%
  \ifnum #1<17\tiny\else \ifnum #1<20\small\else
  \ifnum #1<24\normalsize\else \ifnum #1<29\large\else
  \ifnum #1<34\Large\else \ifnum #1<41\LARGE\else
     \huge\fi\fi\fi\fi\fi\fi
  \csname #3\endcsname}%
\else
\gdef\SetFigFont#1#2#3{\begingroup
  \count@#1\relax \ifnum 25<\count@\count@25\fi
  \def\x{\endgroup\@setsize\SetFigFont{#2pt}}%
  \expandafter\x
    \csname \romannumeral\the\count@ pt\expandafter\endcsname
    \csname @\romannumeral\the\count@ pt\endcsname
  \csname #3\endcsname}%
\fi
\fi\endgroup
\begin{picture}(6623,3543)(901,-3883)
\put(901,-1036){\makebox(0,0)[lb]{\smash{\SetFigFont{12}{14.4}{rm}$y_+$}}}
\put(901,-2086){\makebox(0,0)[lb]{\smash{\SetFigFont{12}{14.4}{rm}$y_i$}}}
\put(901,-3136){\makebox(0,0)[lb]{\smash{\SetFigFont{12}{14.4}{rm}$y_-$}}}
\put(3151,-3811){\makebox(0,0)[lb]{\smash{\SetFigFont{12}{14.4}{rm}$j-1$}}}
\put(5776,-3811){\makebox(0,0)[lb]{\smash{\SetFigFont{12}{14.4}{rm}$j$}}}
\end{picture}

\]

\bs

\ni Set \ $Z(y_+)=\eta_j(Z(y_-))$.

\bs ${\underline{\mbox{Case (iv)}}}$. $y_i$ is the
$y$-coordinate of a maximum, with numbering and shading as below:
\[
	\begin{picture}(0,0)%
\epsfig{file=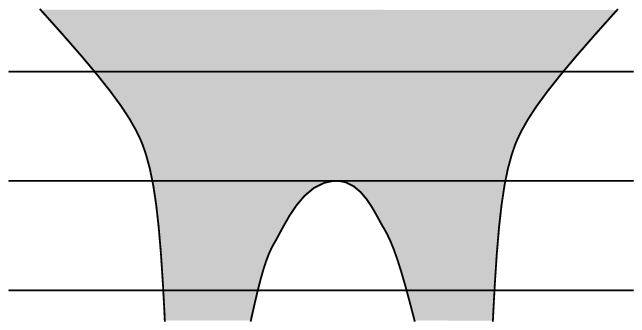}%
\end{picture}%
\setlength{\unitlength}{0.00041700in}%
\begingroup\makeatletter\ifx\SetFigFont\undefined
\def\x#1#2#3#4#5#6#7\relax{\def\x{#1#2#3#4#5#6}}%
\expandafter\x\fmtname xxxxxx\relax \def\y{splain}%
\ifx\x\y   
\gdef\SetFigFont#1#2#3{%
  \ifnum #1<17\tiny\else \ifnum #1<20\small\else
  \ifnum #1<24\normalsize\else \ifnum #1<29\large\else
  \ifnum #1<34\Large\else \ifnum #1<41\LARGE\else
     \huge\fi\fi\fi\fi\fi\fi
  \csname #3\endcsname}%
\else
\gdef\SetFigFont#1#2#3{\begingroup
  \count@#1\relax \ifnum 25<\count@\count@25\fi
  \def\x{\endgroup\@setsize\SetFigFont{#2pt}}%
  \expandafter\x
    \csname \romannumeral\the\count@ pt\expandafter\endcsname
    \csname @\romannumeral\the\count@ pt\endcsname
  \csname #3\endcsname}%
\fi
\fi\endgroup
\begin{picture}(6623,3544)(901,-3883)
\put(901,-1036){\makebox(0,0)[lb]{\smash{\SetFigFont{12}{14.4}{rm}$y_+$}}}
\put(901,-2086){\makebox(0,0)[lb]{\smash{\SetFigFont{12}{14.4}{rm}$y_i$}}}
\put(901,-3136){\makebox(0,0)[lb]{\smash{\SetFigFont{12}{14.4}{rm}$y_-$}}}
\put(5776,-3811){\makebox(0,0)[lb]{\smash{\SetFigFont{12}{14.4}{rm}$j+1$}}}
\put(3226,-3811){\makebox(0,0)[lb]{\smash{\SetFigFont{12}{14.4}{rm}$j$}}}
\end{picture}

\]

\bs

\ni Set \ $Z(y_+)=\mu_{j+1}(Z(y_-))$.

\bs ${\underline{\mbox{Case (v)}}}$. $y_i$ is the
$y$-coordinate of a maximum  with shading as below:
\[
	\begin{picture}(0,0)%
\epsfig{file=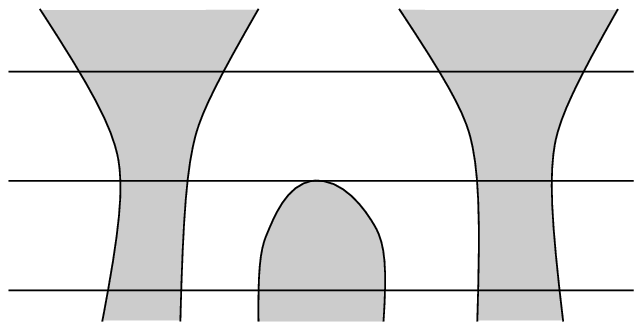}%
\end{picture}%
\setlength{\unitlength}{0.00041700in}%
\begingroup\makeatletter\ifx\SetFigFont\undefined
\def\x#1#2#3#4#5#6#7\relax{\def\x{#1#2#3#4#5#6}}%
\expandafter\x\fmtname xxxxxx\relax \def\y{splain}%
\ifx\x\y   
\gdef\SetFigFont#1#2#3{%
  \ifnum #1<17\tiny\else \ifnum #1<20\small\else
  \ifnum #1<24\normalsize\else \ifnum #1<29\large\else
  \ifnum #1<34\Large\else \ifnum #1<41\LARGE\else
     \huge\fi\fi\fi\fi\fi\fi
  \csname #3\endcsname}%
\else
\gdef\SetFigFont#1#2#3{\begingroup
  \count@#1\relax \ifnum 25<\count@\count@25\fi
  \def\x{\endgroup\@setsize\SetFigFont{#2pt}}%
  \expandafter\x
    \csname \romannumeral\the\count@ pt\expandafter\endcsname
    \csname @\romannumeral\the\count@ pt\endcsname
  \csname #3\endcsname}%
\fi
\fi\endgroup
\begin{picture}(6623,3544)(901,-3883)
\put(901,-1036){\makebox(0,0)[lb]{\smash{\SetFigFont{12}{14.4}{rm}$y_+$}}}
\put(901,-2086){\makebox(0,0)[lb]{\smash{\SetFigFont{12}{14.4}{rm}$y_i$}}}
\put(901,-3136){\makebox(0,0)[lb]{\smash{\SetFigFont{12}{14.4}{rm}$y_-$}}}
\put(2626,-3811){\makebox(0,0)[lb]{\smash{\SetFigFont{12}{14.4}{rm}$j-1$}}}
\put(4351,-3811){\makebox(0,0)[lb]{\smash{\SetFigFont{12}{14.4}{rm}$j$}}}
\put(6226,-3811){\makebox(0,0)[lb]{\smash{\SetFigFont{12}{14.4}{rm}$j+1$}}}
\end{picture}

\]

\bs

\ni Set \ $Z(y_+)=a_j(Z(y_-))$.

\medskip
Finally, we define \ $Z_{\Theta}=\th(Z(y))$ for $y > y_c$,
also written $Z_1(\Theta)$.
Our main job is now to prove that $Z_{\Theta}$ is unchanged if $\Theta$
is changed by isotopy to another standard picture $\Theta'$, with
labels transported by the isotopy. If the isotopy passes only
through standard pictures, critical points can never change order
or be annihilated or created. The pattern of connected components of
the intersection of the black regions with horizontal lines cannot
be changed either as such a change would have to involve two
maxima, minima or cusps having the same $y$ coordinate. So
isotopies through standard pictures do not change $Z_{\Theta}$.

We next argue that if $\phi_t$, $0\leq t\leq 1$ is an isotopy which
preserves the standard form of each cusp, then
$Z_{\Theta}=Z_{\phi_1(\Theta)}$. For now $Z_{\phi_t(\Theta)}$ can
only change if the singularities of the $y$-coordinate function
change. By putting the isotopy in general position we see that this
can be supposed to happen in only two ways (see e.g. [Tu], or
note that this argument can be made quite combinatorial by using
piecewise linear strings).

\begin{verse}
(1) The $y$-coordinates of two of the singularities coincide then
change order while the $x$-coordinates remain distinct.

(2) The $y$-coordinate along some string has a point of
inflection and the picture, before and after, looks locally like
one of the following
\end{verse}
\[
        \begin{picture}(0,0)%
\epsfig{file=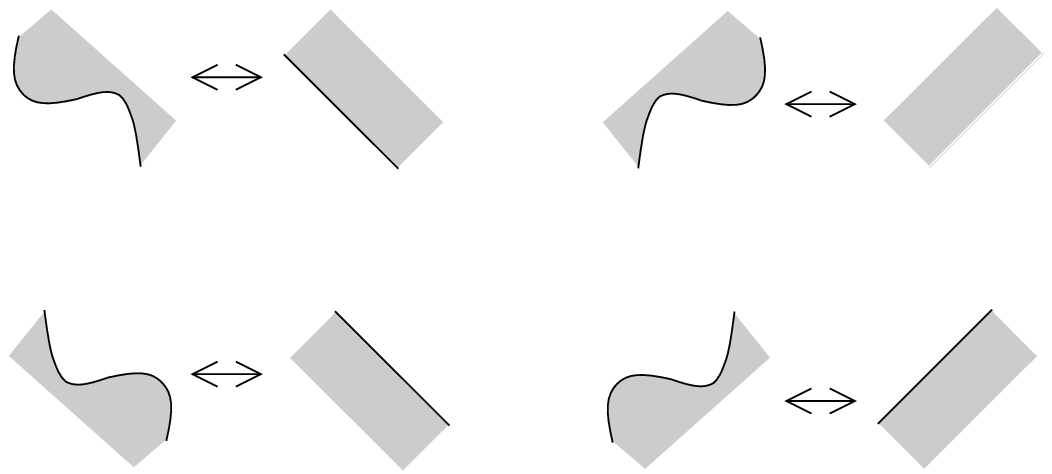}%
\end{picture}%
\setlength{\unitlength}{0.00041700in}%
\begingroup\makeatletter\ifx\SetFigFont\undefined
\def\x#1#2#3#4#5#6#7\relax{\def\x{#1#2#3#4#5#6}}%
\expandafter\x\fmtname xxxxxx\relax \def\y{splain}%
\ifx\x\y   
\gdef\SetFigFont#1#2#3{%
  \ifnum #1<17\tiny\else \ifnum #1<20\small\else
  \ifnum #1<24\normalsize\else \ifnum #1<29\large\else
  \ifnum #1<34\Large\else \ifnum #1<41\LARGE\else
     \huge\fi\fi\fi\fi\fi\fi
  \csname #3\endcsname}%
\else
\gdef\SetFigFont#1#2#3{\begingroup
  \count@#1\relax \ifnum 25<\count@\count@25\fi
  \def\x{\endgroup\@setsize\SetFigFont{#2pt}}%
  \expandafter\x
    \csname \romannumeral\the\count@ pt\expandafter\endcsname
    \csname @\romannumeral\the\count@ pt\endcsname
  \csname #3\endcsname}%
\fi
\fi\endgroup
\begin{picture}(10768,4496)(376,-4221)
\put(451,-436){\makebox(0,0)[lb]{\smash{\SetFigFont{12}{14.4}{rm}(A)}}}
\put(376,-3436){\makebox(0,0)[lb]{\smash{\SetFigFont{12}{14.4}{rm}(B)}}}
\put(6226,-511){\makebox(0,0)[lb]{\smash{\SetFigFont{12}{14.4}{rm}(C)}}}
\put(6226,-3511){\makebox(0,0)[lb]{\smash{\SetFigFont{12}{14.4}{rm}(D)}}}
\end{picture}

\]
In case 2), invariance of $Z_{\Theta}$ is guaranteed by (v) of 4.1.20
(in the order $C,D,B,A$) and in case 1), (i)$\to$(iv) of 4.1.20 is a
systematic enumeration of all 25 possibilities. For instance, the
case
\[
        \begin{picture}(0,0)%
\epsfig{file=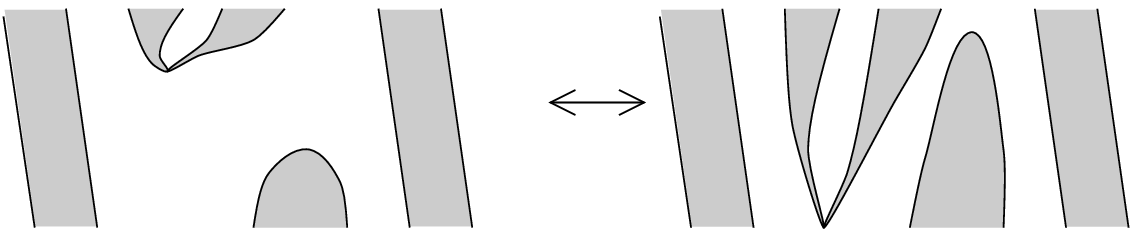}%
\end{picture}%
\setlength{\unitlength}{0.00041700in}%
\begingroup\makeatletter\ifx\SetFigFont\undefined
\def\x#1#2#3#4#5#6#7\relax{\def\x{#1#2#3#4#5#6}}%
\expandafter\x\fmtname xxxxxx\relax \def\y{splain}%
\ifx\x\y   
\gdef\SetFigFont#1#2#3{%
  \ifnum #1<17\tiny\else \ifnum #1<20\small\else
  \ifnum #1<24\normalsize\else \ifnum #1<29\large\else
  \ifnum #1<34\Large\else \ifnum #1<41\LARGE\else
     \huge\fi\fi\fi\fi\fi\fi
  \csname #3\endcsname}%
\else
\gdef\SetFigFont#1#2#3{\begingroup
  \count@#1\relax \ifnum 25<\count@\count@25\fi
  \def\x{\endgroup\@setsize\SetFigFont{#2pt}}%
  \expandafter\x
    \csname \romannumeral\the\count@ pt\expandafter\endcsname
    \csname @\romannumeral\the\count@ pt\endcsname
  \csname #3\endcsname}%
\fi
\fi\endgroup
\begin{picture}(10845,2146)(879,-3085)
\put(5401,-2011){\makebox(0,0)[lb]{\smash{\SetFigFont{12}{14.4}{rm}$\cdots$}}}
\put(11701,-2011){\makebox(0,0)[lb]{\smash{\SetFigFont{12}{14.4}{rm}$\cdots$}}}
\end{picture}

\]
\ni is covered by the first equation of (iii) of 4.1.20 with $j=i$.
Thus $Z_{\Theta}$ is invariant under isotopies preserving
standardness of the cusps.

Now we argue that a general isotopy $\phi_t$ may be replaced by a
$\tilde\phi_t$ for which $\tilde\phi_t$ preserves the standard form
of cusps for all $t$, without changing
$Z_{\Theta}$. To see this, construct a small disc around each cusp of
$\Theta$, sufficiently small so that the $y$-coordinates of all
points in a given disc are distinct from those in any other disc
and distinct from any maxima or minima of $y$ on the strings, and
such that the same is true for the images of these discs under
$\phi_1$. (Remember that $\phi_1(\Theta)$ is also a standard picture.)
Now in each disc $D$, construct a smaller disc $D_0$ inside $D$,
centered at the cusp, sufficiently small so that one can construct
a new isotopy $\tilde\phi_t$ having the properties

\begin{verse}
(i) $\tilde\phi_t$ restricted to each $D_0$ is just
translation in the plane\newline (and $\tilde\phi_t$ (a cusp)
$=\phi_t$ (that cusp)).

(ii) $\tilde\phi_t=\phi_t$ on the complement of the discs
$D$.
\end{verse}
Thus inside $\tilde\phi_t(D_0)$, the cusp remains standard and
$\tilde\phi_t$ is extended somehow to the annular region between
$\tilde\phi_t(D_0)$ and $\phi_t(D)$.  But the mapping class group
of diffeomorphisms of the annulus that are the identity on the
boundary is generated by a Dehn twist of 360${}^\circ$. So in a
neighborhood of each cusp point, $\phi_1$ and $\tilde\phi_1$ differ
only by some integer power of a single full twist. Figure~\ref{pic94}
illustrates  how $\phi_1(\Theta)$ and $\tilde\phi_1(\Theta)$ would
differ if the twist incurred were a single clockwise twist.
\[
        \begin{picture}(0,0)%
\epsfig{file=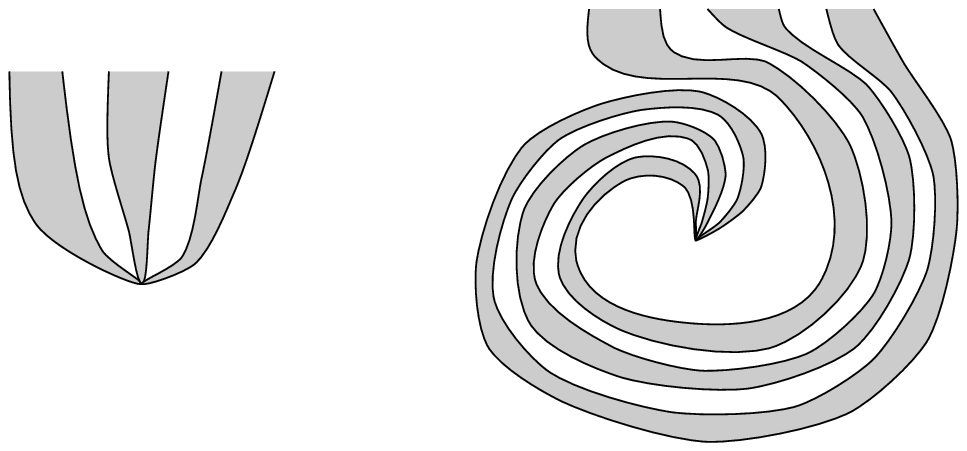}%
\end{picture}%
\setlength{\unitlength}{0.00041700in}%
\begingroup\makeatletter\ifx\SetFigFont\undefined
\def\x#1#2#3#4#5#6#7\relax{\def\x{#1#2#3#4#5#6}}%
\expandafter\x\fmtname xxxxxx\relax \def\y{splain}%
\ifx\x\y   
\gdef\SetFigFont#1#2#3{%
  \ifnum #1<17\tiny\else \ifnum #1<20\small\else
  \ifnum #1<24\normalsize\else \ifnum #1<29\large\else
  \ifnum #1<34\Large\else \ifnum #1<41\LARGE\else
     \huge\fi\fi\fi\fi\fi\fi
  \csname #3\endcsname}%
\else
\gdef\SetFigFont#1#2#3{\begingroup
  \count@#1\relax \ifnum 25<\count@\count@25\fi
  \def\x{\endgroup\@setsize\SetFigFont{#2pt}}%
  \expandafter\x
    \csname \romannumeral\the\count@ pt\expandafter\endcsname
    \csname @\romannumeral\the\count@ pt\endcsname
  \csname #3\endcsname}%
\fi
\fi\endgroup
\begin{picture}(9196,4811)(826,-4558)
\put(826,-4486){\makebox(0,0)[lb]{\smash{\SetFigFont{12}{14.4}{rm}$\phi_1(\Theta )$, near cusp}}}
\put(5776,-4486){\makebox(0,0)[lb]{\smash{\SetFigFont{12}{14.4}{rm}${\tilde\phi_1}(\Theta )$, near cusp.}}}
\end{picture}

\]
\begin{center}
	Figure 4.2.5
\end{center}
\bs We want to show, first in this case and then in the case of an
arbitrary integral power of a full twist, that $Z_{\Theta}$ is
unchanged. For this, consider Figure 4.2.6 which is supposed to be
part of a standard labelled picture in which the maxima and minima
of $y$, in the figure, and its cusps occur as an uninterrupted
sequence in
$\frak S(\Theta)$:
\[
        \begin{picture}(0,0)%
\epsfig{file=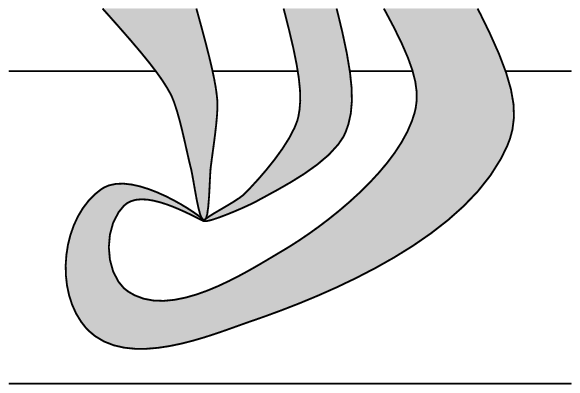}%
\end{picture}%
\setlength{\unitlength}{0.00041700in}%
\begingroup\makeatletter\ifx\SetFigFont\undefined
\def\x#1#2#3#4#5#6#7\relax{\def\x{#1#2#3#4#5#6}}%
\expandafter\x\fmtname xxxxxx\relax \def\y{splain}%
\ifx\x\y   
\gdef\SetFigFont#1#2#3{%
  \ifnum #1<17\tiny\else \ifnum #1<20\small\else
  \ifnum #1<24\normalsize\else \ifnum #1<29\large\else
  \ifnum #1<34\Large\else \ifnum #1<41\LARGE\else
     \huge\fi\fi\fi\fi\fi\fi
  \csname #3\endcsname}%
\else
\gdef\SetFigFont#1#2#3{\begingroup
  \count@#1\relax \ifnum 25<\count@\count@25\fi
  \def\x{\endgroup\@setsize\SetFigFont{#2pt}}%
  \expandafter\x
    \csname \romannumeral\the\count@ pt\expandafter\endcsname
    \csname @\romannumeral\the\count@ pt\endcsname
  \csname #3\endcsname}%
\fi
\fi\endgroup
\begin{picture}(6097,3837)(2026,-4183)
\put(2026,-4111){\makebox(0,0)[lb]{\smash{\SetFigFont{12}{14.4}{rm}$y_-$}}}
\put(2026,-1036){\makebox(0,0)[lb]{\smash{\SetFigFont{12}{14.4}{rm}$y_+$}}}
\end{picture}

\]
\begin{center}
	Figure 4.2.6
\end{center}

\bs\ni From the definition of $Z(y)$ we see that
\begin{eqnarray*}
Z(y_+) & = & a_{{}_j}\mu_{{}_{j+1}} \a_{{}_{j+1,c}} \kappa_{{}_j}\eta_{{}_j} (Z(y_{{}_-})) \\ 
& = &\a_{{}_j,\rho(c)} (Z(y_{{}_-})) \quad \mbox{(by 4.1.13 and 4.1.14)}
\end{eqnarray*}
(Here the cusp is labelled by $c\in N'\cap M_{n-1}$ where
$n=3$ in Figure 4.2.6.) Similarly we see that if the cusp is
surrounded by a full 360${}^\circ$ twist,
\begin{eqnarray*}
Z(y_+) & = & \a_{j,\rho^n(b)} (Z(y_-)) \\
& = & \a_j, b (Z(y_-)) \qquad \mbox{(by Theorem  4.1.18)}
\end{eqnarray*}
Thus if a cusp  (or any part of a picture that is just a scaled
down standard labelled picture) is surrounded by a single clockwise
full twist, the effect on $Z_{\Theta}$ is as if the twist were not
there. If there were anticlockwise full twists around a cusp,
surround it further by the same number of clockwise twists. This
does not change
$Z_{\Theta}$, but the cancelling of the positive and negative
twists involves only isotopies that are the identity near the
cusps. Thus by our previous argument clockwise full twists around
cusps do not change
$Z_{\Theta}$ either. We conclude
\[
Z_{\phi_1(\Theta)} = Z_{\tilde\phi_1(\Theta)} =Z_{\Theta} \ .
\]

We have established that $Z_{\Theta}$ may be assigned to a standard
labelled picture by a product of elementary maps $\a, \kappa,
\eta,a,\mu$ in such a way that $Z_{\Theta}$ is unchanged by isotopies
of $\Theta$. We can now formally see how this makes $\coprod_n V_n$
a planar algebra according to $\S$1. The labelling set will be
$\coprod_n V_n$ itself.  The first step will be to associate a
labelled picture $\beta(T)$ with a labelled tangle $T$. To do this,
shrink all the internal boxes of $T$ to points and  isotope the
standard $k$-box to $[0,1]\x [0,1]$ with all the marked boundary
points going to points in $[0,1]\x 1$. Then distort the shrunk
boxes of $T$ to standard cusps, by isotopy, so that the string
attached to the first boundary point of the box becomes the first
string (from the left) attached to the cusp. The procedure near a
4-box of $T$ is illustrated in Figure~\ref{pic97} :
\[
        \input{xfig/pic96}
\]
\begin{center}
	Figure 4.2.7
\end{center}
\bs\ni The label associated to the cusp is just $\th^{-1}$ of the
label associated to the box but we may reasonably suppress
$\th^{-1}$.

Figure~4.2.8 illustrates a labelled tangle $T$ and a labelled
standard picture $\beta(T)$:
\[
        \begin{picture}(0,0)%
\epsfig{file=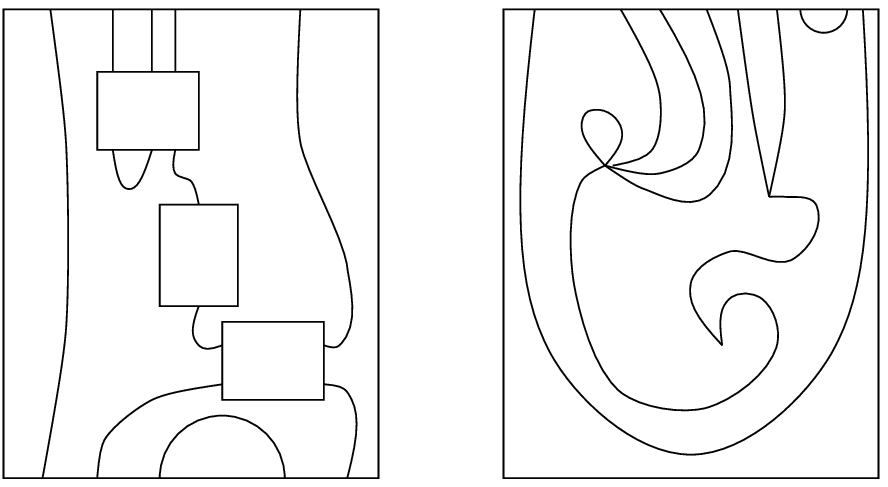}%
\end{picture}%
\setlength{\unitlength}{0.00041700in}%
\begingroup\makeatletter\ifx\SetFigFont\undefined
\def\x#1#2#3#4#5#6#7\relax{\def\x{#1#2#3#4#5#6}}%
\expandafter\x\fmtname xxxxxx\relax \def\y{splain}%
\ifx\x\y   
\gdef\SetFigFont#1#2#3{%
  \ifnum #1<17\tiny\else \ifnum #1<20\small\else
  \ifnum #1<24\normalsize\else \ifnum #1<29\large\else
  \ifnum #1<34\Large\else \ifnum #1<41\LARGE\else
     \huge\fi\fi\fi\fi\fi\fi
  \csname #3\endcsname}%
\else
\gdef\SetFigFont#1#2#3{\begingroup
  \count@#1\relax \ifnum 25<\count@\count@25\fi
  \def\x{\endgroup\@setsize\SetFigFont{#2pt}}%
  \expandafter\x
    \csname \romannumeral\the\count@ pt\expandafter\endcsname
    \csname @\romannumeral\the\count@ pt\endcsname
  \csname #3\endcsname}%
\fi
\fi\endgroup
\begin{picture}(8444,5194)(1779,-5833)
\put(3526,-1486){\makebox(0,0)[lb]{\smash{
\SetFigFont{12}{14.4}{rm}$R$
}}}
\put(4201,-3811){\makebox(0,0)[lb]{\smash{
\SetFigFont{12}{14.4}{rm}$P$
}}}
\put(3901,-2911){\makebox(0,0)[lb]{\smash{
\SetFigFont{12}{14.4}{rm}$Q$
}}}
\put(7501,-2611){\makebox(0,0)[lb]{\smash{\SetFigFont{12}{14.4}{rm}$R$}}}
\put(9076,-2836){\makebox(0,0)[lb]{\smash{\SetFigFont{12}{14.4}{rm}$P$}}}
\put(8551,-4261){\makebox(0,0)[lb]{\smash{\SetFigFont{12}{14.4}{rm}$Q$}}}
\put(2626,-5761){\makebox(0,0)[lb]{\smash{\SetFigFont{12}{14.4}{rm}The tangle $T$}}}
\put(6451,-5761){\makebox(0,0)[lb]{\smash{\SetFigFont{12}{14.4}{rm}A standard picture $\beta (T)$}}}
\end{picture}

\]
\begin{center}
	Figure 4.2.8
\end{center}

\bs\ni Note that $\beta(T)$ is not well defined, but two different
choices of $\beta(T)$ for a given $T$ will differ by an isotopy so
the map $\Phi$, \ $\Phi(T)=Z_{\beta(T)}$, gives a well defined
linear map from the universal planar algebra on $\coprod_k N'\cap
M_{k-1}$.

We now check that $\Phi$ makes $(P^{N\subset M})$ into a
connected spherical $C^*$-planar algebra. The first thing to check
is that
$\Phi$ is a homomorphism of filtered algebras. But if $T_1$ and
$T_2$ are labelled $k$-tangles, a choice of $\beta(T_1T_2)$ is
shown below $(k=3)$.
\[
        \begin{picture}(0,0)%
\epsfig{file=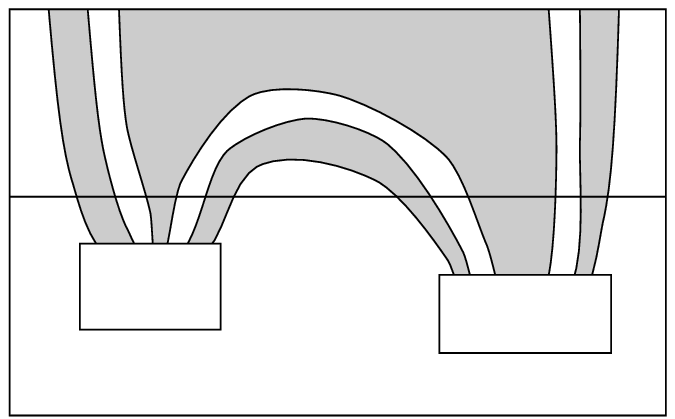}%
\end{picture}%
\setlength{\unitlength}{0.00041700in}%
\begingroup\makeatletter\ifx\SetFigFont\undefined
\def\x#1#2#3#4#5#6#7\relax{\def\x{#1#2#3#4#5#6}}%
\expandafter\x\fmtname xxxxxx\relax \def\y{splain}%
\ifx\x\y   
\gdef\SetFigFont#1#2#3{%
  \ifnum #1<17\tiny\else \ifnum #1<20\small\else
  \ifnum #1<24\normalsize\else \ifnum #1<29\large\else
  \ifnum #1<34\Large\else \ifnum #1<41\LARGE\else
     \huge\fi\fi\fi\fi\fi\fi
  \csname #3\endcsname}%
\else
\gdef\SetFigFont#1#2#3{\begingroup
  \count@#1\relax \ifnum 25<\count@\count@25\fi
  \def\x{\endgroup\@setsize\SetFigFont{#2pt}}%
  \expandafter\x
    \csname \romannumeral\the\count@ pt\expandafter\endcsname
    \csname @\romannumeral\the\count@ pt\endcsname
  \csname #3\endcsname}%
\fi
\fi\endgroup
\begin{picture}(6772,3944)(451,-4583)
\put(1951,-3436){\makebox(0,0)[lb]{\smash{\SetFigFont{12}{14.4}{rm}$\beta (T_1)$}}}
\put(5401,-3661){\makebox(0,0)[lb]{\smash{\SetFigFont{12}{14.4}{rm}$\beta (T_2)$}}}
\put(451,-2611){\makebox(0,0)[lb]{\smash{\SetFigFont{12}{14.4}{rm}$y$}}}
\end{picture}

\]
\ni If $y$ is as marked, from the definition of $Z$,
$Z(y)=Z(\beta(T_1))\otimes Z(\beta(T_2))$. Moreover each pair of
maximal $\bnd$ \hskip8pt  contributes a factor $\chi_i$
$(=a_i\mu_{i+1})$ to $Z(y)$ as $y$ increases, so if $k$ is even,
part (1) of 4.1.24 gives $\Phi(T_1T_2)=\Phi(T_1)\Phi(T_2)$ and if
$k$ is odd (as in the figure) the last maximum has the black
region above so contribute a factor $\mu$ and part (ii) of 4.1.24
applies.

That $\Phi$ is compatible with the filtrations amounts to showing
that
\[
\begin{picture}(0,0)%
\epsfig{file=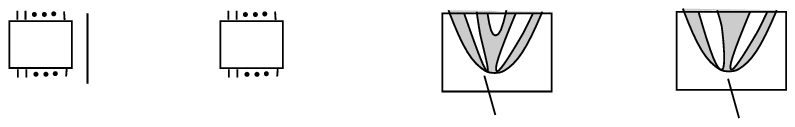}%
\end{picture}%
\setlength{\unitlength}{0.00041700in}%
\begingroup\makeatletter\ifx\SetFigFont\undefined
\def\x#1#2#3#4#5#6#7\relax{\def\x{#1#2#3#4#5#6}}%
\expandafter\x\fmtname xxxxxx\relax \def\y{splain}%
\ifx\x\y   
\gdef\SetFigFont#1#2#3{%
  \ifnum #1<17\tiny\else \ifnum #1<20\small\else
  \ifnum #1<24\normalsize\else \ifnum #1<29\large\else
  \ifnum #1<34\Large\else \ifnum #1<41\LARGE\else
     \huge\fi\fi\fi\fi\fi\fi
  \csname #3\endcsname}%
\else
\gdef\SetFigFont#1#2#3{\begingroup
  \count@#1\relax \ifnum 25<\count@\count@25\fi
  \def\x{\endgroup\@setsize\SetFigFont{#2pt}}%
  \expandafter\x
    \csname \romannumeral\the\count@ pt\expandafter\endcsname
    \csname @\romannumeral\the\count@ pt\endcsname
  \csname #3\endcsname}%
\fi
\fi\endgroup
\begin{picture}(8217,1457)(1276,-2284)
\put(2026,-1241){\makebox(0,0)[lb]{\smash{\SetFigFont{12}{14.4}{rm}$R$}}}
\put(4052,-1242){\makebox(0,0)[lb]{\smash{\SetFigFont{12}{14.4}{rm}$R$}}}
\put(6361,-2191){\makebox(0,0)[lb]{\smash{\SetFigFont{12}{14.4}{rm}$R$}}}
\put(8701,-2221){\makebox(0,0)[lb]{\smash{\SetFigFont{12}{14.4}{rm}$R$}}}
\put(2701,-1246){\makebox(0,0)[lb]{\smash{\SetFigFont{12}{14.4}{rm}$) = \Phi ($}}}
\put(1276,-1231){\makebox(0,0)[lb]{\smash{\SetFigFont{12}{14.4}{rm}$\Phi ($}}}
\put(5053,-1261){\makebox(0,0)[lb]{\smash{\SetFigFont{12}{14.4}{rm}or $Z($}}}
\put(4603,-1261){\makebox(0,0)[lb]{\smash{\SetFigFont{12}{14.4}{rm}$)$}}}
\put(9493,-1336){\makebox(0,0)[lb]{\smash{\SetFigFont{12}{14.4}{rm}$)$}}}
\put(7228,-1291){\makebox(0,0)[lb]{\smash{\SetFigFont{12}{14.4}{rm}$)=Z($}}}
\end{picture}

\]
If $k$ is even this follows from 4.1.25 and if $k$ is odd it
follows from 4.1.26 (together with $\sum_{b\in B} bv_kb^*=\d
E_kE_{k-1}\dots E_2$  to take care of the factor $\kappa$
introduced by the minimum in the picture). So $\Phi$ is a
homomorphism of filtered algebras.

Annular invariance (and indeed the whole operadic picture) is easy.
If $T$ is an element of ${\cal P}(\coprod_k N'\cap M_{k-1})$ (linear
combination of tangles) with $\Phi(T)=0$, then if $T$ is surrounded
by an annular labelled tangle $A$, then we may choose
$\beta(\pi_A(T))$ to look like
\[
        \begin{picture}(0,0)%
\epsfig{file=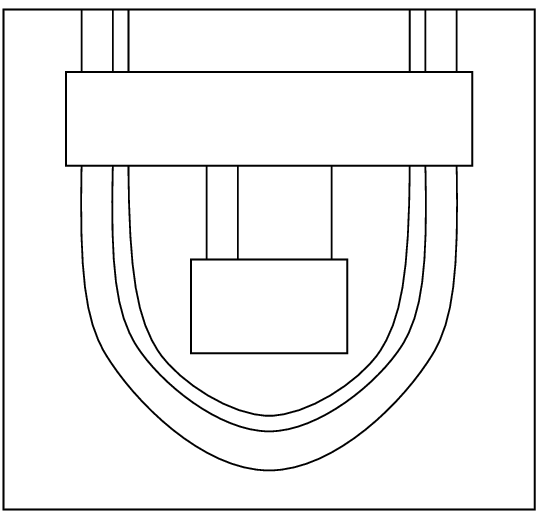}%
\end{picture}%
\setlength{\unitlength}{0.00041700in}%
\begingroup\makeatletter\ifx\SetFigFont\undefined
\def\x#1#2#3#4#5#6#7\relax{\def\x{#1#2#3#4#5#6}}%
\expandafter\x\fmtname xxxxxx\relax \def\y{splain}%
\ifx\x\y   
\gdef\SetFigFont#1#2#3{%
  \ifnum #1<17\tiny\else \ifnum #1<20\small\else
  \ifnum #1<24\normalsize\else \ifnum #1<29\large\else
  \ifnum #1<34\Large\else \ifnum #1<41\LARGE\else
     \huge\fi\fi\fi\fi\fi\fi
  \csname #3\endcsname}%
\else
\gdef\SetFigFont#1#2#3{\begingroup
  \count@#1\relax \ifnum 25<\count@\count@25\fi
  \def\x{\endgroup\@setsize\SetFigFont{#2pt}}%
  \expandafter\x
    \csname \romannumeral\the\count@ pt\expandafter\endcsname
    \csname @\romannumeral\the\count@ pt\endcsname
  \csname #3\endcsname}%
\fi
\fi\endgroup
\begin{picture}(5144,4844)(1179,-4883)
\put(3601,-2161){\makebox(0,0)[lb]{\smash{\SetFigFont{12}{14.4}{rm}$\cdots$}}}
\put(3301,-1186){\makebox(0,0)[lb]{\smash{\SetFigFont{12}{14.4}{rm}$\beta (A)$}}}
\put(3301,-2986){\makebox(0,0)[lb]{\smash{\SetFigFont{12}{14.4}{rm}$\beta (T)$}}}
\end{picture}

\]
\ni (Strictly speaking, one needs to consider such a picture for
each tangle in the linear combination forming $T$, and add.)
Clearly if $Z_{\beta(T)}=0$, so is $Z_{\beta(\pi_A(T))}$, since the
map $\a_{j,\th^{-1}(Z_{\beta(T)})}$ is applied in forming
$Z_{\beta(\pi_A(T))}$.

We now turn to planarity. By definition $V_0=\Bbb C$ so we only
need to show \ dim $V_{1,1}=1$.  A basis element of
${\cal P}_{1,1}(\coprod^\i_{k=1} N'\cap M_{k-1})$ is a 1-tangle $T$
with a vertical straight line and planar networks to the left and
right. We may choose $\beta(T)$ to be as depicted below
\[
        \begin{picture}(0,0)%
\epsfig{file=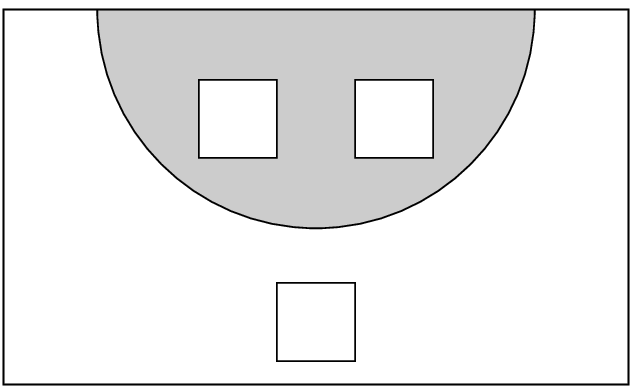}%
\end{picture}%
\setlength{\unitlength}{0.00041700in}%
\begingroup\makeatletter\ifx\SetFigFont\undefined
\def\x#1#2#3#4#5#6#7\relax{\def\x{#1#2#3#4#5#6}}%
\expandafter\x\fmtname xxxxxx\relax \def\y{splain}%
\ifx\x\y   
\gdef\SetFigFont#1#2#3{%
  \ifnum #1<17\tiny\else \ifnum #1<20\small\else
  \ifnum #1<24\normalsize\else \ifnum #1<29\large\else
  \ifnum #1<34\Large\else \ifnum #1<41\LARGE\else
     \huge\fi\fi\fi\fi\fi\fi
  \csname #3\endcsname}%
\else
\gdef\SetFigFont#1#2#3{\begingroup
  \count@#1\relax \ifnum 25<\count@\count@25\fi
  \def\x{\endgroup\@setsize\SetFigFont{#2pt}}%
  \expandafter\x
    \csname \romannumeral\the\count@ pt\expandafter\endcsname
    \csname @\romannumeral\the\count@ pt\endcsname
  \csname #3\endcsname}%
\fi
\fi\endgroup
\begin{picture}(6044,3644)(1179,-3983)
\put(3226,-1561){\makebox(0,0)[lb]{\smash{\SetFigFont{12}{14.4}{rm}$T_1$}}}
\put(4801,-1561){\makebox(0,0)[lb]{\smash{\SetFigFont{12}{14.4}{rm}$T_2$}}}
\put(4051,-3511){\makebox(0,0)[lb]{\smash{\SetFigFont{12}{14.4}{rm}$T_3$}}}
\end{picture}

\]
\ni where there are 0-pictures inside
 the regions $T_1,T_2,\dots$. It is a simple consequence of our
formalism that a closed picture surrounded by a {\it white} region
simply contributes a scalar in a multiplicative way. This is
because one may first isotope the big picture so that all the
maxima and cusps in the 0-picture have $y$ coordinates in an
uninterrupted sequence in $\frak S$, and the last singularity must
be a maximum, shaded below, the first being a minimum, shaded
above. The final map will be an $a_j$ and will send the
contribution of the 0-picture to an element of $N'\cap N=\Bbb
C$. Thus we only need to see that 0-pictures inside a black
region contribute a scalar in a multiplicative way. But we may
isotope the big picture so that the singular $y$-values of the
0-picture occur in uninterrupted succession, and near the
0-picture the situation is as below:
\[
        \begin{picture}(0,0)%
\epsfig{file=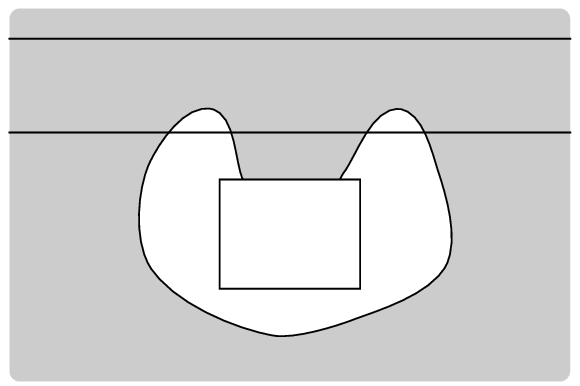}%
\end{picture}%
\setlength{\unitlength}{0.00041700in}%
\begingroup\makeatletter\ifx\SetFigFont\undefined
\def\x#1#2#3#4#5#6#7\relax{\def\x{#1#2#3#4#5#6}}%
\expandafter\x\fmtname xxxxxx\relax \def\y{splain}%
\ifx\x\y   
\gdef\SetFigFont#1#2#3{%
  \ifnum #1<17\tiny\else \ifnum #1<20\small\else
  \ifnum #1<24\normalsize\else \ifnum #1<29\large\else
  \ifnum #1<34\Large\else \ifnum #1<41\LARGE\else
     \huge\fi\fi\fi\fi\fi\fi
  \csname #3\endcsname}%
\else
\gdef\SetFigFont#1#2#3{\begingroup
  \count@#1\relax \ifnum 25<\count@\count@25\fi
  \def\x{\endgroup\@setsize\SetFigFont{#2pt}}%
  \expandafter\x
    \csname \romannumeral\the\count@ pt\expandafter\endcsname
    \csname @\romannumeral\the\count@ pt\endcsname
  \csname #3\endcsname}%
\fi
\fi\endgroup
\begin{picture}(5947,3644)(676,-3683)
\put(676,-1411){\makebox(0,0)[lb]{\smash{\SetFigFont{12}{14.4}{rm}$y_1$}}}
\put(676,-511){\makebox(0,0)[lb]{\smash{\SetFigFont{12}{14.4}{rm}$y_2$}}}
\put(3377,-2591){\makebox(0,0)[lb]{\smash{\SetFigFont{12}{14.4}{rm}picture}}}
\put(3377,-2087){\makebox(0,0)[lb]{\smash{\SetFigFont{12}{14.4}{rm}some}}}
\end{picture}

\]
\ni With $y_1$ as marked, $Z(y_1)$ will be $\sum_{b\in B}
b\otimes x\otimes b^*$ for some element $x$ in $N'\cap M$.
But then $Z(y_2)$ will be $\sum_{b\in B} bxb^*\in M'\cap M=\Bbb C$,
by 4.1.5.  So $P^{N\subset M}$ is a planar algebra.

The spherical property is easy: comparing
\bs
\[
        \begin{picture}(0,0)%
\epsfig{file=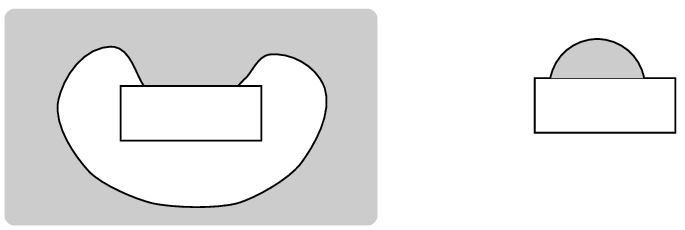}%
\end{picture}%
\setlength{\unitlength}{0.00041700in}%
\begingroup\makeatletter\ifx\SetFigFont\undefined
\def\x#1#2#3#4#5#6#7\relax{\def\x{#1#2#3#4#5#6}}%
\expandafter\x\fmtname xxxxxx\relax \def\y{splain}%
\ifx\x\y   
\gdef\SetFigFont#1#2#3{%
  \ifnum #1<17\tiny\else \ifnum #1<20\small\else
  \ifnum #1<24\normalsize\else \ifnum #1<29\large\else
  \ifnum #1<34\Large\else \ifnum #1<41\LARGE\else
     \huge\fi\fi\fi\fi\fi\fi
  \csname #3\endcsname}%
\else
\gdef\SetFigFont#1#2#3{\begingroup
  \count@#1\relax \ifnum 25<\count@\count@25\fi
  \def\x{\endgroup\@setsize\SetFigFont{#2pt}}%
  \expandafter\x
    \csname \romannumeral\the\count@ pt\expandafter\endcsname
    \csname @\romannumeral\the\count@ pt\endcsname
  \csname #3\endcsname}%
\fi
\fi\endgroup
\begin{picture}(6547,2144)(2079,-3083)
\put(3376,-2086){\makebox(0,0)[lb]{\smash{\SetFigFont{12}{14.4}{rm}picture}}}
\put(6151,-2011){\makebox(0,0)[lb]{\smash{\SetFigFont{12}{14.4}{rm}and}}}
\put(7351,-2011){\makebox(0,0)[lb]{\smash{\SetFigFont{12}{14.4}{rm}picture}}}
\put(8626,-2086){\makebox(0,0)[lb]{\smash{\SetFigFont{12}{14.4}{rm},}}}
\end{picture}

\]
\bs\ni we see that the partition
functions are the same since one gives $E_{M'}$ applied to an
element of $N'\cap M$, and the other gives $E_N$ applied to the
same element with the correct powers of $\d$ contributed from the
minimum in the first picture and the maximum in the second. Either
way we get the trace by extremality.

For the $*$-structure, observe first that
$$
\th(x_1\otimes x_2\dots\otimes x_k)^* =
\th(x_k^*\otimes x_{k-1}^*\otimes\dots\otimes x_1^*)
$$
so
$\th(\a_{j,c}(x))^*=\th(\a_{k-j+2,c^*}(x^*))$, \
$\th(a_j(x))^*=\th(a_{k-j+1}(x^*))$, \
$\th(\mu_j(x))^*=\th(\mu_{n-j+2}(x^*))$, \
$\th(\eta_j(x))^*=\th(\eta_{k-j+2}(x^*))$ and
$\th(\kappa_j(x))^*=\th(\kappa_{k-j+1}(x^*))$.
Moreover if $T$ is a labelled tangle, $\beta(T^*)$ is $\beta(T)$
reflected in the line $x=\frac 12$ and with labels replaced by
their adjoints (via $\th$). We conclude that
$(Z_{\beta(T)})^*=Z_{\beta(T^*)}$ by applying the relations above
at each of the $y$-values in $\frak S(\beta(T))=\frak
S(\beta(T^*))$.

The $C^*$-property is just the positive definiteness of the
partition function. But if $T$ is a labelled $k$-tangle,
\[
\beta(\tangler{T}\ ) = \dtangle{\beta (T)} \ \ \ \ \mbox{so}
Z_{\Phi}(\tangler{T}\ )  =\d E_{M_{k-2}}(\Phi(\tangle{T}\ ))
\]
by 4.1.22.  Applying this $k$ times we get
\[
Z(\boxrloop{T} ) =
\d^k \ {\op{tr}}(\Phi(\tangle{T}\ ))
\]
so the positive definiteness of $Z$ follows from that of tr.

We now verify (i)$\to$(iv) in the statement of the theorem.

\smallskip
(i) Since $\Phi$ is a homomorphism of filtered algebras, it
suffices to prove the formula when $k=i$.
\begin{eqnarray*}
{\mbox{If}} \ k=2p+1, & \ Z_{\beta( | \ |\dots \quad )}= \th(Z_1( \butau )) \\
& = \d^{-P}\th\big( \sum_{b_1,b_2,\dots ,b_p} b_1\otimes b_2\otimes\dots \otimes b_p\otimes 1\otimes 1\otimes b^*_p \otimes\dots \otimes b^*_2\otimes b^*_1\big) \\
& = E_{2p+1} \qquad {\mbox{by Lemma 4.1.27.}}\vspace{2\jot}\\
\end{eqnarray*}
\begin{eqnarray*}
{\mbox {If  $k=2p$}},& \ Z_{\beta( | \ |\dots \quad )} =  \th(Z_1(\butad )) \\
&= \d^{-p-1}\th\big(\sum_{b_1\dots b_{p+1}} b_1\otimes b_2\otimes\dots \otimes b_p\otimes b^*_p b_{p+1}\otimes b^*_{p+1}\otimes b^*_{p-1} \otimes\dots \otimes b^*_1\big) \\
&= E_{2p} \qquad {\mbox{by  Corollary 4.1.28}}.
\end{eqnarray*}

(ii) The first formula follows from 4.1.5 and we showed the second
when we proved the positive defIniteness of the partition function.

(iii) This is just the filtered algebra property.

(iv) We also showed this in the  positive definiteness  proof.

All that remains is to prove the uniqueness of the planar algebra
structure. First observe that, as in Proposition 1.14, a labelled
tangle may be arranged by isotopy so that all of its boxes occur in
a vertical stack. After further isotopy and the introduction of
kinks or redundant loops, one may obtain the picture below for the
tangle (in ${\cal P}_k$)
\[
        \begin{picture}(0,0)%
\epsfig{file=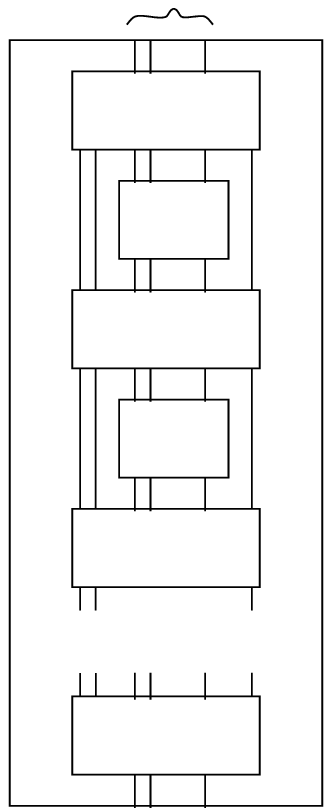}%
\end{picture}%
\setlength{\unitlength}{0.00041700in}%
\begingroup\makeatletter\ifx\SetFigFont\undefined
\def\x#1#2#3#4#5#6#7\relax{\def\x{#1#2#3#4#5#6}}%
\expandafter\x\fmtname xxxxxx\relax \def\y{splain}%
\ifx\x\y   
\gdef\SetFigFont#1#2#3{%
  \ifnum #1<17\tiny\else \ifnum #1<20\small\else
  \ifnum #1<24\normalsize\else \ifnum #1<29\large\else
  \ifnum #1<34\Large\else \ifnum #1<41\LARGE\else
     \huge\fi\fi\fi\fi\fi\fi
  \csname #3\endcsname}%
\else
\gdef\SetFigFont#1#2#3{\begingroup
  \count@#1\relax \ifnum 25<\count@\count@25\fi
  \def\x{\endgroup\@setsize\SetFigFont{#2pt}}%
  \expandafter\x
    \csname \romannumeral\the\count@ pt\expandafter\endcsname
    \csname @\romannumeral\the\count@ pt\endcsname
  \csname #3\endcsname}%
\fi
\fi\endgroup
\begin{picture}(3922,8131)(2101,-8054)
\put(4426,-6286){\makebox(0,0)[lb]{\smash{\SetFigFont{12}{14.4}{rm}$\vdots$}}}
\put(3751,-211){\makebox(0,0)[lb]{\smash{\SetFigFont{12}{14.4}{rm}$k$ boundary point}}}
\put(2101,-4036){\makebox(0,0)[lb]{\smash{\SetFigFont{12}{14.4}{rm}$T =$}}}
\put(4426,-2536){\makebox(0,0)[lb]{\smash{\SetFigFont{12}{14.4}{rm}$R_1$}}}
\put(4426,-4636){\makebox(0,0)[lb]{\smash{\SetFigFont{12}{14.4}{rm}$R_2$}}}
\put(4351,-1411){\makebox(0,0)[lb]{\smash{\SetFigFont{12}{14.4}{rm}$\sigma_1$}}}
\put(4351,-3511){\makebox(0,0)[lb]{\smash{\SetFigFont{12}{14.4}{rm}$\sigma_2$}}}
\put(4351,-7411){\makebox(0,0)[lb]{\smash{\SetFigFont{12}{14.4}{rm}$\sigma_n$}}}
\end{picture}

\]

\ni where the regions marked $\s_1,\dots ,\s_n$ contain only
strings and $\s_2,\dots ,\s_{n-1}$ have a fixed number $p$ of
boundary strings top and bottom, $\geq k$. Clearly $p-k$ is
even  so we conclude that, if $\Phi_1$ is some other planar
algebra structure satisfying (i) and (ii) then
$$
\Phi_1(T) =\d^{-\frac{p-k}{2}} \ \Phi_1
\big(\bigtan{T}\big) \  ,
$$
where we have introduced $\frac{p-k}{2}$ maxima and minima.
(To see this just apply the second formula of (ii) $p-k$ times.)
Thus we find that it suffices to prove that $\Phi'=\Phi$ on a
product of Temperley-Lieb tangles and tangles of the form
\[
        \begin{picture}(0,0)%
\epsfig{file=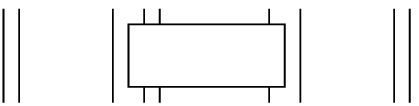}%
\end{picture}%
\setlength{\unitlength}{0.00041700in}%
\begingroup\makeatletter\ifx\SetFigFont\undefined
\def\x#1#2#3#4#5#6#7\relax{\def\x{#1#2#3#4#5#6}}%
\expandafter\x\fmtname xxxxxx\relax \def\y{splain}%
\ifx\x\y   
\gdef\SetFigFont#1#2#3{%
  \ifnum #1<17\tiny\else \ifnum #1<20\small\else
  \ifnum #1<24\normalsize\else \ifnum #1<29\large\else
  \ifnum #1<34\Large\else \ifnum #1<41\LARGE\else
     \huge\fi\fi\fi\fi\fi\fi
  \csname #3\endcsname}%
\else
\gdef\SetFigFont#1#2#3{\begingroup
  \count@#1\relax \ifnum 25<\count@\count@25\fi
  \def\x{\endgroup\@setsize\SetFigFont{#2pt}}%
  \expandafter\x
    \csname \romannumeral\the\count@ pt\expandafter\endcsname
    \csname @\romannumeral\the\count@ pt\endcsname
  \csname #3\endcsname}%
\fi
\fi\endgroup
\begin{picture}(3944,944)(2379,-1283)
\put(4201,-961){\makebox(0,0)[lb]{\smash{\SetFigFont{12}{14.4}{rm}$x$}}}
\put(2701,-961){\makebox(0,0)[lb]{\smash{\SetFigFont{12}{14.4}{rm}$\cdots$}}}
\put(5401,-961){\makebox(0,0)[lb]{\smash{\SetFigFont{12}{14.4}{rm}$\cdots$}}}
\end{picture}

\]
\ni But the Temperley-Lieb algebra is known to be generated by
$\{\|\dots {}^{\cup}\!\!\!\!{}_{\cap_{i \ i+1}}\|\}$ whose images
are the $E_i$'s by (i). By condition (ii), we see that it suffices
to show that\newline $\Phi'(||\tangle{x})=
\Phi( || \tangle{x} )$.  
To this end we begin by showing
$$
\Phi'(\vln\vln\tanglel{a})= \Phi(\vln\vln\tanglel{a})
\quad {\text{for}} \ a\in N'\cap M_k \ .
$$
This follows from the picture below
\[
        \begin{picture}(0,0)%
\epsfig{file=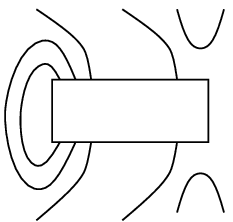}%
\end{picture}%
\setlength{\unitlength}{0.00041700in}%
\begingroup\makeatletter\ifx\SetFigFont\undefined
\def\x#1#2#3#4#5#6#7\relax{\def\x{#1#2#3#4#5#6}}%
\expandafter\x\fmtname xxxxxx\relax \def\y{splain}%
\ifx\x\y   
\gdef\SetFigFont#1#2#3{%
  \ifnum #1<17\tiny\else \ifnum #1<20\small\else
  \ifnum #1<24\normalsize\else \ifnum #1<29\large\else
  \ifnum #1<34\Large\else \ifnum #1<41\LARGE\else
     \huge\fi\fi\fi\fi\fi\fi
  \csname #3\endcsname}%
\else
\gdef\SetFigFont#1#2#3{\begingroup
  \count@#1\relax \ifnum 25<\count@\count@25\fi
  \def\x{\endgroup\@setsize\SetFigFont{#2pt}}%
  \expandafter\x
    \csname \romannumeral\the\count@ pt\expandafter\endcsname
    \csname @\romannumeral\the\count@ pt\endcsname
  \csname #3\endcsname}%
\fi
\fi\endgroup
\begin{picture}(2137,2055)(2529,-1726)
\put(3676,-811){\makebox(0,0)[lb]{\smash{\SetFigFont{12}{14.4}{rm}$a$}}}
\put(3451,-211){\makebox(0,0)[lb]{\smash{\SetFigFont{12}{14.4}{rm}$\cdots$}}}
\put(3451,-1411){\makebox(0,0)[lb]{\smash{\SetFigFont{12}{14.4}{rm}$\cdots$}}}
\end{picture}

\]
\ni for we know that
\[
\Phi'(|\tanglel{a})=
E_{M'}(a)=\Phi(|\tanglel{a}).
\]
\bs\ni
Now to show that \ $X=\Phi'( ||\tangle{x} )=
\Phi(||\tangle{x})=Y$  it suffices to
show that tr$(aX)={\op{tr}}(aY)$ for all
$a\in N'\cap M_k$. But up to powers of $\d$,
\[
	\begin{picture}(0,0)%
\epsfig{file=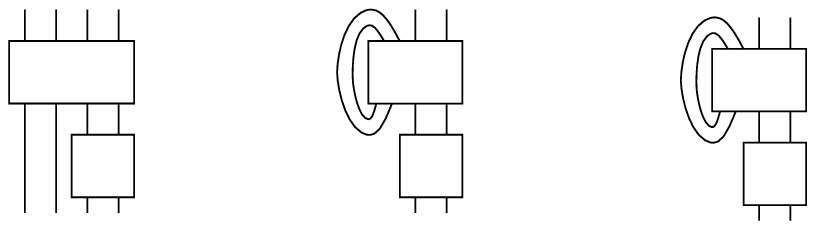}%
\end{picture}%
\setlength{\unitlength}{0.00041700in}%
\begingroup\makeatletter\ifx\SetFigFont\undefined
\def\x#1#2#3#4#5#6#7\relax{\def\x{#1#2#3#4#5#6}}%
\expandafter\x\fmtname xxxxxx\relax \def\y{splain}%
\ifx\x\y   
\gdef\SetFigFont#1#2#3{%
  \ifnum #1<17\tiny\else \ifnum #1<20\small\else
  \ifnum #1<24\normalsize\else \ifnum #1<29\large\else
  \ifnum #1<34\Large\else \ifnum #1<41\LARGE\else
     \huge\fi\fi\fi\fi\fi\fi
  \csname #3\endcsname}%
\else
\gdef\SetFigFont#1#2#3{\begingroup
  \count@#1\relax \ifnum 25<\count@\count@25\fi
  \def\x{\endgroup\@setsize\SetFigFont{#2pt}}%
  \expandafter\x
    \csname \romannumeral\the\count@ pt\expandafter\endcsname
    \csname @\romannumeral\the\count@ pt\endcsname
  \csname #3\endcsname}%
\fi
\fi\endgroup
\begin{picture}(10545,2069)(751,-5108)
\put(3977,-3737){\makebox(0,0)[lb]{\smash{\SetFigFont{12}{14.4}{rm}$a$}}}
\put(4277,-4637){\makebox(0,0)[lb]{\smash{\SetFigFont{12}{14.4}{rm}$x$}}}
\put(7427,-4637){\makebox(0,0)[lb]{\smash{\SetFigFont{12}{14.4}{rm}$x$}}}
\put(7202,-3737){\makebox(0,0)[lb]{\smash{\SetFigFont{12}{14.4}{rm}$a$}}}
\put(10727,-4712){\makebox(0,0)[lb]{\smash{\SetFigFont{12}{14.4}{rm}$x$}}}
\put(10502,-3812){\makebox(0,0)[lb]{\smash{\SetFigFont{12}{14.4}{rm}$a$}}}
\put(2401,-3736){\makebox(0,0)[lb]{\smash{\SetFigFont{12}{14.4}{rm}tr($\Phi'($}}}
\put(4726,-3736){\makebox(0,0)[lb]{\smash{\SetFigFont{12}{14.4}{rm}$)) = $tr$(\Phi'($}}}
\put(7876,-3736){\makebox(0,0)[lb]{\smash{\SetFigFont{12}{14.4}{rm}$)) = $tr$(\Phi ($}}}
\put(11251,-3811){\makebox(0,0)[lb]{\smash{\SetFigFont{12}{14.4}{rm}$)) = $tr$(aY)$}}}
\put(751,-3736){\makebox(0,0)[lb]{\smash{\SetFigFont{12}{14.4}{rm}tr($aX)=$}}}
\end{picture}

\]
\bs\ni
and we are done. \qed

\bs\ni{\bf Definition 4.2.8.} The {\it annular Temperley Lieb}
algebra $AT(n,\d)$, for $n$ even, will be the $*$-algebra with
presentation:
\[
        \input{xfig/pic106}
\]
\bs{\bf Remark.} Since $(F_1F_2\dots F_k)(F_1F_2\dots F_k)^*=\d F_1$
and $(F_1F_2\dots F_k)^*(F_1F_2\dots F_k)=\d F_k$, if $\cal H$
is a Hilbert space carrying a $*$-representation of $AT(n,\d)$,
dim$(F_i\cal H)$ is independent of $i$.

\bigskip
\ni{\bf Corollary 4.2.9} If $N\subset M$ is an extremal subfactor
of index $\d^{-2} >4$, each $N'\cap M_{k-1}$ is a Hilbert space
carrying a $*$-representation of $AT(2k,\d)$. And
${\op{dim}}(F_i(N'\cap M_{k-1}))={\op{dim}}(N'\cap M_{k-2})=
{\op{dim}}(M'\cap M_{k-1})$.

\bigskip
{\sc Proof.} Let $F_i$ be the elements of $\cal A(\phi)$ defined
as follows

\newpage

\ni The relations are easily checked, the Hilbert space structure
on $N'\cap M_{k-1}$ being given by the trace
$-\<a,b\>={\op{tr}}(b^*a)$. Since
$F_k(N'\cap M_{k-1})=N'\cap M_{k-2}$ and $F_{2k}
(N'\cap M_{k-1})=M'\cap M_{k-1}$ (by (ii) of 4.2.1).
We are through. \qed

\bigskip
\ni{\bf Lemma 4.2.10} If $\cal H$ carries an irreducible
 $*$-representation of $AT(n,\d)$ for $\d >2$, $n>4$, and
${\op{dim}}(F_i\cal H)=1$, then ${\op{dim}} \ \cal H=n$
(remember $n$ is even).

\bigskip
{\sc Proof.} Let $v_i$ be a unit vector in $F_i\cal H$ for each
$i$. Then $F_jv_i$ is a multiple of $v_j$ so the linear span of the
$v_i$'s is invariant, thus equal to $\cal H$ by irreducibility.
Here dim $\cal H\leq n$ (this does not require $\d >2$). Moreover
the commutation relations imply $|\<v_i,v_{i+1}\>|=\d^{-1}$ and
$\<v_i,v_j\>=0$ or 1 if $i\neq j\pm 1$. The case
$\<v_i,v_j\>=1$ forces $i=j$ so it only happens if $n=4$.
So by changing the $v_i$'s by phases we may assume that the matrix
$\d\<v_i,v_j\>$ is
\[
\D_n(\o) =\left(
\begin{array}{clcr}
\d & 1 & 0 & 0  \dots  \o \\
 1 & \d & 1 & 0  \dots  0 \\
 0 & 1 & \d & 1  \dots  0 \\
\vdots & {} & {} & {}    \vdots \\
\bar\o & 0 & \dots & {}   {1 \ \ \d}
\end{array}\right).
\]
It is easy to check that det$(\D_n(\o))=P_{2n}(\d)-P_{2n-2}(\d)-2
{\op{Re}}(\o)$ where $P_n(\d)$ are Tchebychev polynomials.
Thus det$(\D_n(\o))$ is smallest, for fixed $\d$, when $\o=1$.
But then \newline $\|\D_n(1)-\d \ {\op{id}}\|=2$ by Perron-Frobenius
so det $\D_n(1) >0$ for $\d >2$. \qed

\bs\ni
{\bf Corollary 4.2.11} Suppose the principal graph of the
subfactor $N\subset M$, $[M:N]>4$, has an initial segment equal to
the Coxeter-Dynkin diagram $D_{n+2}$ with $*$ as shown:
\[
        \begin{picture}(0,0)%
\epsfig{file=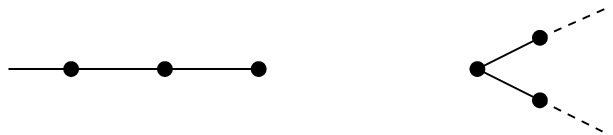}%
\end{picture}%
\setlength{\unitlength}{0.00041700in}%
\begingroup\makeatletter\ifx\SetFigFont\undefined
\def\x#1#2#3#4#5#6#7\relax{\def\x{#1#2#3#4#5#6}}%
\expandafter\x\fmtname xxxxxx\relax \def\y{splain}%
\ifx\x\y   
\gdef\SetFigFont#1#2#3{%
  \ifnum #1<17\tiny\else \ifnum #1<20\small\else
  \ifnum #1<24\normalsize\else \ifnum #1<29\large\else
  \ifnum #1<34\Large\else \ifnum #1<41\LARGE\else
     \huge\fi\fi\fi\fi\fi\fi
  \csname #3\endcsname}%
\else
\gdef\SetFigFont#1#2#3{\begingroup
  \count@#1\relax \ifnum 25<\count@\count@25\fi
  \def\x{\endgroup\@setsize\SetFigFont{#2pt}}%
  \expandafter\x
    \csname \romannumeral\the\count@ pt\expandafter\endcsname
    \csname @\romannumeral\the\count@ pt\endcsname
  \csname #3\endcsname}%
\fi
\fi\endgroup
\begin{picture}(6172,1244)(826,-1283)
\put(5626,-436){\makebox(0,0)[lb]{\smash{\SetFigFont{12}{14.4}{rm}$n$}}}
\put(3826,-736){\makebox(0,0)[lb]{\smash{\SetFigFont{12}{14.4}{rm}$\cdots$}}}
\put(826,-736){\makebox(0,0)[lb]{\smash{\SetFigFont{12}{14.4}{rm}$*$}}}
\end{picture}

\]
\ni Then there are at least two edges of the principal graph
connecting the two points at distance $n+1$ from $*$ to points of
distance $n+2$.

\bs
{\sc Proof.} Since $\d >2$, the Temperley-Lieb algebra generated
by $\{1,e_1,e_2,\dots ,e_k\}$ in $N'\cap M_k$ has dimension
$\frac{1}{k+2}{2k+4 \choose k+1}$. The information on
the principal graph then gives the following Bratteli diagram for
the inclusions $N'\cap M_{n-1}\subset N'\cap M_n\subset
N'\cap M_{n+1}$:
\begin{eqnarray*}
N'&\cap & M_{n+1} \\
&\cup &\\
N'&\cap & M_n\\
&\cup&\\
N'&\cap & M_{n-1}
\end{eqnarray*}
Here we have shown only that part of the Bratteli diagram relevant
to the proof. Only the two 1's in the middle row can be connected
to anything in the top row other than vertices corresponding to the
ideal generated by $e_{n+1}$. We have to show that it is impossible
for just one of these 1's to be connected, with multiplicity one,
to a new principal graph vertex. By contradiction, suppose this
were the case. Then we would have
$$
{\op{dim}}(N'\cap M_{n+1})=\frac{1}{n+3}
{2(n+2)\choose n+2} -(n+1)^2+(n+2)^2
$$
since the only difference between the $N'\cap M_{n+1}$ level
of the Bratteli diagram and the Temperley-Lieb Bratteli diagram
(see [GHJ]) is that the ``$n+1$" in Temperley-Lieb has become
``$n+2$".  Thus
$$
{\op{dim}}(N'\cap M_{n+1})=\frac{1}{n+3}
{2(n+2)\choose n+2} +2n+3 \ .
$$
But consider $N'\cap M_{n+1}$ as a module over $AT(2n+4,\d)$.
The Temperley-Lieb subalgebra is invariant and so therefore is its
orthogonal complement $TL^{\perp}$ of dimension $2n+3$. But
consider the image of $F_{n+2}=\d E_{N'\cap M_n}$. It is
$\dsize{\frac{1}{n+2}{2(n+1)\choose n+1} +1}$ because of the
single extra vertex on the principal graph at distance $n+1$ from $*$.

But, by pictures, the image of $F_{n+2}$ restricted to the
Temperley-Lieb subalgebra of  $N'\cap M_{n+1}$ is
$\dsize{\frac{1}{n+2}{2(n+1)\choose n+1}}$.
Hence on $TL^{\perp}$, \
dim$(F_1(TL^{\perp})=1$. So by Lemma 4.2.10,
dim$(TL^{\perp})\geq 2n+4$, a contradiction. \qed

We have obtained far more powerful results than the following by a
study of the representation theory of $AT(n,\d)$. These results
will be presented in a future paper of this series. We gave the
result here because it was announced some time ago. It is a version
of the ``triple point obstruction" of Haagerup and Ocneanu (see
[Ha]) but proved by a rather different method!

If we apply the argument we have just given when $\d \leq 2$
we obtain nontrivial but known results. The argument is very
simple so we present it.

{\bf Definition 4.2.12} The {\it critical depth} of a planar 
algebra $P$ will be the smallest $k$ for which there is an element
in $P_k$ which is not in the Temperley-Lieb subalgebra $TL_k$.

In the $C^*$- case, if $\d < 2$ ,the norm of the principal graph is less than $2$ so
as in [GHJ] it follows that the principal graph is an $A$,$D$ or
$E$ Coxeter graph with $*$ as far as possible from a vertex of 
valence $3$. In particular if $k$ is the critical depth, the 
dimension of the quotient $\frac {P_k}{V_k}$ is at most $1$ so the
rotation acts on it by multiplication by a $k$-th root of unity.
We will use the term ``chirality" for this root of unity in an
appropriate planar algebra.

The restrictions on the principal graph in the following theorem were first obtained by
Ocneanu.

{\bf Theorem 4.2.13} If $P$ is a $C^*$-planar algebra with $\d < 2$
then the principal graph can be neither $D_n$ with $n$ odd nor $E_7$.
If the principal graph is $D_{2n}$ the chirality is $-1$, if it
is $E_6$ the chirality is $e^{\pm 2\pi i/3}$ and if it is $E_8$ the 
chirality is $e^{\pm 2\pi i/5}$.

\bs
{\sc Proof.} If $k$ is the critical depth, by drawing diagrams one 
sees that the $\o$ in the $(2k+2)$x$(2k+2)$ matrix $\D_n(\o)$ is the chirality. Also if
$\d = z + z^{-1}$, we have  det$(\D_n(\o))= z^{2k+2} + z^{-(2k+2)} - \o - \o ^{-1}$.
If $\kappa$ is the Coxeter number of the principal graph we have 
$z = e^{\pm \pi i/{\kappa}}$. 
 
  On the other hand, by the argument of Corollary 4.2.11, the dimension
of $P_{k+1}$ would be too great if the determinant were non-zero. 
Thus we have, whatever the Coxeter graph may be,  

$$
{ e^{(2k+2)\pi i /\kappa} + e^{-(2k+2)\pi i /\kappa} = \o + \o ^{-1}}
$$

for a $k$-th root of unity $\o$.

The critical depth for $D_m$ would 
be $m-2$ so the left hand side of the equation is $-2$ so that $\o$ has to
be $-1$. But if $m$ is odd, $-1$ is not an $m-2$th root of unity so $D_m$ cannot
be a principal graph. If $m$ is even we conclude that the chirality is $-1$.

The critcal depth for $E_7$ is $4$ and the Coxeter number is 18. The above equation
clearly has no solution $\o$ which is a fourth root of unity.

For $E_6$ the critical depth is $3$ and the Coxeter number is $12$ so $\o = e^{\pm 2\pi i/3}$.
For $E_8$ the critical depth is $5$ and the Coxeter number is $30$ so $\o = e^{\pm 2\pi i/5}$.

\qed

  The above analysis may also be carried out for $\d = 2$ where the Coxeter graphs are
replaced by the extended Coxeter graphs. In the $D$ case the presence of the Fuss Catalan 
algebra of example $2.3$ makes it appropriate to replace the notion of critical depth
by the first integer such that $P_k$ is bigger than the Fuss Catalan algebra. One obtains 
then that the chirality, together with the principal graph, is a complete invariant for
$C^*$-planar algebras with $\d = 2$ (see [EK] p. 586).
  
We end this section by giving more details of the planar structure
on $N\subset M$. In particular we give the subfactor
interpretations of duality, reduction, cabling and tensor product.
The free product for subfactors is less straightforward.

\bs\ni
{\bf Corollary 4.2.12} If $N\subset M$ is an extremal $II_1$
subfactor then, with the notation of $\S$3.2,
$\l_1(P^{N\subset M})=P^{M\subset M_1}$, as planar algebras.
(Note that $P^{M\subset M_1}_k=M'\cap M_{k}$ which is a subset of
$P^{N\subset M}_{k+1}$. We are saying that the identity map is an
isomorphism of planar algebras.)

\bs
{\sc Proof.} Equality of $P^{N\subset M}_k$ and
$\l_1(P^{N\subset M})$ as sets follows immediately from (ii) of
4.2.1. To show equality of the planar algebra structure we use the
uniqueness part of 4.2.1.

By definition of $\l_1(\Phi)$, for $x\in M'\cap M_k\subset N'\cap
M_k$,
\[
x=\Phi^{M\subset M_1}(\boxe{x})=
\Phi^{N\cap M}(\tangle{x}\ )=
\l_1(\Phi)(\boxe{x}\ ).
\]
Properties (i) and (iv) are straightforward as is the second
equation of (ii). So we only need to check the first equation of
(ii), i.e.,
\[
\l_1(\Phi^{N\subset M})( \vln\ \tanglel{x}\ )=\d E_{M_1}(x) 
\mbox{for}  x\in M'\cap M_k.
\]
But by definition
\[
\l_1(\Phi^{N\subset M})(\vln\ \tanglel{x}\ )
=\frac 1\d \Phi(\vln\  \vln\ \tanglel{x}\ ).
\]
If we define $g:\otimes^k_N M\to \otimes^k_N M$ and
$f:\otimes^k_N M\to\otimes^{k-2}_N M$ by
$g(y)=\sum_{b\in B} byb^*$ and
$f(x_1\otimes\dots\otimes x_k)=E(x_1)x_2\otimes\dots\otimes
x_{k-1} E(x_k)$, and if $\th(y)=x$ for $x\in M'\cap M_k$, by the
definition of $\Phi^{N\subset M}$ in Theorem 4.2.1,
\[
\Phi^{N\subset M}(\vln\ \vln\  \tanglel{x}\ )=\th
(\sum_{b\in B} b\otimes f(g(y))\otimes b^*).
\]
But since $\th$
is an $M-M$ bimodule map and $x$ commutes with $M$, \
$g(y)=\d^2y$. But by definition of $\th$ and 4.1.8,
\[
\th(\sum_{b\in B} b\otimes f(y)\otimes b^*)=
\d^{-2}\sum_{b\in B} bE_1 xE_1b^* \ .
\]
Since $\{bE_1\}$ is a basis for $M_1$ over $M$, we are done by
4.1.5. \qed

\bs Iterating, we see that $\l_n(P^{N\subset M})$ is the planar
algebra for the subfactor $M_{n-1}\subset M_n$.

For cabling we have the following

\bigskip
\ni{\bf Corollary 4.2.13} If $N\subset M$ is an extremal $II_1$
subfactor, the cabled planar algebra ${\cal C}_n(P^{N\subset M})$ of
$\S3.3$ is isomorphic to $P^{N\subset M_{n-1}}$.

\bigskip
{\sc Proof.} We will again use the uniqueness part of 4.2.1.
It follows from [PP2] that if we define $E^n_i$ to be
$$
(E_{ni}E_{ni-1}E_{ni-2}\dots E_{n(i-1)+1})
(E_{ni+1}E_{ni}\dots E_{n(i-1)+2})\dots
(E_{n(i+1)-1}E_{n(i+1)}\dots E_{ni})
$$
(a product of $n$ products of $E$'s with indices decreasing by
one), and $v^n_i$ with respect to $E^n_i$ as in 4.1.9, then the map
$x_1\otimes_N x_2\dots \otimes_N x_k\to x_1v^n_1 x_2v^n_2\dots
v^n_kx_k$ establishes, via the appropriate $\th$'s, a $*$-algebra
isomorphism between the $(k-1)$-th algebra in the tower for
$N\subset M_{n-1}$ and $M_{kn-1}$, hence between
$P_k^{N\subset M_{n-1}}$ and $P_{nk}^{N\subset M}=\cal C_n(P^{N\subset M})$.

If $\g$ is the inverse of this map, $\g\circ\cal C_n(\Phi)$ thus
defines a spherical planar algebra structure on
$P^{N\subset M_{n-1}}$. The labelling set is identified with
$P^{N\subset M_{n-1}}$ via $\g$, so
\[
\g\circ\Bbb C_n(\Phi^{N\subset M})(\tangle{x}\ )=
\Phi^{N\subset M})(\tangle{x})=x.
\]
Condition (i) of 4.1 for $\g\circ\Bbb C_n(\Phi)$ follows by
observing that\\ 
$\Bbb C_n(\Phi)(\boxj )=E^n_i$.
Condition (ii) follows from 4.2.12, and condition (iv) is clear.
\qed

\bs Reduction is a little more difficult to prove.

\bigskip
\ni{\bf  Corollary 4.2.14} Let $N\subset M$ be an extremal $II_1$
subfactor and $p$ a projection in $N'\cap M$. The planar algebra
$p(P^{N\subset M})p$ of $\S$3.3 is naturally isomorphic to the
planar algebra of the reduced subfactor
$pN\subset pMp$.

\bigskip
{\sc Proof.}  We first claim that the tautological map
$$
\a:\bigotimes^k_{pN} pMp\to\bigotimes^k_N M \ , \qquad
\a(\bigotimes^k_{i=1} x_i)=\bigotimes^k_{i=1} x_i
$$
is an injective homomorphism of $*$-algebras when both domain and
range of $\a$ are equipped with their algebra structures via the
respective maps $\th$ as in 4.1.10. To see this observe that the
conditional expectation $E_{pN}: pMp\to pN$ is just
$\frac{1}{\op{tr}(p)} E_N$. It follows that
$\a a_j=a_j\a$ and clearly $\a\mu_j=\mu_j\a$, so by 4.1.24,  $\a$
is a $*$-algebra homomorphism.

The next thing to show is that $\a$ takes $pN$-central vectors to
$\th^{-1}(p_k(N'\cap M_{k-1})p_k)$, with $p_k$ as in 3.3, using
the planar algebra structure on $N'\cap M_k$. Let
$\pi:\otimes^k_NM\to \otimes^k_N M$ be the map
$\pi(x_1\otimes x_2\dots\otimes x_k)=px_1p\otimes px_2p
\otimes \dots\otimes px_kp$, which is well defined since $p$
commutes with $N$. Then a diagram shows that, for
$x\in N'\cap M_{k-1}$,
$p_kxp_k=\th(\pi(\th(x)))$, and conversely, if
$\th^{-1}(x)$ is in the image of $\a$,
$\pi(\th^{-1}(x))=\th^{-1}(x)$ so $p_kxp_k=x$. Hence $\a$ induces
a $*$-algebra isomorphism between $P_k^{pN\subset pMp}$ and
$p_k(N'\cap M_{k-1})p_k$. To check that this map induces the
right planar algebra structure, we first observe that $\a$
commutes suitably with the maps $\eta$ and $\kappa$ of 4.1.19.
For $\eta$ we have $\a\circ\eta=\pi\circ\eta\circ\a$ by
definition. For $\kappa$, note that if we perform the basic
construction of [J1] on $L^2(M)$, \
$\frac{1}{\op{tr}(p)} pE_Np$ is the basic construction projection
for $pN\subset pMp$ on $L^2(pMp)$. Thus if $\{b\}$ is a basis for
$pMp$ over $pN$ we have
$\sum_b bE_N b^*={\op{tr}}(p)p$. So if
$$
x=x_1\otimes \dots\otimes x_k\in \bigotimes^k_{pN} pMp \ ,
$$
$$
\a(\kappa_j(x))=\frac{1}{\op{tr}(p)\d} \sum_b
x_1\otimes \dots\otimes x_j
b\otimes b^*\otimes x_{j+1}\otimes\dots\otimes x_k
$$
and
$$
\kappa_j(\a(x))=\frac 1\d \sum_c
x_1\otimes \dots\otimes x_j
c\otimes c^*\otimes x_{j+1}\otimes\dots\otimes x_k \ .
$$
Applying $\th$ we see that $\a \kappa_j=\kappa_j\a$.

With these two commutation results it is an easy matter to check
that $\a$ defines planar algebra isomorphism between
$P^{Np\subset pMp}$ and $p(P^{N\subset M})p$, using the uniqueness
part of 4.2.1 or otherwise.

\bigskip
\ni{\bf Corollary 4.2.15} If $N_1\subset M_1$ and
$N_2\subset M_2$ are extremal finite index subfactors, then
$P^{N_1\otimes N_2\subset M_1\otimes M_2}$ is naturally isomorphic
to the tensor product $P^{N_1\subset M_1}\otimes
P^{N_2\subset M_2}$ of $\S$3.4.

\bigskip
{\sc Proof.} We leave the details to the reader. \qed

\bs\bs\ni{\bf 4.3 \ Planar algebras give subfactors}

\smallskip
The following theorem relies heavily on a result of Popa [Po2].

\bigskip
\ni{\bf Theorem 4.3.1} Let $(P,\Phi)$ be a spherical $C^*$-planar
algebra with invariant $Z$ and trace ${\op{tr}}$. Then there is a
subfactor $N\subseteq M$ and isomorphisms $\Omega: N'\cap M_i\to
P_i$ with
\begin{verse}
(i) $\Omega$ is compatible with inclusions

(ii) ${\op{tr}}(\Omega(x))={\op{tr}}(x)$

(iii) $\Omega(M'\cap M_i)=P_{1,i}$
 (linear span of tangles with vertical first string)

(iv) $\Omega(e_i)=\frac{1}{\d}(\Phi(\bigtank ) )$

(v) $[M:N]=\d^2$. If $x\in N'\cap M_i$, choose $T\in\P$
with $\Phi(T)=\Omega(x)$.

(vi) $\Omega(E_{M'}(x))=\frac{1}{\d}\Phi (\boxk{T})$

(vii) $\Omega(E_{M_{i-1}}(x))=\frac{1}{\d}\Phi (\boxl{T})$
\end{verse}

\bigskip
{\sc Proof.} By theorem 3.1 of [Po2], the pair $N\subseteq M$
exists given a system $(A_{ij})$, $0 \leq i < j < \i$ of finite
dimensional $C^*$-algebras with $A_{i,j}\subset A_{k,\ell}$ if
$k\leq i$, $j\leq\ell$ and a faithful trace on $\bigcup^\i_{n=0}
A_{0n}$ satisfying 1.1.1, 1.1.2, 1.3.3$'$ and 2.1.1 of
[Po2]. We set $A_{i,j}=P_{i,j}$ (Definition 1.20). Then
$A_{ii}=\Bbb C$ since $Z$ is multiplicative and non-degenerate. The
conditions of 
[Po2]  involve $e_i$'s and conditional
expectations
$E_{A_{ij}}$. We define the $e_i$'s in $P$ to be what we have called
$\Omega(e_i)$ (note our $e_i$ is
Popa's ``$e_{i+1}$"). The map $E_{A_{ij}}$ is defined by the
relation tr$(xE_{A_{ij}}(y))={\op{tr}}(xy)$ for $x$ in $A_{ij}$ and
$y$ arbitrary. Since $Z$ is an $S^2$ invariant one easily checks
that $E_{A_{ij}}$ is given by the element of $\cal A_{kj}(\emptyset)$ 
given in the figure below (for $x\in A_{0,k}$) 
\[
	\begin{picture}(0,0)%
\epsfig{file=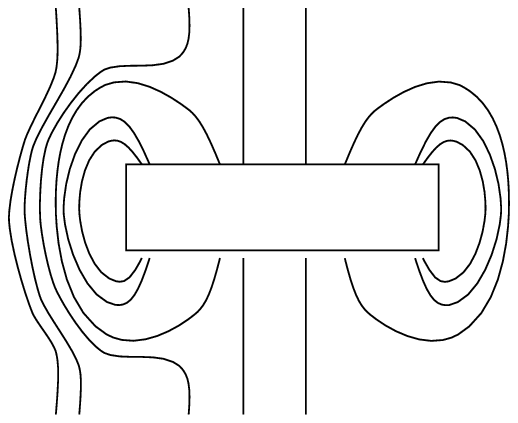}%
\end{picture}%
\setlength{\unitlength}{0.00041700in}%
\begingroup\makeatletter\ifx\SetFigFont\undefined
\def\x#1#2#3#4#5#6#7\relax{\def\x{#1#2#3#4#5#6}}%
\expandafter\x\fmtname xxxxxx\relax \def\y{splain}%
\ifx\x\y   
\gdef\SetFigFont#1#2#3{%
  \ifnum #1<17\tiny\else \ifnum #1<20\small\else
  \ifnum #1<24\normalsize\else \ifnum #1<29\large\else
  \ifnum #1<34\Large\else \ifnum #1<41\LARGE\else
     \huge\fi\fi\fi\fi\fi\fi
  \csname #3\endcsname}%
\else
\gdef\SetFigFont#1#2#3{\begingroup
  \count@#1\relax \ifnum 25<\count@\count@25\fi
  \def\x{\endgroup\@setsize\SetFigFont{#2pt}}%
  \expandafter\x
    \csname \romannumeral\the\count@ pt\expandafter\endcsname
    \csname @\romannumeral\the\count@ pt\endcsname
  \csname #3\endcsname}%
\fi
\fi\endgroup
\begin{picture}(7069,4369)(226,-4408)
\put(3901,-1486){\makebox(0,0)[lb]{\smash{\SetFigFont{12}{14.4}{rm}$\cdots$}}}
\put(4801,-1486){\makebox(0,0)[lb]{\smash{\SetFigFont{12}{14.4}{rm}$\cdots$}}}
\put(4801,-3736){\makebox(0,0)[lb]{\smash{\SetFigFont{12}{14.4}{rm}$\cdots$}}}
\put(3451,-3961){\makebox(0,0)[lb]{\smash{\SetFigFont{12}{14.4}{rm}$\cdots$}}}
\put(2851,-4336){\makebox(0,0)[lb]{\smash{\SetFigFont{12}{14.4}{rm}1}}}
\put(3376,-361){\makebox(0,0)[lb]{\smash{\SetFigFont{12}{14.4}{rm}$\cdots$}}}
\put(3151,-4336){\makebox(0,0)[lb]{\smash{\SetFigFont{12}{14.4}{rm}2}}}
\put(4126,-4336){\makebox(0,0)[lb]{\smash{\SetFigFont{12}{14.4}{rm}$i$}}}
\put(226,-2086){\makebox(0,0)[lb]{\smash{\SetFigFont{12}{14.4}{rm}$E_{A_{ij}}:\frac{1}{\delta^{i+k-j}}$}}}
\put(5851,-1486){\makebox(0,0)[lb]{\smash{\SetFigFont{12}{14.4}{rm}$\cdots$}}}
\put(5401,-4336){\makebox(0,0)[lb]{\smash{\SetFigFont{12}{14.4}{rm}$j$}}}
\put(4426,-4336){\makebox(0,0)[lb]{\smash{\SetFigFont{12}{14.4}{rm}$i+1$}}}
\end{picture}

\]
\ni Popa's  (1.1.1) and (1.1.2) and b)$'$ of (1.3.3) follow
immediately from pictures (note that the power of $\frac 1\d$ is
checked by applying $E_{A_{ij}}$ to 1). Condition a)$'$ of 1.3.3$'$
is dim $A_{ij}=
{\op{dim}} \ A_{i,j+1}e_j={\op{dim}} \ A_{i-1,j+1}$. But it is easy
from pictures that $E_{A_{ij}}$ defines a linear map from
$A_{i,j+1}e_j$ onto $A_{i,j}$, whose inverse is to embed in
$A_{i,j+1}$ and multiply on the right by $e_j$. Moreover th element
of $\cal A(\emptyset)$ (illustrated for $i\!=\!0$) in Figure 1
defines a linear isomorphism from $A_{i,j}$ to $A_{i+1,j+1}$ ---
the inverse is a similar picture.
\[
	\begin{picture}(0,0)%
\epsfig{file=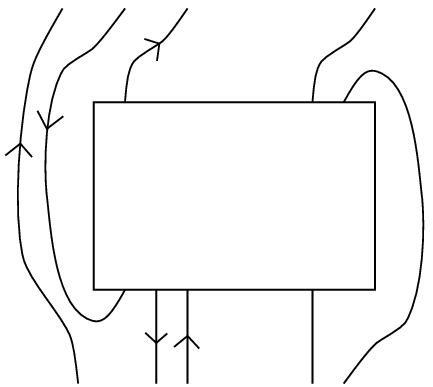}%
\end{picture}%
\setlength{\unitlength}{0.00041700in}%
\begingroup\makeatletter\ifx\SetFigFont\undefined
\def\x#1#2#3#4#5#6#7\relax{\def\x{#1#2#3#4#5#6}}%
\expandafter\x\fmtname xxxxxx\relax \def\y{splain}%
\ifx\x\y   
\gdef\SetFigFont#1#2#3{%
  \ifnum #1<17\tiny\else \ifnum #1<20\small\else
  \ifnum #1<24\normalsize\else \ifnum #1<29\large\else
  \ifnum #1<34\Large\else \ifnum #1<41\LARGE\else
     \huge\fi\fi\fi\fi\fi\fi
  \csname #3\endcsname}%
\else
\gdef\SetFigFont#1#2#3{\begingroup
  \count@#1\relax \ifnum 25<\count@\count@25\fi
  \def\x{\endgroup\@setsize\SetFigFont{#2pt}}%
  \expandafter\x
    \csname \romannumeral\the\count@ pt\expandafter\endcsname
    \csname @\romannumeral\the\count@ pt\endcsname
  \csname #3\endcsname}%
\fi
\fi\endgroup
\begin{picture}(4056,3637)(1231,-3983)
\put(3226,-961){\makebox(0,0)[lb]{\smash{\SetFigFont{12}{14.4}{rm}$\cdots$}}}
\put(3301,-3661){\makebox(0,0)[lb]{\smash{\SetFigFont{12}{14.4}{rm}$\cdots$}}}
\end{picture}

\]
\ni Finally  the commutation relations 2.1.1,
$[A_{ij},A_{k\ell}]=0$ for $i\leq j\leq k \leq\ell$ are trivial
since they are true in $\P$, involving non-overlapping strings.

It is standard theory for external subfactors that $E_{M'}$
restricted to $N'\cap M_i$ is $E_{A_{1,i}}$ and
$E_{M_{i-1}}=E_{A_{1,i-1}}$. So (vi) and (vii) are clear. \qed

\bigskip
\ni{\bf Corollary 4.3.2} If $(P,\Phi)$ is as before, the
Poincar\'e series $\sum^\i_{n=0}\dim (P_n)z^n$ has radius of
convergence $\geq\frac{1}{\d^2}$.

\bigskip
This result could be proved without the full strength of Theorem
4.3.1 (as pointed out by D.Bisch), using the principal graph and
the trace. The Poincar\'e series of planar subfactors enjoy many
special properties as we shall explore in future papers.

\newpage
\centerline{\bf References}

\begin{verse}
[A] V. Arnold, Remarks on the enumeration of plane curves.
{\it Amer. Math. Soc. Transl.} {\bf 173}(1996), 17--32.

[Bi] D.Bisch, Bimodules, higher relative commutants and the fusion algebra associated
to a subfactor,
{\it The Fields Institute for Research in Mathematical 
Sciences Communications Series} {\bf 13} (1997), 13--63.
 
[BW] J.Barrett and B.Westbury, ``Spherical Categories", hep-th/9310164.

[BiW] J.Birman and H.Wenzl, Braids, link polynomials and a new algebra.
{\it Trans. Amer. Math. Soc.} {\bf 163} (1989), 249--273.

[Ba] R. Baxter, {\it Exactly Solved Models in Statistical
Mechanics},   Academic Press, New York, 1982.

[Ban] T. Banica, Hopf algebras and subfactors associated to vertex models,
math.QA/9804016

[BH] D.Bisch and U.Haagerup, Composition of subfactors: new examples of infinite depth subfactors,
{\it Ann. scient. {\'E}c. Norm. Sup.} {\bf 29} (1996),329--383.

[BJ1] D.Bisch and V.Jones, ``Singly generated planar
algebras of small dimension",  preprint 1998.

[BJ2] D.Bisch and V.Jones, Algebras associated to
intermediate subfactors, {\it Invent. Math.} {\bf 128}, 89--157.
(1997).

[BJ3] D.Bisch and V.Jones, A note on free composition of
subfactors, Geometry and physics
 (Aarhus, 1995), 339--361 in {\it Lecture Notes in Pure and Appl.
Math.} {\bf 184}, Dekker, New York, 1997

[Co1] A.Connes, {\it Noncommutative Geometry}, Academic
Press (1994).

[dlH] P. de la Harpe,
Spin models for link polynomials, strongly regular graphs and Jaeger's Higman-Sims model.
{\it Pacific J. Math} {\bf 162}
(1994), 57--96.

[Dr] V.Drinfeld, Quantum groups, {\it Proceedings ICM
1986}, vol.~{\bf 1}, pp.798--820.

[EK] D.Evans and Y.Kawahigashi, {\it Quantum symmetries on operator algebras},
Oxford University Press (1998).

[F+] P.Freyd, D.Yetter, J.Hoste, W.Lickorish, K.Millett,
and A.Ocneanu, A new polynomial invariant of knots and links.
{\it Bull. AMS} {\bf 12} (1985), 183--190.

[FY] P.Freyd and D.Yetter,
 Braided compact closed categories with applications to low
dimensional topology, {\it Adv. in Math.} {\bf 77}
(1989), 156--182.

[FRS] K.~Fredenhagen, K-H.~Rehren and B.~Schroer,
Superselection sectors with braid group statistics and exchange
algebras, {\it Comm.~Math.~Phys.} {\bf 125} (1989), 201--226.

[GHJ] F.Goodman, P. de la Harpe and V.F.R.Jones, {\it
Coxeter Graphs and  Towers of Algebras}, Springer Verlag, MSRI
publications (1989).

[GL] J. Graham and G. Lehrer , ``The representation theory
of affine Temperley-Lieb algebras", preprint, University of Sydney
(1997).

[Gn] S. Gnerre, Free composition of paragroups,
Thesis, UC Berkeley, (1997).

[Ha] U.Haagerup,
Principal graphs of subfactors in the index range $4 < [M:N] <
3 + \sqrt{2}$
{\it Subfactors}, World Scientific, Singapore-New
Jersey-London-Hong Kong (1994), 1--39.

[JMN] F.Jaeger, M.Matsumoto and  K.Nomura, ``Bose-Mesner
algebras related with type  II matrices and spin models," RIMS
preprint (1995).

[Ja] F.JaegerStrongly regular graphs and spin models for the Kauffman polynomial, 
{\it Geom. Dedicata } {\bf 44} (1992), 23--52

[J1] V.Jones,  Index for subfactors, {\it Invent.
Math.} {\bf 72}  (1983), 1--25.

[J2]  V.Jones, A polynomial invariant for knots via
von Neumann algebras.
       {\it Bulletin of the Amer. Math. Soc.} {\bf 12} (1985),
103--112.

[J3] V.Jones,  Hecke algebra representations of
braid groups and link polynomials, {\it Ann. Math.} {\bf 126}
(1987), 335--388.

[J4] V.Jones. On a certain value of the Kauffman polynomial,
{\it Comm. Math. Phys.} {\bf 125} (1989), 459-467

[J5] V.Jones Index for subrings of rings.
{\it Contemporary Math.} {/bf 43} (1985), 181--190.

[JS] V.F.R.Jones and V.Sunder, Introduction to
Subfactors, LMS lecture note
        series no.~{\bf 234} (1997) 162 pages.

[Ka1] L.Kauffman,
 State models and the Jones polynomial, {\it Topology} {\bf 26}
(1987), 395--407.

[Ka2] L.Kauffman,
An Invariant of regular isotopy, {\it Trans. Amer. Math. Soc.}
{\bf 318} (1990), 417--471.

[KS] Krishnan amd V.Sunder, On biunitary permutation matrices and some subfactors of index $9$, 
 {\it  Trans. Amer. Math. Soc} {\bf 348} (1996), 4691--4736.

[Ku] G.Kuperberg, Spiders for rank 2 Lie algebras, {\it
Commun. Math. Phys.} {\bf 180} (1996), 109--151, q-alg/9712003.

[La] Z. Landau, Thesis, UC Berkeley, 1998

[Lo] R.Longo, Index of subfactors and statistics of
quantum fields, I.
 {\it Comm. Math. Phys.} {\bf 126} (1989), 217--247.

[LM] W.Lickorish and K.Millett, Some evaluations of link 
polynomials.
{\it Comment. Math. Helv.} {\bf 61} (1986), 349--359.

[LR] R.Longo and J. Roberts, A Theory of Dimension,
{\it Journal of K-theory} {\bf 11} (1997), 103--159, funct-an/9604008.

[Ma] J.May, Definitions: operads, algebras and
modules, {\it Contemporary
              Mathematics } {\bf 202} (1997) 1--7.

[MvN] F.Murray and J.von Neumann, On rings of operators,
IV, {\it Annals Math.} {\bf 44} (1943)
      716--808.

[Mu] J.Murakami The Kauffman polynomial of links and representation theory,
{\it Osaka Journal of Mathematics} {\bf 24} (1987), 745--758.

[O1] A.Ocneanu,
 Quantum symmetry, differential geometry of finite graphs and
classification of subfactors,  University of Tokyo Seminary Notes
(notes recorded by Y. Kawahigashi).

[Pe] R.Penrose, Application of negative dimensional
tensors, in  {\it Combinatorial
    Mathematics and Its Applications}, Academic Press (1971)
221--244.

[Po1] S. Popa, Classification of subfactors and their
endomorphisms,  {\it CBMS Lecture Notes}, {\bf 86} (1995).

[Po2] S. Popa, An axiomatization of the lattice of higher
relative commutants of a subfactor, {\it Invent. Math.} {\bf120}
(1995), 427--445.

[PP] M.Pimsner and S.Popa, Entropy and index for
subfactors,  {\it Ann.
Sci. Ecole Norm. Sup.} {\bf19} (1986), 57--106.

[TL] H. Temperley and E.Lieb, Relations between the percolation....,
{\it Proc. Royal Soc. London} {\bf A322} (1971), 251--280.

[Tu] V.Turaev , {\it Quantum Invariants for knots and $3-$manifolds},
de Gruyter Berlin (1994).

[TV] V.Turaev and O.Viro, State sum invariants of
3-manifolds and quantum 6-j symbols. {\it
     Topology} {\bf 31} (1992) 86--902.

[Tut] W.Tutte, On dichromatic polynomials, 
{\it J. Combinatorial Theory} {\bf 2} (1967), 301--320.

[V] D.Voiculescu, Multiplication of certain non-commuting
random variables. {\it J. Operator Theory} {\bf 18} (1987),223--235.

[Wa] A.Wassermann, Operator algebras and conformal field theory, III:
Fusion of positive energy representations of LSU(N) using bounded
operators,
{\it Inventiones Math.} {\bf 133} (1998) 467--538, math.OA/9806031.

[We1]  H.~Wenzl,  Hecke algebras of type $A_n$ and
subfactors, {\it Invent.~Math.} {\bf 92} (1988), 249--383.

[We2]  H.~Wenzl,  Quantum groups and subfactors of type B, C
  and D, {\it Comm. Math. Phys.} {\bf 133} (1990), 383--432.

[We3]  H.~Wenzl, C* tensor categories from quantum groups,
{\ J. Amer. Math. Soc.} {\bf 11} (1998), 261--282. 

[X] F.Xu Standard $\lambda$-lattices from quantum groups,
{\it Inventiones Math.} {\bf 134} (1998),
\end{verse}

\enddocument

\end{document}